\documentclass[11pt]{article}
\usepackage[T1]{fontenc}
\usepackage{lmodern}
\usepackage{amsmath,amssymb,amsthm,mathtools}
\usepackage[margin=1in]{geometry}
\usepackage{booktabs,array,longtable,tabularx,tikz}
\usepackage[dvipsnames]{xcolor}
\usepackage[unicode,hypertexnames=false,colorlinks=true,linkcolor=MidnightBlue,citecolor=MidnightBlue,urlcolor=MidnightBlue]{hyperref}
\allowdisplaybreaks[2]
\newcommand{\Szego}{\texorpdfstring{Szegő}{Szego}}
\newcommand{\latintext}{\normalfont}

\newtheorem{theorem}{Theorem}[section]
\newtheorem{lemma}[theorem]{Lemma}
\newtheorem{proposition}[theorem]{Proposition}
\newtheorem{corollary}[theorem]{Corollary}
\newtheorem*{theoremA}{Theorem A}
\newcounter{turnpart}[theorem]

\newtheorem{turnlemma}[turnpart]{Lemma}
\theoremstyle{definition}

\newtheorem{remark}[theorem]{Remark}
\newtheorem{question}[theorem]{Question}
\theoremstyle{remark}

\hypersetup{pdftitle={The P\'olya--Szeg\H{o} inequalities for polygons},pdfauthor={Hedong Wang}}
\title{The P\'olya--\Szego{} inequalities for polygons}
\author{Hedong Wang\\Georgetown University}
\date{}
\begin{document}
\maketitle

\begin{abstract}
For $1\le q\le2$, let $\lambda_q^{\mathrm{S}}(P)$ be the infimum of $\int_P|\nabla u|^2/(\int_P|u|^q)^{2/q}$ over nonzero $u\in H_0^1(P)$. Among simple polygons of prescribed area with at most $n$ sides, $n\ge3$, the regular $n$-gon is the unique minimizer up to congruence. The cases $q=1$ and $q=2$ give the polygonal torsional rigidity and first Dirichlet eigenvalue inequalities conjectured by P\'olya and Szeg\H{o}. We partition the polygon into convex cells, each with a Dirichlet base on a distinct side of the polygon and Neumann conditions elsewhere. Rearrangement on a strip compares each cell with an isosceles triangle. For every $1\le q\le2$, the normalized mixed quotient of this triangle is strictly convex in its opening angle on $(0,\pi)$. At fixed area $A$, the normalized deficit is bounded below by $c_n\min\{1,d^2/A\}$ with $c_n>0$ independent of $q$, where $d$ is the Hausdorff distance to the nearest congruent copy of the equal-area regular $n$-gon. We give explicit constants for torsion and the first eigenvalue. At fixed perimeter the quadratic bound is uniform in $n$. The exponent two is optimal in both estimates. We also treat polygons with unrestricted numbers of reflex vertices. At fixed area the equilateral triangle and the square minimize $\lambda_q^{\mathrm{S}}$ among polygons with at most three and at most four convex vertices, respectively. If $n\ge5$ and $P$ has at most $n-1$ convex vertices, then $\lambda_q^{\mathrm{S}}(P)>(1+1/(1000qn^4))\lambda_q^{\mathrm{S}}(R_n)$ for the equal-area regular $n$-gon $R_n$. Among polygons of fixed area with at most $n$ convex vertices, $R_n$ is a strict local minimizer in Hausdorff distance modulo rigid motions.
\end{abstract}



\part{The polygonal torsion and eigenvalue conjectures}\label{part:endpoints}
\section{Introduction}\label{sec:introduction}
The polygonal conjectures of P\'olya and \Szego{} assert that among polygons of prescribed area and side count the regular polygon maximizes torsional rigidity and minimizes the first Dirichlet eigenvalue \cite{ps1951}. The triangle and quadrilateral cases are classical \cite{polya1948}, \cite[p.~158]{ps1951}. Beyond four sides, Steiner symmetrization may increase the number of sides. Solynin proved that the regular polygon maximizes the conformal radius \cite{solynin1988}, and Solynin and Zalgaller proved that it minimizes logarithmic capacity \cite{sz2004}. Their discussion leaves open the torsion and eigenvalue inequalities \cite[p.~278]{sz2004}.

Let $P$ be the bounded open region enclosed by a simple polygonal curve and let $R_n$ be the regular $n$-gon with the same area. We count sides after removing vertices at which adjacent sides are collinear. For a bounded open set $\Omega$ and $1\le q\le2$ put
\[
 \lambda_q^{\rm S}(\Omega)=\inf_{0\ne v\in H^1_0(\Omega)}
 \frac{\int_\Omega |\nabla v|^2}{(\int_\Omega |v|^q)^{2/q}} .
\]
Then $\lambda_1^{\rm S}=1/T$, where $T(\Omega)=\int_\Omega u$ for $-\Delta u=1$ with zero boundary values, and $\lambda_2^{\rm S}=\lambda_1$ is the first Dirichlet eigenvalue. The quantity $|\Omega|^{2/q}\lambda_q^{\rm S}(\Omega)$ is invariant under dilation.

\begin{theorem}\label{thm:main}
Let $n\ge3$ and $1\le q\le2$. Every simple polygon $P$ with at most $n$ sides satisfies
\[
 |P|^{2/q}\lambda_q^{\rm S}(P)\ge|R_n|^{2/q}\lambda_q^{\rm S}(R_n).
\]
Equality holds exactly when $P$ is congruent to $R_n$. In particular $T(P)\le T(R_n)$ and $\lambda_1(P)\ge\lambda_1(R_n)$ at equal area.
\end{theorem}

The comparison is quantitative. Theorems~\ref{thm:stability} and~\ref{S:thm:full-stability} give quadratic Hausdorff stability at fixed area, for every $n\ge3$, with explicit constants for torsion and the first eigenvalue and a non-explicit constant independent of $q$ for the intermediate exponents. Theorem~\ref{thm:perimeter-stability} gives a bound at fixed perimeter that is uniform in the number of sides. The exponent two is optimal in each case.

\paragraph{Counting convex vertices.} Part~\ref{part:convex-vertices} bounds the number of convex vertices instead of the number of sides and allows any number of reflex vertices. With three or four convex vertices the equilateral triangle, respectively the square, minimizes $|P|^{2/q}\lambda_q^{\rm S}(P)$ for every $1\le q\le2$ (Theorems~\ref{C:thm:three-convex-turns-sharp} and~\ref{C:thm:four-convex-sharp}). If $n\ge5$ and $P$ has at most $n-1$ convex vertices, its normalized Sobolev quotient exceeds that of $R_n$ by an explicit factor for every $q$ (Theorem~\ref{C:thm:convex-vertices} and Corollary~\ref{C:cor:all-q-gap}). Among polygons with at most $n$ convex vertices, $R_n$ is a strict local minimizer of $|P|^{2/q}\lambda_q^{\rm S}(P)$ for every $q$ (Theorem~\ref{C:thm:local}), and the least torsion deficit agrees with that of $R_n$ up to $O(n^{-5})$ (Theorem~\ref{C:thm:sharp-asymptotic}). Whether $R_n$ is optimal in the whole class is open for $n\ge5$ (Question~\ref{C:que:sharp}).

\paragraph{The three ingredients.}
The positive minimizer on a regular $n$-gon is symmetric across each radius joining the center to a vertex. On each of the resulting triangles it satisfies a Dirichlet condition on the base and Neumann conditions on the other two sides. For an isosceles triangle of opening $\alpha$, write $F_q(\alpha)$ for its normalized mixed quotient. We call $\alpha\mapsto F_q(\alpha)$ the profile. In the regular polygon the opening is $a=2\pi/n$. Reflection gives the normalized quotient $nF_q(a)$ for the polygon. Three facts allow a comparison with a general polygon.

\emph{A cell comparison.} On a convex cell $Q$ with a single Dirichlet base, let $\beta,\beta'$ be the base angles. When $\beta+\beta'<\pi$, the opening is $\alpha=\pi-\beta-\beta'$. Lemma~\ref{lem:mixed-cell} shows that the isosceles triangle of opening $\alpha$ has the least normalized mixed quotient, for every $1\le q\le2$. In strip coordinates we express the area density of $Q$ as an average of translates of explicit densities. Two rearrangements of a test function then show that, at fixed Dirichlet energy, the largest possible area average of $|u|^q$ occurs for the isosceles triangle.

\emph{A partition.} To sum the cell inequalities, we need a lower bound on each opening in terms of both the area of the cell and its share of the $L^q$ mass. Let $u$ be the positive minimizer with $\int_Pu^q=1$ and assign a cell the probability $\nu(Q)=\tfrac12(|Q|/|P|+\int_Qu^q)$. Write the tangent of $F_q$ at the regular opening as $s(\alpha-b)$. Lemma~\ref{lem:mass-cap-partition} constructs, inside an arbitrary simple polygon, a partition into convex cells on distinct sides such that, if $P$ has $N$ sides,
\[
 \alpha_i-b\ge(2\pi-Nb)\,\nu(Q_i).
\]

After each cut, let $m$ be the number of original sides still assigned to a piece $C$. At each stage we require $C$ to have an opening of at least $mb+(2\pi-Nb)\nu(C)$. The masses of the pieces add, and each original side retained after a cut is assigned to just one piece. When a straight segment to a reflex vertex would leave the polygon, the construction can use the initial segment of a shortest path inside the piece and vary its starting point until the opening inequality holds on both sides of the cut.

The construction works for any probability density that is positive almost everywhere. To see why the choice of density matters, put $p_i=|Q_i|/|P|$ and $p_i r_i=\int_{Q_i}u^q$. At $q=2$ area weights alone give $\sum_i p_i r_i=1$. With the mixed weights above, Proposition~\ref{prop:geometric-reduction} gives
\[
 \sum_i\nu(Q_i)r_i^{2/q}\ge 1+\frac1q\sum_i p_i(r_i-1)^2 .
\]
At $q=2$ this is an identity. In the stability proof, a small deficit forces the sum on the right to be small. Thus a cell carrying a large fraction of the area must carry nearly the same fraction of the $L^q$ mass.

\emph{Angular convexity.} The supporting line exists because the triangle profile is strictly convex.
\begin{theoremA}
For every $1\le q\le2$ the profile $F_q$ is twice differentiable on $(0,\pi)$ and $F_q''>0$ there.
\end{theoremA}
A dual variational formula for torsion in terms of fluxes gives concavity of $\alpha^2/F_1(\alpha)$
(Lemma~\ref{lem:torsion-curvature}). For the first eigenvalue, the polygonal argument
uses the tangent at $\pi/2$. For $a=2\pi/n\le\pi/2$, convexity on $(0,\pi/2]$ shows that the tangent at $a$ has slope at most $F_2'(\pi/2)$ and value at $\pi/2$ at most $F_2(\pi/2)$. It therefore lies below the tangent at $\pi/2$ throughout $[\pi/2,\pi)$. Strict convexity gives strict support at acute openings. At obtuse openings the same conclusion follows from the ground state and the spectral gap at the opening $\pi/2$, giving \eqref{eq:T90}.

To prove Theorem~A we express $F_q''$ in terms of the minimizer at the given opening. The term we need to bound measures how the minimizer changes as the triangle opens. For $q\le q_\ast=5801/4000$, comparison with a disk covers every opening (Proposition~\ref{S:prop:lowq}). For larger $q$ we divide the openings into three ranges.

For $0<\alpha\le\pi/2$, a conformal map from a circular sector to the triangle gives coordinates on a half-cylinder, with one periodic angular variable. With $s=\alpha/(2\pi)$ and the map normalized to have derivative one at the sector vertex, the image triangle has area $\pi sZ(s)$. The positive coefficients of the odd powers $s^3,s^5,\ldots$ in $\log Z(s)$ give a positive contribution to the curvature. We bound the change in the minimizer by estimating each angular Fourier mode separately (Section~\ref{sec:acute-real}).

For $\pi/2\le\alpha\le2\pi/3$, put $z=\cot^2(\alpha/2)$ and $\Lambda_q(z)=\sqrt z\,F_q(\alpha)$. Starting at $\alpha=\pi/2$, we integrate a differential inequality for the logarithmic derivative $\beta=d\log\Lambda_q/d\log z$. For the higher exponents, we obtain this inequality using the first two eigenvalues of the linear problem at $q=2$ (Section~\ref{S:sec:middle}).
For $2\pi/3\le\alpha<\pi$, we compare the horizontal energy on each slice of the triangle with the corresponding interval problem. This suffices for $q\le27/16$. For larger $q$ we also bound these two eigenvalues on the wide-opening interval (Section~\ref{S:sec:angle-replacement}).

At sufficiently small openings the expansion on the half-cylinder also shows that $F_q''/F_q$ strictly decreases with $q$ (p.~\pageref{par:small-opening-order}). A separate comparison holds at every opening: for $1<q\le2$, fix the positive minimizer $u$ of unit $L^q$ norm at one opening. As the opening varies, the first eigenvalue with this fixed weight $u^{q-2}$ bounds $F_q$ from below, with equality at the chosen opening (Proposition~\ref{prop:frozen-linear}).

\subsection{Relation to earlier work}
Solynin and Zalgaller used proportional triangular coverings in their proof for logarithmic capacity \cite[Theorem~2]{sz2004}. Bogosel, Bucur and Fragal\`a introduced the weighted side cells used in Section~\ref{sec:convex-faces} in their September paper \cite[Section~2.1]{bbf2026}. Their Theorem~7 proves the torsion assertion of Proposition~\ref{prop:mixed-triangles} for triangles of fixed area and apex angle $0<\alpha\le\pi/2$. Their Theorem~8 proves that $g(\alpha)=\alpha^2/F_1(\alpha)$ has $g''(\alpha)<0$ for every $0<\alpha<\pi$. Lemma~\ref{lem:torsion-curvature} gives an explicit lower bound for $-g''$.

The first version of the present paper gave the cell comparison by rearrangement on a conformal strip for every opening $0<\alpha<\pi$ and every $1\le q<\infty$ under the base-angle hypothesis of Lemma~\ref{lem:mixed-cell} \cite[Lemma~2.1 and Section~2.2]{wangv1}. The present paper uses $1\le q\le2$. The first version also stated both triangle comparisons of Proposition~\ref{prop:mixed-triangles} for every such opening \cite[Proposition~2.3]{wangv1}.

The first version, submitted on 17 November 2026, proved the polygonal torsion
and first Dirichlet eigenvalue inequalities for $n\ge5$, including their equality
cases \cite[Theorem~1.1]{wangv1}. It also proved the comparison of Sobolev
quotients for $1\le q\le2$ and quadratic stability at the two endpoints
\cite[Corollary~1.4 and Section~5]{wangv1}. Section~2 of that version
prescribed cell areas in a convex polygon. Sections~3--4 estimated the torsion
deficit from the signed exterior turns and used the Kohler--Jobin inequality
for the nonconvex comparison. Here Section~\ref{sec:convex-faces} gives an interior partition for the polygonal comparison in both cases. A revised form of Sections~3--4 appears in Part~\ref{part:convex-vertices}, where only convex vertices are counted. Theorem~S1.1
of the first version proved strict convexity through $2\pi/5$ and global
support by the tangent at each opening in that range. At $q=2$ it proved the same support property up to $\pi/2$. Theorem~A here gives strictly positive curvature for every
opening and every exponent. Bogosel, Bucur and Fragal\`a submitted an independent proof of the torsion inequality on 18 November \cite{bbf-global}, and Li submitted a proof of the first eigenvalue inequality, together with the optimal linear correction for fewer sides, on 20 November \cite{li2026}.

Other results treat particular classes of competitors. For torsion, Fragal\`a, Gazzola and Lamboley proved the inequality for convex polygons whose perimeter reaches a threshold relative to the equal-area regular polygon \cite{fgl2013}. Gui, Hu and Li proved a quantitative inequality when every side is tangent to a common incircle \cite{ghl2026}. Bogosel gave a validated Hessian computation for $5\le n\le25$ \cite{bogosel2026}. For the first eigenvalue, Indrei proved sharp quantitative stability for triangles and optimality on an explicit class of competitors for sufficiently large side counts \cite{indrei2024}. Bogosel and Bucur reduced the problem to finitely many certified computations \cite[Theorem~7.6]{bb2024} and proved local minimality for the regular pentagon and hexagon by validated computation \cite{bb2026}. Cheng, Gui, Hu and Li proved the inequality for convex polygons with sufficiently many sides \cite{cghl2026}. Dahne, G\'omez-Serrano and Pech-Alberich proved that $\lambda_1(R_n)$ decreases in $n$ at fixed area \cite{dgp2026}; for $n\ge3$ this also follows from the equality case of Theorem~\ref{thm:main}, since $R_n$ is a competitor with at most $n+1$ sides. Bucur and Fragal\`a proved Brunn--Minkowski inequalities for mixed energies of convex sectors \cite{bfmixed2026}.

\paragraph{Reading the proofs for torsion and the first eigenvalue.}
Lemma~\ref{lem:mixed-cell} compares each cell with a triangle.
Lemma~\ref{lem:mass-cap-partition} supplies the partition. Given a supporting
tangent at the regular opening, Proposition~\ref{prop:geometric-reduction}
sums the cell inequalities; reflection gives the sharp constant in
\eqref{eq:regular-reflection}. For triangles and quadrilaterals the direct
Steiner argument on p.~\pageref{par:small-polygons} also gives uniqueness.
For $n\ge5$, Sections~\ref{sec:convex}--\ref{sec:convex-faces} give the common geometric argument. To complete the torsion proof, read Section~\ref{sec:torsion-concavity}. For the first eigenvalue, read Sections~\ref{ang:sec:q2-start}--\ref{ang:sec:square-tangent}. Both proofs are completed within Part~\ref{part:endpoints}.
\begin{center}\small
\begin{tabularx}{\linewidth}{@{}p{.24\linewidth}X@{}}
\toprule
Conclusion & Estimates used after the polygon is partitioned and its cells compared \\
\midrule
Torsion, $q=1$ & The flux proof in Section~\ref{sec:torsion-concavity} gives strict angular convexity.\\[3pt]
First eigenvalue, $q=2$ & The estimates in Sections~\ref{ang:sec:q2-start}--\ref{ang:sec:spectral-curvature}, which start from the disk and from the opening $\pi/2$, give strict convexity at acute openings. Inequality \eqref{eq:T90} places the tangent at $\alpha=\pi/2$ below the profile at obtuse openings.\\[3pt]
Sobolev quotients, $1<q<2$ & Theorem~A gives the tangent at $2\pi/n$ for every $n$. Table~\ref{tab:acute-coverage} lists the estimates for acute openings.\\[3pt]
Quadratic stability & Section~\ref{sec:stability} estimates the tangent remainder, the difference between the two base angles of each cell and the resulting Hausdorff distance. The explicit constants for $q=1$ and $q=2$ use \eqref{eq:spectral-quadratic-gap} and Section~\ref{sec:torsion-concavity}. For $n=3,4$ see also p.~\pageref{par:small-n-area}.\\
\bottomrule
\end{tabularx}\end{center}
The torsion and eigenvalue cases of Theorem~\ref{thm:main} are proved in Part~\ref{part:endpoints} alone; Part~\ref{part:intermediate} is needed only for $1<q<2$ and for some of the stability constants. The stability estimate for all exponents in Section~\ref{S:sec:sob-stability} also uses the uniform angular estimates from Part~\ref{part:intermediate}. Appendix~\ref{app:estimate-register}
records the numerical margins used in the proofs.

\section{Comparison of a mixed cell with a triangle}\label{sec:convex}
Let $I_\alpha$ be the isosceles triangle of height one and apex
angle $\alpha$, with Dirichlet condition on its base and Neumann boundary
conditions on its equal sides. The Dirichlet base is a full side of
the comparison cell. For the weighted partition of a convex polygon, this base is a whole edge of the original polygon.

Let \(U\) be a polygon with one side marked as its Dirichlet side, and
let \(H^1_D(U)\) be the subspace of \(H^1(U)\)
with zero trace on that side. The mixed torsional rigidity and the
first mixed eigenvalue of \(U\) are
\[
 T_{\mathrm{DN}}(U)=\max_{0\ne v\in H^1_D(U)}
   \frac{\bigl(\int_Uv\bigr)^2}{\int_U|\nabla v|^2},
 \qquad
 \lambda_{\mathrm{DN}}(U)=\min_{0\ne v\in H^1_D(U)}
   \frac{\int_U|\nabla v|^2}{\int_Uv^2}.
\]
Put
\[
 \mathcal M_q(U)=
 \sup_{\substack{v\in H^1_D(U)\\\int_U|\nabla v|^2=1}}
       \frac1{|U|}\int_U|v|^q .
\]
Thus \(\mathcal M_q(U)\) is the largest mean of \(|v|^q\) obtainable
with one unit of Dirichlet energy, and
\[
 \mathcal M_1(U)^2=\frac{T_{\mathrm{DN}}(U)}{|U|^2},\qquad
 \mathcal M_2(U)=\frac1{|U|\,\lambda_{\mathrm{DN}}(U)} .
\]
The number \(\mathcal M_q(U)\) does not change when \(U\) is scaled.

\begin{lemma}[The cell comparison]\label{lem:mixed-cell}
Let \(Q\) be a bounded convex polygon with Dirichlet data on one
side and Neumann data elsewhere. If the angles at the ends of the
Dirichlet side are \(\pi a,\pi b\) with \(a+b<1\), set
\(\alpha=\pi(1-a-b)\). Then
\[
 \mathcal M_q(Q)\le\mathcal M_q(I_\alpha)\qquad(1\le q\le 2).
\]
Equality for any \(q\) holds exactly when \(Q\) is similar to
\(I_\alpha\), with corresponding Dirichlet side.
\end{lemma}

\subsection{Comparing a cell with an isosceles triangle}\label{sec:strip}

\begin{proof}[Proof of Lemma~\ref{lem:mixed-cell}]
Map the cell and its comparison triangle to the same horizontal strip.
In these coordinates a test function keeps its energy. What
differs is the area assigned to each part of the strip. We show that
at fixed energy the largest weighted integral of \(u^q\) occurs for
the area density of the isosceles triangle. The proof has four steps.
Step~1 passes to the strip, where the area of \(Q\) has density
\(p\). Step~2 writes \(p\) as an average of translates of members of a
two-parameter family \(p_{\delta,\kappa}\), in which \(p_{0,1}\)
is the density of \(I_\alpha\). Step~3 rearranges the test function. Step~4
proves the inequality \eqref{eq:strip-goal} for each
\(p_{\delta,\kappa}\) from three monotonicity statements.

\emph{Step 1: area in strip coordinates.}
Choose a conformal map from the upper half-plane onto \(Q\) that
sends the positive real axis to the Dirichlet side. Map \(0\) and \(\infty\) to the vertices with angles \(\pi a\)
and \(\pi b\), respectively, and put \(C=a+b\).
The points mapping to the other vertices
lie on the negative real axis. Write them as \(-e^{x_j}\). By the Schwarz--Christoffel formula
\cite[Theorem~2.1]{dt2002}, the derivative in the half-plane is a constant
times \(z^{a-1}\prod_j(z+e^{x_j})^{-c_j}\). Use \(z=e^{x+iy}\) and
\(|z+e^{x_j}|^2=2e^{x+x_j}(\cosh(x-x_j)+\cos y)\). This gives the
area density on \(S=\mathbb R\times(0,\pi)\), up to a constant:
\[
 J=e^{\delta x}\prod_j(\cosh(x-x_j)+\cos y)^{-c_j},
 \qquad c_j>0,\quad \sum_j c_j=C,\quad\delta=a-b.
\]
Here \(\pi c_j\) is the exterior angle at the \(j\)-th of the other
vertices, after removal of flat vertices. The exterior angles of \(Q\)
sum to \(2\pi\), which gives \(\sum_jc_j=a+b\). Let \(\Phi:S\to Q\) be the resulting map. After normalization,
\(p=J/\int_S J=|\Phi'|^2/|Q|\) is the fraction of cell area per
unit strip area. Thus a test \(w\) on \(Q\) and \(u=w\circ\Phi\)
satisfy
\[
 \frac{1}{|Q|}\int_Q |w|^q=\int_S p|u|^q,\qquad
 \int_Q|\nabla w|^2=\int_S|\nabla u|^2 .
\]
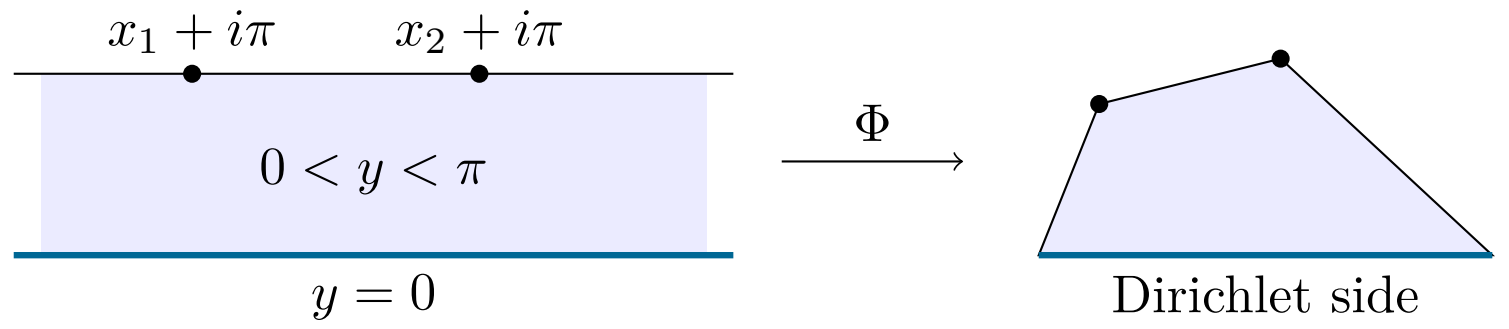
\begin{figure}[htbp]
\centering
\begin{tikzpicture}[font=\small]

\fill[blue!8] (0,0) rectangle (4.4,1.2);
\draw[very thick,MidnightBlue] (-.18,0)--(4.58,0);
\draw (-.18,1.2)--(4.58,1.2);
\foreach \x/\n in {1/1,2.9/2} {\fill (\x,1.2) circle (1.7pt) node[above]{\(x_{\n}+i\pi\)};}
\node at (2.2,.55){\(0<y<\pi\)};
\node[below] at (2.2,0){\(y=0\)};
\draw[->] (4.9,.62)--node[above]{\(\Phi\)}(6.1,.62);
\filldraw[fill=blue!8] (6.6,0)--(7,1.0)--(8.2,1.3)--(9.6,0)--cycle;
\draw[very thick,MidnightBlue] (6.6,0)--(9.6,0);
\fill (7,1) circle(1.7pt);\fill (8.2,1.3) circle(1.7pt);
\node[below] at (8.1,0){Dirichlet side};
\end{tikzpicture}
\caption{The strip extends indefinitely to the left and right. Its lower edge maps to the base of the cell, and marked points on its upper edge map to the other corners, in boundary order.}\label{fig:strip-boundary}
\end{figure}

The base corresponds to \(y=0\), so \(u\) has zero trace there.
To check that the change of variables gives the same \(H^1\) test
space in both directions, first integrate the vertical derivative
from that edge. The Cauchy--Schwarz inequality and integration in \(x,y\) give
\[
 |u(x,y)|^2\le y\int_0^y|\partial_y u(x,s)|^2\,ds,\qquad
 \|u\|_{L^2(S)}^2\le \frac{\pi^2}{2}\|\partial_yu\|_{L^2(S)}^2 .
\]
Conversely, \(p\in L^r(S)\) for some \(r>1\). The tails decay
exponentially since \(|\delta|<C\), and at \((x_j,\pi)\) the density
is bounded by a constant times the distance to that point raised to the power \(-2c_j\), with \(c_j<1\).
With \(r'=r/(r-1)\), Hölder's inequality and the Sobolev
embedding \(H^1(S)\subset L^{2r'}(S)\) yield
\[
 \frac1{|Q|}\int_Q |u\circ\Phi^{-1}|^2
 =\int_S p|u|^2
 \le \|p\|_{L^r(S)}\,\|u\|_{L^{2r'}(S)}^2
 \le C_r\|p\|_{L^r(S)}\,\|u\|_{H^1(S)}^2 .
\]
The change of variables preserves zero trace locally on the open base. Since its two endpoints have zero boundary length, the trace vanishes almost everywhere on the whole base. Hence
\[
 \mathcal M_q(Q)=M_q(p):=
 \sup_{\substack{u\in H^1(S),\ u\ge0,\ u|_{y=0}=0\\
                  \int_S|\nabla u|^2=1}}\int_S p\,u^q .
\]
For example, choose \(r\) close enough to one that \(qr'\ge2\).
For every unit-energy test, the same estimates give the finite bound
\(\int_S p\,u^q\le\|p\|_{L^r(S)}\|u\|_{L^{qr'}(S)}^q
\le C_{q,r}\|p\|_{L^r(S)}\). To obtain the sharp comparison we
must still use the shape of \(p\). The Poincaré inequality and
compact embedding on the bounded polygon give a maximizer.

\emph{Step 2: reduction to a two-parameter family.}
Let \(p_{\delta,\kappa}\) be the density proportional to
\(e^{\delta x}(\cosh x+\kappa\cos y)^{-C}\), normalized to have
integral one. Its parameters have geometric meanings. The factor
\(e^{\delta x}\) records the difference of the two base angles.
For a triangle there is only one other vertex, so a translation
in \(x\) moves it to \(x_1=0\) and gives \(\kappa=1\). In particular \(p_{0,1}\) is
the density of \(I_\alpha\).

For a general cell, write \(d_j=\cosh(x-x_j)+\cos y\). The product in
\(J\) is an average of powers of single denominators:
\[
 \prod_j d_j^{-c_j}
 =\int_{\Delta}\Big(\sum_j t_jd_j\Big)^{-C}d\gamma(t),
 \quad
 d\gamma(t)=\frac{\Gamma(C)}{\prod_j\Gamma(c_j)}
                  \prod_j t_j^{c_j-1}\,dt ,
\]
where \(\Delta=\{t_j\ge0:\sum_j t_j=1\}\) and \(dt\) means integration
in all but the last coordinate. The identity follows by writing each
factor as \(d_j^{-c_j}=\Gamma(c_j)^{-1}\int_0^\infty
\tau_j^{c_j-1}e^{-\tau_jd_j}\,d\tau_j\) and using the sum of the
integration variables as the radial variable. For a single factor take a point
mass. Write
\[
 B_\pm=\sum_jt_je^{\pm x_j},\qquad
 \bar x=\tfrac12\log(B_+/B_-),\qquad K=\sqrt{B_+B_-}.
\]
Then \(\sum_jt_j\cosh(x-x_j)=\tfrac12(e^{-x}B_++e^{x}B_-)
=K\cosh(x-\bar x)\), and the
Cauchy--Schwarz inequality gives \(K\ge\sum_jt_j=1\). Hence
\[
 \sum_jt_jd_j=K\bigl(\cosh(x-\bar x)+\kappa\cos y\bigr),
 \qquad \kappa=\kappa(t)=1/K\in(0,1].
\]
Thus each term, once normalized, is a horizontal translate of
\(p_{\delta,\kappa(t)}\) by \(\bar x(t)\).
If there are two or more distinct \(x_j\), then \(K>1\) throughout
the interior of \(\Delta\). Absorb the positive integral of each
term into the averaging measure and normalize it to a probability
measure \(\mu\) on \(\Delta\).

For any unit-energy test \(u\), put \(u_t(x,y)=u(x+\bar x(t),y)\).
Translation preserves the bottom trace and the energy. Consequently
\[
 \int_S p\,u^q
 =\int_\Delta\!\left(\int_S p_{\delta,\kappa(t)}\,u_t^q\right)d\mu(t).
\]
Since \(\mu\) is a probability measure, the lemma follows once we show
that
\begin{equation}\label{eq:strip-goal}
 \int_S p_{\delta,\kappa}\,u^q\ \le\ M_q(p_{0,1})
\end{equation}
for all \(|\delta|<C\), all \(0<\kappa\le1\), and every nonnegative
\(u\in H^1(S)\) with zero trace on \(y=0\) and
\(\int_S|\nabla u|^2=1\). We prove \eqref{eq:strip-goal} in the form
\begin{equation}\label{eq:strip-chain}
 \underbrace{\int_S p_{\delta,\kappa}u^q
       \ \le\ \int_S p_{0,\kappa}v^q}_{\text{equalize the base angles}}
       \ \le\
 \underbrace{\int_S p_{0,1}v^q}_{\text{the isosceles triangle}}
 \ \le\ M_q(p_{0,1}),
\end{equation}
where \(v\) is a rearrangement of \(u\) with no larger energy.

\emph{Step 3: rearranging the test function.}
Rearrange \(u\) twice. First, for each \(x\), replace the function
\(u(x,\cdot)\) on \((0,\pi)\) by the nondecreasing function
\(u^\sharp(x,\cdot)\) with the same distribution of values:
\(|\{y:u^\sharp(x,y)>t\}|=|\{y:u(x,y)>t\}|\) for every \(t\). Second,
for each \(y\), replace \(u^\sharp(\cdot,y)\) by the even function
\(v(\cdot,y)\) that is nonincreasing on \((0,\infty)\) and satisfies
\(|\{x:v(x,y)>t\}|=|\{x:u^\sharp(x,y)>t\}|\) for every \(t\). Then
\begin{equation}\label{eq:rearranged-test}
 v \text{ is nondecreasing in } y \text{ and nonincreasing in } |x|,
 \qquad v|_{y=0}=0,\qquad
 \int_S|\nabla v|^2\le\int_S|\nabla u|^2 .
\end{equation}
The energy inequality for the second step is the P\'olya--Szeg\H{o}
inequality for Steiner symmetrization, and for the first step it is
its analogue for the monotone rearrangement
\cite[Ch.~II]{kawohl1985}. The first step can also be reduced to the
second: reflect \(u\) evenly across \(y=\pi\), extend by zero, and
symmetrize about the line \(y=\pi\). In both steps the derivative in
the rearranged direction does not increase in \(L^2\), and the other
derivative is controlled by
\(\|f^*-g^*\|_2\le\|f-g\|_2\) \cite[Theorem~3.5]{ll2001}, applied
to neighboring slices and then to their difference quotients. The vertical rearrangement has zero trace at the bottom. Horizontal rearrangement preserves this trace because it preserves the \(L^2\) norm on each horizontal line. It also preserves the pointwise order of the horizontal slices, so the result remains nondecreasing in \(y\).

We use the rearrangement inequality of Hardy and Littlewood
\cite[Theorem~3.4]{ll2001} in two forms. If \(f\ge0\) and
\(\rho\ge0\) on \(\mathbb R\), then
\(\int_{\mathbb R}\rho f\le\int_{\mathbb R}\rho^*f^*\), where the star
denotes the even rearrangement that is nonincreasing on
\((0,\infty)\). If \(f\ge0\) and \(\rho\ge0\) on \((0,\pi)\) and
\(\rho\) is nondecreasing, then
\(\int_0^\pi\rho f\le\int_0^\pi\rho f^\sharp\); this follows from the
first form after even reflection across \(y=\pi\).

\emph{Step 4a: the vertical step and an expansion of the density.}
On a vertical line the area density increases towards the top:
\[
 \partial_y\log p_{\delta,\kappa}(x,y)
 =\frac{C\kappa\sin y}{\cosh x+\kappa\cos y}\ge0
 \qquad(0<y<\pi).
\]
By the second form of the Hardy--Littlewood inequality, applied on
each vertical line,
\(\int_Sp_{\delta,\kappa}u^q\le\int_Sp_{\delta,\kappa}(u^\sharp)^q\).
For the horizontal step, expand the density as a mixture of products
of probability densities:
\begin{equation}\label{eq:mixture}
 p_{\delta,\kappa}=\sum_{k\ge0}w_{k,\delta,\kappa}\rho_{C+k,\delta}(x)\eta_{k,\kappa}(y),
 \qquad
 \rho_{b,\delta}\propto e^{\delta x}(\cosh x+1)^{-b},
 \qquad
 \eta_{k,\kappa}\propto(1-\kappa\cos y)^k .
\end{equation}
This is the binomial series
\((1-\xi)^{-C}=\sum_k\binom{C+k-1}{k}\xi^k\) with
\(\xi=(1-\kappa\cos y)/(\cosh x+1)\in[0,1)\), applied to
\(\cosh x+\kappa\cos y=(\cosh x+1)(1-\xi)\). Its terms are
nonnegative, and the weights are
\begin{equation}\label{eq:mixture-weights}
 w_{k,\delta,\kappa}=\frac1{N_{\delta,\kappa}}\binom{C+k-1}{k}\,
 2^{C+k}\mathrm B(C+k+\delta,C+k-\delta)\,Z_k(\kappa),
 \qquad Z_k(\kappa)=\int_0^\pi(1-\kappa\cos y)^k\,dy ,
\end{equation}
where \(2^b\mathrm B(b+\delta,b-\delta)
=\int_{\mathbb R}e^{\delta x}(\cosh x+1)^{-b}dx\) and
\(N_{\delta,\kappa}\) makes the weights sum to one.
The first form of the Hardy--Littlewood inequality, applied in \(x\)
for each \(y\) and each \(k\), gives
\begin{equation}\label{eq:after-HL}
 \int_Sp_{\delta,\kappa}u^q\ \le\
 \sum_{k\ge0}w_{k,\delta,\kappa}
 \int_S\rho_{C+k,\delta}^*(x)\,\eta_{k,\kappa}(y)\,v^q\,dx\,dy .
\end{equation}

\emph{Step 4b: the factor \(e^{\delta x}\) spreads the density out.}
We claim that for \(b>|\delta|\) and every \(s>0\)
\begin{equation}\label{eq:centered-intervals}
 \int_{-s}^s \rho_{b,\delta}^*(x)\,dx
 \le \int_{-s}^s \rho_{b,0}(x)\,dx .
\end{equation}
The left side is the largest \(\rho_{b,\delta}\)-mass of an interval
of length \(2s\). Indeed
\((\log\rho_{b,\delta})'(x)=\delta-b\tanh(x/2)\) decreases and changes
sign, so every set \(\{\rho_{b,\delta}>c\}\) is an interval
\((m-\ell,m+\ell)\), and \(\rho^*_{b,\delta}(\ell)=c\). Let
\[
 L=(\log\rho_{b,\delta})'(m-\ell)>0,\qquad
 R=-(\log\rho_{b,\delta})'(m+\ell)>0 .
\]
As \(\log c\) rises, the left and right endpoints move at the rates
\(1/L\) and \(-1/R\). The half-length \(\ell\) therefore changes at
the rate \(-\tfrac12(1/L+1/R)\), and
\(-\frac{d}{d\ell}\log\rho^*_{b,\delta}(\ell)=2LR/(L+R)\).
With \(r=\tanh(\ell/2)\) and \(z=\tanh(m/2)\),
\[
 \frac{2LR}{L+R}\le\frac{L+R}{2}
  =\frac b2\Bigl(\tanh\frac{m+\ell}{2}-\tanh\frac{m-\ell}{2}\Bigr)
  =br\frac{1-z^2}{1-z^2r^2}
  \le br=-\frac{d}{d\ell}\log\rho_{b,0}(\ell).
\]
Thus \(\rho_{b,\delta}^*/\rho_{b,0}\) is nondecreasing on
\((0,\infty)\). Both densities have integral \(1/2\) over
\((0,\infty)\), so their difference changes sign at
most once there, from negative to positive. This proves
\eqref{eq:centered-intervals}.
If \(\delta\ne0\), then \(m\ne0\) for every \(\ell>0\), and the
inequality is strict for every \(0<s<\infty\). For fixed \(y\)
the sets \(\{x:v(x,y)^q>t\}\) are centered intervals up to sets of
measure zero. Writing \(v^q=\int_0^\infty1_{\{v^q>t\}}\,dt\) and
using \eqref{eq:centered-intervals} for each \(t\) gives
\[
 \int_{\mathbb R}\rho_{b,\delta}^*(x)\,v(x,y)^q\,dx
 \le\int_{\mathbb R}\rho_{b,0}(x)\,v(x,y)^q\,dx ,
\]
strictly if \(\delta\ne0\) and \(v(\cdot,y)\) is not identically zero.
With the notation
\[
 a_k(\kappa)=\int_S \rho_{C+k,0}(x)\eta_{k,\kappa}(y)v^q\,dx\,dy
\]
we obtain from \eqref{eq:after-HL}
\begin{equation}\label{eq:after-centering}
 \int_Sp_{\delta,\kappa}u^q\ \le\ \sum_{k\ge0}w_{k,\delta,\kappa}\,a_k(\kappa),
\end{equation}
with strict inequality if \(\delta\ne0\) and \(u\ne0\).

\emph{Step 4c: three monotonicity statements.}
The three comparisons below use one observation about averages.
Let \(d\mu\) be a probability measure on an ordered set. If
\(\varphi\ge0\) and \(f\) are monotone in the same direction, then
\[
 \int f\varphi\,d\mu-\int f\,d\mu\int\varphi\,d\mu
 =\tfrac12\iint (f(s)-f(t))(\varphi(s)-\varphi(t))\,d\mu(s)d\mu(t)
 \ge0.
\]
For bounded functions this follows by expanding the product;
truncation gives the same inequality when the three means are finite.
We use it in two forms. (W1) Multiplying a density by \(\varphi\) and
normalizing gives a mean of \(f\) at least as large as before.
(W2) For positive probability weights \(w_k,w'_k\), if \(w'_k/w_k\)
and \(a_k\ge0\) are both nondecreasing, then
\(\sum_kw_ka_k\le\sum_kw'_ka_k\). Take \(\mu(\{k\})=w_k\),
\(\varphi(k)=w'_k/w_k\), and \(f(k)=a_k\) in (W1).

\emph{(M1) \(a_k(\kappa)\le a_{k+1}(\kappa)\).} Raising \(k\) by one
multiplies the two densities, up to normalization, by
\(1/(\cosh x+1)\) and \(1-\kappa\cos y\), respectively. The first
factor decreases in \(|x|\) and the second increases in \(y\), while
\(v^q\) is nonincreasing in \(|x|\) and nondecreasing in \(y\). Write
\[
 a_{k+1}-a_k=\int_S(\rho_{C+k+1,0}-\rho_{C+k,0})\,\eta_{k+1,\kappa}\,v^q
 +\int_S\rho_{C+k,0}\,(\eta_{k+1,\kappa}-\eta_{k,\kappa})\,v^q .
\]
By (W1), applied on each horizontal half-line \(x>0\) in the first
integral and on each vertical line in the second, both integrals are
nonnegative.

\emph{(M2) \(a_k(\kappa)\le a_k(1)\) for \(0<\kappa\le1\).} Let
\(X=-\cos y\). The ratio
\(\eta_{k,1}/\eta_{k,\kappa}\propto((1+X)/(1+\kappa X))^k\) is
nondecreasing in \(y\), and \(v^q\) is nondecreasing in \(y\). The
claim follows from (W1) on each vertical line. The
inequality is strict for \(k=1\), \(\kappa<1\) and \(v\ne0\), because
then the ratio is strictly increasing, and the nonzero vertical slices
of \(v\) have zero bottom trace and are not constant.

\emph{(M3) The ratios \(w_{k,0,\kappa}/w_{k,\delta,\kappa}\) and
\(w_{k,0,1}/w_{k,0,\kappa}\) are nondecreasing in \(k\).} By
\eqref{eq:mixture-weights} and
\(\mathrm B(x+1,y+1)=\mathrm B(x,y)\,xy/((x+y)(x+y+1))\), passing from
\(k\) to \(k+1\) multiplies \(w_{k,\delta,\kappa}/w_{k,0,\kappa}\) by
\(1-\delta^2/(C+k)^2\le1\). For the second ratio,
\[
 \frac{Z_{k+1}(\kappa)}{Z_k(\kappa)}=1+\kappa\,\mathbb E_{k,\kappa}X,
\]
where the expectation is for the density \(\eta_{k,\kappa}\).
The expectation is nonnegative by (W1), because \(X\) has mean zero
under the uniform density on \((0,\pi)\) and both \(X\) and
\((1+\kappa X)^k\) increase with \(y\). It is nondecreasing in
\(\kappa\) by (W1) and the comparison of densities used in (M2). Hence
\(Z_{k+1}(1)/Z_k(1)\ge Z_{k+1}(\kappa)/Z_k(\kappa)\), which is the
claim.

\emph{Step 4d: conclusion.}
Start with \eqref{eq:after-centering}. To change the weights in the first and last sums below, apply (W2): the weight ratios increase by (M3), and the terms being averaged increase by (M1). Between these changes use (M2) to increase \(a_k(\kappa)\) to \(a_k(1)\). This gives
\[
 \int_S p_{\delta,\kappa}u^q
 \le \sum_{k\ge0}w_{k,\delta,\kappa}a_k(\kappa)
 \le \sum_{k\ge0}w_{k,0,\kappa}a_k(\kappa)
 \le \sum_{k\ge0}w_{k,0,\kappa}a_k(1)
 \le \sum_{k\ge0}w_{k,0,1}a_k(1).
\]
The comparisons of nonnegative increasing sequences follow first for their bounded truncations and then by monotone convergence.
By \eqref{eq:mixture} for \(\delta=0\), the second sum is
\(\int_Sp_{0,\kappa}v^q\) and the last is \(\int_Sp_{0,1}v^q\). This
gives the first two inequalities in \eqref{eq:strip-chain}. Finally,
by the definition of \(M_q\) and \eqref{eq:rearranged-test},
\[
 \int_S p_{0,1}v^q
 \le M_q(p_{0,1})\Big(\int_S|\nabla v|^2\Big)^{q/2}
  \le M_q(p_{0,1}).
\]
This proves \eqref{eq:strip-goal}. Applying it to each translated test
\(u_t\) in Step~2 and taking the supremum over \(u\) gives
\(\mathcal M_q(Q)=M_q(p)\le M_q(p_{0,1})=\mathcal M_q(I_\alpha)\).

\emph{Equality.}
Take a maximizer on \(Q\). Its strip test \(u\) is nonzero, and so is
every translate \(u_t\). The first inequality in
\eqref{eq:strip-chain} is strict if \(\delta\ne0\), by Step~4b. The
second is strict if \(\kappa<1\), by (M2) with \(k=1\). With more than one
of the points \(-e^{x_j}\), \(\kappa(t)<1\) on the interior of the simplex,
which has full \(\mu\)-measure. Equality therefore requires
\(\delta=0\) and a single such point, that is, a triangle with
equal base angles. These conditions characterize
\(I_\alpha\) up to similarity, and this triangle gives equality.
\end{proof}

\begin{proposition}[Triangles of prescribed area and apex angle]\label{prop:mixed-triangles}
Fix \(0<\alpha<\pi\) and an area. Among triangles with apex angle
\(\alpha\), with Dirichlet condition on the opposite side and Neumann
conditions on the two sides meeting at the apex, the isosceles
triangle uniquely maximizes mixed torsional rigidity and uniquely
minimizes the first mixed eigenvalue, up to congruence.
\end{proposition}

\begin{proof}
The two angles at the ends of the Dirichlet side add to
\(\pi-\alpha\), so Lemma~\ref{lem:mixed-cell} applies with
\(a+b=1-\alpha/\pi\). The areas are equal and the equality condition
in the lemma gives congruence.
\end{proof}

The torsion assertion for \(0<\alpha\le\pi/2\) is \cite[Theorem~7]{bbf2026}.
For all \(0<\alpha<\pi\), both assertions appeared in \cite[Proposition~2.3]{wangv1}. The September proof uses monotonicity of the mixed torsion function in the direction of the inward
normal to the base for \(\alpha\le\pi/2\), from \cite[Theorem~1.3]{ly2026}.
Here we map the domains to a strip, rearrange the test functions, and compare the area densities of the domains.
The strip argument requires only \(|a-b|<a+b<1\).

\section{Reduction to the triangle profiles}\label{sec:sobolev-consequence}

Write $F_q(\alpha)=\mathcal M_q(I_\alpha)^{-2/q}$,
$F=F_2$ and $g(\alpha)=\alpha^2/F_1(\alpha)$.
Thus \(F_2(\alpha)=|I_\alpha|\lambda_{\rm DN}(I_\alpha)\) and
\(F_1(\alpha)=|I_\alpha|^2/T_{\rm DN}(I_\alpha)\). For every \(q\),
\begin{equation}\label{eq:Fq-quotient}
 F_q(\alpha)=|I_\alpha|^{2/q}\inf_{0\ne v\in H^1_D(I_\alpha)}
 \frac{\int_{I_\alpha}|\nabla v|^2}
      {\bigl(\int_{I_\alpha}|v|^q\bigr)^{2/q}} .
\end{equation}
We record four properties of \(F_q\). The first three are elementary.
The supporting-line property will follow from the angular estimates.

\emph{(F1) Symmetric tests suffice.}
Let \(u\ge0\) be a test function on \(I_\alpha\) and let \(u^*\) be its
reflection in the bisector. The symmetric function
\(w=((u^2+(u^*)^2)/2)^{1/2}\) satisfies
\[
 |\nabla w|^2\le (|\nabla u|^2+|\nabla u^*|^2)/2,\qquad
 w^q\ge (u^q+(u^*)^q)/2.
\]
These are the Cauchy--Schwarz inequality and the concavity of
\(t^{q/2}\), respectively. Integration shows that the quotient of \(w\)
in \eqref{eq:Fq-quotient} is at most that of \(u\). A minimizer exists
because \(H^1(I_\alpha)\) embeds compactly in \(L^q(I_\alpha)\).
Hence there is a minimizer that is symmetric in the bisector.

\emph{(F2) Two-sided bounds.}
For \(0<\alpha<\pi\),
\begin{equation}\label{eq:Fq-bounds}
 \frac{\pi^2}{4}\tan\frac{\alpha}{2}
 \le F_q(\alpha)
 \le \left(\frac{(q+1)(q+2)}{2}\right)^{2/q}\tan\frac{\alpha}{2}.
\end{equation}
For the lower bound, let \(v\) vanish on the base. On a segment
perpendicular to the base, of length \(L\le1\), the one-dimensional
inequality \(\int_0^L v^2\le(4L^2/\pi^2)\int_0^L(v')^2\) holds for
functions vanishing at one end. Integrating over these segments gives
\(\int_{I_\alpha}v^2\le(4/\pi^2)\int_{I_\alpha}|\nabla v|^2\). By
H\"older's inequality,
\[
 \frac1{|I_\alpha|}\int_{I_\alpha}|v|^q
 \le\Bigl(\frac1{|I_\alpha|}\int_{I_\alpha}v^2\Bigr)^{q/2}
 \le\Bigl(\frac{4}{\pi^2|I_\alpha|}\int_{I_\alpha}|\nabla v|^2\Bigr)^{q/2},
\]
and \(|I_\alpha|=\tan(\alpha/2)\). For the upper bound, test with the
distance to the base. Its gradient has length one, and its mean
\(q\)-th power is \(2\int_0^1 y^q(1-y)\,dy=2/((q+1)(q+2))\).
In particular \(F_q(0+)=0\), \(F_q(\pi-)=+\infty\), and
\(0<g(\alpha)\le4\alpha^2/(\pi^2\tan(\alpha/2))\to0\) as
\(\alpha\downarrow0\).

\begin{lemma}[The affine family of triangles]\label{lem:affine-family}
Let \(t=\tan(\alpha/2)\), \(z=t^{-2}\), and
\(D=\{(x,y):0<y<x<1\}\). Equip \(D\) with the probability measure
\(d\mu=2\,dx\,dy\); all tests have zero trace on \(x=1\).
Write \(E_x(u)=\int_D u_x^2\,d\mu\), \(E_y(u)=\int_D u_y^2\,d\mu\), and
\(Q_z=E_x+zE_y\). Use the same symbols for their bilinear forms. For \(1\le q\le2\),
\begin{equation}\label{eq:Fq-affine}
 F_q(\alpha)=t\,\Lambda_q(z),\qquad
 \Lambda_q(z)=\inf_{u\ne0}\frac{Q_z(u)}{(\int_D|u|^q\,d\mu)^{2/q}} .
\end{equation}
The function \(\Lambda_q\) is nondecreasing and concave.
\end{lemma}
\begin{proof}
By (F1) it suffices to take symmetric tests. The map
\((x,y)\mapsto(x,ty)\) sends \(D\) to the upper half of the triangle
of height one. Reflection accounts for the factor two in \(d\mu\).
The Jacobian is \(t\), and the transverse derivative is multiplied by
\(1/t\). Substitution in \eqref{eq:Fq-quotient} gives the formula.
Every test gives an affine function of \(z\) with nonnegative slope.
Their infimum is therefore concave and nondecreasing.
\end{proof}
\begin{figure}[htbp]
\centering
\begin{tikzpicture}[font=\small]

\foreach \shift/\height/\top/\name in {0/1.8/{(1,1)}/D,5.7/1.05/{(1,t)}/{I_\alpha^+}}{
 \begin{scope}[shift={(\shift,0)}]
 \fill[blue!8] (0,0)--(2.4,0)--(2.4,\height)--cycle;
 \draw (0,0)node[below]{\(0\)}--(2.4,0)node[below]{\((1,0)\)}--(2.4,\height)node[above]{\(\top\)}--cycle;
 \draw[very thick,MidnightBlue] (2.4,0)--(2.4,\height);
 \node at (1.7,.26){\(\name\)};
 \end{scope}}
\draw[->] (2.85,.55)--node[above]{\((x,y)\mapsto(x,ty)\)}(5.2,.55);
\end{tikzpicture}
\caption{The horizontal coordinate stays fixed. Stretching the vertical coordinate gives the half triangle at opening \(\alpha\). The heavy edge carries the zero trace. Coordinates are labeled in each panel.}\label{fig:affine}
\end{figure}
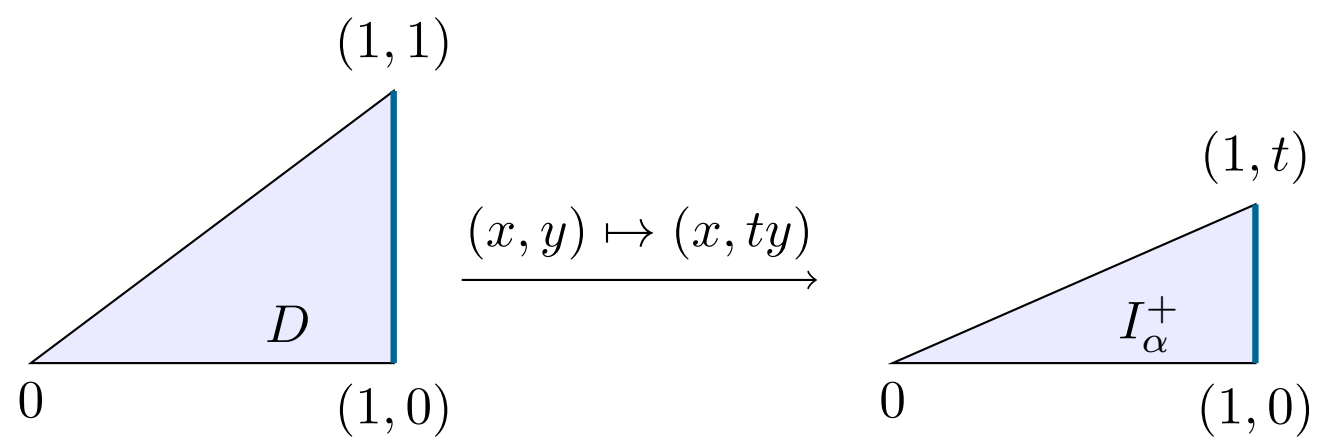

\paragraph{The regular polygon as an equality example.}
On each radial triangle of the regular polygon the two base angles equal $(\pi-a)/2$. Uniqueness of the normalized positive minimizer makes it invariant under rotation by $2\pi/n$. Each radial triangle therefore carries $1/n$ of the area and $1/n$ of its $L^q$ mass. Hence $\nu_i=1/n$ and $\alpha_i-b=(2\pi-nb)\nu_i$, giving equality both in the supporting-line inequality at $a$ and in the bound on each cell opening.

Put $a=2\pi/n$ and $A=|R_n|$. Write $e_i$ for the energies of a test function on the $n$
radial triangles. Its mean $q$-mass on triangle $i$ is bounded by
$(e_i/F_q(a))^{q/2}$. Averaging over $i$ and using concavity of
$x^{q/2}$ bounds its mean on $R_n$ by
$(\sum_i e_i/(nF_q(a)))^{q/2}$. For the reverse bound paste $n$
copies of a bisector-symmetric mixed minimizer. Its traces agree on
each cut, including the last one, and its trace on the polygon boundary vanishes. Hence
\begin{equation}\label{eq:regular-reflection}
 A^{2/q}\lambda_q^{\rm S}(R_n)=nF_q(a).
\end{equation}

\subsection{Supporting tangents at the regular opening}\label{sec:angular-input}

Angles and derivatives in $\alpha$ are in radians. We use $a=2\pi/n$ for the reference angle, $t=\tan(\alpha/2)$ and $z=t^{-2}$ for affine stretching. In the sector formulas $s=\alpha/(2\pi)$, so $d/d\alpha=(2\pi)^{-1}d/ds$. A slope of a supporting line is always a derivative in $\alpha$.

\emph{(F3) Monotonicity.}
By Lemma~\ref{lem:affine-family}, \(F_q(\alpha)/t=\Lambda_q(t^{-2})\)
is nonincreasing in \(\alpha\), and is continuous by finiteness and concavity
of \(\Lambda_q\) on \((0,\infty)\).

\emph{(F4) The supporting line.}
For \(1\le q\le2\), fix \(n\ge3\), put \(a=2\pi/n\), and let
\[
 \ell_a(\alpha)=F_q(a)+F_q'(a)(\alpha-a)=F_q'(a)\,\alpha+c_a,
 \qquad c_a=F_q(a)-aF_q'(a),
\]
be the tangent line to $F_q$ at $a$. We need the strict inequality
\begin{equation}\label{eq:support}
 F_q(\alpha)>\ell_a(\alpha)\qquad(0<\alpha<\pi,\ \alpha\ne a).
\end{equation}
Theorem~A implies \eqref{eq:support} for every exponent. The proofs for torsion and for the first eigenvalue use, respectively,
the flux inequality of Section~\ref{sec:torsion-concavity} at $q=1$ and, at $q=2$
with $a\le\pi/2$, Propositions~\ref{prop:spectral-endpoint-curvature}
and~\ref{prop:square-tangent}. For $q=2$ and $n\ge3$ the regular openings satisfy the stronger bound
\begin{equation}\label{eq:spectral-quadratic-gap}
 F_2(\alpha)-\ell_a(\alpha)\ge c\,(\alpha-a)^2/n^2
 \quad(0\le\alpha<\pi),\qquad
 c=\begin{cases}1/1369,&n=3,4,\\3/4,&n=5,\\1,&n\ge6.\end{cases}
\end{equation}
Section~\ref{ang:sec:square-tangent} obtains \eqref{eq:spectral-quadratic-gap} by integrating curvature bounds. For $n\ge5$ the tangent at $\alpha=\pi/2$ suffices at obtuse openings, by \eqref{eq:T90}; for $n=3,4$ we use the curvature estimates for $q=2$ at obtuse openings.
We need two consequences:
\begin{equation}\label{eq:line-constants}
 c_a<0,\qquad 0<F_q'(a)\le \frac{F_q(a)}{\sin a}.
\end{equation}
Strict convexity and $F_q(0)=0$ give $c_a<0$, hence $aF_q'(a)=F_q(a)-c_a>0$. Finally (F3) gives
\(F_q'(a)/F_q(a)\le\frac{d}{d\alpha}\log\tan(\alpha/2)\big|_{\alpha=a}
=1/\sin a\).

There is a sharp correction for torsional rigidity itself. Write
$\tau_j=T(R_j)/|R_j|^2$. Reflection gives $4\pi^2\tau_j=jg(2\pi/j)$.
Since $g(0)=0$ and $g$ is strictly concave by Lemma~\ref{lem:torsion-curvature},
$\tau_j$ strictly increases (see also \cite[Theorem~1.5]{ghlz2026}).
Jensen's inequality for $-g$ at the two neighboring openings, with weights
$j-1,j+1$, gives $2\tau_j>\tau_{j-1}+\tau_{j+1}$ for $j\ge4$. Thus for every
simple $k$-gon with $3\le k\le m<n$,
\[
 \frac{T(P)}{|P|^2}\le \tau_n-(n-k)\frac{\tau_n-\tau_m}{n-m}.
\]
Indeed Theorem~\ref{thm:main} gives $T(P)/|P|^2\le\tau_k$, and averaging the
decreasing differences proves the displayed inequality. Equality holds exactly
for a regular $m$-gon, so the coefficient is optimal.

\section{Mixed-cell partitions of polygons}\label{sec:convex-faces}
\subsection{Weighted side cells of a convex polygon}
Let $P$ be a convex $m$-gon of positive area with its flat vertices removed. Write $E_i=[V_i,V_{i+1}]$, let $\theta_i$ be the interior angle at $V_i$, and let $d_i$ be the affine inward distance to the line of $E_i$. For positive weights $r_i$ the side cells are
\[
 Q_i=\{x\in\overline P:r_i d_i(x)\le r_j d_j(x)\text{ for every }j\}.
\]
These convex polygons cover $\overline P$ and have disjoint interiors. Each has positive area, since the inequality against every other side is strict in the relative interior of its base $E_i$. Let $\beta_i,\beta_i'$ be its angles at $V_i,V_{i+1}$, respectively, and put $\alpha_i=\pi-\beta_i-\beta_i'$.

\begin{lemma}[Prescribing the openings of weighted side cells]\label{lem:balance}
On a convex polygon let $\nu$ be a finite positive measure with density positive almost everywhere. Mark either all edges (in cyclic order) or one nonempty consecutive chain. Edge $i$ receives the cell on which $e^{t_i}d_i$ is least among the marked supporting distances; write $\alpha_i$ for its actual opening. Here a marked edge is a whole side of the polygon, with consecutive collinear boundary segments counted as one.
For a proper consecutive block $I$ of marked sides put
\[
 D(I)=\pi|I|-\sum_{v\ {\rm internal\ to}\ I}\sigma_v
       -\sum_{v\ {\rm terminal\ and\ met\ by}\ I}\sigma_v ,
\]
where $\sigma_v$ is the polygon angle. The terminal sum counts only endpoints of an open marked chain. Put $\Theta=k\pi-\sum_v\sigma_v$, summing over all vertices of the marked chain or cycle.
Given nonnegative costs $b_i$ and $K>0$ with $\sum b_i+K\nu(P)=\Theta$, the equations
\[
               \alpha_i=b_i+K\nu(Q_i)
\]
have a solution precisely when $\sum_{i\in I} b_i<D(I)$ for every nonempty proper consecutive block. Adding a common constant to all the $t_i$ gives the only freedom in the solution.
\end{lemma}
\begin{proof}
Where two marked edges meet at interior angle $\sigma$, the portion $\beta$ on side $i$ satisfies $e^{t_i-t_j}\sin\beta=\sin(\sigma-\beta)$. It decreases strictly from $\sigma$ to zero as $t_i-t_j$ increases. A terminal marked side takes the entire terminal polygon angle. Two weighted distances coincide only on a set of area zero, so the cell masses depend continuously on the weights. The equation for $\beta$ at each junction gives the same continuity for the angles. Every marked cell has positive area, since its weighted distance is strictly least near any point in the relative interior of its base. Summing openings gives $\Theta$; on a proper block $I$ the sum is $D(I)$ minus the positive angles taken at the interfaces with its complement. This proves necessity.

For sufficiency set $R_i(t)=\alpha_i-b_i-K\nu(Q_i)$. Each residual is bounded and continuous, and $\sum_i R_i=0$.
For $\varepsilon>0$, Brouwer applied to $t\mapsto-R(t)/\varepsilon$ on a sufficiently large ball in $\{\sum t_i=0\}$ gives $R(t)+\varepsilon t=0$.

If these solutions are unbounded as $\varepsilon\downarrow0$, pass to a subsequence along which every $t_i-\min_j t_j$ has a limit in $[0,\infty]$. The indices with finite limit form a nonempty proper set $J$, and all $t_j$ there tend to $-\infty$. Put $S=J^c$. At every interior point the distances are positive, so $\nu(\bigcup_{i\in S}Q_i)\to0$. At a junction between the two sets, the angle assigned to $S$ tends to zero. Consequently the sum of its openings tends to $\sum_I D(I)$ over its maximal consecutive blocks. Thus $\lim\sum_S R_i=\sum_I D(I)-\sum_S b_i>0$.
This contradicts $\sum_S R_i=\varepsilon\sum_J t_i<0$. A convergent subsequence yields a zero of $R$.

Finally compare two solutions $t,t'$ on the set $S$ where $t_i-t_i'$ is largest. If $S$ is proper, its union of cells can only shrink on passing from $t'$ to $t$, while the equation for $\beta$ at each junction shows that no opening in $S$ decreases and at least one increases strictly. Summing the equations $\alpha_i=b_i+K\nu(Q_i)$ over $S$ gives a contradiction. Thus the two solutions differ by a common shift.
\end{proof}

\begin{remark}
With all sides marked, the graph of $\min_i r_i d_i$ is a concave polyhedral roof over $P$ whose faces project to the weighted cells. Its height vanishes on the boundary of $P$.
\end{remark}

\subsection{Cutting a simple polygon into mixed cells}\label{sec:geometric-reduction}
At the opening $a=2\pi/n$ of the regular polygon, write the supporting line as
$s(\alpha-b)$. On a cell of area $S$, the lower bound
$\alpha\ge b+(2\pi-nb)S/|P|$ gives a uniform lower bound for
$F_q(\alpha)/S$. We prove a partition theorem with an arbitrary
positive density, which will also allow the use of the eigenfunction
in the stability argument.

\begin{lemma}[A measure and a cost for each base]\label{lem:mass-cap-partition}
Let $P$ be a simple polygon with $N$ effective sides, $d\nu=f\,dx$ a probability with $f>0$ almost everywhere, and assign each side a cost $b_i\ge0$, with $B=\sum_{i=1}^N b_i<\pi$. There is a partition into convex polygons with one Dirichlet base each on distinct original sides. If the base of $Q_i$ lies on side $i$, its opening satisfies
\[
                  \alpha_i\ge b_i+(2\pi-B)\nu(Q_i).
\]
The bases may be subsegments of the original sides. Below we use the case $b_i=b$ for every side.
\end{lemma}

Put $h=2\pi-B$. A piece of the construction, called a \emph{cap},
has a chain of $k\ge1$ Dirichlet edges on distinct original sides
and a complementary Neumann chain. We call an edge of the Neumann chain a \emph{seam}. The Neumann chain must turn through
nonnegative exterior angles at its internal vertices. For an interior
angle $\sigma_v$ at a vertex $v$, set
\[
 \Theta(C)=\pi-\sum_{v\text{ internal to the Neumann chain}}(\pi-\sigma_v),
 \qquad M(C)=\sum_{i\ {\rm in}\ C} b_i+h\nu(C).
\]
We call a cap admissible if $\Theta(C)\ge M(C)$.
A straight Neumann chain gives opening $\pi$. Two Neumann segments meeting at $W$
give opening $\sigma_W$. By the polygon angle sum,
$\Theta(C)=k\pi-\sum_{j=0}^k\sigma_j$, where the sum runs over the
vertices of the Dirichlet chain, including its two ends. For one base
this is the actual opening used in Lemma~\ref{lem:mixed-cell}.
Releasing an edge means assigning it Neumann data and dropping its
side label. We may also do this on the boundary of $P$: a function that vanishes there is still an admissible test when the zero-trace condition is removed.
Figure~\ref{fig:cap-cuts} shows the two cuts used in the recursion.

\begin{figure}[htbp]
\centering
\begin{tikzpicture}[
  x=1.13cm,y=1.13cm,font=\small,line join=round,
  dir/.style={MidnightBlue,line width=1.15pt},
  seam/.style={black,line width=.65pt},
  cut/.style={black,dash pattern=on 4pt off 2.4pt,line width=.8pt}
]
\begin{scope}
  \coordinate (A) at (0,0); \coordinate (V) at (1.35,0); \coordinate (X) at (4.05,0);
  \fill[blue!9] (0,0) rectangle (4.05,1.1);
  \fill[orange!14] (1.35,-1.1) rectangle (4.05,0);
  \draw[dir] (V)--(1.35,-1.1)--(4.05,-1.1)--(4.05,1.1)--(0,1.1)--(A);
  \draw[seam] (A)--(V); \draw[cut] (V)--(X);
  \node at (1.8,.62) {$\nu_+=5/7$};
  \node at (2.6,-.57) {$\nu_-=2/7$};
  \draw[seam] (4.05,-1.1)--(X);
  \draw (3.80,.0)--(3.80,.24)--(4.05,.24);
  \fill (A) circle(1.2pt) node[below left] {$A$};
  \fill (V) circle(1.2pt) node[below left] {$V$};
  \fill (X) circle(1.2pt) node[right] {$X$};
\end{scope}
\begin{scope}[shift={(6,-.65)}]
  \coordinate (A) at (0,0);
  \coordinate (B) at (4.6,0);
  \coordinate (W) at (.55,0);
  \coordinate (U) at (2,.52);
  \coordinate (V) at (3.58,1.94);
  \fill[MidnightBlue!4]
    (A)--(B)--(4.6,2.16)--(4.08,2.55)--(V)--(2.84,2.58)
    --(U)--(1.38,2.22)--(0,1.83)--cycle;
  \draw[dir] (B)--(4.6,2.16)--(4.08,2.55)--(V)--(2.84,2.58)
    --(U)--(1.38,2.22)--(0,1.83)--(A);
  \draw[seam] (A)--(B);
  \draw[cut] (W)--(U);
  \draw[densely dotted,line width=.8pt] (U)--(V);
  \node[below,inner sep=3pt] at (A) {$A$};
  \node[below,inner sep=3pt] at (B) {$B$};
  \fill (W) circle (1.3pt) node[below,inner sep=3pt] {$W$};
  \fill (U) circle (1.3pt) node[below right,inner sep=2pt] {$U$};
  \fill (V) circle (1.3pt) node[below right,inner sep=2pt] {$V$};
\end{scope}
\end{tikzpicture}
\caption{The heavy boundary is the Dirichlet chain. Left, the Neumann
seam extends through the reflex endpoint $V$ to $X$. Above $X$ the right edge remains Dirichlet. Below $X$ it becomes Neumann.
Right, the shortest path from $W\in AB$ to $V$ first bends at $U$.
Only the initial segment $WU$ is cut. The dotted continuation $UV$ lies
in one child.}
\label{fig:cap-cuts}
\end{figure}

\paragraph{A nonconvex example with unequal mass and area.}
In the left panel of Figure~\ref{fig:cap-cuts}, take vertices
$(0,0),(1,0),(1,-1),(3,-1),(3,1),(0,1)$ in order and use the seam from $(0,0)$ to $(1,0)$. Set every side cost to zero. Use the density $1+y/2$, normalized to give the parent mass one. Extending the seam splits the polygon into rectangles of areas $3$ and $2$, with unnormalized masses $15/4$ and $3/2$. The parent cap $C$ has opening $\pi$, so with this normalization its local coefficient is $k_C=\Theta(C)/\nu(C)=\pi$. Before releasing a side both children have opening $\pi$. Changing either portion of the right edge to Neumann decreases the opening of its rectangle by $\pi/2$. The excess of opening over required mass is $2\pi/7$ above the cut and $5\pi/7$ below it. Only the lower rectangle can accommodate this decrease. We therefore release the lower edge. The upper child's opening per unit mass is then $7\pi/5$, while the lower child's is $7\pi/4$. Both exceed the parent's $\pi$.

\begin{lemma}[Starting the partition]\label{lem:cap-start}
There is a partition of $P$ into admissible caps with disjoint sets of
side labels.
\end{lemma}
\begin{proof}
Triangulate $P$ at its vertices. Give a triangle $K$ mass $h\nu(K)$
plus the costs of any sides it shares with $P$. The masses total $2\pi$.
If a diagonal bisects this mass, each piece has mass and opening $\pi$, as required. Otherwise follow
branches of mass greater than $\pi$ in the dual tree. No step can be
reversed, so there is a triangle $K$ with every branch of mass less than
$\pi$. For an interior side $i$, call the polygonal piece across it a wing and denote its mass by $L_i$. If side $i$ lies on the boundary of $P$, set $L_i$ equal to its cost. Then $L_i<\pi$ and $\sum_iL_i+h\nu(K)=2\pi$.

Let $v_i$ be the angle opposite side $i$ of $K$. If
$v_i\ge L_i+h\nu(K)$, give the whole triangle to that wing, retaining
side $i$ as the base if the wing is empty. The new cap has opening
$v_i$. Each other nonempty wing has a single straight seam and hence opening $\pi$.
If a new seam is a side of $P$, give it Neumann data.

Otherwise apply Lemma~\ref{lem:balance} with all three sides of $K$ marked, with side costs $L_i$ and coefficient $h$. For a single side the hypothesis is $L_i<\pi$. For two sides it is $L_j+L_k<2\pi-v_i$, or equivalently $v_i<L_i+h\nu(K)$. The sum condition is the preceding mass identity. The three weighted distances agree at a unique interior point. Indeed, $|E_i|d_i/(2|K|)$ are barycentric coordinates and equality of the three weighted distances fixes a positive triple with sum one. The segments joining that point to the vertices separate the three cells. The triangular cell on side $i$ now has opening $L_i+h\nu(K_i)$. Joining it to its wing gives the required admissible cap.
\end{proof}

\paragraph{Choosing a cut so that the opening bound holds on both sides.}
For a cap of positive mass put
\[
 k_C=\frac{\Theta(C)-\sum_{i\in C}b_i}{\nu(C)},\qquad
 D_K(C)=\Theta(C)-\sum_{i\in C}b_i-K\nu(C).
\]
Every cut below can be made with $k_{C'}\ge k_C$ in every child $C'$.
For a single cut, fix $K=k_C$.
If $m_i$ children provisionally carry parent label $i$, measure additivity gives
\[
 \sum_{C'}D_K(C')=D_K(C)+\sum_{C'}\Theta(C')-\Theta(C)
                           +\sum_i(1-m_i)b_i .
\]
When terminal edge $i$ becomes a Neumann edge, its cost $b_i$ drops out and the opening decreases by its turn $\omega$. Thus $D_K$ changes by $b_i-\omega$ on that piece. The next two lemmas choose a cut for which every piece has $D_K\ge0$, first by extending a Neumann edge and then by varying the starting point of a shortest path. For the next cut we take $K$ to be the quotient of the piece being divided.

\begin{lemma}[Extending a seam]\label{lem:cap-extend}
If a junction of the two chains in an admissible cap is flat or
reflex, the cap can be replaced by admissible caps with disjoint
side labels, each having a smaller pair $(k,\#\mathrm{reflex})$
in lexicographic order.
\end{lemma}
\begin{proof}
In this proof use $M(C')=\sum_{i\in C'}b_i+k_C\nu(C')$ for pieces of the parent $C$. Merge adjacent collinear edges of the same type. At a flat junction
release the terminal Dirichlet edge. This decreases $k$ and preserves
$\Theta$. Such a junction requires $k>1$, since for a single base the formula $\Theta(C)=\pi-\sigma_0-\sigma_1$ would give $\Theta(C)<0$.

At a reflex junction $V$, extend the terminal Neumann edge straight
through $V$ until it first meets the boundary again at $X$. The ray initially lies inside
the cap. The Neumann edge directions turn by a total of
$\pi-\Theta<\pi$, so all have positive projection on one direction.
Along the extension the projection exceeds its value at the end of the Neumann chain, so the ray can meet only the Dirichlet chain. One child extends the old seam and
keeps $\Theta$, the other has opening $\pi$. Each has a subset of the parent labels and smaller measure, hence satisfies its opening bound.

A hit at a vertex creates no common label. A hit in the interior of
an edge gives both children that label provisionally. Let $b_X$ be the duplicated cost. Their two provisional surpluses $D_K$ are nonnegative and sum to $\pi-b_X$ by the preceding identity for the sum of $D_K$. If the two interior angles at $X$ are $y_1,y_2$, their sum is $\pi$. For each child, compute the value of $D_K$ it would have if its portion of that edge became Neumann. The two values sum to $(\pi-b_X)-(\pi-2b_X)=b_X\ge0$. Make this change on a side where the computed value is nonnegative, and keep the edge Dirichlet in the other child. Its new Neumann turn $\pi-y$ is
nonnegative. Its last base cannot be released this way: for a single
base with other angle $\sigma'>0$, substituting into the angle formula after release would give $-\sigma'<0$, contrary to the opening bound. In both children the reflex at $V$ disappears.
Splitting the angle at $X$ creates no reflex unless $X$ was already
reflex, and then at most one remains. This proves the decrease.
\end{proof}

\begin{lemma}[Cutting along the first segment of a shortest path]\label{lem:cap-sweep}
An admissible nonconvex cap whose two junctions are strictly convex
can be cut into admissible caps with disjoint side labels, each having
fewer Dirichlet edges.
\end{lemma}

When a cut passes through a vertex, one of the labeled edges can become part of the cut. The next lemma lets us keep track of its boundary measure in the limit.

\begin{lemma}[Persistence of the component containing a fixed point]\label{lem:grazing}
Let $X$ be a compact polygonal region. Let $S_j,S$ be nonempty closed subsets of $X$ such that every limit of a convergent subsequence $x_{j_k}\in S_{j_k}$ belongs to $S$.
Fix $o\in X\setminus S$. Let $U$ be its path component in $X\setminus S$ and, for all sufficiently large $j$, let $U_j$ be its path component in $X\setminus S_j$. Then for any finite Borel measure $\eta$ on $X$,
\[
                    \eta(U)\le\liminf_{j\to\infty}\eta(U_j).
\]
\end{lemma}
\begin{proof}
Every compact set $L\subset X\setminus S$ eventually misses $S_j$. Otherwise a subsequence of points in $L\cap S_j$ converges in $L$ to a point of $S$. Apply this first to $\{o\}$ and then, for each $x\in U$, to the image of a path in $U$ from $o$ to $x$. Hence $x\in U_j$ eventually. The complements are relatively open in a locally path-connected polygonal region, so their path components are Borel. Fatou's lemma applied to $\mathbf1_U\le\liminf_j\mathbf1_{U_j}$ proves the claim, including when $\eta(S)>0$.
\end{proof}

\begin{proof}[Proof of Lemma~\ref{lem:cap-sweep}]
Fix $K=k_C$ throughout this proof, so the interior measure is $K\nu$. Orient the cap counterclockwise, with Neumann chain from $A$ to $B$.
Fix an internal reflex vertex $V$ on the Dirichlet chain. For $W$ on
the Neumann chain take the shortest path in the closed cap from $W$
to $V$. Such a path exists by compactness and bends only at reflex
vertices by local shortening. It is polygonal, since it cannot revisit
a vertex and only finitely many vertices can be bends. If two minimizing paths differed, two of their subpaths would form a simple closed polygonal curve between a point where they separate and a point where they meet again. Its interior lies in
the cap by simple connectivity. Every simple polygon has at least three strictly convex corners. One of them
must be an interior point of one of the two subpaths, where cutting across the corner would shorten the path. This proves uniqueness.

\emph{Which segment to cut.}
Let $S(W)$ be the whole initial straight segment ending where the path first changes direction, or at $V$ if the path is straight. In particular the cut continues past any collinear
boundary contacts. For $W$ in the open chain, its first boundary
contact after $W$ is an internal Dirichlet vertex: a straight path
cannot graze a strictly convex vertex or touch an open side
transversely and continue, and starting along the Neumann boundary
would encounter a strictly convex vertex first. At $A$ or $B$ it may
follow a terminal Dirichlet edge, which is released.
There is a direct description valid at every grazing contact. Close each component arc of $\partial C\setminus S(W)$ by the straight segment between its endpoints in $S(W)$, keeping the old boundary orientation. Each resulting curve is simple and closed. Its open boundary arc misses the closing segment, and the two endpoints are distinct. Its interior lies in $C$ since the exterior of $\overline C$ is connected, and its orientation is positive, as is seen on the old boundary arc. Off $S(W)$ their winding numbers sum to that of $\partial C$: the difference is a closed polygonal chain supported on a single segment. Thus these regions cover the cap with disjoint interiors. For each region, the portions of $S(W)$ beyond its closing chord cannot enter its interior. Such a portion ends on $\partial C$ and meets the boundary of the region only at the end of the chord. The closed regions with the cut removed are therefore the components of $\overline C\setminus S(W)$. The two main pieces adjoin the old Neumann subarcs $AW,WB$; any extra piece has just one straight seam and hence opening $\pi$. A boundary portion contained in the cut is released. Since the contacts on the Dirichlet chain are vertices or whole edges, no label is duplicated. The first internal contact splits the chain, or a terminal edge is released, so every descendant has fewer labels. For each extra piece the sum of its side costs and its interior measure is at most the corresponding sum for the parent, hence at most $\pi$.

\emph{What varies continuously.}
Lift the first and last Neumann edge directions to $d_0,d_1$ with
$d_1-d_0=\pi-\Theta$. For $W$ in the open chain, lift the outgoing cut direction
$t$ into the local interior sector. The main openings are
\[
 \Theta_A=\pi+d_0-t,\qquad \Theta_B=t-d_1,\qquad
 \Theta_A+\Theta_B=\Theta .
\]
The first chain turns from $d_0$ to $t$ and the second from $t-\pi$
to $d_1$. This also holds at a Neumann corner. For $W_j\to W$, the distance to $V$ changes by at most the length of the chain between $W_j$ and $W$. The constant-speed minimizing paths are therefore equicontinuous. Every uniform subsequential limit is the minimizing path from $W$ to $V$, with constant speed: lower semicontinuity of length gives minimality, and its Lipschitz bound on $[0,1]$ equals the endpoint distance. Uniqueness identifies the limit.

The far endpoint of $S(W_j)$ takes values in the finite set consisting of $V$ and the possible first bends, all disjoint from the Neumann chain. On a subsequence it is a fixed vertex $Z$. The paths up to $Z$ are straight and converge to the initial segment $[W,Z]$ of the limiting path (the traversal times converge to $|Z-W|/\operatorname{dist}_C(W,V)$). Hence $[W,Z]\subset S(W)$, including at a grazing contact. Thus the initial direction $t$ varies continuously, and any limit of points on the cuts $S(W_j)$ belongs to $S(W)$, as required by Lemma~\ref{lem:grazing}. Figure~\ref{fig:grazing} shows how the inclusion can be strict.

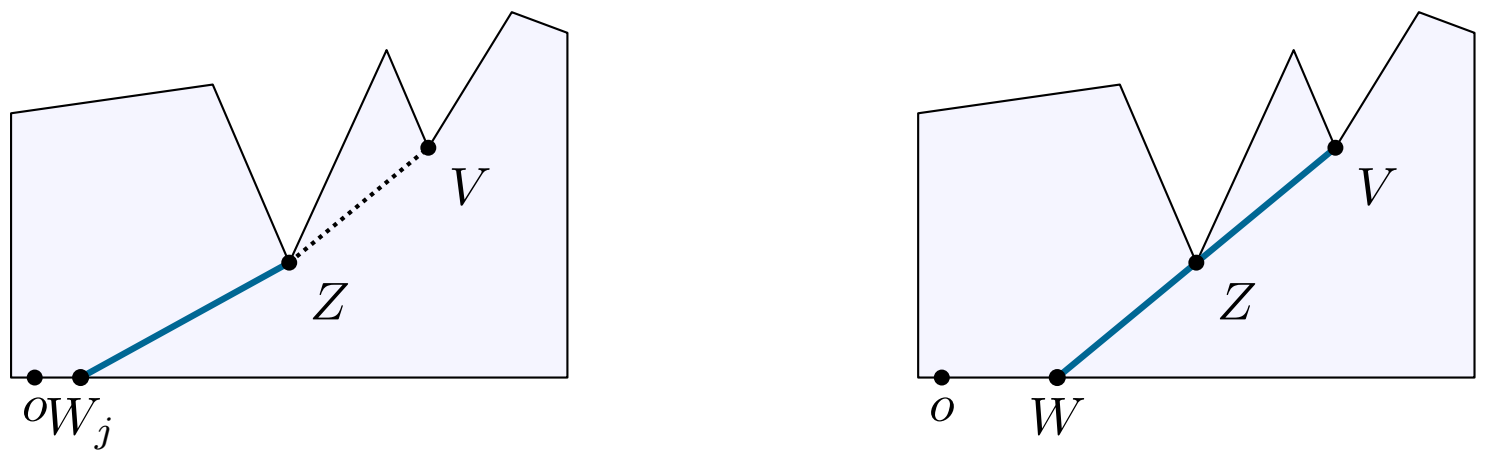
\begin{figure}[htbp]
\centering
\begin{tikzpicture}[x=.92cm,y=.76cm,font=\small,line join=round]
\foreach \offset/\start/\lab in {0/.5/{W_j},6/1/{W}}{
  \begin{scope}[xshift=\offset cm]
    \filldraw[fill=blue!4] (0,0)--(4,0)--(4,3)--(3.6,3.18)--(3,2)--(2.7,2.85)--(2,1)--(1.45,2.55)--(0,2.3)--cycle;
    \ifnum\offset=0
       \draw[densely dotted,thick] (2,1)--(3,2);
       \draw[MidnightBlue,very thick] (\start,0)--(2,1);
    \else
       \draw[MidnightBlue,very thick] (\start,0)--(3,2);
    \fi
    \fill (\start,0) circle(1.6pt) node[below] {$\lab$};
    \fill (2,1) circle(1.5pt) node[below right] {$Z$};
    \fill (3,2) circle(1.5pt) node[below right] {$V$};
    \fill (.17,0) circle(1.5pt) node[below] {$o$};
  \end{scope}}
\end{tikzpicture}
\caption{The solid blue cut stops at the first bend. As $W_j$ approaches $W$ from the left, the bend at $Z$ disappears. The limiting cut reaches $V$ and contains the limit $[W,Z]$ of the earlier cuts. The point $o$ marks the same boundary subarc in both panels.}\label{fig:grazing}
\end{figure}

\emph{Lower semicontinuity of the measure on either side of the cut.}
On the closed cap put
\[
 \eta=K\nu+\sum_E b_E\frac{\mathcal H^1|_E}{|E|},
\]
summing over its current Dirichlet edges. For each starting point $W$, either the relative interior of a given Dirichlet edge $E$ lies in one component of $\overline C\setminus S(W)$, or $E\subset S(W)$ and its label is released. The endpoints have zero arclength. Thus the sum of the side costs and the interior measure for a main piece is $\eta(U)$, where $U$ is the corresponding component after removing $S(W)$ from the closed cap. In particular the
edge measure counts no released labels. For the piece adjoining $AW$, choose $o$ in the relative interior of an old Neumann edge on the positive-length subarc $AW$, away from $W$ and from the cut. For $W_j$ close to $W$, that point still lies on $AW_j$. Its path component in $\overline C\setminus S(W_j)$ is therefore the main piece on the same side. Lemma~\ref{lem:grazing} gives lower semicontinuity of $\eta(U)$ as the starting point varies. The argument for $WB$ uses a point on $WB$. It applies at an endpoint to the piece that retains the positive-length arc. In particular, a boundary edge absorbed by the limiting cut drops out of the limiting component and contributes nothing to its $\eta$-measure.

For each main piece we require $\Theta\ge\eta(U)$. At an interior starting point $W$, the two openings sum to $\Theta(C)$ and the two measures sum to at most $\sum b_i+K\nu(C)=\Theta(C)$. Thus at least one piece satisfies its inequality. At $A$ the piece adjoining $AW$ is empty or has a single straight seam and opening $\pi$, so its inequality holds there. The opening on the $B$ side equals its one-sided limit. Its measure is lower semicontinuous there because that piece retains the open Neumann chain and we can choose the fixed point $o$ away from $A$. If the cut follows a terminal Dirichlet edge,
releasing it makes the reversed cut the incoming seam at $A$. The turn
from this seam to the old chain gives the same opening $t-d_1$ as the
one-sided limit. The convention at $B$ is symmetric. On each side, the set of starting points for which the opening inequality holds is closed. Together the two sets cover the chain from $A$ to $B$, and each contains its respective endpoint. Connectedness gives a point in their intersection. This cut has the required opening bounds.
All new Neumann turns are nonnegative, as the cut leaves $W$ inside
its local interior sector and its other portions form straight seams.
\end{proof}

\begin{lemma}[Assigning one base to each convex cell]\label{lem:cap-finish}
Every admissible convex cap can be partitioned as in
Lemma~\ref{lem:mass-cap-partition} using its surviving side labels.
For each chain left by the terminal releases below, the logarithms of
the side weights prescribing its cell openings are unique up to a common
additive constant.
\end{lemma}
\begin{proof}
List the marked bases along the Dirichlet chain. Release any proper prefix or suffix $I$ for which $D(I)\le \sum_{i\in I}b_i$, repeating until no such block remains. Here $D(I)$ is precisely the sum of the exterior turns converted into internal Neumann turns at a release, so admissibility persists and at least one base remains.
For the resulting chain set $K=(\Theta-\sum b_i)/\nu(C)\ge h$. Since no such block remains, every proper terminal block satisfies the strict inequality in Lemma~\ref{lem:balance}. Each internal block has $D(I)=\pi+\sum_{v\,{\rm internal}}(\pi-\sigma_v)\ge\pi>B$. Each remaining marked edge is a whole side of the cap. Consecutive marked edges meet at original, nonflat vertices. At a flat terminal junction we could still release the end edge (with a single base, such a junction would force $\Theta\le0$). Lemma~\ref{lem:balance} now gives openings $b_i+K\nu(Q_i)$, and the logarithms of the weights are unique up to a common additive constant. This also holds for one base.
\end{proof}

\begin{proof}[Proof of Lemma~\ref{lem:mass-cap-partition}]
Start with Lemma~\ref{lem:cap-start}. Each initial cap has $k_C\ge h$. At each cut, the quotient for each new piece is at least that of the piece being cut. It is therefore at least $h$ in every final cell. Apply Lemmas~\ref{lem:cap-extend}
and~\ref{lem:cap-sweep} to nonconvex pieces. Each branch decreases
$(k,\#\mathrm{reflex})$ lexicographically, so the finitely branching
subdivision terminates. All surviving labels stay distinct.
Finish with Lemma~\ref{lem:cap-finish} on each convex cap.
\end{proof}

\begin{proposition}[The fixed-side comparison]\label{prop:geometric-reduction}
Fix $1\le q\le2$, $n\ge3$ and $a=2\pi/n$. Suppose that $F_q$ has a supporting line
$s(\alpha-b)$ with $s>0$, $0<b<\pi/n$, and equality precisely at $a$.
Every simple polygon $P$ with $N\le n$ effective sides then satisfies
$|P|^{2/q}\lambda_q^{\rm S}(P)\ge nF_q(a)$, with equality exactly for a regular $n$-gon.
\end{proposition}
\begin{proof}
Write $A=|P|$, $h=2\pi-Nb$, $c=2/q$, and normalize the positive
minimizer by $\int_P u^q=1$. Use Lemma~\ref{lem:mass-cap-partition}
with the probability measure $d\nu=(A^{-1}+u^q)\,dx/2$.
For each cell put $p_i=|Q_i|/A$ and $r_i=p_i^{-1}\int_{Q_i}u^q$.
Thus $\nu_i=\nu(Q_i)=p_i(1+r_i)/2$ and $\sum p_i=\sum p_i r_i=1$.
Convexity gives $r^c\ge1+c(r-1)$ for $r\ge0$. Multiplication by $(1+r)/2$
and summation yield
\[
 \sum_i \nu_i r_i^c\ge 1+\frac c2\sum_i p_i(r_i-1)^2.
\]
The mixed-cell comparison bounds the energy on each cell, multiplied by $A^c$, below by
$r_i^c F_q(\alpha_i)$. The supporting line and the partition give
$F_q(\alpha_i)\ge s(\alpha_i-b)\ge sh\nu_i$. Since
$sh-nF_q(a)=sb(n-N)$, we obtain
\[
 A^c\lambda_q^{\rm S}(P)-nF_q(a)
 \ge sb(n-N)+\frac{sh}{q}\sum_i p_i(r_i-1)^2
       +\sum_i r_i^c\bigl(F_q(\alpha_i)-sh\nu_i\bigr).
\]
All terms are nonnegative. If equality holds, the first two terms on the right give $N=n$ and $r_i=1$ on every cell. The remaining terms then give $F_q(\alpha_i)=s(\alpha_i-b)=sh\nu_i$ separately. Strict support gives $\alpha_i=a$; division by $sh$ gives $\nu_i=p_i=(a-b)/(2\pi-nb)=1/n$. Since the masses sum to one there are exactly $n$ cells.
Each cell comparison
is sharp as well, so the cells are $n$ isosceles triangles of equal area.
On an open portion of a Neumann side the
positive minimizer on the cell extends by even reflection and is positive by
the strong minimum principle. Since the polygonal minimizer has zero trace on the boundary of $P$, the bases must cover that boundary. There are $n$ distinct bases and $n$ true
sides, so each base is a whole side of the polygon. At each original vertex the two
base angles sum to $\pi-a$; their sum over all vertices is already
$(n-2)\pi$. Every interior angle is therefore $\pi-a$, and the
equal triangle areas give equal side lengths. This forces regularity.
\end{proof}

\subsection{Proof of Theorem~\ref{thm:main}}\label{sec:main-proof}
\emph{Three and four sides.}\label{par:small-polygons}
Steiner symmetrization preserves the area and the $L^q$ mass and decreases
the Dirichlet energy. Use the classical geometric construction
\cite[Section~6.3.1, proof of Theorem~6.2]{S:laugesen-siudeja2017}. The same construction
works for all the quotients considered here. We record why a step is
strict whenever it changes the shape.

In the symmetrization direction write the sections as intervals, with transverse
coordinate $y$. If their centers vary, choose a closed subinterval between
vertex heights on which the center has nonzero slope. The positive minimizer
vanishes on the two straight edges there, with nonzero inward derivative by
the Hopf lemma. For all small $t>0$ its superlevel section has exactly two
endpoints $a(y,t),b(y,t)$ on this subinterval. As $t\downarrow0$ their
$y$-derivatives converge locally uniformly to the edge slopes by boundary
$C^1$ regularity and the implicit function theorem. Put $A=a_t>0$, $B=-b_t>0$.
In the coarea proof of the Steiner energy inequality, the energy density
minus its symmetrized value on this set is
\[
 \frac{1+a_y^2}{A}+\frac{1+b_y^2}{B}
  -\frac{4+(b_y-a_y)^2}{A+B}
 =\frac{(A-B)^2+(Ba_y+Ab_y)^2}{AB(A+B)} .
\]
It is strictly positive there for small $t$ since $a_y+b_y$ stays away from
zero. Elsewhere the coarea and Cauchy--Schwarz inequalities give a
nonnegative difference. This proves strict decrease of the quotient.
Constant section centers make the two domains translates.

A minimizing triangle exists at fixed area $S$. A longest side $d$ puts
the triangle in a strip of width $2S/d$, where one-dimensional Poincar\'e
and H\"older give $\lambda_q^{\rm S}\ge S^{1-2/q}(\pi d/(2S))^2$.
The remaining triangles form a compact family up to rigid motions, and
affine pullback gives continuity. Choose a side as the base and center every chord parallel to the base on the perpendicular bisector of the base.
This gives an isosceles triangle with the same base and area.
By the strictness criterion above, the original chords of a minimizing triangle
must already have a common center coordinate. The base fixes this coordinate
at its midpoint, so its perpendicular bisector is an axis of symmetry. This characterizes the
equilateral triangle.

For any simple quadrilateral choose an interior diagonal. Its two adjacent
triangles lie on opposite sides of the diagonal, so symmetrizing parallel
to it gives a kite, including when the original quadrilateral is nonconvex.
Symmetrizing parallel to the other diagonal gives a rhombus, and parallel
to a rhombus side gives a rectangle. The sections in these steps are
intervals with piecewise affine centers. A triangle is included by inserting
a point on one side and joining it to the opposite vertex. Among rectangles
of area $S$ a minimizer exists by the same strip bound and affine continuity.
If it had unequal sides, the first diagonal symmetrization in the construction
would strictly improve it, since its diagonal is not an axis, and the next
two steps would give a rectangle of smaller quotient. It follows that the
square minimizes among all competitors with at most four sides. Equality
at every step forces congruence to the square.

\emph{Five or more allowed sides.} Set $a=2\pi/n$.
Section~\ref{sec:angular-input} gives the strict support \eqref{eq:support}
and the slope bounds \eqref{eq:line-constants}. For $q=1$ these follow from the torsion flux inequality. For $q=2$ we use the eigenvalue comparison at acute openings and \eqref{eq:T90}. For intermediate exponents we use Theorem~A. With $b=a-F_q(a)/F_q'(a)$ these bounds yield
$0<b\le a-\sin a<a/2$.
Proposition~\ref{prop:geometric-reduction} applies to every simple polygon
with at most $n$ sides, including convex ones, and gives the stated equality case.
The regular polygon attains equality by \eqref{eq:regular-reflection}.\phantomsection\label{par:geometry-route-end}

\section{Angular estimates for torsion and the first eigenvalue}
\label{sec:spectral-angle}\label{sec:sobolev-angles}
\begingroup
\providecommand{\dd}{\,\mathrm d}
\providecommand{\dn}{\mathrm{DN}}
\subsection{Concavity of the normalized mixed torsional rigidity}\label{sec:torsion-concavity}

Concavity of $g(\alpha)=\alpha^2/F_1(\alpha)$ in the opening was proved in \cite[Theorem~8]{bbf2026}. The zeroth and first transverse moments used below also appear, as integrals against polynomial test functions, in \cite[Section~6]{bbf2026}. The following curvature bound gives the constant used in Section~\ref{sec:tangent-gaps}. Its proof holds the minimizing flux fixed as the opening varies.
\begin{lemma}\label{lem:torsion-curvature}
For \(0<\alpha<\pi\) and the half-width \(t\) of Lemma~\ref{lem:affine-family},
\[
         -g''(\alpha)\ge\frac{5t^3}{6(t^4+6t^2+3)}>0.\tag{TC}
\]
In particular the bound \(-g''\ge t^3/[6(1+t^2)]\)
holds for \(t\le4/5\).
\end{lemma}

Strict concavity of $g$ and strict support of every tangent to $F_1$ follow directly
from (TC). Indeed, \(g>0\), and differentiating \(F_1(\alpha)=\alpha^2/g(\alpha)\) twice gives
\begin{equation*}
 F_1''(\alpha)=\frac{2\bigl(g(\alpha)-\alpha g'(\alpha)\bigr)^2}{g(\alpha)^3}
                -\frac{\alpha^2 g''(\alpha)}{g(\alpha)^2}>0 .
\end{equation*}
Thus \(F_1'\) is strictly increasing. Integrating between \(a\) and \(\alpha\) gives \(F_1(\alpha)>F_1(a)+F_1'(a)(\alpha-a)\) whenever \(\alpha\ne a\).

\begingroup
\def\dd{\mathrm d}
\begin{proof}
Put $\theta=\alpha/2$, $t=\tan\theta$, $T=t^2$ and $f=2\theta^2/t$.
On $D=\{0<y<x<1\}$ call a flux $J=(p,q)$ admissible if
$\int_D J\cdot\nabla v=\int_Dv$ for all tests vanishing at $x=1$.
Write $X=\int_Dp^2$, $Y=\int_Dq^2$. Lemma~\ref{lem:affine-family} and
Cauchy--Schwarz give the following variational formula for the flux
\[
      g(2\theta)/4=f(\theta)\min_J(X+TY).
\]
Indeed, on the triangle of height one the area is $t$ and its torsional rigidity is $2t\int_D u$ by reflection and the affine change of variables. Equality in Cauchy--Schwarz is attained by the flux
$(u_x,T^{-1}u_y)$ of the weak torsion solution on $D$.  The set of
admissible fluxes is independent of $\theta$. We will prove that
the minimizing flux has
\begin{equation*}
        X\ge1/16,\qquad (T^2+6T+3)Y\le X.                 \tag{*}
\end{equation*}

There are two elementary admissible fluxes:
$J_0=(-x/2,-y/2)$ and $J_1=(y-x,0)$. Both have divergence $-1$
and zero normal component on the two free edges.
Let $J=(p,q)$ minimize $X+TY$, set $r=p+x/2$, and put $h=y-x/2$.
The horizontal flux through the vertical section at $x$ is $-x^2/2$: the divergence is $-1$, the triangle to the left has area $x^2/2$, and the flux across its two free sides vanishes. Thus $\int_0^x r(x,y)\,dy=0$ for almost every $x$, and $X=1/16+\|r\|^2$, with norm and pairing in $L^2(D)$.
Orthogonality of $J$ to the variations $J_0-J$ and $J_1-J$ in the
weighted norm yields
\[
 \ell:=\|r\|^2+TY=\langle h,r\rangle=-T\langle y/2,q\rangle .
\]
Here $\langle x,h\rangle=0$ by integration on $0<y<x$.
To obtain the first moment on the same section, apply the flux balance to $yJ$. Its divergence is $q-y$, so
\[
 \int_0^x y p(x,y)\,dy=-x^3/6+\int_0^x\int_0^s q(s,y)\,dy\,ds.
\]
In the weak formulation, test with a smooth approximation to the indicator of the subtriangle, first on its own and then multiplied by $y$. This gives the same two identities. Integrating in $x$ and using the zero slice mean of $r$ yields
$\ell=\langle h,r\rangle=\langle y,r\rangle=1/48+\langle1-x,q\rangle$.
Put $M=48\ell$. Since
$\|h\|^2=\|y/2\|^2=1/48$ and $\|x-1\|^2=1/12$, the last three
pairings give
\[
  M^2\le48\|r\|^2,\qquad M^2\le48T^2Y,\qquad
  (1-M)^2\le192Y .
\]
In particular $48TY\le M(1-M)$. These inequalities and the
pairings exclude $M=0,1$ and imply
$T/(T+4)\le M\le T/(T+1)$. Consequently
\[
 \frac X Y\ge T\frac{3+M^2}{M(1-M)}
   =T\left(\frac3M+\frac4{1-M}-1\right)\ge T^2+6T+3 .
\]
This proves (*).

Now fix $\theta$ and keep its minimizing flux fixed at nearby angles.
The function $fX+(fT)Y$ is an upper bound for $g(2\theta)/4$, with
equality at the chosen angle. Thus $g''(2\theta)\le f''X+(fT)''Y$,
where primes on $f$ and $T$ mean $\theta$ derivatives.
The weak torsion solution is smooth in $T>0$ since its defining coercive
form is affine in $T^{-1}$. As $(fT)''=(2\theta^2\tan\theta)''>0$,
it suffices, by (*), to bound $-f''-(fT)''/(T^2+6T+3)$ from below.

Write $B=T^2+6T+3$, $d=1-\theta/t$. Differentiation gives
\[
 \frac t4[-Bf''-(fT)'']
 =3T^2+(1+T)d\{2T(1+T)-(2T^2+6T+3)d\}.
\]
Indeed $t f''/4=-T+(1+T)d^2$ and
$(fT)''/(4t)=1+2(1+T)(1-d)+T(1+T)(1-d)^2$.
Since $d=\int_0^1 Tu^2/(1+Tu^2)\,du$, comparison with
$Tu^2/(1+T)$ and concavity in $u^2$ give
$T/[3(1+T)]\le d\le T/(T+3)$.
The expression in braces is at least $T(2T+3)/(T+3)\ge T$.
The display is therefore at least $10T^2/3$. Using (*) now yields
\[
          -g''(2\theta)\ge\frac{5t^3}{6(t^4+6t^2+3)}.
\]
This is (TC). For $t\le4/5$, the remaining assertion follows from $5(1+t^2)\ge t^4+6t^2+3$.
\end{proof}
\endgroup
\label{end:torsion}

\subsection{Curvature of the first eigenvalue profile for \texorpdfstring{$0<\alpha\le\pi/3$}{openings up to 60 degrees}}\label{ang:sec:q2-start}\label{ang:sec:acute-curvature}
\begin{proposition}\label{prop:acute-curvature}
For $s=\alpha/(2\pi)$ and $0<s\le1/6$, the mixed eigenvalue profile satisfies
\[
 F_2''(2\pi s)>4s^2.
\]
\end{proposition}

\paragraph{Sector coordinates.}
For $s=1/n$, take the Schwarz--Christoffel map from the disk onto the regular $n$-gon and restrict it to one sector \cite[Eq.~(5)]{dgp2026}, \cite[Ch.~4]{dt2002}. For continuous $s$, take the primitive $f(0)=0$ with $f'(\xi)=(1-\xi^{1/s})^{-2s}$ on the sector branch $0<\arg\xi<2\pi s$. This map sends $0<\arg\xi<2\pi s$, $|\xi|<1$, onto an isosceles triangle. The two straight radii have image length $s\mathrm B(s,1-2s)$, and on $\xi=e^{is\theta}$ the identity $\arg(1-e^{i\theta})=\theta/2-\pi/2$ gives the tangent angle $\pi/2+\pi s$. These three segments form a simple boundary, so the argument principle gives univalence. The positive mixed ground state is symmetric across the angle bisector. Thus its traces agree on the free rays, and imposing periodicity there preserves the first eigenvalue.

Use $\tau=-\log|\xi|$, $\theta=\arg\xi/s$ and measure $d\tau\,d\theta/(2\pi)$ on $(0,\infty)\times(\mathbb R/2\pi\mathbb Z)$. Tests vanish at $\tau=0$. Put
\[
 E_s(u)=\int(|u_\tau|^2+s^{-2}|u_\theta|^2),\qquad
 I_s(u)=\int e^{-2\tau}|1-e^{-\tau/s+i\theta}|^{-4s}|u|^2,\qquad
 L(s)=\inf E_s/I_s .
\]
The profile is $F_2(2\pi s)=\pi s Z(s)L(s)$. The area is $\pi s Z(s)$, where $Z(s)=2\int e^{-2\tau}|1-e^{-\tau/s+i\theta}|^{-4s}$. The side length and beta reflection formula give
\[
 Z(s)=\frac{\Gamma(1+s)^2\Gamma(1-2s)}{\Gamma(1-s)^2\Gamma(1+2s)},\qquad
 \log Z(s)=\sum_{m\ge3,\,m\ {\rm odd}}\frac{(2^{m+1}-4)\zeta(m)}m s^m ,
 \qquad |s|<1/2.
\]
For $j=j_{0,1}$ the normalized radial ground state for the reference mass $I_0(u,v)=\int e^{-2\tau}u\bar v$ is
$\psi(\tau)=\sqrt2\,J_0(je^{-\tau})/J_1(j)$.
It satisfies $I_0(\psi,\psi)=1$. Since $(rJ_1(r))'=rJ_0(r)\ge0$ for $0\le r\le j$ and $\psi(0)=0$, it also satisfies $|\psi(\tau)|\le\sqrt2 j\tau$.

Extend radial Dirichlet test functions on the disk by zero to $(-1,1)^2$. The extensions are even in both variables. Comparing with the first two square eigenvalues in this symmetry class gives the lower bounds below. Testing on the disk with $(1-r^2)(3-r^2)$ gives the upper bound. Here $\mu$ is the second radial eigenvalue.
\begin{equation}\label{eq:acute-bessel}
 \mu\ge5\pi^2/2>24,\qquad \pi^2/2\le j^2\le110/19<29/5.
\end{equation}
We use $6/5<\zeta(3)\le1+1/8+2/25$, by summing the first 16 terms
for the lower bound and using the midpoint integral after two terms
for the upper bound.

Write $W_s=|1-e^{-\tau/s+\mathrm i\theta}|^{-4s}$ and $K_s=I_s-I_0$.

\paragraph{Projection onto the orthogonal complement of the disk ground state.}\label{par:independent-real-projection}
Let $H=\psi^{\perp_{I_0}}$ and equip $H$ with the energy norm $E_s^{1/2}$.
Write $b_s=K_s(\psi,\cdot)|_H$, $d_s=K_s(\psi,\psi)$ and $\gamma=j^2/L-1$.
The circular mean of $W_s$ is at least one, so the radial trial gives $\gamma\ge0$.
Normalize the positive eigenfunction to $\psi+v$ with $v\in H$. Projection of its equation gives
\[
 C v=b_s,\qquad \gamma=d_s+b_s(v),\qquad C=(1+\gamma)E_s/j^2-I_0-K_s\quad\hbox{on }H .
\]
For $h=1/6$ and $t=s/h$, the bounds below use $K_j=s^j\partial_s^j K_s$, $b_j=s^j\partial_s^j b_s$, $d_j=s^j\partial_s^j d_s$. In particular $K_0=K_s$.
The indices refer to ordinary derivatives taken with $\tau$ fixed.
All decimal bounds in this proof denote rational numbers.

\begin{lemma}[Bounds for the projected equation and its derivatives]\label{lem:acute-eliminate}
For $0<s\le h$, in the energy norm and its dual,
\[
 C\ge .1168 E_s,\quad \|K_1\|<.13,\quad \|K_2\|<.36,
 \quad (\|b_0\|,\|b_1\|,\|b_2\|)<t^{3/2}(.009,.029,.095).
\]
Also $0\le d_0\le15s^5$, $|d_1|\le75s^5$ and $|d_2|\le300s^5$.
\end{lemma}
\begin{proof}
For $c_i(a)=(a)_i/i!$ put $D_i(a)=\sum_{m=0}^{i-1}a/(a+m)$ and
$U_i(a)=\sum_{m=0}^{i-1}am/(a+m)^2$. Empty sums vanish.
The coefficient indexed by $(i,l)$ in $W_s-1$ is
$c_i(2s)c_l(2s)e^{-kx}e^{\mathrm i(i-l)\theta}$, where $k=i+l>0$ and $x=\tau/s$.
For derivatives of order $0,1,2$ its real multipliers, with $D=D_i+D_l$, $U=U_i+U_l$, are
\[
 P_0(y)=1,\quad P_1(y)=D+y,\quad
 P_2(y)=D(D-1)+U+2(D-1)y+y^2,\qquad y=kx.
\]
These polynomials have nonnegative coefficients, since $D\ge1$.
The quantity $s\partial_s D_i=U_i$ increases with $s$ here (each nonzero summand has $m\ge1>2s$).

On $H$, the radial gap in \eqref{eq:acute-bessel} and the circle Poincar\'e inequality give $I_0\le E_s/24$. The absolute series for $W_s-1$ is bounded by $4s/(e^{\tau/s}-1)$ because $c_k(4s)\le4s$. We use the following ground-state inequality. Suppose $\varphi>0$ satisfies $-\operatorname{div}(A\nabla\varphi)=V\varphi$, with $A$ symmetric positive definite. For a smooth real test $v$ supported away from the Dirichlet boundary, with zero conormal derivative of $\varphi$ on any Neumann boundary met by the support, testing with $v^2/\varphi$ gives the ground-state identity
\[
 \int A\nabla v\cdot\nabla v-\int Vv^2
 =\int \varphi^2 A\nabla(v/\varphi)\cdot\nabla(v/\varphi)\ge0 .
\]
For $V\ge0$ the inequality extends by density and Fatou's lemma to the energy space with zero Dirichlet trace. Applied on the half-line to $\varphi(\tau)=1-e^{-a\tau}$, it gives $\int |v|^2/(e^{a\tau}-1)\le a^{-2}\int|v'|^2$. Since
$e^{-2\tau}/(e^{a\tau}-1)\le(1+2/a)/(e^{(a+2)\tau}-1)$, this gives $\|K_0\|\le4s^3/(1+2s)$.
In particular $1/j^2-1/24-4h^3/(1+2h)>.1168$.
For the two derivatives, sum the positive majorants at $\theta=0$. At degree $k$ their multipliers have the same formulas with $D=D_k(4s)$ and $U=U_k(4s)$. The bound $\|v(\tau,\cdot)\|_2\le\sqrt\tau\,\|v\|_{E_s}$ gives integrals of $1,kx,(kx)^2$ against $xe^{-kx}$ proportional to $1,2,6$.
Bound $D\le D_k(2/3)\le D_*:=1+(2/3)\log k$ by midpoint integration of $1/(x+2/3)$, and $U\le D$. Dropping $e^{-2\tau}$ bounds the norms by $\sum_{k\ge1}(D_*+2)/(54k^2)$ and $\sum_{k\ge1}(D_*^2+4D_*+2)/(54k^2)$.
Both summands decrease; the first term plus the integral from one gives the stated bounds.

For the radial estimates the circle average retains $i=l=n$. Write
$J_n=\int e^{-2\tau}\psi^2e^{-2n\tau/s}\,d\tau\le j^2s^3/(2n^3)$.
The variable $y=2n\tau/s$ has, in this integral, density proportional to $y^2e^{-y}$ times a decreasing function: $\psi$ is concave since $-\psi''=j^2e^{-2\tau}\psi$. The polynomials $y-3$ and $y^2-2y-6$ change sign once and have zero mean against $y^2e^{-y}$. Hence $0\le sJ_n'\le3J_n$ and $-J_n\le s^2J_n''\le6J_n$.
The recurrence gives $c_n(2s)\le2s/\sqrt n$ and $D_n(2s)\le\sqrt n$ (use $1/i\le3(\sqrt{i+1}-\sqrt i)$).
For $c_n^2$ the ratios of its first two scaled derivatives to itself are $2D_n$ and at most $4D_n^2-2$, by taking two logarithmic derivatives and retaining the factors at index zero. Its derivatives are nonnegative. The product rule now bounds the absolute scaled derivatives of $c_n^2J_n$ by itself times $1,5n,20n$. Summing gives the last three assertions since $2j^2\zeta(3)<15$.

For the forcing terms scale the test as $v(\tau,\theta)=\sqrt s\,w(x,\theta)$ so that $E_s(v)=E_1(w)$. The inequality $\psi\le\sqrt2 j\tau$ supplies the common factor $H_s=\sqrt2j s^{5/2}$. In mode $n=|i-l|$, pair $w$ with the density
$c_i c_l\,x e^{-(k+2s)x}P_j(kx)$. The norm of $x^m e^{-Bx}$ on this mode is bounded by
\[
 A_m(B,n)=
 \begin{cases}
 \Gamma(m+3/2)/B^{m+3/2},& n=0,\\
 \sqrt{(2m)!}/[n(2B)^{m+1/2}],& n>0 .
 \end{cases}
\]
For $n=0$, use the preceding bound on the trace at each $x$. For $n>0$, apply Cauchy--Schwarz in $L^2$ and use the angular energy of the mode.
For the first two degrees there are three pairs of degree and absolute Fourier frequency, $(k,n)=(1,1),(2,0),(2,2)$.
For each representative $(i,l)=(1,0),(1,1),(2,0)$ form
\[
 L_{j,i,l}=c_i c_l\sum_{m=0}^j [y^m]P_j(y)\,k^m A_{m+1}(k+2s,n).
\]
Distinct signed Fourier modes are orthogonal; the three displayed representatives have multiplicities $2,1,2$.
After division by $t^{3/2}$ each upper bound is largest at $s=h$. Indeed, $H_s/t^{3/2}$ contributes one power of $s$ and $c_ic_l$ contributes at least one more. The logarithmic derivative $s\partial_s\log (k+2s)^p$ of the denominator is at most $p/4\le9/8$.
Thus for these degrees use $H_h(2L_{j,1,0}^2+L_{j,1,1}^2+2L_{j,2,0}^2)^{1/2}$ at $h$.

For $k\ge3$ the coefficient vector has squared norm at most $c_k(2s)c_k(4s)$: its largest coefficient is $c_k(2s)$ by log-convexity in the index, and its sum is $c_k(4s)$. Both factors decrease with $k$. Midpoint integration gives $D_i(1/3)\le1+\tfrac13\log(6i/5)$ for $i\ge1$, whence $D\le d(x):=2+\tfrac23\log x$ at $x=k$ (also when one index vanishes). Apply the trace bound at each $x$ before summing over $k$. The successive moments of $kx$ under $x^{3/2}e^{-kx}$ are $1,5/2,35/4$. Since $U\le D$, the three bounds have integrands $k^{-5/2}$ times $1,d(k)+5/2,d(k)^2+5d(k)+15/4$. All three functions are convex on $[5/2,\infty)$ by differentiation. Their sums are bounded by the integrals from $w=5/2$. With $m=d(w)+4/9$ and
\[
 T=H_h\Gamma(5/2)\sqrt{c_3(1/3)c_3(2/3)}\,\frac{2}{3w^{3/2}},
\]
the tail bounds divided by $t^{3/2}$ are $T$ times $1,m+5/2,m^2+(4/9)^2+5m+15/4$. For clarity, the upper bounds before adding are
\[
\begin{array}{c|ccc}
 j&0&1&2\\ \hline
 k=1,2&.0061&.0142&.0186\\
 k\ge3&.0026&.0142&.0725
\end{array}
\]
Substitution uses $j^2<29/5$, $\pi<22/7$, $\log(5/2)<.917$. These estimates also give locally uniform convergence for the twice differentiated forms by the same summable bounds.
\end{proof}

\begin{proof}[Proof of Proposition~\ref{prop:acute-curvature}]\label{par:spectral-acute-close}
Let $g=.1168$, $V=.009/g$ and $J=4.9<j^2$.
The projected equation gives $\|v\|\le Vt^{3/2}$ and $0\le\gamma\le(15h^5+.009V)t^3<.003t^3$. All forms are $C^2$ locally for $s>0$. Substitute $v=C^{-1}b_s$ into $\gamma=d_s+b_s(v)$. At fixed $s$, the derivative of the left side minus the right side with respect to $\gamma$ is $1+E_s(v)/j^2$, so the real implicit function theorem applies.
Put $p=s\gamma'$, $z=s^2\gamma''$. Stationarity gives exactly
\[
 (1+E_s(v)/j^2)p=d_1+2b_1(v)+K_1(v,v)-(1+\gamma)s E_s'(v,v)/j^2.
\]
Since $\|s E_s'\|\le2$, set
\[
 P=75h^5+2(.029)V+[.13+2(1.003)/J]V^2 .
\]
Then $|p|\le Pt^3$. The derivative of $Cv=b_s$ gives
$C(sv')=b_1-[pE_s/j^2+(1+\gamma)s E_s'/j^2-K_1]v$.
Upon differentiating the stationary equation a second time at fixed tests, the terms containing $v'$ combine to $2C(sv',sv')$. After moving $z E_s(v)/j^2$ to the left, the other terms are
\[
 d_2+2b_2(v)+K_2(v,v)
 -(1+\gamma)s^2 E_s''(v,v)/j^2-2p\,sE_s'(v,v)/j^2 .
\]
Here $E_s''\ge0$. Coercivity of $C$ bounds $2C(sv',sv')$ by twice the squared dual norm of $C(sv')$ divided by $g$. Consequently $z\le Qt^3$, with
\[
 Q=300h^5+2(.095)V+(.36+4P/J)V^2
       +2[.029+(.13+(2.006+P)/J)V]^2/g .
\]
In particular $P<.0174$ and $Q<.142$, so $(Q+13P/6)/h^3<39$. This is the only consequence of the two bounds used below.
Writing $A=\log Z$, its positive series gives $(sZ)''/Z\ge48\zeta(3)s^2$ and
$sA'\le\zeta(3)s^3[12+64s^2/(1-4s^2)]<1/12$.
Twice differentiating $sZ/(1+\gamma)$, dropping the positive square and using the preceding bounds gives
\[
 F_2''(2\pi s)\ge \frac{j^2}{4\pi}\frac{Z}{1+\gamma}
         (48\zeta(3)-39)s^2>4s^2.
\]
Indeed $Z\ge1$, $1+\gamma<1.003$ and $\zeta(3)>6/5$. This proves the required bound on the entire interval.
\end{proof}

\subsection{Curvature of the first eigenvalue profile for \texorpdfstring{$\pi/3\le\alpha\le\pi/2$}{openings from 60 to 90 degrees}}\label{ang:sec:spectral-curvature}

\begin{proposition}[A lower bound for $F_2''$]
\label{prop:endpoint-A5}
Let $F_2(\alpha)=|I_\alpha|\lambda_{\mathrm{DN}}(I_\alpha)$ for the
isosceles triangle of opening $\alpha$, with Dirichlet base and Neumann
sides. With derivatives in radians,
\[
 F_2''(\alpha)>\tfrac14
 \qquad\bigl(\pi/3\le\alpha\le\pi/2\bigr).
\]
\end{proposition}
\begin{proof}
Put $t=\tan(\alpha/2)$, $z=t^{-2}=1+\delta\in[1,3]$,
and write $F_2=t\Lambda$, where $\Lambda$ is the lowest eigenvalue
of $Q_z=X+zY$ on $D=\{0<y<x<1\}$, by Lemma~\ref{lem:affine-family}.
Here $d\mu=2\,dx\,dy$, $X(v)=\int_Dv_x^2\,d\mu$, $Y(v)=\int_Dv_y^2\,d\mu$,
and $v|_{x=1}=0$. Reflection across the diagonal
and separation on the square give the unit ground state
$e=2\cos(\pi x/2)\cos(\pi y/2)$ of $Q_1$, its eigenvalue
$\lambda=\pi^2/2$, and $Q_1\ge5\lambda I$ on $e^\perp$.
Direct integration gives $b=Y(e,e)=\lambda/2-2<23/49$, using $\pi<22/7$. The one-dimensional Poincar\'e inequality on horizontal slices gives $X\ge(\lambda/2)I$. For the positive unit
ground state $u$ at $z$, simplicity and the common form domain give
$\Lambda'=Y(u,u)$ and $\Lambda''=-2S(u',u')$, where $S=Q_z-\Lambda I$
and primes on $u,\Lambda$ denote $z$ derivatives. Monotonicity,
concavity and the trial $e$ give $\lambda\le\Lambda\le\lambda+\delta b$
and $0\le\Lambda'\le b$. Differentiating $t\Lambda$ twice gives
\begin{equation}\label{ang:e:curv}
 \frac{F_2''}{t(1+t^2)/2}=\Lambda+z(z-1)\Lambda'-4R S(u',u'),
 \qquad R=z^2(1+z).
\end{equation}

At the chosen $z$ put $a=\Lambda'$, $h=z-1$ and $\nu=\Lambda-ha$.
Subtract $a(z-1)I$ from $Q_z$ and keep $a$ fixed when differentiating. The resulting lowest eigenvalue has derivative zero at the chosen $z$. In the inner product $Q_1-\nu I$ on
$e^\perp$ let $B_z$ represent $Y-aI$, $f$ represent $Y(e,\cdot)$,
and $C=(I+hB_z)^{-1}$. Projection of the eigenvalue equation gives
$u=c(e-hCf)$, $0<c\le1$. Differentiating that projection at fixed $a$
gives $S(u',u')=c^2\langle C^3 f,f\rangle$: the derivative of the
stationary eigenvalue vanishes, and $(hC)'=C^2$ there.
Testing against $e$ gives the first identity below; its derivative
uses the same stationarity:
\[
 \Lambda=\lambda+hb-h^2\langle Cf,f\rangle,\qquad
 a=b-h\langle(C+C^2)f,f\rangle .
\]
For these derivatives we fix the inner product at its value for
the chosen $z$.

Put $d=a/(5\lambda-\nu)$, $D=1/(40-h)$. The gap,
$a\le b<\lambda/10$ and $\nu\le\lambda+hb$ give $d\le D$.
Cauchy--Schwarz for $Y$ yields $f=(B_z+dI)^{1/2}\xi$, $\|\xi\|^2\le b$
by extending the functional from the range of the square root.
The horizontal Poincar\'e inequality gives $B_z\le(9\lambda/2-a)/(5\lambda-\nu)$,
since $(x-\lambda/2-a)/(x-\nu)$ is nonincreasing for $x\ge5\lambda$: indeed $\nu=Q_1(u,u)=X(u,u)+a\ge\lambda/2+a$ by horizontal Poincar\'e.
Thus each spectral value $y$ of $B_z+dI$ lies in $[0,45D]$.
For the exact multiplier $m=(1+h(y-d))^{-1}$ of $C$ we need
\[
 G=y[h^2((1+z)m+zm^2)+4z^2(1+z)m^3]-h(1+z)\le9 .
\]
Indeed inserting the two identities in \eqref{ang:e:curv} and integrating
against the spectral measure of $\xi$ then gives
$F_2''/[t(1+t^2)/2]\ge\lambda-9b$.
We use one comparison on the full spectral interval. Since $\nu=\Lambda-ha$ and $\Lambda\le\lambda+hb$, we have
$d=a/(5\lambda-\nu)\le a/(4\lambda-hb+ha)\le b/(4\lambda)<1/40$,
since the last ratio with variable $a$ increases up to $a=b$.
Put $d_*=1/40$,
$M=1+h(y-d_*)$, $H=9+h(2+h)$. Thus the actual multiplier of $C$ is bounded by $M^{-1}$. The desired bound follows from
\[
 P=HM^3-y\{h^2[(2+h)M^2+(1+h)M]+4(1+h)^2(2+h)\}\ge0 .
\]
The upper bound on $y$ gives $L=9-(40-h)y/5\ge0$.
We subtract $(1-h^2/4)L$, which equals $L$ at $h=0$ and vanishes at $h=2$. The remaining polynomial satisfies
\[
 P-(1-h^2/4)L\ \ge h^2\left(\frac{18}{5}-\frac{39}{2}y+27y^2\right).
\]
The quadratic equals $27(y-13/36)^2+19/240$.

Here are the coefficients of the comparison. Write $d=d_*$ for this
calculation, with $0<h\le2$. In
$[P-(1-h^2/4)L]/h^2$ the cubic coefficient is $9h$ and the
quadratic coefficient is
$27+h[1-d(h^2+2h+27)]\ge27$.
The constant coefficient is at least
\[
 \tfrac94+(2-27d)/h+1-6d-3dh\ge18/5,
\]
by $(1-dh)^3\ge1-3dh$. The linear coefficient is at least
\[
 \frac{34}{5h}-15-54d-(3+7d-1/20)h-3dh^2\ge-39/2 .
\]
Indeed dropping $d^2 h(2h^2+4h+27)$ gives the first inequality;
each remaining term is minimized at $h=2$. At $h=0$ one has $P=L$.

Thus $F_2''\ge[t(1+t^2)/2](\lambda-9b)>(3/8)(4-7b)>15/56>1/4$,
using $b<23/49$.
\end{proof}

\begin{proposition}[A lower bound for $F_2''$ at acute openings]\label{prop:spectral-endpoint-curvature}
For $0<\alpha\le\pi/2$, one has $F_2''(\alpha)>(\alpha/\pi)^2$.
On $\pi/3\le\alpha\le\pi/2$ one also has $F_2''(\alpha)>1/4$.
\end{proposition}
\begin{proof}
Apply Proposition~\ref{prop:acute-curvature} up to $\pi/3$.
For $\pi/3\le\alpha\le\pi/2$, Proposition~\ref{prop:endpoint-A5}
gives $F_2''(\alpha)>1/4\ge(\alpha/\pi)^2$.
\end{proof}
\subsection{The tangent at opening \texorpdfstring{$\pi/2$}{90 degrees} and its quadratic remainder}\label{ang:sec:square-tangent}

\begin{proposition}[Strict support by the tangent at $\pi/2$]\label{prop:square-tangent}
One has \(F_2(\pi/2)=\pi^2/2\), \(F_2'(\pi/2)=4\), and
\begin{equation}\tag{T90}\label{eq:T90}
 F_2(\alpha)>\frac{\pi^2}{2}+4(\alpha-\pi/2)
 \qquad(0<\alpha<\pi,\ \alpha\ne\pi/2).
\end{equation}
\end{proposition}
\begin{proof}
On $D=\{0<y<x<1\}$ with $d\mu=2\,dx\,dy$ and zero trace on $x=1$, write
$X(v)=\int_Dv_x^2\,d\mu$ and $Y(v)=\int_Dv_y^2\,d\mu$; by Lemma~\ref{lem:affine-family}
$F_2(\alpha)$ is the least value of $K_r=rX+r^{-1}Y$ on unit $L^2(d\mu)$ tests, $r=\tan(\alpha/2)$.
Put $P=\pi^2/4$. At $\alpha=\pi/2$, where $r=1$, the unit ground state $e$ of $K_1=X+Y$ has $X(e)=P+2$, $Y(e)=P-2$,
so $F_2(\pi/2)=2P$ and $F'_2(\pi/2)=4$. For $\alpha<\pi/2$ the claim follows from the
strict convexity of $F_2$ on $(0,\pi/2]$ in Proposition~\ref{prop:spectral-endpoint-curvature}.
For $r>1$ put $u=(r-1)/(r+1)$. Since $\alpha-\pi/2=2\arctan u$, the tangent at $\pi/2$ has value $2P+8\arctan u\le T:=2P+8u$.
We have $T\le rP+3P/r$: the equivalent bound $8/P\le r+4/(r-1)-3/r$ follows since
this right side exceeds $(r-1)+3/(r-1)\ge2\sqrt3$.

Let $A=X-PI$, $B=Y$. Both forms are nonnegative. After $2P$ the next eigenvalue of $X+Y$ is $10P$, as one sees by reflecting in the diagonal to a square. Hence $A+B\ge9PI$ on functions $v$ orthogonal to $e$ in $L^2(d\mu)$, and the ground-state equation gives $B(e,v)=-A(e,v)=c$.
Therefore
\[
 (K_r-TI)(v)\ge[r-1/(3r)]A(v)+2B(v)/(3r)\ge c^2L,\qquad
 L=\tfrac12\Bigl(r+\frac{6-P}{3(P-2)r}\Bigr),
\]
by Cauchy--Schwarz and $A(e)=2$, $B(e)=P-2$. On $ae+v$ the cross term is
$2ac(r^{-1}-r)$ and the coefficient of $a^2$ is $(r-r^{-1})^2D$, $D=r(P+2u)/(r+1)^2$.
For $a\ne0$, completing the square in $c$ gives strict positivity if $DL>1$. For $a=0$ and $v\ne0$, strict positivity follows from the positive coefficients of $A(v)$ and $B(v)$ and the lower bound on $A+B$. Now $(6-P)/(3(P-2))\ge(P+1)/(P-1)$ because
$2<P<5/2$. With $w=2u/P<u$ this gives
\[
 DL\ge\frac{P^2}{4(P-1)}(1+w)(1-w+u^2)\ge(1+w)(1-w+w^2)=1+w^3>1.
\]
Attainment of the first eigenvalue gives $F_2(\alpha)>T$.
\end{proof}

\paragraph{The quadratic tangent remainder.}
Put $F=F_2$, $a=2\pi/n$, $b=\pi/2$ and $d=\alpha-a$, with $n\ge5$.
For $0<\alpha\le b$ integrate $F''(t)>(t/\pi)^2$ from
Proposition~\ref{prop:spectral-endpoint-curvature}. Since $a+ud\ge(1-u)a$,
\[
 F(\alpha)-\ell_a(\alpha)
 =d^2\int_0^1(1-u)F''(a+ud)\,du
 \ge 4(d/n)^2\int_0^1(1-u)^3\,du=(d/n)^2 .
\]
Continuity includes $\alpha=0$. Beyond $b$ the tangent at $\pi/2$ lies below $F$
by \eqref{eq:T90}. The difference $\ell_b(\alpha)-\ell_a(\alpha)-c(d/n)^2$
is concave in $\alpha$ for $c>0$, so it suffices to compare the two tangents
at $b$ and $\pi$. The bound at $b$ was just proved. At $\pi$ use
$\partial_t\ell_t(\pi)=(\pi-t)F''(t)$ from $a$ to $b$. This gives
\[
 \frac{n^2[\ell_b(\pi)-\ell_a(\pi)]}{(\pi-a)^2}
 \ge 4\int_0^{(n-4)/(2(n-2))}
              (1-u)\big(1+(n-2)u/2\big)^2\,du .
\]
For $n\ge6$ the integral covers $[0,1/4]$ and the integrand there is at
least $(1-u)(1+2u)^2\ge1$. For $n=5$ it covers $[0,1/6]$, where
$(1-u)(1+3u/2)^2\ge(1-u)(1+3u)\ge1+3u/2$.
Integration gives $c=1$ and $c=3/4$, respectively, in \eqref{eq:spectral-quadratic-gap}.
For $n=3,4$ use the unweighted estimates in Proposition~\ref{S:prop:wide-linear} and on p.~\pageref{par:spectral-obtuse}.
At obtuse openings they give $\mathcal C_2\ge1/1369$; both $k_\alpha$ and $F_2(\alpha)$ are at least one there by \eqref{eq:Fq-bounds}.
Together with the acute estimate this gives $F_2''(x)\ge (x/\pi)^2/1369$ on $(0,\pi)$.
Integration along $x=a+ud\ge(1-u)a$ as above gives the stated coefficient.
\phantomsection\label{par:spectral-route-end}

\endgroup

\section{Quadratic stability}\label{sec:stability}

\subsection{Prescribed area and optimizing sequences}

Write \(A=|P|\) and \(d(P)=\inf_g d_H(\overline P,g\overline R_n)\), where \(g\)
runs over plane isometries and \(R_n\) has area \(A\). The distance
\(d(P)/\sqrt A\) is invariant under common dilation.

\begin{theorem}[Prescribed area]\label{thm:stability}
Let $n\ge3$ and let $P$ be a simple polygon with at most $n$ sides. One may take
$c_n=10^{-18}$ (not optimized) for $n=3,4$ and $c_n=(18432\pi^6n^5)^{-1}$ for $n\ge5$ in the inequality
\[
 \min\left\{\frac{T(R_n)-T(P)}{A^2},
                A\bigl(\lambda_1(P)-\lambda_1(R_n)\bigr)\right\}
 \ge c_n\min\{1,d(P)^2/A\}.
\]
The exponent two is optimal. For convex polygons the normalized eigenvalue deficit
\(A(\lambda_1(P)-\lambda_1(R_n))\) is at least an explicit
positive constant times \(d(P)^2/A\).
\end{theorem}

For instance, if either normalized deficit is at most
\(\varepsilon<c_n\), then \(d(P)/\sqrt A\le\sqrt{\varepsilon/c_n}\).
For $n\ge5$ we control how far each cell opening $\alpha_i$ is from $a=2\pi/n$ and how far its two base angles are from one another. For the first eigenvalue, we multiply the pointwise bound \eqref{eq:spectral-quadratic-gap} by the mean $r_i$ of $u^2$ on each cell and sum. For torsion at unit area, we first apply \eqref{eq:stability-angle-torsion} to $v_i=2\pi|Q_i|$, which sum to $2\pi$, and then bound $\alpha_i-v_i$. The strip comparison controls the difference between the base angles. Sections~\ref{sec:tangent-gaps}--\ref{sec:stability-cell-geometry} combine these bounds with $d(P)\le16ne$ to obtain $c_n$. For $n=3,4$ see p.~\pageref{par:small-n-area}.

For $n=3,4$ the cell estimates on p.~\pageref{par:small-n-area}, and for $n\ge5$ the small-deficit argument in Section~\ref{sec:small-deficit-convexity} force an optimizing sequence for either functional to consist eventually of convex polygons with \(n\) effective sides. Theorem~\ref{thm:stability}
gives Hausdorff convergence to \(R_n\) up to rigid motions. To see convergence of their vertices in
cyclic order, pass from any subsequence to one along which all
ordered vertices converge. Their convex hull is \(R_n\), so the \(n\) distinct extreme points
of \(R_n\) are precisely the limiting vertices. Use a consistent
orientation of the cyclic order and choose the starting index to match
the regular polygon. This gives convergence of the vertices. In particular, no side length tends to zero and no interior angle tends to $\pi$.

The minimum with $1$ is necessary for torsion when $n\ge3$ and for the first eigenvalue when $n\ge4$. At unit area the torsion deficit is bounded above by $T(R_n)$. Right triangles of unit area with one leg tending to infinity have $d(P)\to\infty$ and are admissible for every $n\ge3$.

For the eigenvalue, take the quadrilateral with vertices
\[
 (0,-L),\quad(L^{-2},0),\quad(1,0),\quad(0,1),\qquad L>2.
\]
Its area is $(1+L^{-1})/2$, its diameter tends to infinity, and it contains the fixed triangle with vertices $(1/3,1/4),(1/2,1/4),(1/3,1/3)$. After normalization to unit area its first eigenvalue stays bounded by domain monotonicity, while $d(P)\to\infty$.

\subsection{Prescribed perimeter}\label{sec:perimeter}

Quantitative polygonal isoperimetry controls the side and radius
variances and the loss of convexity \cite{indrei2015}. For fixed \(n\), the convex
Hausdorff estimate also follows from \cite[Proposition~2.1]{cm2016}.
We use a direct
Hausdorff estimate with an explicit constant for the closed regions,
then transfer it to both functionals by their scaling laws.

Fix \(n\ge3\) and \(L>0\). Let \(R=R_{n,L}\) be the regular \(n\)-gon
of perimeter \(L\), and let \(A_*=|R|=L^2/[4n\tan(\pi/n)]\). For a simple polygon \(P\) with
the same perimeter and at most \(n\) sides, write
\(\delta=1-|P|/A_*\) for the relative area deficit of $P$. Here \(d(P)\) is measured relative to \(R\).

\begin{theorem}[Quadratic stability at fixed perimeter]\label{thm:perimeter-stability}
With these conventions,
\[
 \frac{T(P)}{T(R)}\le (1-\delta)^2,\qquad
 \frac{\lambda_1(P)}{\lambda_1(R)}\ge (1-\delta)^{-1},
\]
and
\[
 \min\left\{1-\frac{T(P)}{T(R)},\frac{\lambda_1(P)}{\lambda_1(R)}-1\right\}
 \ge \frac{d(P)^2}{4L^2}.
\]
The constant is uniform for all \(n\ge3\).
The regular \(n\)-gon is the unique optimizer for either functional
up to isometry. The power two is optimal.
\end{theorem}

The geometric assertion used here is
\begin{equation}\label{eq:perimeter-hausdorff}
 |P|\le A_*,\qquad d(P)\le 2 L\sqrt{\delta}.
\end{equation}
Section~\ref{sec:perimeter-geometry} proves this estimate, including
for polygons with reflex vertices. To see how area loss enters the two
functionals, put \(r=\sqrt{|P|/A_*}=\sqrt{1-\delta}\).
The regular polygon \(rR\) has exactly the area of \(P\). Thus
Theorem~\ref{thm:main} and scaling give
\[
 T(P)\le T(rR)=r^4T(R),\qquad
 \lambda_1(P)\ge\lambda_1(rR)=r^{-2}\lambda_1(R).
\] Each relative deficit is at least
\(\delta\), which by \eqref{eq:perimeter-hausdorff} is at least
\(d(P)^2/(4L^2)\). If either functional attains its value on \(R\),
then \(\delta=0\), and Theorem~\ref{thm:main} gives \(P\cong R\).

\paragraph{Sharpness in both normalizations.}
Let \(P_\varepsilon=c_\varepsilon D_\varepsilon R\), where
\(D_\varepsilon=\operatorname{diag}(e^\varepsilon,e^{-\varepsilon})\)
for \(\varepsilon>0\). Choose the first axis parallel to a diametral segment of \(R\).
At fixed area take \(c_\varepsilon=1\). At fixed perimeter the
rescaling factor satisfies \(c_\varepsilon=1+O(\varepsilon^2)\),
since the sum of the first variations of the side lengths vanishes
by symmetry. The chosen segment stretches by $e^\varepsilon$, the
operator norm of $D_\varepsilon$, so
$\operatorname{diam}(P_\varepsilon)=c_\varepsilon e^\varepsilon
\operatorname{diam}(R)$. For any two compact sets the difference of their
diameters is at most twice their Hausdorff distance. Since rigid motions
preserve diameter, this gives
\[
 d(P_\varepsilon)\ge \tfrac12\operatorname{diam}(R)
 (c_\varepsilon e^\varepsilon-1)
 =\tfrac12\operatorname{diam}(R)\,\varepsilon+O(\varepsilon^2).
\]

For $1\le q\le2$ let $u$ be the normalized positive Sobolev minimizer on $R$.
Uniqueness and rotational invariance make $\int_R\nabla u\,\nabla u^{\mathsf T}$ a scalar matrix.
Since $\det D_\varepsilon=1$, the stretch preserves $\int |u|^q$. The energy of
$u\circ D_\varepsilon^{-1}$ on the stretched polygon is the energy of $u$ on $R$ multiplied by $\cosh(2\varepsilon)$. Dilation therefore gives
\[
 \lambda_q^{\rm S}(P_\varepsilon)
 \le c_\varepsilon^{-4/q}\cosh(2\varepsilon)\lambda_q^{\rm S}(R).
\]
For $q=1$ take reciprocals to obtain
$T(P_\varepsilon)\ge c_\varepsilon^4T(R)/\cosh(2\varepsilon)$.
Thus every deficit is $O(\varepsilon^2)$ while $d(P_\varepsilon)$ is bounded below by a
positive multiple of $\varepsilon$. This proves sharpness for each exponent at fixed area and fixed perimeter.

\subsection{Geometric estimates for stability}\label{sec:area-quantitative}

For $n\ge3$ set $|P|=1$, $a=2\pi/n$, $t=\tan(a/2)$ and $\eta=a/8$.
Write $F=F_2$ and $g(\alpha)=\alpha^2/F_1(\alpha)$. The parameter
$\delta_i=(\beta_i-\beta_i')/\pi$ of a convex mixed cell measures the difference of its two base angles.

\subsubsection{The angular errors}\label{sec:tangent-gaps}

Append zeros to the list of openings until it has length $n$. If their sum is $2\pi$, then $\sum_i(\alpha_i-a)=0$, so the linear terms in the tangent bounds cancel. Since $F(0)=g(0)=0$, \eqref{eq:spectral-quadratic-gap} gives
\begin{equation}\label{eq:stability-angle-lambda}
 \sum_i F(\alpha_i)-nF(a)\ge K_\lambda\sum_i (\alpha_i-a)^2,
 \qquad K_\lambda=\begin{cases}1/(1369n^2),&n=3,4,\\1/(2n^2),&n\ge5.\end{cases}
\end{equation}
For torsion take any $n\ge3$. If $|\alpha-a|\le a/8$, then
$2a/5<7a/16\le v=\tan(\alpha/2)\le\tan(3\pi/8)=1+\sqrt2$.
For $v\le1$, (TC) gives $-g''\ge5(2a/5)^3/60=2a^3/375$.
For $v\ge1$, the denominator in (TC) satisfies
$(v^4+6v^2+3)/v^3=v+6/v+3/v^3<12$.
Thus $-g''>5/72>2a^3/375$, using $a\le2\pi/3<21/10$.
In particular $-g''\ge a^3/256$ throughout this neighborhood for every $n\ge3$.
Integration on $[a-\eta,a+\eta]$ gives a quadratic tangent remainder there; outside that interval the tangent gap cannot decrease as one moves away from $a$. It follows that
\begin{equation}\label{eq:stability-angle-torsion}
 \frac{ng(a)-\sum_i g(\alpha_i)}{4\pi^2}\ge
 K_T\sum_i\min\{(\alpha_i-a)^2,\eta^2\},
 \qquad K_T=\frac{a^3}{2048\pi^2}=\frac{\pi}{256n^3}.
\end{equation}

\subsubsection{The quadratic loss in the base-angle difference}\label{sec:stability-tilt}

Use the strip $S=\mathbb R\times(0,\pi)$, with Dirichlet trace on $y=0$.
For a probability density $h$ on $S$, let
$M_q(h)=\sup_{\|\nabla u\|_2\le1}\int_S h|u|^q$, with this trace.
The comparison below uses the binomial mixture \eqref{eq:mixture} and the rearrangement inequalities (M1)--(M3).
Write $0<C<1$, $|\delta|<C$, $0<\kappa\le1$ and
$p_{\delta,\kappa}\propto e^{\delta x}(\cosh x+\kappa\cos y)^{-C}$, normalized on $S$.
For $b>0$, let $\rho_b(x)\propto(1+\cosh x)^{-b}$ be normalized on $\mathbb R$. Put
\[
 p=p_{0,1},\quad h_0(x,y)=\rho_C(x)/\pi,\quad
 h_1(x,y)=\rho_{C+1}(x)(1-\cos y)/\pi .
\]
Let $P_k$ be the weights of the centered mixture for $p$, so that $h_0,h_1$ are its first two product densities.

\paragraph{Effect of unequal base angles.}
For $q=1,2$,
\begin{equation}\label{sr:mixturegain}
 M_q(p)-M_q(p_{\delta,\kappa})\ \ge\
 \frac{P_0(1-P_0)}{C^2}\,b_q M_q(h_0)\,\delta^2,
 \qquad b_1=\frac{2}{\pi^2},\quad b_2=\frac12.
\end{equation}
In particular, if $n\ge3$ and $\max\{2/5,1-3/n\}\le C\le1-1/n$, the right side is at least \(k\delta^2\), uniformly in \(\kappa\), where \(k=b_1/(12\pi n)\).

\paragraph{A discrete comparison.}
Let $w_k$ be the weights of $p_{\delta,\kappa}$ and $r_k=w_k/P_k$. By (M3), $r_k$ is nonincreasing and $r_0\ge1$. The beta recurrence and $\int_0^\pi(1-\kappa\cos y)\,dy=\pi$ give $r_1/r_0=1-\delta^2/C^2$. Write $s=\sum_{k\ge1}w_k/(1-P_0)$. Then $r_0-s\ge r_0\delta^2/C^2$ and $1-s=P_0(r_0-s)\ge P_0\delta^2/C^2=:e$. Each tail ratio $\sum_{k\ge j}w_k/\sum_{k\ge j}P_k$, $j\ge1$, is the $P$-average of the decreasing sequence $r_k$ over a terminal interval. It is therefore at most $s$. Summation by parts gives, for every nonnegative nondecreasing sequence $A_k$,
\[
 \sum_k w_k A_k\le e A_0+(1-e)\sum_k P_k A_k.
\]
Here and below infinite sums follow by monotone truncation. For a unit-energy test, apply the rearrangements of Section~\ref{sec:strip} to obtain $v$, with energy at most one. Let $A_k$ be the integral of $|v|^q$ against the $k$th centered product density. These integrals are nondecreasing by (M1). Centering the horizontal factors and using (M2) gives $\int p_{\delta,\kappa}|u|^q\le\sum_k w_k A_k$. Hence
\[
 M_q(p_{\delta,\kappa})\le e M_q(h_0)+(1-e)M_q(p).
\]
It remains to estimate the gain from the first two centered densities.

\paragraph{Comparison of the first two densities for $q=1,2$.}
For $q=2$, a maximizer for $h_0$ has the form $u(x,y)=f(x)\sin(y/2)$ with $f\ge0$ even and nonincreasing on $[0,\infty)$. Indeed, expansion in the vertical Dirichlet--Neumann modes diagonalizes the energy and numerator. The lowest mode has the smallest energy for each horizontal coefficient. A one-dimensional maximizer exists by local compactness and decay of $\rho_C$, and symmetric decreasing rearrangement gives the stated order. Passing from $\rho_C$ to $\rho_{C+1}$ increases the mean of $f^2$ because the likelihood ratio decreases with $|x|$. The mean of $\sin^2(y/2)$ is $1/2$ under $dy/\pi$ and $3/4$ under $(1-\cos y)\,dy/\pi$. Consequently, for the energy-normalized maximizer, $\int(h_1-h_0)u^2\ge M_2(h_0)/2$.

For \(q=1\), put \(\rho=\rho_C\) and solve \(-\Delta w=h_0=\rho/\pi\)
on the mixed strip. Then \(m=M_1(h_0)=\|\nabla w\|_2\), and \(w/m\)
is a maximizer. The rearrangements in \eqref{eq:rearranged-test}
decrease the energy of \(w\) and increase its integral against \(h_0\).
Uniqueness in the quadratic variational principle therefore gives
\(w\) nondecreasing in \(y\) and nonincreasing in \(|x|\).
As in the case $q=2$, replacing $\rho_C$ by $\rho_{C+1}$ in the horizontal variable can only increase the integral against \(1-\cos y\), so it suffices to
bound \(\int_S h_0(-\cos y)w\) from below.

Write \(\gamma_j=j+1/2\), \(\phi_j(y)=\sqrt{2/\pi}\sin(\gamma_j y)\),
and \(D_j=\langle\rho,(-\partial_x^2+\gamma_j^2)^{-1}\rho\rangle\).
The \(D_j\) decrease with \(j\). Expansion in the \(\phi_j\) gives
\[
 \pi^2m^2=\sum_{j\ge0}\langle1,\phi_j\rangle^2D_j\le\pi D_0,\qquad
 \pi^2\int_S h_0(-\cos y)w=\sum_{j\ge0}d_j D_j,
 \quad d_j=\frac{-2}{\pi(\gamma_j^2-1)} .
\]
Indeed \(d_j=\langle1,\phi_j\rangle\langle-\cos y,\phi_j\rangle\).
Their sum is \(\langle1,-\cos y\rangle=0\), and only
\(d_0=8/(3\pi)\) is positive. The last sum is therefore at least
\(d_0(D_0-D_1)\).

The horizontal variation of \(\rho\) is bounded by
\(\|\rho'\|_2^2=-\frac12\int(\log\rho)''\rho^2\le\|\rho\|_2^2/4\),
since \(-(\log\rho)''=C/(2\cosh^2(x/2))\le1/2\).
Testing the first resolvent with \(\rho\) gives
\[
 D_0\ge\frac{\|\rho\|_2^4}{\|\rho'\|_2^2+\|\rho\|_2^2/4}
 \ge2\|\rho\|_2^2,\qquad
 D_1\le\frac49\|\rho\|_2^2\le\frac14D_0 .
\]
It follows that \(\int_S h_0(-\cos y)w\ge
(3d_0/(4\pi))m^2=(2/\pi^2)m^2\). Division by \(m\) gives the
bound for $q=1$ with \(b_1=2/\pi^2\).
For $q=1$ or $q=2$, let $u$ denote this maximizer for $h_0$. By (M1), the integral of $|u|^q$ against each density with index $j\ge1$ is at least its integral against $h_1$. Using this maximizer as a test for $p$ gives \(M_q(p)-M_q(h_0)\ge(1-P_0)b_qM_q(h_0)\), which proves \eqref{sr:mixturegain}.

\paragraph{A spherical coordinate on the strip.}
Set \(\psi(x,y)=\sin(y/2)/\cosh(x/2)\). Under the product density \(h_0\),
the variables \(U=\tanh(|x|/2)\) and \(\theta=y/2\) are independent.
The latter is uniform on \([0,\pi/2]\), and the former has density
proportional to \((1-u^2)^{C-1}\) on \([0,1]\). Integration of
\((u(1-u^2)^C)'\) gives \(\mathbb E U^2=1/(2C+1)\), hence
\[
 \mathbb E\psi^2=\mathbb E(1-U^2)\,\mathbb E\sin^2\theta
 =\frac{C}{2C+1}.
\]
The constant term of the centered mixture has mass
\(P_0=(\mathbb E(1-\psi^2)^{-C})^{-1}\), since
\(\cosh x+\cos y=(1+\cosh x)(1-\psi^2)\).
For \(C\le1\) the density of \(U\) is increasing, so \(U\)
stochastically dominates a uniform variable on \([0,1]\).
Replacing \(U\) by that uniform variable increases \(\psi\)
in stochastic order. In this case
\((U,\sqrt{1-U^2}\cos\theta,\sqrt{1-U^2}\sin\theta)\) is uniform on the
positive octant of the unit sphere, so its last coordinate is uniform
on \([0,1]\) as well. The tangent bound \((1-s)^{-C}\ge1+Cs\)
and the spherical comparison give
\[
 1+\frac{C^2}{2C+1}\le P_0^{-1}
 \le\int_0^1(1-v^2)^{-C}\,dv\le\frac1{1-C}.
\]
Thus \(P_0(1-P_0)/C^2\ge(1-C)/(1+C)^2\ge1/(4n)\) in the
stated range.

With \(f(x)=1/\cosh(x/2)\), the same \(\psi\) satisfies \(-\Delta\psi=f^2\psi/2\) and the mixed boundary conditions. The substitution \(u=\tanh(x/2)\) gives
\(\int_{\mathbb R}f^2=4\) and \(\int_{\mathbb R}(f')^2=1/3\).
Consequently \(\|\nabla\psi\|_2^2=(\pi/2)(1/3+1)=2\pi/3\).
Since \(0\le\psi\le1\), \(C\ge2/5\), and this energy is greater than one,
the moments for $q=1$ and for $q=2$ both satisfy
\[
 M_q(h_0)\ge
 \frac{\mathbb E\psi^2}{\|\nabla\psi\|_2^2}
 \ge\frac{2/9}{2\pi/3}=\frac1{3\pi}=:\mu.
\]
Put \(b=2/\pi^2\) and \(k=b\mu/(4n)\). Since \(b_q\ge b\),
\eqref{sr:mixturegain} has coefficient at least \(k\).
Translation and averaging in the conformal strip
transfer the estimate to every convex mixed cell, by Section~\ref{sec:strip}.

\paragraph{Application to a mixed cell.}
For $n\ge5$, if $|\alpha-a|\le a/4$, then $1-3/n\le C=1-\alpha/\pi\le1-1/n$. The uniform strip bound and $M_1(p)\ge M_1(h_0)\ge\mu$ give
\begin{align}
 \frac{T_{\rm DN}(I_\alpha)}{|I_\alpha|^2}-\frac{T_{\rm DN}(Q)}{|Q|^2}
     &\ge k\mu\,\delta^2,\label{eq:tilt-torsion}\\
 |Q|\lambda_{\rm DN}(Q)-F(\alpha)&\ge k F(\alpha)^2\delta^2.\label{eq:tilt-lambda}
\end{align}
Indeed the torsion quotient squares the two moments and the eigenvalue quotient takes their reciprocals. The difference of the reciprocals is at least $F(\alpha)^2$ times the difference of the moments.

\subsubsection{Small deficits and convexity}\label{sec:small-deficit-convexity}
\begin{lemma}[An angle bound for a reflex corner]\label{lem:partition-rigidity}
Suppose the partition of Lemma~\ref{lem:mass-cap-partition} of a polygon
with at most $n$ sides has $n$ cells. Let $\beta_*$ be the largest of
their base angles and $E=\sum_i\alpha_i-2\pi$. If a final base is a
proper subsegment of its original side, then $E\ge\pi-\beta_*$.
If $P$ is nonconvex, then $E\ge\pi-2\beta_*$.
\end{lemma}
\begin{proof}
There are exactly $n$ original labels, and each survives in a distinct
cell. First consider any intermediate cap $C$ whose $k$ current bases
survive in full in its final cells. For each of the $k+1$ vertices $v$
of its Dirichlet chain, including its endpoints, write $\sigma_v$ for
the angle of $C$ and $b_v$ for the sum of the final base angles at $v$.
These angles occupy disjoint sectors inside $C$, up to their rays.
Subtracting the formula
$\Theta(C)=k\pi-\sum_v\sigma_v$ from the sum of the final openings gives
\[
     \sum_{Q_i\subset C}\alpha_i-\Theta(C)
       =\sum_v(\sigma_v-b_v)\ge0 .
\]
In particular the left side is at least the unused angle at either endpoint of the chain, where $b_v\le\beta_*$.

If a base has been shortened, choose the last shortening in the
construction. The initial triangles and the cuts of Lemma~\ref{lem:cap-sweep} cut the Dirichlet
chain only at vertices or release whole edges. Thus the shortening is
an extension hitting an open side at $X$. In the child releasing its
portion of that side let the angle at $X$ be $x$. Before this portion of the side is released, the two openings sum to the original opening plus $\pi$. Releasing it reduces their sum by $\pi-x$. The net increase is $x$. Just before this cut the current
caps (including finished cells) have total opening at least $2\pi$:
they are admissible and all $n$ labels survive. Writing $B=\sum_i b_i$ for the total of the original side costs, the prescribed minimum openings sum to $B+(2\pi-B)\nu(P)=2\pi$. Thus the openings just after the cut sum
to at least $2\pi+x$. From this stage every current base survives
in full. In the retaining child $X$ is a terminal endpoint with
$\sigma_X=\pi-x$. The displayed identity over all current caps gives
$E\ge x+(\pi-x-\beta_*)=\pi-\beta_*$.

If all bases are full original sides, let $\sigma_v$ be the interior
angle of $P$ and $b_v$ the sum of the two final base angles at $v$.
The polygon angle sum gives $E=\sum_v(\sigma_v-b_v)$, again with
nonnegative summands. A reflex vertex contributes more than
$\pi-2\beta_*$. This proves both assertions.
\end{proof}

For the following explicit deficit estimates assume $n\ge5$; the cases $n=3,4$ are treated on p.~\pageref{par:small-n-area}. At area one write $\Delta_T=T(R_n)-T(P)$, $\Delta_\lambda=\lambda_1(P)-\lambda_1(R_n)$ and let $\Delta$ denote whichever of these two deficits is being considered. Set
\[
 L_0=\frac{k\mu}{4\pi^2 n^2}=\frac1{72\pi^6n^3},\qquad
 c=\frac{L_0}{(16n)^2},\qquad e=\sqrt{\Delta/L_0}.
\]
For $0<\Delta\le c$ we have $e\le1/(16n)$; equality at zero was characterized above.
Use Lemma~\ref{lem:mass-cap-partition} on $P$, of any shape, with the measures below.
Write $p$ for the number of positive-area cells, $S_i=|Q_i|$ and $h=2\pi-nb$. Since $P$ has at most $n$ sides, the partition guarantees $\alpha_i\ge b+h\nu(Q_i)$.
In both cases we show that all $n$ cells have positive area and that
\begin{equation}\label{eq:one-cap-data}
 0\le E:=\sum_i\alpha_i-2\pi\le e/100,\quad
 |\alpha_i-a|\le e/10,\quad \pi|\delta_i|\le e,\quad S_i\ge1/n-e/4 .
\end{equation}

For torsion set $b=0$ and $\mathrm d\nu=\mathrm dx$. Put
$v_i=2\pi S_i\le\alpha_i$, and pad this list with zeros to length
$n$. Concavity of $g$ and $g(0)=0$ give
$g(\alpha_i)/\alpha_i\le g(v_i)/v_i$. The mixed torsion bounds on
the positive-area cells therefore give
\[
 T(P)\le\sum_{i=1}^p T_{\rm DN}(Q_i)
 \le\sum_{i=1}^p S_i^2\,\frac{g(\alpha_i)}{\alpha_i^2}
 \le\frac{1}{4\pi^2}\sum_{i=1}^p g(v_i)\frac{v_i}{\alpha_i}.
\]
The first inequality is obtained by restricting the torsion test to
each cell in the variational formula $T(P)=\sup_{v\in H^1_0(P)}(2\int_Pv-\int_P|\nabla v|^2)$. Since $\sum_i v_i=2\pi$,
\eqref{eq:stability-angle-torsion} yields
\[
 \Delta_T\ge K_T\sum_{i=1}^n\min\{(v_i-a)^2,\eta^2\}
       +\frac1{4\pi^2}\sum_{i=1}^p
          g(v_i)\frac{\alpha_i-v_i}{\alpha_i}.
\]
By \eqref{eq:tilt-torsion}, we may add to the right side the term
$k\mu S_i^2\delta_i^2$ for each $|\alpha_i-a|\le a/4$.
As $\sqrt{\Delta/K_T}=e\sqrt{256/(72\pi^7)}<e/25<\eta$,
all $n$ cells have positive area and $|v_i-a|<e/25$.
Thus $S_i=v_i/(2\pi)\ge1/n-e/4\ge1/(2n)$ and $v_i<\pi/2$.
By \eqref{eq:Fq-bounds}, $g(v_i)\ge v_i/9>1/(2n)$.
Consequently $E=\sum_i(\alpha_i-v_i)\le8\pi^3 n\Delta<e/100$.
All $|\alpha_i-a|<e/10<a/4$, and the bound on the base-angle differences gives
$\pi|\delta_i|\le2\pi n\sqrt{\Delta/(k\mu)}=e$.

For the eigenvalue normalize the positive first eigenfunction by
$\int_Pu^2=1$ and use $\mathrm d\nu=(1+u^2)\,\mathrm dx/2$.
Write $L=nF_2(a)=\lambda_1(R_n)$, $s=F'_2(a)$ and
$b=a-F_2(a)/s$, so $L=s(2\pi-nb)$.
The slope bound in \eqref{eq:line-constants} gives
$0<b\le a-\sin a<a/2$. Put
$r_i=S_i^{-1}\int_{Q_i}u^2$ and $\nu_i=S_i(1+r_i)/2$.
Restriction gives
$\int_{Q_i}|\nabla u|^2\ge r_iF_2(\alpha_i)$.
The supporting line satisfies $s(\alpha_i-b)\ge L\nu_i$.
As $\sum_i S_i=\sum_i S_ir_i=1$, the identity
$\sum_i r_i\nu_i=1+\tfrac12\sum_i S_i(r_i-1)^2$
and \eqref{eq:spectral-quadratic-gap} imply
\[
 \Delta_\lambda\ge \frac L2\sum_i S_i(r_i-1)^2
       +K_\lambda\sum_i r_i(\alpha_i-a)^2.
\]
The additional local term from \eqref{eq:tilt-lambda} is
$kr_iF_2(\alpha_i)^2\delta_i^2$ wherever $|\alpha_i-a|\le a/4$.

Retain the $j$ cells with $r_i\ge1/2$. The discarded area and also its $\nu$-mass are at most $8\Delta/L$. The preceding angular bound gives
\[
 \left|\sum_{\rm retained}(\alpha_i-a)\right|
 \le\sum_{\rm retained}|\alpha_i-a|
 \le\sqrt{2n\Delta/K_\lambda}=\frac{e}{\sqrt{18}\,\pi^3}<e/100.
\]
The partition bounds give a lower bound $(n-j)(a-b)-8h\Delta/L$ for the signed sum. Here $a-b>\pi/n$, $L\ge n(\pi^2/4)\tan(\pi/n)>1$ and $8h\Delta/L\le16\pi c<1/(16n)$. Hence $j=n$, $0\le E<e/100$, and each $|\alpha_i-a|<e/100$.
The term from \eqref{eq:tilt-lambda} and $F_2(\alpha_i)\ge1/n$ give $\pi|\delta_i|\le\pi n\sqrt{2\Delta/k}=e\sqrt{\mu/2}<e$.
Finally the nonnegative residuals $\alpha_i-b-h\nu_i$ sum to $E$, hence
$|\nu_i-1/n|\le(|\alpha_i-a|+E)/h<e/100$, using $h>\pi$.
The variance estimate gives $|\nu_i-S_i|=S_i|r_i-1|/2\le\sqrt{\Delta/(2L)}<e/100$.
This proves \eqref{eq:one-cap-data}.

\subsubsection{Geometry of the cells}\label{sec:stability-cell-geometry}

The following implication holds for every $n\ge3$. Suppose $P$ has area one,
all $n$ cells of the partition have positive area, and \eqref{eq:one-cap-data}
holds for some $0<e\le1/(16n)$. We show $d(P)\le16ne$.

Each actual base angle differs from $\beta_{\mathrm{reg}}=\pi/2-a/2$ by at most $11e/20$.
For the three weighted side cells of a triangle the bases are full and the polygon is convex. For the general partition, $a>2e$ and Lemma~\ref{lem:partition-rigidity} exclude both a shortened base and a reflex vertex. At each vertex, the unused sector between the two incident base angles is nonnegative; their total is $E$. The exterior turns $\omega_i$ consequently satisfy $|\omega_i-a|\le11e/10+E<2e$.

Write $r=(nt)^{-1/2}$ and $s=2rt$ for the regular inradius and side length. With $B=n/4$ we have $S_i\ge(1-Be)/n$. For $\rho=r(1-He)$, let $W$ be the polygon circumscribed about a circle of radius $\rho$ with the outward normals of $P$. We choose $H$ to make its side lengths at most those of $P$. For $n\ge5$ put $H=B+4/t$ and $E_*=\min\{a/8,t/4,(2H)^{-1}\}$.
The bound $nt>3$ gives $H<19n/12$. Hence $e\le E_*$ and $\rho\ge r/2$.

Extend the two sides of each cell that meet its base. Together with the base, they enclose the convex cell in a triangle of area
$L_i^2/[2(\cot\beta_i+\cot\beta'_i)]$, where $L_i$ is the full side length.
For $n\ge5$ the base angles are in $[\pi/4,\pi/2)$ and $|(\cot)'|\le2$ there. Since $s^2=4t/n$,
\[
 L_i^2\ge s^2(1-Be)(1-2e/t),\qquad L_i\ge s[1-(B+2/t)e].
\]
As $|\omega_i-a|\le2e$, $0<\omega_i\le\pi/2$ and $(\tan(\omega/2))'\le1$,
the side of $W$ with the same outward normal as side $i$ of $P$ has length at most $2\rho(t+2e)\le s[1-(B+2/t)e]\le L_i$.

For $n=3,4$, take $H=3$. Put $U=21/10$ for $n=3$ and $U=3/2$ for $n=4$.
With $y=a/2$ and $e\le1/(16n)$ one has $\sec(y+e)\le U$, since $\cos(y+e)\ge\cos y-e$.
The base-angle errors are at most $.55e$ and $|\omega_i-a|\le1.11e$.
Thus $|(\cot)'|\le U^2$ on the base-angle intervals, and
$|\tan(\omega_i/2)-t|\le .56U^2e$. The enclosing triangles used above give
\[
 L_i\ge s\sqrt{(1-Be)(1-.55U^2 e/t)}
 \ge s(1-me),\qquad m=\max\{B,.55U^2/t\}.
\]
Since $t=\sqrt3$ or $1$, respectively, $m+.56U^2/t<3$. Hence the corresponding side of $W$ has length at most $2\rho(t+.56U^2e)\le s(1-me)\le L_i$.

For $n=3,4$ we also have $\rho=r(1-3e)\ge r/2$ since $3e\le1/16$. For all $n\ge3$ we now have $\rho\ge r/2$ and the side-length comparison.
The differences of matched oriented edges form a closed sequence with nonnegative lengths in cyclic angular order, hence bound a convex polygon $K$ (possibly a point or segment). Indeed every vertex lies to the left of each directed edge, since projecting the cyclic increments to its left normal gives positive and then negative increments with zero total. For each outward normal of $W$, the exposed face of $W+K$ is the sum of the exposed faces of $W$ and $K$ in that direction. The face of $K$ may be a vertex. Since every side direction of $K$ comes from $W$, the oriented side vectors of $W+K$ are exactly those of $P$. They determine the polygon up to translation.
Take $0$ to be one endpoint of a diameter $[0,v]$ of $K$. The area of $W+[0,v]$ exceeds $|W|$ by at least $2\rho |v|$, and Jensen gives $|W|=\rho^2\sum_i\tan(\omega_i/2)\ge(1-He)^2$. Thus
\[
 d_H(P,W)\le\operatorname{diam}K\le (1-|W|)/(2\rho)\le 2He/r .
\]
Place $W$ and the regular polygon so that their incircles have the same center and one pair of outward normals agrees. Each vertex lies on the bisector of its two adjacent normal directions, with radii $\rho\sec(\omega_i/2)$ and $r\sec(a/2)$. The angle between each pair of matched vertices, measured at the common center, is at most $(2n+1)e$, for all $n$. For $n\ge5$ we have $|(\sec)'|\le\sqrt2$ on $[0,\pi/4]$. Their Euclidean distances are therefore at most $\sqrt2 r(H+2n+2)e$, which also bounds Hausdorff distance by convexity. Since $3<nt<4$, $1/2<r<1/\sqrt3$ and $n\ge5$, this proves
\[
 d(P)\le [2H/r+\sqrt2 r(H+2n+2)]e<10ne\le16ne.
\]
For $n=3,4$, the bisector directions differ by at most $1.12ne$ by summing the turn errors from the aligned normals. The matched radii differ by at most $re(3U+.56U^2)$, since $|(\sec)'|\le U^2$ on the intervening half-angle interval. The regular radius is at most $rU$. Adding the diameter bound for $K$ gives
\[
 d(P)\le e\,[6/r+r(3U+.56U^2+1.12nU)].
\]
For $n=3$ use $2/5<r<1/2$, $U=21/10$ to bound the bracket by $24=8n$. For $n=4$ use $r=1/2$, $U=3/2$ to bound it by $20=5n$. This proves $d(P)\le16ne$ throughout the asserted range.

For $n\ge5$ this implication with $e=\sqrt{\Delta/L_0}$ gives the stability bound. For $\Delta>c$ the bound follows from $\min\{1,d(P)^2\}\le1$.

To bound $\Delta_\lambda$ below by a multiple of $d(P)^2$ for convex polygons, put $d_R=\operatorname{diam}R_n$ and $\Lambda_0=\pi^2 n t$. The diameter estimate below gives $D^2\le 4(\Delta_\lambda+\Lambda_0)/\pi^2$ for $D=\operatorname{diam}P$. Translating the two polygons to share a point gives $d(P)^2\le D^2+d_R^2$. For $\Delta_\lambda\ge c$, the quotient of $\Delta_\lambda$ by $4(\Delta_\lambda+\Lambda_0)/\pi^2+d_R^2$ is increasing. The small- and large-deficit estimates therefore give the constant
\[
 \min\left\{c,\frac{c}{4(c+\Lambda_0)/\pi^2+d_R^2}\right\}.
\]
Rescaling gives both area-normalized statements. The affine perturbations in Section~\ref{sec:perimeter} show that the exponent two is optimal.

\paragraph{The diameter bound.}
The inscribed disk contains a square of side \(\sqrt2\,r\), so domain
monotonicity gives \(\lambda_1(R_n)\le\Lambda_0=\pi^2/r^2\).
A convex polygon of area one and diameter \(D\) has width at most
\(2/D\) perpendicular to a diameter, by adding the areas of the
triangles on either side. One-dimensional Dirichlet Poincaré on the
enclosing strip gives \(D^2\le4\lambda_1(P)/\pi^2\).

\subsubsection{The geometric estimate at fixed perimeter}\label{sec:perimeter-geometry}

Put \(s_n=n\tan(\pi/n)\). We first compare the convex hull to a polygon
with the same side normals and an incircle. The tangent remainder for
\(\tan(\theta/2)\) then compares those normals with the regular ones.
We then bound the distance from the polygon to its hull in terms of their perimeters.
\begin{proof}[Proof of \eqref{eq:perimeter-hausdorff}]
First let \(K\) be a convex polygon with at most \(n\) sides and
perimeter \(\ell\). Let \(\nu_i\) be its successive outer normals and \(\theta_i\) their angular gaps.
Repeat normals if needed to pad the list to length \(n\), using zero
gaps for the additional entries. Set
\[
 W=\bigcap_i\{x:x\cdot\nu_i\le1\},\qquad
 s=\sum_i\tan(\theta_i/2),\qquad r=\frac{\ell}{2s}.
\]
The polygon $W$ has area $s$ and perimeter $2s$. The support function $p$ of a polygon satisfies $p+p''=\sum_i L_i\delta_{\nu_i}$ on the unit circle: between successive normals it is a vertex projection, and at a normal its derivative jumps by the side length $L_i$. Hence perimeter is $\int p$ and area is $\tfrac12\int(p^2-p'^2)$, by the signed triangle area formula. Set $H=h_K-h_{rW}$. Equal perimeters give $\int H=0$. The cross term $\int h_{rW}(h_K+h_K'')=r\ell=2|rW|$, since $h_{rW}=r$ on the normals of $K$. Thus
\[
 \Delta:=sr^2-|K|=\tfrac12\int (H'^2-H^2).
\]
Translate $K$ to remove the first harmonic of $H$. For its remaining Fourier coefficients $a_j,b_j$, Parseval gives $2\Delta=\pi\sum_{j\ge2}(j^2-1)(a_j^2+b_j^2)$. Cauchy--Schwarz and $\sum_{j\ge2}(j^2-1)^{-1}=3/4$ give
\[
 d_H(K,rW)^2=\|H\|_\infty^2\le 3\Delta/(2\pi).
\]
In particular $|K|\le\ell^2/(4s)\le A_{*,\ell}:=\ell^2/(4s_n)$ by Jensen. Put $\delta_K=1-|K|/A_{*,\ell}$.
We next compare $W$ with the regular polygon $W_0$ of inradius one.
Put $a=2\pi/n$, $y=a/2$, and $f(\theta)=\tan(\theta/2)$.
The minimum on $[0,\pi)$ of the tangent remainder divided by $(\theta-a)^2$, with the ratio extended continuously at $\theta=a$, is
\[
 c_a=\frac{a f'(a)-f(a)}{a^2}
     =\frac{a\sec^2(a/2)/2-\tan(a/2)}{a^2}.
\]
Indeed the remainder divided by $(\theta-a)^2$ is
$\int_0^1(1-u)f''(a+u(\theta-a))\,du$, also at $\theta=a$ by continuity.
Since $f''$ increases, this integral is smallest at $\theta=0$. Summing over the gaps and using $\sum_i(\theta_i-a)=0$ gives $\sum_i(\theta_i-a)^2\le(s-s_n)/c_a$. We may rotate $W_0$ before comparing its normals with those of $W$.
Write $g_i=\theta_i-a$ and $S_j=\sum_{i=1}^j g_i$, with $S_0=S_n=0$.
For any block $J$ of $k$ successive indices, its sum is
$\langle g,\mathbf 1_J-(k/n)\mathbf 1\rangle$.
The second vector has squared norm $k(1-k/n)\le n/4$.
Hence the range of the partial sums has length at most
$\sqrt n\,\|g\|_2/2$. Choose the rotation that centers this range
at zero. Every pair of matched normals is then at angular distance
at most $\sqrt n\,\|g\|_2/4$, including repeated normals.
For any ray from the origin, its angular distances $d,d_0$ to the two sets of normals satisfy
the same bound, so
\[
 |d-d_0|\le \tfrac14\sqrt{n(s-s_n)/c_a}.
\]
The radial functions of the two polygons on this ray are $\sec d$
and $\sec d_0$. Since $d_0\le y$, replacing $d$ by $\min(d,y)$
bounds their difference by $\sec y\tan y\,|d-d_0|$, apart from
the excess $\sec x-\sec y$ at the largest half-gap
$x=\tfrac12\max_i\theta_i$. This last error is controlled by the
same tangent remainder. For $y\le v\le x$,
\[
 \frac{d}{dv}(\sec v-\sec y)^2
 \le 2\tan x\,(\sec^2v-\sec^2y).
\]
Integrating, using $\tan x\le s$ and bounding the largest gap's
remainder by the sum of all the remainders, gives
$(\sec x-\sec y)^2\le 2s(s-s_n)$.
Matching points on rays and using $r=\ell/(2s)$,
$\tan y=s_n/n$ and $(s-s_n)/s\le\delta_K$ now yields
\[
 d_H(rW,rW_0)\le \ell\sqrt{\delta_K}
 \left(\frac1{\sqrt2}+\frac{\sec y}{8}
                 \sqrt{\frac{\tan y}{c_a}}\right).
\]
For general \(P\), write \(K=\operatorname{conv}\overline P\), and let
\(\ell\) and \(L\) be the perimeters of \(K\) and \(P\), respectively.
The vertices of \(K\) occur in their cyclic order on \(\partial P\),
since disjoint arcs in a disk cannot join alternating boundary points.
Each intervening arc has length \(b\) at least the length \(d\) of its
chord, so \(\ell\le L\), \(|P|\le |K|\le A_{*,\ell}\le A_*\),
and \(\delta_K\le\delta\).

The same chords control the distance to the hull. Pair the point
\(p\) at arclength \(tb\) on an arc from \(A\) to \(B\) with
\(c=(1-t)A+tB\). The identity
\[
 |p-c|^2=(1-t)|p-A|^2+t|p-B|^2-t(1-t)d^2
 \le t(1-t)(b^2-d^2)
\]
gives \(|p-c|\le \tfrac12\sqrt{L^2-\ell^2}\), since
\(b-d\le L-\ell\) and \(b+d\le L+\ell\).
The straight segments from \(p\) to \(c\) give a homotopy of the two
boundaries within that distance of \(\partial P\). Any point in the
interior of \(K\) outside this neighborhood has the same winding number
for both curves and thus belongs to \(P\). Taking closures gives
\[
 d_H(K,\overline P)\le\tfrac12\sqrt{L^2-\ell^2}
 \le \tfrac L2\sqrt{\delta},
\]
where the last inequality uses \(|P|\le A_{*,\ell}=\ell^2/(4s_n)\).
The target polygon of perimeter $L$ is $L W_0/(2s_n)$.
The area bound for $K$ gives
\[
 1-\delta\le (\ell/L)^2s_n/s\le (\ell/L)s_n/s
             =r\big/\bigl(L/(2s_n)\bigr)\le 1 .
\]
Thus the Hausdorff distance from $rW_0$ to its rescaling of perimeter $L$ is at most $L\sec(\pi/n)\delta/(2s_n)$. Since $4s\Delta=\ell^2-4s|K|$ and $s|K|\ge s_n|P|$, Cauchy--Schwarz gives
\[
 d_H(\overline P,K)+d_H(K,rW)
 \le \sqrt{\tfrac14+\tfrac{3}{8\pi s}}\,L\sqrt\delta .
\]
Adding the three comparisons and using $\delta\le\sqrt\delta$, $s\ge s_n$, and $\ell\le L$ gives $d(P)\le \kappa_nL\sqrt\delta$, where
\[
 \kappa_n=\sqrt{\tfrac14+\tfrac{3}{8\pi s_n}}+\frac1{\sqrt2}
       +\frac{\sec y}{8}\sqrt{\frac{\tan y}{c_a}}
       +\frac{1}{2n\sin y}.
\]
For $n=3,4$, substitution of $y=\pi/3,\pi/4$ gives $c_a>.55,.23$, respectively. Using $3.14<\pi<22/7$ in the four terms gives
\[
 \kappa_3<.524+.708+.445+.193=1.870,\qquad
 \kappa_4<.530+.708+.370+.177=1.785 .
\]
For $n\ge5$, $y\le\pi/5<.63$. The third term equals
$\tfrac y4\sqrt{\tan y/(y-\sin y\cos y)}$.
The fifth-order Taylor bound for $\sin(2y)$ gives
$y-\sin y\cos y\ge (2/3)y^3(1-y^2/5)$, and $\tan y/y\le1/\cos y<5/4$.
Thus this term is less than $.358$. The fourth term is at most $1/(10\sin(\pi/5))<.172$, since $\sin y/y$ decreases. With $s_n>\pi>3$ the first term is less than $.541$. Hence $\kappa_n<1.779$ for $n\ge5$, and the required bound $d(P)\le2L\sqrt\delta$ holds for all $n\ge3$.

\end{proof}

\subsection{Hausdorff stability for all exponents}\label{S:sec:sob-stability}
Write $d(P,R)=\inf_g d_H(\overline P,g\overline R)$, with $g$ ranging over plane isometries.
\begin{theorem}\label{S:thm:full-stability}
Let $n\ge3$, $1\le q\le2$, and let $P$ be a simple polygon with at most $n$ sides, area $A$ and perimeter $\ell$. Let $R$ and $R_\ell$ be regular $n$-gons of the same area and the same perimeter as $P$, respectively. Some $c_n>0$, independent of $q$, satisfies
\[
 A^{2/q}(\lambda_q^{\rm S}(P)-\lambda_q^{\rm S}(R))\ge
 c_n\min\{1,d(P,R)^2/A\}.
\]
At fixed perimeter,
\[
 \frac{\lambda_q^{\rm S}(P)}{\lambda_q^{\rm S}(R_\ell)}-1
 \ge \frac{d(P,R_\ell)^2}{2q\,\ell^2}.
\]
The power two is optimal in both estimates for each $q$.
\end{theorem}

\begin{proof}
Scale to $A=1$, and put $a=2\pi/n$, $\gamma=2/q$, $s=F_q'(a)$,
$b=a-F_q(a)/s$ and $L=nF_q(a)=s(2\pi-nb)$. All comparison constants in the fixed-area proof may depend on $n$. The variational bounds in Section~\ref{S:sec:response} bound $L$ above and below by positive constants uniformly in $q$.

\emph{A lower bound in terms of the base-angle difference.}
On any compact interval of openings about $a$ contained in $(0,\pi)$, a convex cell of opening $\alpha$ and base angles differing by $\pi\delta$ has normalized mixed quotient at least $F_q(\alpha)+k\delta^2$ for some $k>0$ independent of $q$.
To see this use the strip representation of Section~\ref{sec:convex} and Section~\ref{sec:stability-tilt}. In its notation $C=1-\alpha/\pi$, $p=p_{0,1}$, $h_0=\rho_C/\pi$, and $P_0$ is the weight of the first product density in $p$. Write $M_q(h)=\sup_{\|\nabla v\|_2\le1}\int h|v|^q$ on the strip with zero bottom trace. The discrete comparison there, obtained by summation by parts on the mixture weights, gives for every $q$
\[
 M_q(p)-M_q(p_{\delta,\kappa})\ge
 (P_0\delta^2/C^2)\,[M_q(p)-M_q(h_0)] .
\]
The bracket is strictly positive. Take a nonnegative maximizer for $h_0$ rearranged to increase in $y$ and decrease in $|x|$. Its first nonconstant product density $\rho_{C+1}(x)(1-\cos y)/\pi$ gives a strictly larger integral: the change in the horizontal factor can only increase it, and each nonzero vertical slice increases from zero trace and has strictly positive covariance with $1-\cos y$. The later components give at least this value by monotonicity of their likelihood ratios.

For compact $C$-ranges in $(0,1)$ these suprema depend continuously on $(C,q)$. Indeed $p$ and $h_0$ vary continuously in $L^r$ for some $r>1$ sufficiently close to one. Their tails decay exponentially, and near the singular point on the upper edge the densities are bounded by a constant times the distance to that point raised to the power $-2\max C$. Choose $r$ so $(1-\varepsilon)r'\ge2$ for some $\varepsilon>0$. The Poincar\'e inequality for functions vanishing on the bottom edge bounds the $H^1$ norm of all unit-energy tests, hence their $L^t$ norms for $2\le t\le(2+\varepsilon)r'$. H\"older and the bound on $|v|^q|\log |v||$ by a constant times $|v|^{q-\varepsilon}+|v|^{q+\varepsilon}$ give continuity uniformly over tests. Local compact embedding and the decaying tails also give attainment for $h_0$. The component weight $P_0$ is continuous and positive by the same normalization integrals. Taking the power $M\mapsto M^{-2/q}$ now gives the claimed uniform lower bound; translations and averages preserve the upper bound on $M_q$ for a general cell.

\emph{A lower bound for the deficit from the partition.}
Put $h=2\pi-nb$. For $n\ge4$, \eqref{eq:line-constants} gives $b\le a-\sin a<a/2$ since $0<a\le\pi/2$, so Lemma~\ref{lem:mass-cap-partition} applies.
For $n=3$ use Lemma~\ref{lem:balance} on all three sides of the triangle. Here
$b\le 2\pi/3-\sqrt3/2<\pi/2$ and $h>0$. A single side requires $b<\pi$, and two sides meeting at angle $\sigma$ require $2b<2\pi-\sigma$, which follows from $2b<\pi$.
This gives $\alpha_i=b+h\nu_i$ on all three cells, so the same summation applies.
Use the notation of Proposition~\ref{prop:geometric-reduction} at area one,
so $S_i=|Q_i|$, $r_i=S_i^{-1}\int_{Q_i}u^q$, and $\nu_i=S_i(1+r_i)/2$
for $\int_Pu^q=1$. With $D_q(\alpha)=F_q(\alpha)-s(\alpha-b)$,
let $J$ be a compact neighborhood of $a$ in $(0,\pi)$. Retaining
the term for unequal base angles in the cell comparison of that proposition gives
\[
 \Delta:=\lambda_q^{\rm S}(P)-L\ge \tfrac12 L\gamma\sum_i S_i(r_i-1)^2+
 \sum_i r_i^\gamma\bigl[D_q(\alpha_i)+k\mathbf1_J(\alpha_i)\delta_i^2\bigr].
\]
Here $s(\alpha_i-b)\ge L\nu_i$ because there are at most $n$ sides.

By Theorem~A and Lemma~\ref{S:lem:exponent-continuity}, $D_q$ has a uniform quadratic lower bound on $J$ and a uniform positive lower bound outside $J$, by convexity on either side of $a$. Retain the $j$ cells with $r_i\ge1/2$. Their openings differ from $a$ by $O_n(\sqrt\Delta)$ when $\Delta$ is small, since $r_i^\gamma\ge1/4$. The discarded cells have total area and total $\nu$-mass $O_n(\Delta)$. Summing the opening bounds over retained cells yields
\[
 C_n\sqrt{\Delta}\ \ge\ \sum_{\rm retained}(\alpha_i-a)
 \ \ge\ (n-j)(a-b)-(2\pi-nb)\nu(\text{discarded cells}).
\]
As $a-b\ge\sin a$, sufficiently small $\Delta$ forces $j=n$. The term involving $\delta_i^2$ gives $|\delta_i|=O_n(\sqrt\Delta)$ on every cell.
Also $\nu_i\le(\alpha_i-b)/(2\pi-nb)\le1/n+O_n(\sqrt\Delta)$. Summation gives the matching lower bound. The variance controls $|S_i-\nu_i|=S_i|1-r_i|/2$, so $S_i=1/n+O_n(\sqrt\Delta)$.

\emph{Recovering the polygon.}
Since all $n$ cells have positive area, summing their partition bounds gives
$0\le E=\sum_i\alpha_i-2\pi=O_n(\sqrt\Delta)$.
Thus $E$, $|\alpha_i-a|$, $\pi|\delta_i|$ and $|S_i-1/n|$ are all
$O_n(\sqrt\Delta)$. For every $n\ge3$ these bounds imply \eqref{eq:one-cap-data} with $e=C_n\sqrt\Delta$, where $C_n$ is independent of $q$. Section~\ref{sec:stability-cell-geometry} therefore gives $d(P,R)\le16ne$ for sufficiently small $\Delta$. For larger deficits the estimate follows from $\min\{1,d(P,R)^2\}\le1$, after decreasing $c_n$ if necessary. Dilation restores the original area.

At fixed perimeter put $\delta_A=1-|P|/|R_\ell|$. The fixed-area comparison and scaling give a quotient ratio at least $(1-\delta_A)^{-2/q}$. The geometric inequality $d(P,R_\ell)\le2\ell\sqrt{\delta_A}$ of \eqref{eq:perimeter-hausdorff} and Bernoulli's inequality give the claim.

The area-preserving stretch in Section~\ref{sec:perimeter} proves sharpness for every $q$.
\end{proof}

\paragraph{Explicit fixed-area constants for three and four sides.}\label{par:small-n-area}
We give the constants in Theorem~\ref{thm:stability} using the deficit inequality in the preceding proof for $q=1,2$. Take $J=[a-1/10,a+1/10]$. On $J$ we have $1/2<\tan(\alpha/2)<3$ and $1/4<C=1-\alpha/\pi<3/5$, so \eqref{eq:Fq-bounds} gives $F_q>\pi^2/8>1$.
For $q=1$, since $1/10<a/8$, Section~\ref{sec:tangent-gaps} gives $-g''\ge a^3/256$ on $J$. Also $F_1/\alpha\ge\pi^2/8>1$ by \eqref{eq:Fq-bounds} and $\tan(\alpha/2)\ge\alpha/2$. Hence $F_1''\ge(F_1/\alpha)^2(-g'')>1/20000$ on $J$. Integrating twice and using convexity beyond the endpoints of $J$ gives the bound below. For $q=2$, \eqref{eq:spectral-quadratic-gap} gives it directly since $1369n^2<40000$. Thus for $q=1$ and for $q=2$
\[
 D_q(\alpha)\ge \min\{(\alpha-a)^2,1/100\}/40000 .
\]
In \eqref{sr:mixturegain}, $P_0(1-P_0)/C^2\ge(1-C)/(1+C)^2>1/8$ and the test $\psi$ gives $M_q(h_0)\ge [C/(2C+1)]/(2\pi/3)>1/16$. Thus the moment decreases by at least $\delta_i^2/1000$ for $\alpha_i\in J$, including after translation and averaging. Since $M_q(p)=F_q^{-q/2}<1$, taking the power $-2/q$ gives at least the same gain in the quotient.

Write $x=\sqrt\Delta$, where $\Delta=\lambda_q^{\rm S}(P)-L$ at unit area, and suppose $x\le x_0:=1/12{,}800{,}000$. The variance bound gives $\sum_i S_i(r_i-1)^2\le2\Delta/L$, so the cells with $r_i<1/2$ have total area and total $\nu$-mass at most $8\Delta/L$. Each retained cell has $r_i^{2/q}\ge1/4$, so $\min\{(\alpha_i-a)^2,1/100\}\le160000x^2$. Since $400x_0<1/10$, its opening lies in $J$ and $|\alpha_i-a|\le400x$. If only $j<n$ cells were retained, summing their opening bounds would give
\[
 1600x\ge(n-j)(a-b)-8h\Delta/L,
\]
which is impossible since $a-b\ge\sin a\ge\sqrt3/2$, $L>1$, $h\le2\pi$ and $1600x_0+16\pi x_0^2<1/2$. Thus all $n$ cells occur and $0\le E\le1600x$. The bound for unequal base angles gives $\pi|\delta_i|\le\pi\sqrt{4000}\,x$. The nonnegative remainders $\alpha_i-b-h\nu_i$ sum to $E$, so $|\nu_i-1/n|\le2000x/h$; here $h=n(a-b)>2$. Finally $|S_i-\nu_i|\le x/\sqrt{2L}$. All four errors are therefore at most $2000x$. With $e=200000x$ this gives \eqref{eq:one-cap-data} and $e\le200000x_0=1/64\le1/(16n)$. The implication proved in Section~\ref{sec:stability-cell-geometry} gives $d(P,R_n)\le16ne\le x/x_0$.

For $x>x_0$ we have $\Delta>x_0^2$. Combining the two cases gives $\Delta\ge x_0^2 m$ for both quotients, $q=1$ and $q=2$, where $m=\min\{1,d(P,R_n)^2\}$. At $q=1$ the torsion deficit is $\Delta/[L(L+\Delta)]$, increasing in $\Delta$, and $L\le9n\tan(\pi/n)<72$ by \eqref{eq:Fq-bounds}. Substitution of $x_0^2m$ gives a lower bound of $x_0^2m/(72\cdot73)$ since $x_0^2m<1$. As $x_0^2>6\cdot10^{-15}$ and $72\cdot73<6000$, a constant valid in both cases after rescaling is $c_n=10^{-18}$. For convex polygons the diameter argument of Section~\ref{sec:stability-cell-geometry} then gives an eigenvalue bound proportional to $d(P,R_n)^2$ at unit area, with this value of $c$ in the formula for the constant.

\clearpage
\part{Intermediate exponents and angular convexity}\label{part:intermediate}
\section{Derivatives of the minimum on a fixed triangle}\label{S:sec:response}
In Sections~\ref{S:sec:response}--\ref{sec:acute-real} we write $\lambda_q$ for $\lambda_q^{\rm S}$. For a domain $G$ with a specified Dirichlet edge or arc, $H_D^1(G)$ consists of the $H^1$ functions of zero trace there, and $\lambda_q(G)$ is the infimum of $\int_G|\nabla v|^2/(\int_G|v|^q)^{2/q}$ over nonzero $v\in H^1_D(G)$. Thus $F_q(\alpha)=|I_\alpha|^{2/q}\lambda_q(I_\alpha)$. Put
\[
 z=\cot^2(\alpha/2),\qquad k_\alpha=(1+z)/(2z),\qquad
 \mathcal C_q(z)=F_q''(\alpha)/(k_\alpha F_q(\alpha)).
\]
As the opening $\alpha$ decreases from $\pi$ to zero, $z$ increases from zero to infinity. In these variables Theorem~A asserts $\mathcal C_q>0$.

The fixed triangle is $D=\{(x,y):0<y<x<1\}$ with probability measure
$d\mu=2\,dx\,dy$ and Dirichlet edge $x=1$. Let $X,Y$ denote the
quadratic forms $\int_D v_x^2\,d\mu,\int_D v_y^2\,d\mu$, with the same
notation for their polarizations. Put $Q_z=X+zY$. Taking the
quadratic mean of a nonnegative test and its reflection in the
bisector reduces energy and increases $q$-mass, by Cauchy--Schwarz
and concavity of $t^{q/2}$. It suffices to use symmetric tests.
The map $(x,y)\mapsto(x,y/\sqrt z)$ onto a half triangle gives
\begin{equation}\label{S:eq:aff}
 F_q(\alpha)=\Lambda_q(z)/\sqrt z,\qquad
 \Lambda_q(z)=\inf_{\|u\|_{L^q(d\mu)}=1}Q_z(u).
\end{equation}
Each quotient at a fixed test is affine and nondecreasing in $z$.
Their infimum is therefore concave and nondecreasing.
A horizontal line meets $D$ in an interval of length at most one.
The one-dimensional Poincar\'e inequality on each such interval, followed by H\"older's inequality for the probability measure $d\mu$, gives
$X(v)\ge (\pi^2/4)\|v\|_q^2$. The test function $1-x$ shows
$\Lambda_q(z)\le ((q+1)(q+2)/2)^{2/q}$. Hence $F_q(0+)=0$.

Let $u>0$ be the normalized minimizer, $\lambda=\Lambda_q(z)$,
$p=q-1$, $r=2-q$ and $e=u/\sqrt\lambda$. We take $Q_z$ as the inner product and define $R$ and $W$ by
$Q_z(Rv,w)=zY(v,w)$ and $Q_z(Wv,w)=\int_D\lambda u^{-r}vw\,d\mu$. The weak Euler equation gives
$We=e$, $\|e\|_{Q_z}=1$. The ground-state inequality from the proof of Lemma~\ref{lem:acute-eliminate}, applied to $u$ by testing with $v^2/(u+\epsilon)$ and letting $\epsilon\downarrow0$, gives $0\le W\le I$. The operator $R$ is between $0$ and
$I$. Define
\[
 \beta=z\lambda'/\lambda,\qquad G=-z^2\lambda''/(2\lambda),
 \qquad g=(R-\beta)e.
\]
\begin{lemma}[First and second derivatives with respect to $z$]\label{S:lem:response}
On $e^\perp$ in the energy inner product,
\begin{align}
 G&=\langle g,(I-pW)^{-1}g\rangle,\label{S:eq:response}\\
 \mathcal C_q(z)&=1+(z-1)\beta-4(1+z)G,\qquad
 \frac{d\beta}{d\log z}=\beta-\beta^2-2G.\label{S:eq:curv}
\end{align}
In particular $0\le G\le r^{-1}\|g\|^2$ when $q<2$.
If $X\ge kW$ as forms, with $0<k<1$ and $p\le k$, then
\begin{equation}\label{S:eq:twoform}
 0\le\beta\le1-k,\qquad
 G\le\beta-\frac{\beta^2}{1-k}.
\end{equation}
\end{lemma}
\begin{proof}
Differentiating the normalization gives $\int_Du^{q-1}u'\,d\mu=0$. The Euler equation then gives $u'\perp e$. The envelope identity $\lambda'=Y(u)$
gives $\langle Re,e\rangle=\beta$ and $g\perp e$. Moreover $Y(u)>0$.
Otherwise $u$ depends only on $x$, and the Neumann condition on the diagonal forces $u_x=0$ there for every $0<x<1$. Such a constant function cannot vanish on the base and have unit mass. In particular $\beta>0$. Differentiating the
Euler equation on the tangent space gives
$zu'/\sqrt\lambda=-(I-pW)^{-1}g$. Since $\lambda''=2Y(u,u')$,
this proves \eqref{S:eq:response}. Equation \eqref{S:eq:curv} follows by
differentiating \eqref{S:eq:aff} and $z=\cot^2(\alpha/2)$.
On $e^\perp$ the inverse has norm at most $1/(1-p)$ for $p<1$. For $p=1$ it exists on the same subspace by compactness and simplicity of the first eigenvalue.

The bound~\eqref{S:eq:twoform} is the case $p\le k$
of the Cauchy--Schwarz argument given in Section~\ref{S:sec:compact-curvature};
that argument uses only~\eqref{S:eq:response}.
\end{proof}

\begin{proposition}[Two linear problems at one opening]\label{prop:frozen-linear}
Fix $1<q\le2$ and an opening $\alpha_0$. Let $u$ be its positive minimizer, normalized by $\|u\|_q=1$. Put $d\nu=u^q\,d\mu$, $K_\alpha=z^{-1/2}Q_z$, and
\[
 P_t(\alpha)=\inf_{v\ne0}K_\alpha(v)/\|v/u\|_{L^t(\nu)}^2 .
\]
Then $P_1\ge P_q=F_q\ge P_2$ at every opening, with equality throughout at $\alpha_0$. The Euler equation for $P_1$ has source proportional to $u^{q-1}$, and $P_2$ is the first eigenvalue with weight $u^{q-2}$. In particular, for $q<2$,
\[
 F_q(\alpha)=
 \max_{\int m^{-q/(2-q)}\,d\mu=1}\
 \inf_{v\ne0}\frac{K_\alpha(v)}{\int m v^2\,d\mu},
\]
where $m$ ranges over positive form-bounded weights. If $\alpha_0\le\pi/2$, all three curves have the same tangent there.
\end{proposition}
\begin{proof}
The norm increases with $t$ on the probability space $(D,\nu)$.
At $\alpha_0$, the inequality $W\le I$ gives the lower bound for $P_2$, and $v=u$
attains it in all three quotients. For any admissible $m$, H\"older
gives $\|v\|_q^2\le\int mv^2\,d\mu$. Taking $m=u^{q-2}$ at the
opening in question attains the displayed maximum. Positive tests suffice, so for $t=1$ the norm in the denominator equals $\int_Dvu^{q-1}\,d\mu$ before squaring. For $z\ge1$ the product comparison in
Lemma~\ref{S:lem:exponent-continuity} gives
$\lambda u^{q-2}\le (1+z)(\pi^2/4)
[\cos(\pi x/2)\cos(\pi y/2)]^{q-2}$
(for $q=2$ use the product as a trial). The weight thus lies in
some $L^a$, $a>1$. By compactness of the weighted embedding and simplicity, the derivative of each minimum at $\alpha_0$ is the partial derivative of its quotient with $v=u$. The three derivatives agree.
\end{proof}

\paragraph{Regularity for fixed exponent.}\label{par:fixed-regularity}
For $q=2$ the first eigenvalue of $Q_z$ is simple,
so ordinary perturbation theory applies. For $1\le q<2$ the rescaled minimizer
$U=\lambda^{-1/r}u$ solves $Q_z(U,v)=\int U^p v\,d\mu$.
The comparison argument below proves uniqueness and yields, for all $p<1$, the pointwise bound
\[
 U^{-r}\le B(x,y)=(1-y)^{-2}
 \varphi\big((x-y)/(1-y)\big)^{-r}.
\]
Indeed, in the Euclidean coordinates $(x,y/\sqrt z)$ center a disk of
radius $1-y_0$ at $(y_0,y_0/\sqrt z)$ on the diagonal. Its scaled
radial solution, extended by zero, has zero trace at $x=1$ and
nonpositive outward flux on the bottom and its circle. Its diagonal
flux vanishes. For any positive subsolution $v$, test the difference
with $d=(v-U)_+$. Concavity of the right-hand side and the
ground-state inequality give $Q_z(d)\le pW(d)\le pQ_z(d)$.
This also gives uniqueness by interchanging two positive solutions.
Taking $y_0=y$ gives the displayed bound,
where $\varphi$ is the positive radial disk solution of $-\Delta\varphi=\varphi^p$
with zero boundary value. We use that $\varphi$ vanishes linearly at the boundary.
On a unit slice $|v(s)|^2\le(1-s)\int_s^1 |v'|^2$.
Since $\varphi(s)\ge c_p(1-s)$, rescaling and integration give
$\int Bv^2\,d\mu\le C_p X(v)$, also for $p=0$.

Let $S_h$ be the difference quotient of $t^p$ between $U_z$ and $U_{z+h}$, taking the derivative where the two values agree.
Convexity of $t^{-r}$ gives $0\le S_h\le
(p/2)(U_z^{-r}+U_{z+h}^{-r})\le pB$. The ground-state inequalities at
both parameters imply
$Q_z-S_h\ge rQ_z-p|h|Y/2\ge (r/2)Q_z$ locally. Subtraction gives
$(Q_z-S_h)(U_{z+h}-U_z,\cdot)=-hY(U_{z+h},\cdot)$, where
multipliers denote their bilinear forms. Local energy bounds for the
rescaled minimizers follow from normalization and the bounds for $\Lambda_q$.
Thus $U_z$ is locally Lipschitz in the energy norm. If $0\le V_j,V\le CB$
and $V_j\to V$ in measure, the corresponding operators converge strongly:
apply dominated convergence in measure to $|V_j-V|v^2$ for fixed $v$,
then Cauchy--Schwarz in the second test. This applies to $S_h$
with limit $pU_z^{-r}$. Uniform coercivity and the resolvent identity
give strong convergence of the inverses. The difference equation now
gives $U'_z=-(Q_z-pU_z^{-r})^{-1}Y(U_z,\cdot)$ in the energy norm, and this derivative is continuous in $z$. Normalization gives a continuously differentiable $u$, and
$\lambda'=Y(u)$ shows that $\lambda$ is twice continuously differentiable. At $p=0$ we have $S_h=0$ and normalization uses the integral of the positive
minimizer.
\begin{lemma}[Continuity in the exponent and opening]\label{S:lem:exponent-continuity}
The normalized curvature $\mathcal C_q(z)$ and $F_q''(\alpha)$ are jointly
continuous on $[1,2]\times(0,\infty)$ in the variables $(q,z)$.
\end{lemma}
\begin{proof}
Let $(q_i,z_i)\to(q,z)$ with $d\le z_i,z\le D_0$ for $d>0$.
Let $u_i$ be the positive unit-mass minimizer and $\lambda_i=\Lambda_{q_i}(z_i)$.
Testing with the minimizer at $(q,z)$ gives $\limsup_i\lambda_i\le\lambda$. Bounded energy gives compactness
in every finite $L^k$ on the fixed triangle. In particular the limit has
unit $L^q$ mass by uniform integrability. Lower semicontinuity and uniqueness of
the positive minimizer imply $\lambda_i\to\lambda$ and strong energy convergence
$u_i\to u$. We may pass to almost everywhere convergent subsequences.

Put $m=\max\{1,1/d\}$, $c(x)=\cos(\pi x/2)$ and $\Phi=c(x)c(y)^m$.
For $A_z=-\partial_x^2-z\partial_y^2$, differentiation gives
$A_{z_i}\Phi\le M\Phi$, $M=(\pi^2/4)(1+D_0m)$. The bottom flux
vanishes and the outward diagonal flux has the sign of $1-z_i m\le0$.
For $r_i=2-q_i>0$, the function $M^{-1/r_i}\Phi$ is therefore a weak
subsolution of $A_{z_i}U=U^{1-r_i}$. The comparison argument on
p.~\pageref{par:fixed-regularity} gives $V_i=\lambda_i u_i^{-r_i}\le M\Phi^{-r_i}$.
At $r_i=0$ the potential is constant. On each set $x\le1-\delta$ the
potentials are bounded uniformly and converge almost everywhere to $V=\lambda u^{q-2}$.

For $q<2$ this local control suffices. Functions vanishing near $x=1$
are dense in the mixed form domain $\mathcal H$. Indeed the one-dimensional
Hardy inequality on each horizontal slice gives
$\int v^2/(1-x)^2\,d\mu\le4X(v)$ and justifies the usual cutoff in $1-x$.
On functions supported away from $x=1$, dominated convergence gives
$(V_i-V)v\to0$ in $\mathcal H^*$. The ground-state inequality $W_i\le I$
and uniform equivalence of the energy norms extend strong convergence
to every fixed $v\in\mathcal H$. The Riesz maps of $Q_{z_i}$ converge in norm.
Since $1-(q_i-1)$ stays positive uniformly, the resolvent identity gives
strong convergence of the inverses of $I-(q_i-1)W_i$.

At $q=2$, eventually $4m r_i\le1/2$. The product bound and
$c(x)\ge1-x$ give $\sup_i\|V_i\|_4<\infty$ after discarding finitely many terms, by
integration on the unit square. Convergence in measure now gives
$V_i\to\lambda$ in $L^2$. The embedding $\mathcal H\hookrightarrow L^4$
implies operator-norm convergence of $W_i$. The energy projections $P_i$
onto $u_i$ converge in norm. The limiting $I-W+P$ is invertible by simplicity
and the gap between the first two eigenvalues. Thus $(I-(q_i-1)W_i+P_i)^{-1}$ converges
in norm and equals the required inverse on the energy-orthogonal complement of $u_i$.
In both cases the representatives of $z_iY$ converge in norm, so the
normalized minimizers $e_i$, the numbers $\beta_i$ and the sources $g_i$ converge in
the second-derivative formula~\eqref{S:eq:response}. Equation~\eqref{S:eq:curv}
and the normalization of $\mathcal C$ prove the assertion.
\end{proof}

\begin{lemma}[A horizontal disk comparison]\label{S:lem:shared-disk-barrier}
Let $0\le p\le1/2$ and let $\varphi$ be the positive radial solution of
$-\Delta\varphi=\varphi^p$ on the unit disk, with zero boundary values.
For every $v\in H^1(0,1)$ with $v(1)=0$,
\[
 \int_0^1|v'|^2\,ds\ge k\int_0^1\varphi^{p-1}v^2\,ds
\]
with $k=4/9$ on the full range, and with $k=p_*=1801/4000$ when $p\le p_*$.
Consequently $X\ge kW$ on the fixed triangle at every $z>0$.
\end{lemma}
\begin{proof}
Set $Y=-s\varphi'/\varphi$, $Z=s^2\varphi^{p-1}$ and $R=Y/Z$.
Differentiation in $\log s$ gives $\dot Y=Z+Y^2$ and
$\dot Z=[2+(1-p)Y]Z$. Thus $Y$ increases from zero to infinity, and
\[
 Y(1+YR)R'=R(1-2R+pYR),\qquad R(0+)=1/2 .
\]
For $0<p\le1/2$ take $a=3/5$, $t=pY/(4a)$ and
$S=(1+t)^a/2$. The inequality for $S$ to be a supersolution is equivalent to
\[
 (1+(1-a)t)(1+t)^{-a}
 \le 1+(1-2a)t+\frac{2a(a-p)}{p}t^2 .
\]
The left side has second derivative
$a(1+t)^{-a-2}[(3a-1)-(1-a)^2t]\le a(3a-1)$.
Its Taylor bound at zero gives the inequality since
$a\ge3p/(4-3p)$. The differential equation's right side, after
division by $Y(1+YR)$, decreases in $R$ near $(0,1/2)$, as seen by
differentiation. Thus $(R-S)_+$ is nonincreasing on an initial
interval, and $R\le S$ there by the common initial value.
The ODE comparison principle extends the bound to every $Y>0$.

Put $B=\varphi^{1-p/5}$. Then
\[
 -B''/(\varphi^{p-1}B)=(1-p/5)[1-(1-pY/5)R].
\]
If $pY\ge5$ the bracket is at least one. Otherwise the term subtracted
is bounded by the maximum of $(1-12t/25)(1+t)^{3/5}/2$, attained at
$t=5/32$. This maximum is less than $101/200$, since Taylor's cubic
upper bound gives $(37/32)^{3/5}\le1+3/32-3/1024+7/32768<101\cdot 2/(5\cdot37)$.
Thus the interval ground-state identity gives
$k>(99/200)(1-p/5)$. This exceeds $4/9$ for $p\le1/2$ and is at least
$p_*=1801/4000$ for $p\le p_*$ since $p_*<495/1099$.
At $p=0$ the explicit disk solution $\varphi=(1-s^2)/4$ gives $k=1/2$.
The bound first holds for tests vanishing near $1$, by the
ground-state identity and $B'(0)=0$, and extends by density and Fatou.
On each horizontal interval the disk subsolution on
p.~\pageref{par:fixed-regularity} bounds the potential by
$(1-y)^{-2}\varphi((x-y)/(1-y))^{p-1}$. Rescale the interval inequality
and integrate in $y$.
\end{proof}

\begin{proposition}[The disk estimate at low exponents]\label{S:prop:lowq}
If $1\le q\le q_\ast$, then $F_q''/F_q\ge 1/[20z(1+z)]$ for all
$z=\cot^2(\alpha/2)>0$.
\end{proposition}
\begin{proof}
Put $p=q-1\le p_*=q_\ast-1$, $z=\cot^2(\alpha/2)$ and
$\lambda=\Lambda_q(z)$. Lemma~\ref{S:lem:shared-disk-barrier}
gives $X\ge p_*W$. Applying \eqref{S:eq:twoform} with $k=p_*\ge p$
and using $1/(1-p_*)>20/11$ gives $G\le\beta-(20/11)\beta^2$;
therefore $d\beta/d\log z\ge-\beta+(29/11)\beta^2$.
As $\beta>0$ by Lemma~\ref{S:lem:response}, $(1/\beta-29/11)/z$ is nonincreasing, and so is $1/(z\beta)$.

We next bound $z\beta(z)$ as $z\to\infty$.
Let $\lambda_\infty$ be the quotient restricted to functions of $x$,
so $\lambda\le\lambda_\infty$. The slice mean
$m(x)=x^{-1}\int_0^x u(x,y)\,dy$ satisfies
\[
 \lambda_\infty\|m\|_q^2\le\int_0^1 2x|m'|^2dx
 \le (\sqrt X+\sqrt{Y/3})^2\le (1+1/(3z))\lambda.
\]
Here $m'$ is the mean of $u_x+(y/x)u_y$ and the mean of $(y/x)^2$ is $1/3$.
For $q>1$ slice Poincar\'e bounds $\|u-m\|_q^q\le C_qY^{q/2}$.
The linear term in expanding $|u|^q$ about $m$ has mean zero,
so $\|m\|_q^2\ge1-C_qY^{q/2}$ for large $z$; for $q=1$ the mass is
preserved exactly. Initially $Y\le\lambda_\infty/z$.
Expanding the middle square with $X=\lambda-zY$ gives
$zY=O(\sqrt Y+Y+Y^{q/2})$, hence $Y=O(z^{-2})$.
The last displayed bound now gives $\lambda\to\lambda_\infty$ and
$\limsup z\log(\lambda_\infty/\lambda)\le1/3$ also at $q=1$.
As $z\beta(z)$ increases,
$\log(\lambda_\infty/\lambda(z))=
\int_z^\infty\beta(x)\,dx/x\ge\beta(z)$.
Consequently $\limsup z\beta(z)\le1/3$; the decreasing quantity
$(1/\beta-29/11)/z$ is at least its limit, hence at least $3$.

Write $1/\beta=V+w$, $V=3z+29/11$, $w\ge0$.
Substituting $G\le\beta-(20/11)\beta^2$ into \eqref{S:eq:curv} gives
\[
 \mathcal C_q\ge
 \frac{w^2+(3z+3/11)w+2z/11+126/121}{(V+w)^2}
 \ge\frac1{10(1+z)^2}.
\]
Indeed the numerator exceeds $1+w^2$, and
$(V+w)^2\le(1+V^2)(1+w^2)$ with $1+V^2\le10(1+z)^2$.
Since $F_q''/F_q=(1+z)\mathcal C_q/(2z)$ this gives
$F_q''/F_q\ge1/[20z(1+z)]$ and proves the asserted estimate.
\end{proof}

\section{Spectral and trial-function estimates}
\label{S:sec:compact-curvature}
For the spectral estimates in this section let $1\le q\le2$, $r=2-q$, $p=q-1$,
$P=\pi^2/4$, $V=\Lambda_q(z)u^{-r}$, and define $W$ by $Q_z(Wv,w)=\int_D Vvw\,d\mu$. On the fixed triangle $0<y<x<1$, the side $x=1$ carries the Dirichlet condition. Put $h=1-y$ and
$s=(x-y)/h$. All triangle integrals use $d\mu=2\,dx\,dy$.

\begin{lemma}[A lower bound for the Green operator on sublevel sets]\label{lem:linear-hemisphere-profile}
Put $A_z=-\partial_x^2-z\partial_y^2$ and let $J$ be a positive
function in the form domain, bounded by one, with bounded gradient. Assume $J$ vanishes on the base and has nonpositive outward conormal derivative on the free sides. Suppose for $d=1$ or $2$ and
$M>0$ that $A_zJ\le MJ$ and $|\nabla J|_z^2\le M(1-J^2)/d$.
Define $T=MA_z^{-1}$ and the probability measure
\[
 d\eta_d(t)=d\,t(1-t^2)^{d/2-1}
          \left(\int_0^t(1-v^2)^{-d/2}\,dv\right)dt,\qquad 0<t<1 .
\]
Then $T[J\mathbf1_{\{J\le a\}}]\ge J\eta_d([0,a])$ for $0<a\le1$.
In particular, for every $r=2-q\in[0,1]$,
\[
 V\le M J^{-r}\left(\int_0^1 t^{-r}\,d\eta_d(t)\right)^{-1}.
\]
For $d=1$ the last factor is at most $e^{-\gamma r}$ with
$\gamma=1/9+2/75+8/735>4/27$. For $d=2$ it is at most $1-r/4$.
\end{lemma}
\begin{proof}
Use $H_a(0)=0$ and
$H_a'(t)=[1-((1-a^2)/(1-t^2))^{d/2}]_+$.
The function $H_a$ is nondecreasing and concave with
$d\,tH_a'-(1-t^2)H_a''=d\,t\mathbf1_{\{t<a\}}$ almost everywhere
for $a<1$. The chain rule gives $A_zH_a(J)\le M J\mathbf1_{\{J\le a\}}$
in the weak sense, also across $J=a$ by continuity of $H_a'$.
Thus $T[J\mathbf1_{\{J\le a\}}]\ge H_a(J)\ge JH_a(1)$.
Differentiating $H_a(1)=\int_0^a H_a'(t)dt$ in $a$ gives precisely
the displayed density. At $a=1$ one has $H_1(t)=t$, $TJ\ge J$ and
$H_1(1)=1$. Integrating the sublevel-set bound using the layer-cake formula for $t^{-r}$ gives
$T(J^{1-r})\ge bJ$ for $b=\int t^{-r}\,d\eta_d$.

For $r>0$ the rescaled minimizer $s=(M/\Lambda_q)^{1/r}u_q$
satisfies $s=T(s^{1-r})$. The comparison argument on p.~\pageref{par:fixed-regularity} gives $s\ge J$ since
$A_zJ\le MJ\le MJ^{1-r}$. Every lower bound $s\ge kJ$ improves to
$s\ge b k^{1-r}J$ by positivity of $T$. Starting at $k=1$ yields
$s\ge b^{1/r}J$, proving the potential estimate. At $r=0$ use $J$ in the unweighted Rayleigh quotient.

For $d=1$ write $t=\sin y$; then $d\eta_1=y\sin y\,dy$ on $(0,\pi/2)$.
Jensen gives $\int t^{-r}d\eta_1\ge\exp(r\int(-\log t)d\eta_1)$.
The first three terms of $-\log\sin y=\tfrac12\sum_{j\ge1}\cos^{2j}y/j$
give the lower bound $\gamma$ by integration by parts, since
$\int_0^{\pi/2}y\sin y\cos^{2j}y\,dy=(2j+1)^{-1}\int_0^{\pi/2}\cos^{2j+1}y\,dy$.
For $d=2$ one has $d\eta_2=2t\operatorname{arctanh}(t)\,dt$.
Its power series gives
\[
 \int t^{-r}d\eta_2=\sum_{j\ge0}\frac{b_j}{1-r/(2j+3)},\qquad
 b_j=\frac2{(2j+1)(2j+3)},\qquad \sum b_j=1 .
\]
The weighted harmonic mean is at most the arithmetic mean, and
$\sum b_j/(2j+3)\ge2/9+2/75+2/245>1/4$. This proves the second bound.
\end{proof}

\begin{lemma}[A bound for the horizontal energy from the interval minimizer]\label{lem:interval-product}
Put $r=2-q$, $a=\max\{1,2z/(1+z)\}$ and $b=a/z$. If $(a+b)r\le2$, then
$X\ge W/(2a)$. In particular this holds at every acute opening for every $1\le q\le2$.
\end{lemma}
\begin{proof}
Write $v=ar<2$. Let $\Phi$ be the positive solution of the mixed boundary value problem
$-\Phi''=\ell\Phi^{1-v}$, $\Phi(0)=1$, $\Phi'(0)=\Phi(1)=0$.
It is positive and concave, with bounded derivative: multiply the ODE by
$\Phi'$ to get $|\Phi'|^2=2\ell(1-\Phi^{2-v})/(2-v)$, which also
constructs the interval solution by separation of variables for all $0\le v<2$.
Set $B=\Phi(x)^a\Phi(y)^b$ and $D=-\Phi'/\Phi$. Since $zb=a$,
\[
 \frac{A_zB}{a B}=\ell\big(\Phi(x)^{-v}+\Phi(y)^{-v}\big)
 -(a-1)D(x)^2-(b-1)D(y)^2\le2\ell B^{-r}.
\]
Indeed $D(x)\ge D(y)\ge0$, $a\ge1$, $a+b\ge2$, and $\Phi(x)\le\Phi(y)\le1$.
The conormal derivative of $B$ vanishes on the free sides. Both first derivatives are bounded,
using $a+b\ge2$ when $b<1$. Truncating at the vertices gives the weak
subsolution inequality. The comparison argument on p.~\pageref{par:fixed-regularity} therefore gives $V\le2a\ell B^{-r}$.
For $h=1-y$, $s=(x-y)/h$, concavity gives $B\ge h^{a+b}\Phi(s)^a$.
On each horizontal slice the interval ground-state inequality bounds the energy below by the squared norm with weight $\ell h^{-2}\Phi(s)^{-v}$. This weight is at least $V/(2a)$.
At $r=0$ instead use $X\ge P I$ and the trial $B$ to get $\Lambda_2\le2aP$.
\end{proof}

\paragraph{Two lower bounds for the normalized curvature.}\label{par:two-form}
We bound $G$ in two ways to obtain curvature estimates for all $1\le q\le2$.
Let $0<k<1$, $0\le w\le1$, $X\ge kW$ and
$W|_{e^\perp}\le wI$. Put $u_0=1-pw$, $b=1-k$ and
$D=b\min\{1,u_0/(1-kw)\}$. If $u_0>0$, then
\begin{align}
 \mathcal C_q(z)&\ge 1-\frac{[(1-z)D+4(1+z)b]^2}{16(1+z)D},
 \label{S:cc:twosource}\\
 \mathcal C_q(z)&\ge 1-(1-S)
       \frac{[(1-z)u_0+4(1+z)]^2}{16(1+z)u_0},
 \quad S\le\Lambda_q(0)/\Lambda_q(z).
 \label{S:cc:ranksource}
\end{align}
For $v\perp e$ we have $\langle(I-R-kW)e,v\rangle=-\langle Re,v\rangle$, since $We=e$. Both operators are nonnegative, and their quadratic forms at $e$ have values $b-\beta$ and $\beta$, respectively. Cauchy--Schwarz gives
$|\langle g,v\rangle|^2\le
\beta(b-\beta)\langle(I-kW)v,v\rangle/b$.
The minimum of $(1-pt)/(1-kt)$ on $[0,w]$ is $D/b$.
By the second-derivative formula \eqref{S:eq:response},
\begin{equation}\label{eq:G-weighted-bound}
 G\le\beta(b-\beta)/D .
\end{equation}
For the second estimate use $P_ev=\langle v,e\rangle e$.
The Euler equation gives
$|\langle e,v\rangle|\le\sqrt{\Lambda_q(z)}\|v\|_q$, so
$X\ge SP_e$. Applying the same two Cauchy--Schwarz inequalities to
$R$ and $I-R-SP_e$ gives
$\|g\|^2\le \beta(1-S-\beta)/(1-S)$.
Use $(I-pW)^{-1}\le I/u_0$ on the tangent space. Completing the
square in $1+(z-1)\beta-4(1+z)G$ gives both displays.

Let $\Phi=\Phi_q$ be the positive interval solution with $-\Phi''=\ell_q\Phi^{q-1}$, $\Phi(0)=1$ and $\Phi'(0)=\Phi(1)=0$. Its first integral gives
$d_q=\mathrm B(1/q,1/2)/q$ and $\ell_q=qd_q^2/2$.

\begin{lemma}[Projection onto the first horizontal cosine]\label{S:lem:horizontal-projection}
Put $P=\pi^2/4$ and $\eta_0=z/(1+z)$. For $0<\eta\le\eta_0$ let
$T_2(\eta)$ be the second eigenvalue of the form below on functions $f\in H^1(0,1)$ with $f(0)=0$,
\[
 \eta\left(\int_0^1 |f'|^2-|f(1)|^2\right)
 +P\int_0^1 h^{-2}f^2,\qquad \|f\|_2=1 .
\]
The second unweighted eigenvalue of $Q_z$ is at least $\min\{5P,T_2(\eta)\}$.
\end{lemma}
\begin{proof}
Write $h=1-y$, $s=(x-y)/h$, $c(s)=\sqrt2\cos(\pi s/2)$ and
$w=\sqrt h\,v=A(y)c+b$ with $b\perp c$ on each slice.
With $Tg=g/2-(1-s)g_s$, the energy divided by two is the integral of
$h^{-2}w_s^2+z(w_y+Tw/h)^2$ in $ds\,dy$, and the mass divided by
two is $\int(A^2+\|b\|_2^2)\,dy$. The next horizontal eigenvalue
is $9P$, so $E=\|b_s\|_2^2-P\|b\|_2^2\ge(8/9)\|b_s\|_2^2$.
As $b(1)=0$ and $\langle c,b\rangle=0$, the primitive
$C(s)=\int_0^s c$ satisfies $\int C b_s=0$. Hence
\[
 L=\langle c,Tb\rangle=-\int_0^1[(1-s)c-C/2]b_s,\qquad
 \|(1-s)c-C/2\|_2^2=\tfrac13+\pi^{-2}<\tfrac49 .
\]
Thus $L^2\le E/2$. Since $\langle c,Tc\rangle=c(0)^2/2=1$,
projection of the vertical energy onto $c$, followed by adding
$E/(2h^2)$, gives
\[
 z(A'+A/h+L/h)^2+(L/h)^2\ge\eta_0(A'+A/h)^2 .
\]
The remaining horizontal term is at least
$h^{-2}(PA^2+5P\|b\|_2^2)$. For $f(h)=A(1-h)$,
integration by parts gives $\int(A'+A/h)^2dy=\int|f'|^2dh-|f(1)|^2$.
Since $f/h$ and $f'-f/h$ both belong to $L^2$, we have $f\in H^1$ and
$f(0)=0$, so $f(h)^2/h\to0$. These slice estimates extend by density.
Requiring $f$ to be orthogonal to the first eigenfunction of the one-dimensional form imposes at most one linear condition on $v$. Since $h^{-2}\ge1$ and $\eta\le\eta_0$, min--max proves the claim.
\end{proof}
For the interval solution above,
$|\Phi_q'|^2=(2\ell_q/q)(1-\Phi_q^q)$ and multiplication of its
equation by $\Phi_q$ give $\int_0^1\Phi_q^q=2/(q+2)$.
After $t=e^{-v}$ in $d_q=\int_0^1(1-t^q)^{-1/2}\,dt$, the factor
$q/(1-e^{-qv})$ increases with $q$ by $e^{qv}>1+qv$.
Thus $\ell_q=qd_q^2/2\ge\ell_1=2$.

\begin{lemma}[The horizontal minimizer]\label{S:lem:moving}
Let $1\le q<2$, $r=2-q$, $B_q=2/(q+2)$ and $n=2/r$.
The horizontal quotient has minimizer
$v=h^n\Phi_q((x-y)/h)$, where $h=1-y$, and
\[
 \Lambda_q(0)=\ell_q[(q+2)/r]^{r/q},\qquad
 Y(v)/X(v)\le 1+n^2/\ell_q .
\]
For $q=2$ we have $\Lambda_2(0)=\pi^2/4$.
\end{lemma}
\begin{proof}
The unit-interval horizontal quotient is $\ell_q B_q^{-r/q}$.
Rescale to each horizontal interval of length $h$, then apply H\"older
in $2\,dy$. The transverse integral is $2\int_0^1 h^{(q+2)/r}\,dy=r/2$.
Equality holds for the displayed $v$, whose $q$-mass is $2B_q/(nq+2)$.
With $\sigma=(x-y)/h$ one has
\[
 X(v)=\ell_q B_q/n,\qquad
 Y(v)=n^{-1}\int_0^1[n\Phi_q+(1-\sigma)\Phi_q']^2\,d\sigma .
\]
The cross term has the sign of
$2\int(1-\sigma)\Phi_q\Phi_q'=\int\Phi_q^2-1\le0$.
Now $\int\Phi_q^2\le B_q$ gives the claimed bound.
At $q=2$ the mixed interval inequality gives the lower bound,
and the horizontal quotients of $h^n\cos(\pi\sigma/2)$ converge
to $\pi^2/4$ as $n\to\infty$.
\end{proof}

The same slice integrals for $v=h^n\Phi_q((x-y)/h)$ also give a
uniform rate for the horizontal quotient. Put $s=z^{1/3}$ and
$n=2/(r+2s)$. Then
\[
 \frac{X(v)}{\|v\|_q^2\Lambda_q(0)}
 =\left(\frac{r}{r+2s}\right)^{r/q}(1+s)^{2/q},\qquad
 zY(v)/X(v)\le z+s/2 ,
\]
using $n\le1/s$ and $\ell_q\ge2$. Thus
$1\le\Lambda_q(z)/\Lambda_q(0)\le(1+s)^{2/q}(1+s/2+z)$,
with the first factor in the display interpreted as one at $r=0$.

\begin{proposition}[Uniform rate at the flat angle]\label{S:prop:uniform-flat}
For all sufficiently small $z>0$,
$1-4z^{1/6}<\mathcal C_q(z)\le1$ uniformly for $1\le q\le2$.
For fixed $\eta>0$ $1-\mathcal C_q(z)=O_\eta(z)$ when $q\le2-\eta$.
\end{proposition}
\begin{proof}
The upper bound follows from \eqref{S:eq:curv}. For the lower bound write
$t=z^{1/3}$, $a=61/50$, $u=1+at$, $v=(1+t)^3$, $k=1/u$ and $b=1-k\le at$.
We use the horizontal bound $X\ge kW$ and estimate (W) on p.~\pageref{par:flat-transfer}. The proof of (W) uses the first two unweighted eigenvalues, the upper bound on the potential in Lemma~\ref{S:lem:linear-potential-transfer}, and the mixed Green kernel on the interval; none of these uses the present proposition. Here $\beta\le b$ since $X\ge kW$. With $d=\min\{1,(1-pw)/(1-kw)\}$, maximizing $4\beta(b-\beta)$ in \eqref{eq:G-weighted-bound} gives $4G\le b/d$.
If $r=2-q\ge t^{3/2}/2$, then $d\ge\min\{1,r/b\}$, so
\[
 1-\mathcal C_q\le b+(1+z)\max\{b,b^2/r\}
                 \le(2a^2+o(1))\sqrt t<4\sqrt t .
\]
If $r\le t^{3/2}/2$, then $r\le 2/\nu$ for small $t$, since
$\nu t^{3/2}\to\pi/2$. From (W), $\log w\le-(3-a-3\pi/8)t+O(t^{3/2})$,
with $3-a-3\pi/8>0$. Consequently $1-w\ge c t$
for a fixed $c>0$. Now $d\ge (1-w)/(1-w+b)$ is bounded below, so
$1-\mathcal C_q=O(t)$. Both estimates are uniform in $q$.

Finally, for
$r=2-q\ge\eta$ the horizontal minimizer of Lemma~\ref{S:lem:moving} has
$Y/X\le1+4/(\ell_q\eta^2)$. Hence testing with it gives
$\Lambda_q(z)/\Lambda_q(0)\le1+O_\eta(z)$;
\eqref{S:cc:ranksource} with $w=1$ proves the last assertion.
\end{proof}

We also use the first eigenfunction $E$ at $z=1/3$, defined below and normalized to have maximum one. Its energy is
\[
 Q_z(E)=4\pi^2/81+3/8+z(4\pi^2/27-9/8).
\]
Indeed $Q_{1/3}(E)=(16P/9)\int E^2$, and differentiating its
explicit formula followed by integration in $y$ gives
$Y(E)=(32\pi^2/27)\int_0^1\sin^2(\pi(1-x)/3)
[x/2-\sin(2\pi x)/(4\pi)]\,dx=4\pi^2/27-9/8$.
\paragraph{Computing the mass of the eigenfunction $E$.}
In the fixed triangle the test at $z=1/3$ is
$E=4\sin b(\cos b+\cos\pi y)/(3\sqrt3)$, $b=\pi(1-x)/3$.
This $E$ has maximum one. Successive even reflections give six congruent copies in the full equilateral triangle. There the reflected function is
$8\sin a\sin b\sin c/(3\sqrt3)$, with $a,b,c$ the positive affine
coordinates that sum to $\pi$ (each is $\pi/3$ times the distance to the
corresponding side). Reflection preserves the probability distribution of
$E$. With $A=2a$, $B=2(a+b)$ the triangle becomes $0<A<B<2\pi$, and
$E=|\Delta(1,e^{iA},e^{iB})|/(3\sqrt3)$, where
$\Delta(X,Y,Z)=(X-Y)(X-Z)(Y-Z)$. Symmetry in $A,B$ identifies
the normalized triangle integral with the normalized torus integral
on $[0,2\pi]^2$. The six distinct Fourier monomials of $\Delta$ have
coefficients $\pm1$, hence $\int |\Delta|^2=6$ and $\int_D E^2\,d\mu=2/9$.

\paragraph{A lower bound for the second eigenvalue in the longitudinal variable.}
For $0<z\le1/3$ write $P=\pi^2/4$, $t=z^{1/3}$ and $\eta=z/(1+z)$.
The second eigenvalue of the quadratic form
\[
 \eta\!\int_0^1 |f'|^2+P\!\int_0^1 h^{-2}f^2-\eta |f(1)|^2,
 \qquad f(0)=0,\qquad \|f\|_2=1 ,
\]
is at least $P(1+t)^3$. Indeed, put $v=P/\eta$ and use the tangent
$h^{-2}\ge d^2(3-2d)+2d^3(1-h)$ for $d\ge1$.
Reflect at $h=1$ and extend by zero beyond the interval. After
dilation, the second eigenvalue divided by $\eta$ is bounded below by
\[
 v d^2(3-2d)+(2v)^{2/3}d^2 e_2,
\]
where $e_2$ is the second eigenvalue on the half-line of
$\int(|g'|^2+y g^2)-\sigma|g(0)|^2$, with
$\sigma=[d(2v)^{1/3}]^{-1}<5/12$ since $v\ge3P$.

Impose the single condition $g(3/4)=0$. On the left, the trace bound
$|g(0)|^2\le(3/4)\int |g'|^2$ and mixed interval Poincar\'e give
a lower bound $(1-5/16)(16P/9)=11P/9>29/10$.
On the right translate $y$ by $3/4$. The positive Dirichlet test
$G(y)=y\exp[-\sqrt2\,y^{3/2}/3]$ gives
\[
 y-G''(y)/G(y)=y/2+5/(2\sqrt{2y})\ge3(25/64)^{1/3}>13/6 .
\]
The ground-state inequality adds this bound to $3/4$. Thus $e_2\ge29/10$
by the min--max principle, since $g(3/4)=0$ imposes just one linear condition. Taking $e=29/10$ and $d=1+2^{2/3}e/(3v^{1/3})$
gives the lower bound $P d^3$. Finally $\eta/z\ge3/4$ and $(29/10)^3>9P$ (use $P<5/2$), so $d\ge1+t$.

\subsection{Comparison with a disk or circular sector}\label{S:sec:exponent}
If $B$ is a disk, write $c_q=|B|^{2/q}\lambda_q(B)/(2\pi)$.
A circular sector of opening $\alpha$ with Dirichlet arc has normalized
quotient $\alpha c_q$, by radial rearrangement in the wedge.
Use $\alpha=2\pi$ for the full plane. The relative perimeter constant
is $2\alpha$ both for a wedge of opening $0<\alpha<\pi$ and, under
this convention, for the plane. This gives a single exponent comparison.
\begin{proposition}[Comparison with a disk or circular sector as the exponent varies]\label{S:prop:exponent-round}
Let $G$ be a bounded connected Lipschitz polygon in the plane, or in a
wedge of opening $0<\alpha<\pi$. Impose Dirichlet data on its boundary
inside the ambient open set and Neumann conditions on the wedge rays.
Write $S_A$ for the disk, respectively the round sector centered at the
apex, of the same area $A$. Then
\[
  p\longmapsto \bigl(\lambda_p(G)/\lambda_p(S_A)\bigr)^p
\]
is strictly increasing on $[1,2]$. For a disk in the plane or a circular sector centered at the wedge apex, the ratio is one. In the triangle case write
$\Phi_p=F_p(\alpha)/(\alpha c_p)=\lambda_p(I_\alpha)/\lambda_p(S_A)$.
\end{proposition}
\begin{proof}
Write $I=G$. We use the single-crossing argument of
Carroll and Ratzkin~\cite[Sec.~4]{S:carroll-ratzkin}.
Fix $1<p<2$. Scale $S_A$ to a domain $S$ for which
$\lambda_p(S)=\lambda_p(I)=\lambda$. Its area $b$ is strictly less than $A$. Indeed radial
rearrangement preserves $p$-mass and lowers the energy. Equality would
force almost every positive superlevel set to attain the relative
isoperimetric inequality. Since $0<\alpha<\pi$, the wedge is an open convex cone containing no line.
By Figalli and Indrei~\cite[Thm.~2.2]{figalli-indrei2013}, each such level set agrees up to a null set
with a round sector centered at the apex. See also Pacella and Tralli~\cite[Thm.~3.7]{S:pacella-tralli}
for the corresponding rearrangement inequality. In the plane the equality sets are disks. Taking a
nested sequence of levels tending to zero would show that $I$ agrees almost everywhere with a circular sector in the wedge or a disk in the plane. In the plane, the centers of the nested disks converge because their distance is at most the difference of their radii. This excludes a
polygon. Now use $\lambda_p(kI)=k^{-4/p}\lambda_p(I)$.
Let $u$ and $v$ be positive minimizers on $I$ and $S$, each normalized
to have $p$-mass one, with $v$ radial and decreasing. Denote by $U(s)$
and $V(s)$ their decreasing rearrangements in the area variable.

The relative perimeter inequality $P(E)^2\ge 2\alpha |E|$, applied to the positive superlevel sets after zero extension across the Dirichlet boundary, yields
\[
 2\alpha s(-U'(s))\le \lambda\int_0^s U(t)^{p-1}\,dt ,
 \qquad
 2\alpha s(-V'(s))= \lambda\int_0^s V(t)^{p-1}\,dt .
\]
In fact, for $E=\{u>h\}$ of area $s$, its boundary flux in the
ambient open set is $\lambda\int_Eu^{p-1}$. Cauchy--Schwarz and coarea
give $P(E)^2\le(-ds/dh)\lambda\int_Eu^{p-1}$. The coarea and
inverse-distribution formulas apply at almost every level; positivity
and interior analyticity exclude critical sets of positive area.
The equality for $V$ follows by radial integration.

Since $V(b)=0<U(b)$ and both total $p$-masses equal one, the graphs
cross in $(0,b)$. Let $s_0$ be their last crossing. We claim $V\ge U$
on $(0,s_0)$. Otherwise put $W=\max(U,V)$ up to $s_0$ and $W=V$
on $[s_0,b]$. This is continuous and decreasing, with
$M=\int_0^bW^p>1$. Wherever $W$ follows $U$ before $s_0$, its
preceding values dominate $U$; everywhere that it follows $V$, its
preceding values dominate $V$. It therefore obeys the same
differential inequality as $U$. Its radial representative on $S$ has
zero boundary value on the arc and energy
\[
 2\alpha\int_0^b s|W'|^2\,ds
 \le \lambda\int_0^b (-W'(s))\int_0^s W(t)^{p-1}\,dt\,ds
 =\lambda M .
\]
The rearrangements have finite local Sobolev energy by wedge
symmetrization, so the maximum and gluing above are admissible.
Its quotient is at most $\lambda M^{1-2/p}<\lambda$, a contradiction.
Extend $V$ by zero to $[b,A]$. Thus $V\ge U$ before $s_0$ and
$V\le U$ afterwards. In particular,
\[
 \int_Iu^p\log u<\int_Sv^p\log v .
\]
In fact the same comparison holds, weakly, for any convex function of
the $p$-mass density for which the integrals exist. Subtract a supporting
line at $d=U(s_0)^p$ from the convex function. The remainder decreases
up to $d$ and increases thereafter; at each point it is therefore larger
at $V^p$ than at $U^p$. The line integrates to the same value by mass
normalization. For $x\log x$ strictness holds on $(b,A)$.

For any one of these domains $G$, let $K_G(r)=r\log\lambda_r(G)$.
It is the infimum of the concave functions
$r\log\int_G|\nabla w|^2-2\log\int_G|w|^r$ over nonzero admissible
$w$. It is finite for $0<r<\infty$ by the planar Sobolev embedding
and Poincar\'e inequality, hence locally Lipschitz. At a
differentiability point $p$ the derivative agrees with that of the
function inside the infimum evaluated at any minimizer $w$ with unit $p$-mass,
\[
 K_G'(p)=\log\lambda_p(G)-2\int_G w^p\log w .
\]
Take a common differentiability point of $K_I$ and $K_{S_A}$.
The dilation chosen at this $p$ changes $K_{S_A}(r)$ by a constant
independent of $r$, so $K'_S(p)=K'_{S_A}(p)$. The identity $\lambda_p(S)=\lambda_p(I)$ and the strict inequality $\int_Iu^p\log u<\int_Sv^p\log v$ proved above give
$(K_I-K_{S_A})'(p)>0$. Integration proves strict increase on $[1,2]$.

\end{proof}

The round sector of radius $R$ has torsion function $(R^2-|x|^2)/4$,
so $T(S_A)=A^2/(4\alpha)$ and $c_1=4$. Applying
Proposition~\ref{S:prop:exponent-round} between exponents one and $q$ gives
$R_q(s)^q\ge R_1(s)$, where $R_q(s)=F_q(2\pi s)/(2\pi s\,c_q)$.

\paragraph{The torsion ratio of the disk to the square.}
Reflection of the four mixed cells gives
$R_1(1/4)=T(B)/T(S)$ for $S=(-1,1)^2$ and the disk $B$ of area four;
the four comparison sectors form $B$. On $S$ the field
\[
 J=-\tfrac12(x,y)+u(x^3/3-xy^2,\ y^3/3-x^2y)
\]
has divergence $-1$. Its mean squared norm is
$1/6+4u/45+8u^2/105$, by the even moments $1/(2k+1)$ of each coordinate.
Taking $u=-1/2$ in the dual variational principle for torsion gives $T(S)\le4(89/630)$.
Since $T(B)=2/\pi$ and $\pi<22/7$, this gives
\begin{equation}\label{eq:square-flux-ratio}
 R_1(1/4)\ge315/(89\pi)>9/8 .
\end{equation}

\subsection{Boundary comparisons for the fixed triangle}\label{S:sec:boundary-inputs}
In this subsection, $W(v)=\int_D\lambda u^{q-2}v^2\,d\mu$,
$\lambda=\Lambda_q(z)$ and $\beta=z\partial_z\log\Lambda_q(z)$.
\begin{lemma}[Coordinate monotonicity]\label{S:lem:coordinate-monotonicity}
On $0<z\le1$ the positive minimizer $u$ decreases in both coordinates for $1\le q\le2$.
\end{lemma}
\begin{proof}
Work in the physical triangle $0<h<m\xi$, $\xi<1$, $m=z^{-1/2}$.
First let $0<p=q-1<1$, put $V=\lambda u^{p-1}$, and set
$w=-u_\xi$, $v=-u_h$. Both satisfy $(-\Delta-pV)w=0$,
$(-\Delta-pV)v=0$ in the interior. On the sloping side
$v=mw$, and $\partial_n w+m\partial_n v=0$ by tangential
differentiation of $\partial_n u=0$. On the bottom $v=0$,
$\partial_n w=0$. At the vertical Dirichlet edge $v=0$ and $w\ge0$.
Test by $-w_-$ and $-v_-$ and add: the sloping-edge terms cancel,
and the ground-state inequality $\int Vb^2\le\int|\nabla b|^2$ for zero Dirichlet trace gives
\[
 (1-p)\int(|\nabla w_-|^2+|\nabla v_-|^2)\le0.
\]
Here $w,v\in H^1$: the Poisson source belongs to $L^2$, the apex
is a convex corner with Neumann data on both sides, and even reflection across the Neumann side gives a convex
Dirichlet corner at the upper endpoint and a flat Dirichlet point at the lower.
Local $H^2$ regularity for pure data on convex polygonal corners applies \cite[Chapter~3]{S:grisvard}.
For the tests one may first truncate the negative parts at $M,mM$,
respectively, preserving their sloping-edge ratio. Logarithmic cutoffs at the
vertices have vanishing Dirichlet energy, so the cutoff-gradient terms vanish
for these bounded tests. Then let $M\to\infty$, using the ground-state
bound on each negative part. On open edges the integrations are justified
by reflection across the Neumann edges and by the Dirichlet trace on the vertical edge
(the Hopf derivative there is positive and $V=O((1-\xi)^{p-1})$ locally).
The zero traces force $w_-=v_-=0$. Strong energy convergence as $q\to1,2$
extends the inequalities to the endpoints.
\end{proof}

\begin{lemma}[An upper bound for $\beta_q(1)$]\label{S:lem:square-slope}
For \(1\le q\le2\),
\[
  \beta_q(1)\le \frac34-\frac{\pi^2}{16}<\frac2{15}.
\]
\end{lemma}
\begin{proof}
Reflect the minimizer in the diagonal of \(D=\{0<y<x<1\}\).
On the square it solves \(-\Delta u=\lambda u^p\), \(p=q-1\),
with Neumann data at \(x=0,y=0\) and Dirichlet data at \(x=1,y=1\).
Write \(k=\pi/2,\ P=k^2,\ w=-u_x,\ V=p\lambda u^{p-1}\) and
\(H=-\Delta-V\). For \(0<p<1\) we first prove
\begin{equation}
 w>0,\qquad w_y\le0,\qquad
 w_x\ge k\cot(kx)w,\qquad
 |u_y|/|u_x|\le\tan(ky)/\tan(kx)\quad (y<x).                 \label{S:eq:local-s1}
\end{equation}

We record the boundary estimates used in the comparisons. Even
reflection in a Neumann side and odd reflection in a Dirichlet side
give a weak equation with source
\(\lambda\operatorname{sgn}(u)|u|^p\). Interior elliptic regularity
first gives a locally Lipschitz solution, then \(C^{2,p}\).
At an open Dirichlet side use inward normal and tangential coordinates
\(s,t\), and put \(A(t)=u_s(0,t)>0\) by Hopf's lemma.
There \(u_t\) also belongs to \(C^{2,p}\): its source is
\(p\lambda s^p(u/s)^{p-1}(u_t/s)\), with both quotients \(C^{0,p}\)
and \(u/s\) bounded away from zero. Both \(u\) and \(u_t\) have zero
second normal derivative on the side. Consequently
\[
 u=sA+O(s^{2+p}),\quad u_s=A+O(s^{1+p}),\quad
 u_t=sA'+O(s^{2+p}),\quad u_{st}=A'+O(s^{1+p}).
\]
At the double Dirichlet corner put \(s=1-x,\ t=1-y\) and
\(\rho=(s^2+t^2)^{1/2}\). There \(u_{st}=A_*+O(\rho^p)\) with
\(A_*>0\): on a small corner square compare \(u\) with \(cst\),
choosing \(c>0\) on the two inner sides by Hopf's lemma.
Integration gives \(u=st[A_*+O(\rho^p)]\),
\(u_s=t[A_*+O(\rho^p)]\), \(u_t=s[A_*+O(\rho^p)]\).
In particular \(u\ge c(1-x)(1-y)\) on the square, so \(V\in L^a\)
for some \(a>1\); its form is continuous on \(H^1\).

The positive solution \(u\) gives
\(\int(|\nabla v|^2-Vv^2)\ge(1-p)\int|\nabla v|^2\)
when \(v\) vanishes on the Dirichlet sides, by testing with
\(v^2/u\) (or regularizing \(u\) by \(u+\epsilon\)).
Now \(Hw=0\), \(w=0\) at \(x=0,y=1\), \(w_y=0\) at \(y=0\),
and \(w\ge0\) at \(x=1\) by Hopf's lemma. Testing with its negative
part gives \(w\ge0\), strictly positive inside. Symmetry gives \(u_y<0\).
In addition \(w_x=0\) at \(x=1\), since \(u_{xx}=-u_{yy}-\lambda u^p=0\)
there. The ground-state identity with \(w\) shows that the form of
\(H\) is also nonnegative for tests vanishing on \(x=0,y=1\).

Put \(f=\log u,\ Z=f_{xy},\ r=1-p\). Twice differentiating
\(\Delta f+|\nabla f|^2+\lambda u^{-r}=0\) gives
\[
 -u^{-2}\operatorname{div}(u^2\nabla Z)
 +[(2+r)\lambda u^{-r}+2|\nabla f|^2]Z
 =r^2\lambda u^{-r} f_xf_y\ge0 .
\]
Here \(Z=0\) on the Neumann sides. Testing with the cutoff negative
part gives \(Z\ge0\): the cutoff error is bounded by a constant times
\(\int u^2Z_-^2|\nabla\chi|^2\).
Indeed \(uZ=u_{xy}-u_xu_y/u=O(s^{1+p})\) at an open Dirichlet
side and \(O(\rho^p)\) at the corner. For fixed corner cutoff first
let the edge-strip width tend to zero, then shrink the corner annulus.
Thus \(u_{xy}\ge u_xu_y/u>0\), proving \(w_y\le0\).

Write \(b=k\tan(ky)\), \(d=k\cot(kx)\), \(T=w_y+bw\), \(L=w_x-dw\).
Commutation using \(b'=P+b^2,\ d'=-P-d^2\) gives
\[
 (H+2b')T=V_yw\ge0,\qquad (H-2d')L=V_xw\ge0 .
\]
Both added potentials are at least \(2P\). The function \(T\) vanishes on \(x=0,y=0,y=1\); reflect evenly across \(x=1\), where \(w\) and
\(T\) are even by the odd reflection of \(u\). The differentiated
equation holds weakly there: \(V_y=O(|1-x|^{p-1})\) away from vertices.
The function \(L\) vanishes on \(x=0,x=1,y=1\), with even reflection at
\(y=0\). To justify the negative-part tests, \(T=O(s^{1+p})\) on
approaching \(y=1\), while \(L=O(s)\) there, by the tangential regularity
above. At \(x=1\) we have \(L=\lambda u^p+O(s)>0\) locally, since
\(u_{yy}=O(s)\) and \(d=O(s)\). On \(x=0\), \(w\) is odd and the
odd Taylor expansion gives \(L=O(x^2)\) away from
vertices. At the double Dirichlet corner both \(T\) and \(L\) are
\(O(\rho^p)\). At a mixed corner the even reflection gives the same
flat-edge estimates; for \(L\) near \(x=0,y=1\) one can also use
\(w_x-w/x=O(x^p)\) from \(C^{1,p}\) regularity. Hence an annular
cutoff error is \(O(\rho^{2p})\). Testing against the cutoff negative
parts and using the ground-state inequality for the positive solution \(w\) of \(Hw=0\)
now gives \(T,L\ge0\).

Finally put \(a=k\tan(kx)\) and \(B=a u_y-bu_x\) on \(D\). Then
\[
 (H+2a')B=2(b^2-a^2)(u_{xy}+bu_x)=2(a^2-b^2)T\ge0 .
\]
The traces on \(y=0,y=x\) are zero. At \(x=1\),
\(B=T(1,y)+O(s^{1+p})\), and \(B=O(\rho^p)\) at the double
Dirichlet corner. Its cutoff negative part is therefore admissible for
the ground-state inequality applied to \(u\) on \(D\), proving \(B\ge0\)
and completing \eqref{S:eq:local-s1}. For \(p=1\), \eqref{S:eq:local-s1} follows directly from
\(u=\mathrm{const.}\cos(kx)\cos(ky)\).

It remains to average the pointwise gradient ratio. Under the measure
\(d\nu=\sin^2(kx)\,dx\,dy\) on \(D\), set
\(g=w^2/\sin^2(kx)\) and \(s=\tan(ky)/\tan(kx)\).
By \eqref{S:eq:local-s1}, \(g\) increases in \(x\) and decreases in \(y\).
Conditionally on \(s\), the densities of \(t=\tan(kx)\) and \(v=st\)
are proportional, respectively, to
\[
 \frac{t^3}{(1+t^2)^2(1+s^2t^2)},\qquad
 \frac{v^3}{(s^2+v^2)^2(1+v^2)} .
\]
For \(s_2>s_1\) the ratios of these conditional densities at \(s_2\)
to those at \(s_1\) are respectively decreasing and increasing.
Each ratio crosses one at most once by monotonicity. Since the densities in each pair have the same integral, their difference must take both signs unless they coincide. Integrating the difference therefore orders the distribution functions, and hence their quantiles. At a common quantile
level, \(t\) decreases and \(st\) increases as \(s\) increases.
Consequently the conditional mean of \(g\) decreases in \(s\);
truncating \(g\) first justifies the same statement if it is unbounded.
Chebyshev's integral inequality gives
\(\int s^2g\,d\nu\le\int s^2\,d\nu\,\int g\,d\nu/\nu(D)\).
By \eqref{S:eq:local-s1}, \(Y/X=\int_D u_y^2/\int_D u_x^2\le\int s^2\,d\nu/\nu(D)\).
The two elementary integrals are
\[
 \nu(D)=\int_0^1x\sin^2(kx)\,dx=\tfrac14+\pi^{-2},\qquad
 \int s^2\,d\nu=\int_0^1\cos^2(kx)[\tan(kx)/k-x]\,dx
                =3\pi^{-2}-\tfrac14 .
\]
Thus \(\beta_q(1)=Y/(X+Y)\le3/4-\pi^2/16<2/15\), using
\(\pi>3.1415\). As \(q\downarrow1\), normalized positive minimizers
converge strongly in \(H^1\): compactness gives a normalized limit,
testing with the torsion minimizer gives the matching upper energy
bound, and uniqueness identifies the limit. The inequality passes to
\(q=1\).
\end{proof}

\section{Wide openings}\label{S:sec:angle-replacement}
We call an opening wide if $2\pi/3\le\alpha<\pi$, that is, if $0<z\le1/3$.

\begin{proposition}[Curvature at wide openings for $q=2$]\label{S:prop:wide-linear}
For \(q=2\) and \(0<z\le1/3\), \(\mathcal C_2(z)>1/16\).
\end{proposition}
\begin{proof}
Write \(t=z^{1/3}\), \(P=\pi^2/4\), and let \(\lambda,\lambda_2\)
be the first two unweighted eigenvalues of \(Q_z=X+zY\) in \(L^2(D,\mu)\).
We prove $\lambda/P\le1+61t/50$ and $\lambda_2/P\ge(1+t)^3$.
For the upper bound set $h=1-y$, $s=(x-y)/h$, $m^3=9P/(2z)$,
and test with $A(h)\cos(\pi s/2)$, where
$A=(1+m\tau)e^{-m\tau}$ and $\tau=-\log h$. In the variable $\tau$, both energy integrals are moments of exponentials. Writing
$I=\int_0^\infty A^2\,d\tau$, $J=\int_0^\infty A_\tau^2\,d\tau$ and
$M=\int_0^\infty e^{-2\tau}A^2\,d\tau$, we obtain
\[
 I=\frac5{4m},\qquad J=\frac m4,\qquad M=\frac{5m^2+6m+2}{4(m+1)^3}.
\]
Integration in $s$ gives the quotient $R=[PI+z(J-1/2+(P/3-1/2)I)]/M$.
With $d=1/m<1/3$ and $P<5/2$, this yields
\[
 \frac RP\le\frac{(1+d)^3(5+9d/2-9d^2+15d^3/2)}{5+6d+2d^2}
 \le 1+68d/25 .
\]
For the second inequality, multiply the difference by $(5+6d+2d^2)100/d$.
The result is $10-118d+644d^2-1350d^4-750d^5\ge1+(20d-3)^2$,
since $1350d^2+750d^3<180$. Finally $68d/25<(61/50)t$
because $(136/61)^3<9P/2$ by $\pi>3.14$.
Take $\eta=z/(1+z)$ in Lemma~\ref{S:lem:horizontal-projection}.
The estimate for the longitudinal problem gives $T_2(\eta)\ge P(1+t)^3$.
Since $(1+t)^3<5$ by $t<7/10$, it follows that $\lambda_2\ge P(1+t)^3$.

Write $a=61/50$ and let $u=\lambda/P$, $v=(1+t)^3$, $b=1-1/u$.
Since $Q_z\le Q_{1/3}$, the ground state at $z=1/3$ gives
$u\le\lambda(1/3)/P=16/9$ by Lemma~\ref{S:lem:eqspec}.
We also have $u\le1+at<v$. In \eqref{S:cc:twosource} use $k=1/u$,
$w=u/v$ and $d=(1-w)/(1-kw)=(v-u)/(v-1)$. Expanding the square
in \eqref{S:cc:twosource} and using $d\le1$ gives
\[
 1-\mathcal C_2\le b\left[H+(1+z)\frac{u-1}{v-u}\right],
 \qquad H=\frac{(5+3z)^2}{16(1+z)}\le\frac{27}{16}.
\]
Here $b\le7/16$, so $bH\le189/256$. Since $v-u\ge t(3-a+3t)$,
the second term, including the factor $b$, is at most
\[
 \frac{4a^2}3\frac{t}{(1+at)(3-a+3t)}
 \le\frac{4a^2}{3(\sqrt3+\sqrt{a(3-a)})^2}
 <\frac{39}{200}.
\]
The middle inequality is AM--GM applied to $(3-a)/t+3at$.
For the last bound use $\sqrt3>173/100$ and $\sqrt{a(3-a)}>147/100$.
Hence $\mathcal C_2>1-189/256-39/200>1/16$.
\end{proof}

\begin{lemma}[An upper bound on $V_q$ from the unweighted ground state]\label{S:lem:linear-potential-transfer}
Let $z>0$, $r=2-q$, $P=\pi^2/4$ and $\lambda=\Lambda_2(z)$.
Write $\phi$ for the positive first eigenfunction for $Q_z$, normalized to have maximum one, and put
$\gamma=1/9+2/75+8/735$. Then
\begin{equation}\label{S:eq:intervalu}
 V_q\le\lambda e^{-\gamma r}\phi^{-r}.
\end{equation}
For $0<z\le1$ put $\nu=\pi/(2\arctan\sqrt z)$, $h=1-y$ and $l=(1-x)/h$.
On this range $\phi\ge h^\nu l$.
\end{lemma}
\begin{proof}
In physical coordinates
$\mathcal P=|\nabla\phi|^2+\lambda\phi^2\le\lambda$.
Indeed, away from critical points, take the first coordinate direction
along the gradient at the point of differentiation. Then
\[
 \Delta\mathcal P+
 \frac{2\lambda\phi}{|\nabla\phi|^2}\nabla\phi\cdot\nabla\mathcal P
 =4[(\phi_{11}+\lambda\phi)^2+\phi_{12}^2]\ge0 .
\]
An interior maximum exceeding $\lambda$ would force constancy on its
component of noncritical points, by the strong maximum principle.
Even and odd reflection give the same argument on open Neumann and
Dirichlet sides, respectively, and at the lower right vertex after even
reflection across the bottom side. The gradient tends to zero at the Neumann apex
and upper mixed vertex: in their physical sectors the first nonconstant Neumann
mode and the first mixed mode have degrees $\pi/\arctan(1/\sqrt z)>2$
and $\pi/(2\arctan\sqrt z)>1$, respectively. Separation of the Helmholtz
equation gives these orders for the regular Bessel modes, and the radial
Neumann mode has gradient $O(R)$ at distance $R$ from the vertex.
A component on which $\mathcal P$ is constant at a value above $\lambda$ could
therefore have no relative boundary, contradicting the apex value.

Apply Lemma~\ref{lem:linear-hemisphere-profile} with $d=1$, $M=\lambda$
and $J=\phi$ to get \eqref{S:eq:intervalu} at every opening.

Now assume $z\le1$. Lemma~\ref{S:lem:coordinate-monotonicity} places the
maximum at the apex. To bound $\phi$ from below, write $\theta=\arctan\sqrt z$ and put $g(s)=\operatorname{Im}(1+i\sqrt z\,s)^\nu$.
For $0<s\le1$ and $\psi=\arctan(\sqrt z\,s)$,
\[
 0\le-g''/g=z\nu(\nu-1)
 \frac{\sin((\nu-2)\psi)}{(1+zs^2)\sin(\nu\psi)}
 \le z\nu(\nu-1)\le P .
\]
Here $0\le(\nu-2)\psi\le\nu\psi\le\pi/2$,
$\sin\theta\le\theta$ and convexity on $[0,\pi/4]$ gives
$\sin^2\theta\le2\theta/\pi$.
Thus $g$ is concave, and the Wronskian with $\sin(ks)$, $k=\pi/2$,
gives $l\le g(l)/g(1)\le\sin(kl)$. The function $h^\nu g(l)/g(1)$
is harmonic for $Q_z$, vanishes on the base, and has zero conormal derivative on the diagonal.
On the bottom it is at most $\cos(kx)\le\phi(x,0)$ by the derivative comparison in Lemma~\ref{S:lem:horizontal-sturm}.
Comparison with the superharmonic $\phi$ gives $\phi\ge h^\nu l$.
\end{proof}

\subsection{Horizontal comparison with the first eigenfunction}
Write $c(t)=\cos(\pi t/2)$, $P=\pi^2/4$, $\lambda=\Lambda_2(z)$.
Normalize the unweighted first eigenfunction $\phi$ by $\phi(0,0)=1$.
\begin{lemma}[Horizontal Sturm comparison]\label{S:lem:horizontal-sturm}
For $0<z\le1$, both coordinate derivatives of $\phi/[c(x)c(y)^{1/z}]$ are
nonnegative. For the positive normalized minimizer $u_q$ and its potential
$V_q=\Lambda_q u_q^{q-2}$, one has, simultaneously for $1\le q\le2$,
\[
  \kappa_q(z):=\inf_{v\ne0}\frac{X(v)}{\int_D V_q v^2\,d\mu}\ \ge\ \frac P\lambda .
\]
Equality holds for $q=2$, and the inequality is strict for $q<2$.
\end{lemma}
\begin{proof}
First impose $\partial_x\phi_k+a(1)\phi_k=0$ at the base, where
$0<k<\pi/2$ and $a(t)=k\tan(kt)$. Let $\phi_k$ be the positive unit
$L^2$ ground state with eigenvalue $\lambda_k$. Put
$B=\cos(kx)\cos(ky)^{1/z}$, $m=B^2$, $h=\phi_k/B$ and $A=\operatorname{diag}(1,z)$.
Both $B$ and $\phi_k$ satisfy the same Robin condition on the base and have zero conormal derivative on the other sides.
Thus $h$ has homogeneous Neumann data for
\[
 H=-m^{-1}\operatorname{div}(m A\nabla\,\cdot\,)+M-\lambda_k,\qquad
 M=2k^2-(z^{-1}-1)a(y)^2 .
\]
Since $h>0$ and $Hh=0$, the quadratic form of $H$ is nonnegative.
The first-order terms are $2a(x)\partial_x+2a(y)\partial_y$. Hence
\[
 (H+2a'(x))h_x=0,\qquad (H+2a'(y))h_y=-M_yh\ge0 .
\]
Test by $-(h_x)^-$ and $-z(h_y)^-$ in $L^2(m\,d\mu)$ and add.
On the diagonal $h_x=zh_y$ and $Nh_x+Nh_y=0$ for $N=-\partial_x+z\partial_y$,
by tangential differentiation of $Nh=0$. The boundary terms cancel there.
At the base $h_x=\partial_xh_y=0$, and at the bottom $h_y=\partial_yh_x=0$.
These tests are legitimate since $B$ is bounded away from zero:
in physical coordinates $h$ solves a Neumann Poisson equation with $L^2$
right-hand side on a convex triangle, hence $h\in H^2$
\cite[Chapter~3]{S:grisvard}. One can truncate the negative parts at heights
$zK,K$ and remove the vertices by logarithmic cutoffs, then pass to the limit.
Nonnegativity of the two forms and $a'>0$ give $h_x,h_y\ge0$.

As $k\uparrow\pi/2$ the Rayleigh bound $\lambda_k\le\lambda$ bounds the
energy and forces the squared base trace to be $O(1/a(1))$. Compactness
and the Dirichlet variational principle identify the unit $L^2$ limit as
the ground state with Dirichlet data on the base. On interior compact sets the denominators
converge uniformly to $c(x)c(y)^{1/z}>0$. Both distributional derivative
inequalities pass to the limit, and interior regularity gives the assertion.

For bounded positive nondecreasing \(b=b(x)\), let \(\lambda_D[b]\)
be the first eigenvalue on \(D\) for energy \(Q_z\) and mass
\(\int_D b v^2\,d\mu\). Let \(\lambda_I[b]\) denote the analogous
unit-interval eigenvalue for \(\int_0^1|f'|^2\), with \(f(1)=0\).
For \(0<z\le1\),
\[
              \lambda_D[b]\ \le\ \frac{\lambda}{P}\lambda_I[b].
\]
Compare the mass of $\phi$ along a vertical slice with the interval density $c(x)^2$ by setting
\[
  J(x)=\frac{\int_0^x\phi(x,y)^2\,dy}{c(x)^2}
      =\int_0^x c(y)^{2/z}h(x,y)^2\,dy .
\]
It is strictly increasing by the derivative comparison just proved and the positive
term at its moving upper limit. It is bounded up to \(x=1\), since
\(\phi\) is Lipschitz and vanishes on the base.

Write \(f=cg\) for the positive interval eigenfunction and
\(\eta=\lambda_I[b]\). Its equation gives
\((c^2g')'=(P-\eta b)c^2g\). The flux vanishes at both endpoints,
and \(P-\eta b\) changes sign at most once, from positive to
negative. Hence \(g'\ge0\). Multiplying the equation by \(Jg\)
and integrating gives
\[
 \eta\int_0^1 b c^2Jg^2
  = P\int_0^1 c^2Jg^2+\int_0^1 c^2J(g')^2+\int_0^1 c^2gg'\,dJ
  \ \ge\ P\int_0^1 c^2Jg^2+\int_0^1 c^2J(g')^2 .
\]
The endpoint terms vanish since \(g\) is bounded and \(c^2g'\)
vanishes there. Testing the triangular quotient with \(\phi g(x)\)
gives its energy, divided by two, as
\(\lambda\int c^2Jg^2+\int c^2J(g')^2\). Also \(\lambda>P\)
by horizontal Poincar\'e on the slices of lengths \(1-y\).
Multiplying the preceding bound by \(\lambda/P\) proves the claim.

For \(q<2\), write \(r=2-q,\ p=q-1\) and
\(U=\Lambda_q(z)^{-1/r}u_q\), so
\(-\operatorname{div}(A\nabla U)=U^p\) and \(V_q=U^{-r}\).
Lemma~\ref{S:lem:coordinate-monotonicity}
gives $U_x,U_y\le0$. Both are strictly negative in the interior, since
each negative derivative is a nonnegative supersolution of $-\operatorname{div}(A\nabla\,\cdot\,)$, and an identically vanishing derivative would force a constant minimizer by the diagonal condition.
Put \(b(x)=V_q(x,0)\), an increasing weight. For any cut \(y=a\),
put \(l=1-a\). The part above this cut is a copy of \(D\) under
\((x,y)\mapsto(a+lx,a+ly)\). Under this map, its positive solution at $(a+lx,a+ly)$ equals $l^{2/r}U(x,y)$. The restriction of \(U\) is a
supersolution because its outward derivative at the cut is
nonnegative. Call the solution on the subtriangle $U_a$. Test the difference of the equations against $(U_a-U)_+$. The tangent inequality for $t^p$ at $U$ and the ground-state inequality for $U$ give $U_a\le U$. In terms of $V_q$ this reads
\[
              l^2 V_q(a+lt,a)\le b(t),\qquad 0<t<1 .
\]
Applying the interval quotient on each horizontal slice and integrating
therefore gives \(\kappa_q(z)\ge\lambda_I[b]\).

The comparison with the interval eigenvalue $\lambda_I[b]$, proved above for bounded $b$, extends to this \(b\) by monotone truncation.
Indeed \(b(x)\le C(1-x)^{-r}\) by Hopf's lemma at the bottom of
the base after even reflection. In a strip \(1-x<\delta\), the
weighted mass is bounded by \(C\delta^{2-r}\) times the horizontal
energy, both on the interval and on the triangle. This follows by
integrating \(|v(x,y)|^2\le(1-x)\int_x^1 |v_x(t,y)|^2dt\).
Thus the weighted embeddings are compact, the first eigenvalues are
attained, and the truncated quotients converge.

The first eigenvalue on \(D\) with mass \(V_q\) is exactly one,
by the equation of \(U\) and the ground-state inequality. Since
\(V_q(x,y)>b(x)\) in the interior, testing with the positive
minimizer for \(b\) gives \(\lambda_D[b]>1\). Consequently
\[
 \kappa_q(z)\ \ge\ \lambda_I[b]\ \ge\
       \frac{P}{\lambda}\lambda_D[b]\ >\ \frac{P}{\lambda}.
\]
When \(q=2\), the potential equals the constant \(\lambda\).
Every horizontal slice has first eigenvalue \(P/(1-y)^2\).
This gives the lower bound \(P/\lambda\), and functions with a cosine on each interval and support in \(0<y<\varepsilon\)
approach it as \(\varepsilon\downarrow0\).
\end{proof}

\paragraph{Horizontal comparison at wide openings.}\label{par:horizontal-direct}
Lemma~\ref{S:lem:horizontal-sturm} gives $X\ge kW$ with $k=P/\Lambda_2(z)$.
For $0<z\le1/3$, comparison with $Q_{1/3}$ gives
$\Lambda_2(z)\le16P/9$, so take $k=9/16$, $b=1-k=7/16$.
If $q\le27/16$ the ground-state inequality $w\le1$ implies
$d=\min\{1,(1-pw)/(1-kw)\}\ge5/7$, where $p=q-1$.
Expanding the square in \eqref{S:cc:twosource} and using $d\le1$ gives
\[
 1-\mathcal C_q\le b\left[\frac{1+z}{d}
         +\frac{(1-z)(9+7z)}{16(1+z)}\right]
 \le\frac7{16}\left[\frac75(1+z)+\frac{5-7z+4/(1+z)}{16}\right].
\]
The bracket increases with $z$, since its derivative is at least $7/5-11/16$.
At $z=1/3$ this yields $\mathcal C_q\ge109/3840>1/1000$.

\begin{lemma}[A bound on the mass operator from two unweighted eigenvalues]\label{lem:obtuse-projection}
Let $0<z\le1$, $P=\pi^2/4$, $\tau=4\arctan\sqrt z/\pi$.
Choose $u,v>0$ so that $\Lambda_2(z)/P\le u$ and the second unweighted eigenvalue is at least $Pv$.
When $r=2-q\le\tau$, the norm $w$ of the mass operator $W$ on $e^\perp$ satisfies
\begin{equation*}\tag{W}\label{par:flat-transfer}
 w\le (u/v)\exp\{r[2/5+(\log v)/\tau]\}.
\end{equation*}
\end{lemma}
\begin{proof}
On the unit interval with zero trace at zero the Green kernel is
$L(x,y)=\min(x,y)$ for the measure $dx/x$.
Its operator is Hilbert--Schmidt and positive. The trace of its cube is
$6\int_{0<x<y<z<1}(x/z)\,dx\,dy\,dz=1/3$.
Thus the first eigenvalue of this weighted interval problem is at least $\kappa=3^{1/3}$.
With $h=1-y$, $l=(1-x)/h$, H\"older interpolation with the unweighted
mixed interval estimate and rescaling on each horizontal slice give
$X(b)\ge P(\kappa/P)^\tau\int h^{-2}l^{-\tau}b^2\,d\mu$.
For functions orthogonal in $L^2(D,\mu)$ to the unweighted first eigenfunction, interpolate this inequality with
$Q_z\ge Pv I$, using weights $r/\tau,1-r/\tau$. Since $\nu\tau=2$,
Lemma~\ref{S:lem:linear-potential-transfer} and min--max give
$w\le (u/v)e^{-\gamma r}(P/\kappa)^r v^{r/\tau}$.
Finally $P/\kappa<43/25$ by $P<121/49$, $\kappa>36/25$,
and $e^{2/5+4/27}>43/25$ by the cubic Taylor bound for the exponential.
\end{proof}

\subsection{Curvature at wide openings for all exponents}
\begin{proposition}\label{S:prop:spectral-wide}
For $0<z\le1/3$ and $1\le q\le2$ one has $\mathcal C_q(z)>1/1000$.
\end{proposition}
\begin{proof}
The direct comparison on p.~\pageref{par:horizontal-direct} covers
$q\le27/16$. Fix $q>27/16$ and \(0<z\le1/3\). On \(D=\{0<y<x<1\}\) write
\(Q_z=X+zY\), with zero trace on \(x=1\) and Neumann conditions
on the other two edges. Set \(P=\pi^2/4\), \(\theta=\arctan\sqrt z\),
\(\nu=\pi/(2\theta)\), \(t=z^{1/3}\).
Let \(\lambda,\lambda_2\) be the first two unweighted eigenvalues
of \(Q_z\). The proof of Proposition~\ref{S:prop:wide-linear} gives
\(\lambda/P\le1+(61/50)t\) and \(\lambda_2/P\ge(1+t)^3\).
Testing with the eigenfunction $E$ at $z=1/3$ gives the additional upper
bound \(\lambda/P\le16/9+(z-1/3)[8/3-81/(4\pi^2)]
\le(71+27z)/45\). The lower bound on the second eigenvalue in Lemma~\ref{S:lem:eqspec} and
\(Q_z\ge3zQ_{1/3}\) give \(\lambda_2/P\ge64z/3\).
Thus \(\lambda/P\le u\), \(\lambda_2/P\ge v\) with the
choices $u=\min\{1+at,(71+27z)/45\}$ and $v=\max\{(1+ct)^3,64z/3\}$,
where $a=61/50$, $c=21/20$. Indeed for $z\le2/9$ the longitudinal
test in Section~\ref{S:sec:compact-curvature} gives the coefficient
$[4(29/10)^3/(27P(1+z))]^{1/3}>c$, and $(1+ct)^3<5$ since $t<61/100$.
For $z\ge2/9$, $t>3/5$ gives $(1/t+c)^3<64/3$.

Write $r=2-q$, $p=q-1$, $\tau=2/\nu$, and let $w$ be the norm on $e^\perp$ of the mass operator $W$ with weight $V=\Lambda_q u_q^{-r}$ in the $Q_z$ inner product, as in Lemma~\ref{lem:obtuse-projection}.
Lemma~\ref{lem:obtuse-projection} gives (W) for $r\le\tau$.
Always $w\le1$ by the ground-state inequality applied to the minimizer.
By Lemma~\ref{S:lem:horizontal-sturm} one may take \(k=1/u\)
in \(X\ge kW\), with \(W\) the energy-space mass operator.
Set $\delta=1/1000$ and define
\[
 b=1-k,\quad h=\frac{(1-z)(9+7z)}{16(1+z)},\quad
 H=h+1+z,\quad B_\delta=\frac{b(1+z)}{1-\delta-bh},\quad
 R_\delta=b B_\delta,\quad J_\delta=1-B_\delta .
\]
Here $0<B_\delta<1$ since $b\le7/16$ and $H\le27/16$.
Set \(d=\min\{1,(1-pw)/(1-kw)\}\). Inequality~\eqref{S:cc:twosource} gives
\[
 \mathcal C_q\ge1-
 b\frac{[(1-z)d+4(1+z)]^2}{16(1+z)d}.
\]
On $0<d\le1$, expanding the square gives the upper bound
\(b[(1+z)/d+h]\) on $1-\mathcal C_q$. Thus \(d>B_\delta\) suffices.
For \(r>R_\delta=bB_\delta\), one has
\(d\ge\min\{1,r/b\}>B_\delta\). For \(r\le R_\delta\),
the second inequality in \eqref{S:eq:widechecks} below guarantees
\(R_\delta<\tau\). Also $R_\delta<b\le7/16$, so (W) applies.
The condition \((1-pw)/(1-kw)>B_\delta\) is equivalent to
\(w[1+(R_\delta-r)/J_\delta]<1\). By (W), the logarithm of the left-hand side is at most
\[
 -\log(v/u)+(R_\delta-r)/J_\delta+
       r[2/5+(\nu/2)\log v].
\]
This is affine in \(r\), and \eqref{S:eq:widechecks} makes it negative
at both endpoints $0,R_\delta$.

\paragraph{Checking the bound at $r=0$ and $r=R_\delta$.}
Write $a=61/50$, $t=z^{1/3}$. The conditions on the affine bound at
$r=0,R_\delta$ are
\begin{equation}\label{S:eq:widechecks}
 \log(v/u)>R_\delta/J_\delta,\qquad
 (1-R_\delta/\tau)\log v>\log u+2R_\delta/5 .
\end{equation}
For the first condition,
$R_\delta/J_\delta=b^2(1+z)/(1-\delta-bH)<1$
by $b\le7/16$ and $H\le27/16$. If $z\ge2/9$, the two bounds obtained from $z=1/3$ give
$v/u\ge960z/(71+27z)\ge640/231>11/4>e$.
If $z\le2/9$, then $t<61/100$ and $b\le at/(1+at)$. Hence
\[
 \frac{R_\delta}{tJ_\delta}\le
 \frac{(11/9)a^2 t}{(1+at)[1-\delta-(11/16+\delta)at]}<\frac43 .
\]
The right-hand side increases with $t$, and its value at $61/100$ is below $4/3$.
Meanwhile $\log(v/u)\ge(3-a)\log(1+t)\ge(3-a)t/(1+t/2)>4t/3$.

For the second condition we bound $R_\delta/\tau$ on the whole interval.
With $d=u-1\le at$, put $R=d^2(1+z)/[(1+d)(1+(1-h)d)]$.
This increases with $d$, and $R_\delta/R=(1-bh)/(1-bh-\delta)<501/500$.
Since $1-h\ge7/16+5z/8$, one has
$(1+t^3)(1+t^3/3)\le1+(1-h)at$: indeed
$(13/9)t^2-3t^3/4$ increases to less than $1/2$ at $7/10$, and $7a/16>1/2$.
Thus $\arctan\sqrt z\ge\sqrt z/(1+z/3)$ and AM--GM give
\[
 R_\delta/\tau\le\frac{501}{500}\frac{\pi a^2}4\frac{\sqrt t}{1+at}
 \le\frac{501}{500}\frac{\pi a^{3/2}}8<\frac8{15}.
\]
Also $\tau/t\le2/3^{2/3}<1$: the ratio increases up to $z=1/3$ since
$3x/[2(1+x^2)]\ge x\ge\arctan x$ for $x=\sqrt z\le1/\sqrt3$.
It now suffices that $L(t)=(7/5)\log(1+ct)-\log(1+at)>(16/75)t$.
This function is concave by $a>c$ and $a^2<(7/5)c^2$. At $t=7/10$
it exceeds $(2/5)(11/20)-7/100=3/20$: here
$\log(347/200)>11/20$ by the first two terms of its arctanh series,
and $\log(1+119/1735)<7/100$. Concavity and $L(0)=0$ now give
$L(t)>(3/14)t$, with excess over the required bound greater than $t/1050$.

\end{proof}

\section{Openings from \texorpdfstring{$\pi/2$ to $2\pi/3$}{90 to 120 degrees}}\label{S:sec:middle}

\begin{lemma}[The spectral gap at opening $2\pi/3$]\label{S:lem:eqspec}
For the fixed-triangle form at $z=1/3$, $\nu_1=L=4\pi^2/9$ and
\begin{equation}\label{S:comp:q2:equilateral-even-gap}
 \nu_2(1/3)\ge4L.
\end{equation}
\end{lemma}
\begin{proof}
Reflect the mixed eigenfunction evenly across the two Neumann sides.
The resulting Dirichlet eigenfunction on the equilateral triangle of side
$2\sqrt3$ is invariant under permutations of its vertices. Write these
vertices as $0,u,v$ and take lattice coordinates $tu+sv$. The product
\[
 4\sin(\pi t)\sin(\pi s)\sin(\pi(1-t-s))
   =\sin(2\pi t)+\sin(2\pi s)-\sin(2\pi(t+s))
\]
is positive inside, vanishes on the boundary and has all these symmetries.
Each Fourier mode on the right has eigenvalue $L=4\pi^2/9$.

Extend across the sides by odd reflection to the triangular tiling.
The extension is invariant under $120^\circ$ rotations about $0$
(two side reflections) and about the centroid (vertex symmetry).
Composing the centroid rotation with the inverse rotation about $0$
translates $0$ to the next vertex. The two orientations give the periods
$u$ and $v$. Hence every Fourier frequency has energy
$L(a^2-ab+b^2)$ with $a,b\in\mathbb Z$.
The first eigenvalue $L$ is simple by positivity. The integer
$a^2-ab+b^2\equiv(a+b)^2\pmod3$ cannot equal two. For $a^2-ab+b^2=3$
the only pairs are $\pm(2,1),\pm(1,2),\pm(1,-1)$, by completing the
square. Their frequency vectors are parallel to $u,v,v-u$, respectively.
Each corresponding line through $0$ is a reflection axis of the tiling with odd sign (the third is a rotated side). Its reflection fixes the
frequency and negates its coefficient, so all six coefficients vanish.
This gives $\nu_2\ge4L$.
\end{proof}

\begin{proposition}[Curvature for openings between $\pi/2$ and $2\pi/3$]\label{S:prop:central-sum}
For every $1\le q\le2$ and $1/3\le z\le1$,
\[
 G\le(12/5)\beta(1/2-\beta),\qquad
 \mathcal C_q(z)\ge\frac{37751}{54274560}>\frac1{1500}.
\]
\end{proposition}
\begin{proof}
Write $P=\pi^2/4$, $p=q-1$, $r=1-p$, and let $w$ be the norm on $e^\perp$ of the mass operator $W$ with weight $V=\Lambda_q u_q^{-r}$ in the $Q_z$ inner product, as in Lemma~\ref{lem:obtuse-projection}.
If $X(v)\ge k\int V v^2\,d\mu$ for every test, $k\ge1/2$, and $(6p-5k)w\le1$, the bound \eqref{eq:G-weighted-bound} gives
\[
 G\le\beta\left(1-\frac{\beta}{1-k}\right)
     \max\{1,(1-kw)/(1-pw)\}\le(12/5)\beta(1/2-\beta).
\]
For all $q$ the horizontal Sturm comparison in
Lemma~\ref{S:lem:horizontal-sturm} gives $k=P/\Lambda_2(z)\ge1/2$
using the ground state at $z=1$ as a trial function. If $r\ge5/12$, $w\le1$ now suffices and we can use the differential inequality on p.~\pageref{par:middle-backward}.

\paragraph{A lower bound on the second unweighted eigenvalue.}\label{par:ordinary-middle-gap}
Let $\lambda,\lambda_2$ be the first two unweighted eigenvalues at $z$.
Their sum is concave in $z$, by the variational formula
over orthonormal pairs. Its values at $1/3$ and $1$ are at least
$80P/9$ and $12P$, by Lemma~\ref{S:lem:eqspec} and the eigenvalues $2P$ and $10P$ at $z=1$.
Using the ground state at $z=1$ as a trial function gives
$\lambda/P\le2-(1-8/\pi^2)(1-z)$. Subtract this upper bound from the linear lower bound on $(\lambda+\lambda_2)/P$. Since $1-8/\pi^2>1/6$,
\begin{equation}\label{S:eq:local-c1}
           \lambda_2/P\ge v(z):=(11+9z)/2\qquad(1/3\le z\le1).
\end{equation}
For $z\ge1$, monotonicity gives $\lambda_2(z)\ge10P$.
\paragraph{The mass estimate for $r<5/12$.}\label{par:middle-uniform}
Put $x=z-1/3$, $y=1+z$, and take
$m=\min(16/9+5x/8,(29+3z)/16)\ge\lambda/P$.
The ground state at $z=1/3$ gives the first upper bound since
$8/3-81/(4\pi^2)<5/8$. The ground state at $z=1$ gives the second since $1-8/\pi^2>3/16$. For $r<5/12$ one has
$r\le\tau=4\arctan\sqrt z/\pi$. By (W),
$w\le(m/v)e^{rH}$ with $H=2/5+(\log v)/\tau$.
It suffices to prove $(B-Sr)e^{rH}\le v$ with $S=6m$, $B=S-5$.
The inequality $B\le13v/16$ holds at $z=1/3,5/9,1$. On each interval between these points both sides are affine, so the inequality holds throughout.
The first upper bound also gives $B\le17y/4$.
Concavity of $\tau$ gives $\tau\ge y/2$ by interpolation at $z=1/3$ and $z=1$. The tangent to $\log v$ at
$v=7$ and $\log7=3\log2-\log(8/7)<3(25/36)-2/15=39/20$
by the arctanh series give
\[
 yH\le(2/5)y+2\log v
 \le133/30+(2/5+9/7)x\le(16/9)(y+20/17)\le16yS/(9B).
\]
The third inequality holds at $x=0$ and the right-hand side has the larger slope.
For $t=Sr/B\in[0,1]$ we obtain
\[
 (B-Sr)e^{rH}/v\le(13/16)(1-t)e^{16t/9}
 \le (117/256)e^{7/9}<1 .
\]
The last bound follows from
$\log(256/117)\ge2[139/373+(139/373)^3/3]>7/9$.
For $t>1$ the numerator is negative.
\paragraph{Integrating the differential inequality for $\beta$ downward from $z=1$.}\label{par:middle-backward}
On the whole interval $1/3\le z\le1$, \eqref{S:eq:curv} gives
$d\beta/d\log z\ge(-7\beta+19\beta^2)/5$. Integrating for $1/\beta$ backwards
from $z=1$, where $\beta(1)<2/15$ by Lemma~\ref{S:lem:square-slope}, gives $\beta\le14/(38+67z^{7/5})$. To simplify this upper bound on $\beta$, we bound its denominator from below by a linear function of $z$.
Let $H=105/4+74z$ and $b=14/H$. In fact
\[
 \max_{z\ge0}(74z-67z^{7/5})
 =\tfrac{148}{7}(370/469)^{5/2}
 <\tfrac{148}{7}(8/9)^5<47/4,
\]
so $\beta\le b$. As a function of $\beta$, the lower bound
$1+(z-1)\beta-(48/5)(1+z)\beta(1/2-\beta)$
decreases up to $(29+19z)/[96(1+z)]$. This point exceeds $b$ because
$(29+19z)H-1344(1+z)=127/18+(26857/12)y+1406y^2$
at $y=z-1/3\ge0$. Substituting $z=(14/b-105/4)/74$ now gives
\[
 \mathcal C_q\ge (9168b^2-3901b+416)/1480
 \ge\frac{37751}{54274560}>\frac1{1500}.
\]
\end{proof}
\paragraph{The case $q=2$ on $1/3\le z\le1$.}\label{par:spectral-obtuse}
For $q=2$ the bound $\mathcal C_2\ge1/37^2$ also follows directly from the
unweighted eigenvalues at the openings $\pi/2$ and $2\pi/3$. On $1/3\le z\le1$ the ground
state at $z=1$ gives $\Lambda_2(z)\le P+2+(P-2)z\le2P$ with $P=\pi^2/4$.
By monotonicity and Lemma~\ref{S:lem:eqspec}, $\nu_2(z)\ge64P/9$.
Take $k=1/2$, $w=9/32$ in \eqref{eq:G-weighted-bound}.
Then $D=23/55\ge2/5$, so $G\le(5/2)\beta(1/2-\beta)$ and
$d\beta/d\log z\ge-3\beta/2+4\beta^2$. At $z=1$ we have
$\beta(1)=1/2-4/\pi^2<12/125$. Backward integration for $1/\beta$ gives
$\beta\le12/(32+93z^{3/2})$. The tangent bound
$z^{3/2}\ge z-4/27$ at $4/9$ yields $\beta\le v=4/(6+31z)<1/4$.
The lower bound for $\mathcal C_2$ given by $1-(6+4z)\beta+10(1+z)\beta^2$ decreases up to $1/4$.
At $v$ it equals $[1+2h+465h^2]/(6+31z)^2\ge1/37^2$, with $h=z-1/3$.

\section{Acute openings and the proof of Theorem~A}\label{sec:acute-real}
Proposition~\ref{S:prop:lowq} covers $q\le q_\ast$ at every opening, Proposition~\ref{S:prop:spectral-wide} covers $z\le1/3$ and Proposition~\ref{S:prop:central-sum} covers $1/3\le z\le1$ for every $1\le q\le2$. It remains to treat $q_\ast\le q\le2$ and $z\ge1$, that is, openings $0<\alpha\le\pi/2$. In the cylinder coordinates introduced below, the area factor has positive Taylor coefficients,
and the same coordinates cover the whole range from the limit of small openings up to $\alpha=\pi/2$.
Section~\ref{sec:cylinder-curvature-bound} bounds the change in the minimizer throughout this range. Proposition~\ref{prop:real-acute-curvature}
combines this estimate with the disk bound at low exponents.

\subsection{The variational problem on a fixed cylinder}\label{S:sec:acute-complete}
Put $h=1/4$, $s=\alpha/(2\pi)$, $t=s/h$, $p=q-1$, $r=1-p$, and $c=2/q$.
After rotating the sector if necessary, the map with derivative $(1-\zeta^{1/s})^{-2s}$ sends the sector
of angle $2\pi s$ to the isosceles mixed triangle: the two rays have
straight images, and the tangent to the image of the unit arc has constant direction. The length of each ray image is $s\mathrm B(s,1-2s)$. Since the image of the boundary is a simple curve, the argument principle gives injectivity. In coordinates $\tau=-\log|\zeta|$, with normalized angular
measure $d\theta/(2\pi)$, write an overbar for the angular average and define
\[
 E_s(v)=\int_0^\infty\overline{v_\tau^2+s^{-2}v_\theta^2}\,d\tau,
 \quad I_s(v)=\int_0^\infty e^{-2\tau}\overline{W_s|v|^q}\,d\tau,
 \quad W_s=|1-e^{-\tau/s+i\theta}|^{-4s}.
\]
Tests vanish at $\tau=0$. Periodic tests give the same minimum as the
reflection-symmetric ones: replacing a nonnegative test by the quadratic
mean of it and its reflected copy decreases the energy and increases the
$q$-mass. Set
\[
 L=\inf E_s/I_s^{\,c},\quad I_s(u)=1,\quad U=\overline u,\quad v=u-U,
 \quad Z=2\int e^{-2\tau}\overline{W_s},\quad A=\log Z,\quad \mathcal A=sZ^c.
\]
Conformal invariance gives $F_q(2\pi s)=2\pi\,2^{-c}\mathcal A L$. The area identity
and the gamma reflection formula yield
\begin{equation}
 Z(s)=\frac{\Gamma(1+s)^2\Gamma(1-2s)}
            {\Gamma(1-s)^2\Gamma(1+2s)},\qquad
 A(s)=\sum_{\substack{m\ge3\\m\ {\rm odd}}}
       \frac{(2^{m+1}-4)\zeta(m)}m s^m.
 \label{S:eq:local-b-1}
\end{equation}
All coefficients of $A$ are positive.
\subsection*{The radial minimizer on the disk at small openings}
We first compute how the conformal weight changes the quotient evaluated at the disk minimizer. We then compute how much further the quotient decreases when we minimize over all tests. Both expansions can be differentiated twice with respect to the opening.
\begin{proposition}[Small-opening expansion about the disk minimizer]\label{S:prop:disk-limit}
Let $L_0(q)$ denote the radial value and $\sigma_q=(4L_0(q)/q)^{1/2}$ with the normalization above. Uniformly for $1\le q\le2$,
\[
 \frac{L_q(s)}{L_0(q)}=1-c D_q s^{q+3}+o(s^{q+3}),\qquad
 D_q=2^{1-q}\sigma_q^q\Gamma(q+1)\zeta(q+3).
\]
The remainder may be differentiated $j\le2$ times, giving $o(s^{q+3-j})$. In particular $(sZ^cL_q)''/L_0=48c\zeta(3)s^2+O(s^{q+2})$. Also $L_q(s)/L_0(q)=1+O(|s|^{4-2\eta})$ locally holomorphically in $|q-Q|<\eta$ for each real $Q\in[1,2]$, with an absolute $0<\eta<1/4$ and $s$ in a fixed sector about the positive axis.
For the radial minimizer $g=g_q$ of unit mass on the disk, put $M_s=I_s(g)$.
The relative difference between the quotient at $g$ and its minimum satisfies, uniformly for $1\le q\le2$,
\begin{align*}
 \Delta_q(s)&:=\frac{L_0 M_s^{-2/q}}{L_q(s)}-1
 =L_0\sigma_q^{2q-2}N_q\,s^{2q+3}+o(s^{2q+3}),\\
 N_q&=2^{2-2q}\zeta(2q+3)\left(2\Gamma(2q)/q-\Gamma(q)^2\right)>0 .
\end{align*}
After $j$ derivatives in $s$, the remainder is $o(s^{2q+3-j})$ for $j=1,2$.
\end{proposition}
\begin{proof}
Fix $Q\in[1,2]$. For the radial minimizer $g>0$ on the disk at exponent $Q$, normalized so that $d\mu=e^{-2y}g^Qdy\,d\theta/(2\pi)$ has mass one, write $p=Q-1$. Its equation is $g(t)=L_0\int\min(t,y)e^{-2y}g(y)^pdy$ and multiplying $-g''=L_0e^{-2t}g^p$ by $g'$ gives $g'(0)=\sigma_Q$. The trace bound $|v(t)|\le\sqrt t\,\|v'\|_2$ and the test $1-e^{-2t}$ bound $L_0$ above and below uniformly. Since $g$ increases and $g/y$ decreases, the expectation of the increasing function $y\mapsto\min(t,y)$ under the density proportional to $e^{-2y}g^p$ lies between its expectations under the densities proportional to $e^{-2y}$ and $y^pe^{-2y}$. This gives
\[
 \sigma_Q(1-e^{-2t})/2\le g(t)\le\sigma_Q\min(t,Q/2).
\]
For $|\arg\xi|<\pi/6$ put $R=|\xi|$ and $P_w=e^{-w|D_\theta|}$, and let $\Pi$ denote the average over the circle. For real $\xi=s>0$, $K_s$ is the Dirichlet Green kernel for $-\partial_y^2-s^{-2}\partial_\theta^2$ on the half-cylinder, and $W_s$ is the conformal weight. Both are given by the following formulas, which we also use for complex $\xi$:
\[
 K_\xi(t,y)=\tfrac12\int_{|t-y|}^{t+y}P_{v/\xi}\,dv,\qquad
 W_\xi=\exp(\xi H),\qquad H=2\sum_{n\ge1}\frac{e^{-ny/\xi}}n(e^{in\theta}+e^{-in\theta}).
\]
For $|\arg w|<\pi/6$ we have $|w^2+x^2|\ge (|w|^2+x^2)/2$, so the kernel $w/[\pi(w^2+x^2)]$ has $L^1(\mathbb R)$ norm at most two. Periodization gives the uniform operator bound $\|P_w\|_{L^\infty\to L^\infty}\le2$. The circle kernel itself, relative to $d\theta/(2\pi)$, takes the value $\coth(w/2)$ at $\theta=0$ for real $w>0$, which diverges as $w\downarrow0$. The operator bound for $\Re w\le1$ and the Fourier series for $\Re w\ge1$ give $\|P_w-\Pi\|_{L^\infty\to L^\infty}\le C e^{-\Re w}$. Thus $K_\xi$ is bounded by $C\min(t,y)$ and differs from $K_0=\min(t,y)\Pi$ by at most this bound times $Ce^{-\Re(1/\xi)|t-y|}$. With $H_*=-4\log(1-e^{-y\cos(\pi/6)/R})$ we have $|H|\le H_*$. For each fixed $j\ge1$ and $-1/4\le m\le3$,
\[
 \int_0^\infty y^m H_*^j e^{R H_*}\,dy\le C_jR^{m+1}.
\]
After the substitution $y=Rx$, the integrand near zero is bounded by $O(x^{-3/4}(1+|\log x|^j))$ for $4R\le1/2$, and its tail decays exponentially.

For $q$ near $Q$ define $T_\xi h=(L_0(Q)/g)\int K_\xi e^{-2y}g^{q-1}W_\xi h\,dy$, using $W_0=1$ at $\xi=0$. For real $\xi=s>0$ and $1\le q\le2$, the Euler equation for a positive $u=gh$ takes the integral form $u=\text{const}\cdot\int K_s e^{-2y}W_s u^{q-1}\,dy$. Equivalently, $h$ is proportional to $T_s(h^{q-1})$. For complex parameters we use the same integral equation with $K_\xi$ and $W_\xi$. We therefore look for a fixed point of the map
\[
 h\longmapsto T_\xi(h^{q-1})/\mu[T_\xi(h^{q-1})]
\]
among the functions $h\in L^\infty$ with $\mu h=1$ that are close to $h=1$. At $\xi=0$ and $q=Q$ the equation of $g$ gives $T_01=1$, so $h=1$ is a fixed point. With respect to $\mu$ the operator $T_0$ has the kernel $L_0\min(t,y)/(g(t)g(y))$. This kernel is symmetric, so $T_0$ preserves $\mu$, and it is bounded below by $Q/4$. Hence $T_0$ contracts functions of $\mu$-mean zero in $L^\infty$ by the factor $1-Q/4$. The derivative of the map at $h=1$, acting on functions of $\mu$-mean zero, is $pT_0$, so its norm is at most $p(1-Q/4)\le3/4$.

The map remains a contraction in $L^\infty$ on a ball about $1$ whose radius does not depend on $Q$, whenever $|q-Q|<\eta$ and $R$ is sufficiently small. To see this, compare it with the map at $\xi=0$, $q=Q$ in three steps. First, at a fixed exponent $q$, replacing $K_0$ by $K_\xi$ while still omitting the weight $W_\xi$ changes the operator by $O(R^{1-\eta})$; this follows from the bound for $K_\xi-K_0$ above, after integrating $y^{-\eta}$ near zero. Secondly, inserting the weight $W_\xi$ changes it by $O(R^{\Re q+1})$, by the moment bound with $j=1$. Thirdly, at $\xi=0$, changing the exponent from $Q$ to $q$ changes it by $O(|q-Q|)$, by integration against $(1+|\log\min(y,1)|)\min(y,1)^{-\eta}$. The same bounds hold for the derivatives of the map with respect to $h$ on the ball. Iterating the map therefore gives fixed points $h_{\xi,q}$ that are holomorphic in $(\xi,q)$ for $\xi\ne0$ in the sector. At $\xi=0$ the fixed point is radial and depends holomorphically on $q$. On radial functions $K_\xi$ and $K_0$ act in the same way, so the maps at $\xi$ and at $0$ differ on $h_{0,q}$ only through the weight $W_\xi$. The contraction estimate then gives $\|h_{\xi,q}-h_{0,q}\|_\infty=O(R^{\Re q+1})$.

From now on write $g$ for $g h_{0,q}$, let $\mu_q$ be the complex measure proportional to $e^{-2y}g^qdy\,d\theta/(2\pi)$ of total mass one, and write the solution constructed above, after multiplication by a scalar, as $g(1+b)$ with $\mu_q b=0$. Then $\|b\|_\infty=O(R^{\Re q+1})$. For $i=q-1,q$ put $J_i=\mu_q[W_\xi(1+b)^i]$. Testing the equation of this solution against the radial function $g$ gives
\[
                         L_q(\xi)/L_0(q)=J_q^{1-c}/J_{q-1}.
\]
For real $s>0$ and $1\le q\le2$ the positive solution has finite energy, by the trace estimate and $W_s\le(1-e^{-y/s})^{-4s}$. Testing its equation with $v^2/(u+\epsilon)$ and using H\"older shows it realizes the minimum. The same argument at $\xi=0$ identifies the radial minimizer on the disk for $1\le q\le2$. The displayed ratio gives the holomorphic continuation for complex parameters in the sector.

Since $\overline H=0$, the part of $J_i-1$ independent of $b$ is $O(R^{\Re q+3})$ by the moment bound with $j=2$. Its mixed term is $O(R^{\Re q+2}\|b\|_\infty)$, and the remaining terms contribute $O(\|b\|_\infty^2)$. This gives the asserted complex-exponent bound.

\emph{Rescaling near the boundary.}
For the rest of the proof let $1\le q\le2$ be real, $p=q-1$, $\omega=e^{-2t}$, and choose the normalization $u=g+v$, $\int\omega g^p v=0$. In particular $E_\xi(g,v)=0$ and $\|v/g\|_\infty=O(R^{q+1})$. For $\phi$ orthogonal to $g$ in energy, define
\[
 b_\xi(\phi)=\int\omega(W_\xi-1)g^p\phi .
\]
For $\xi=s>0$ the change $\phi(t,\theta)=\sqrt{s}\,\psi(t/s,\theta)$ preserves the energy, $E_s(\phi)=E_1(\psi)$. With this change of variables, $s^{-q-3/2}b_s$ converges in dual norm to $\sigma_q^p B_q$, where
\[
 H(x,\theta)=-4\log|1-e^{-x+i\theta}|,\qquad B_q(\psi)=\int x^p H(x,\theta)\psi .
\]
Indeed $W_s(sx)=e^{sH}$, $g(sx)/s\to\sigma_q x$ and $\|\psi(x,\cdot)\|_2\le\sqrt x\|\psi\|_{E_1}$. For $s\le1/10$, dominated convergence follows from $C\sqrt x(1+x) h(x)e^{h(x)/10}$ with $h(x)=-4\log(1-e^{-x})$. The convergence is uniform in $q$ by the equation of $g$ and the uniform bounds for $g$ above. Since $\overline H=0$, the functional $B_q$ vanishes on every radial test function. Restricting to tests orthogonal to $g$ therefore does not change this leading term.

\emph{The leading term in $\Delta_q$.}
Fix the real exponent $q$, take $Q=q$ in the preceding contraction, and
write $u=g(1+b)$, $\mu b=0$, $\|b\|_\infty=O(R^{q+1})$.
Set $M=\mu W_\xi$, $a=\mu(W_\xi b)/M$ and $d=\mu(W_\xi b^2)/M$.
Since $\mu b=0$, and since $\mu|W_\xi-1|=O(R^{q+2})$ by the moment bound with $j=1$, we have
$a=\mu[(W_\xi-1)b]/M=O(R^{2q+3})$. The identity $L_q(\xi)/L_0(q)=J_q^{1-c}/J_{q-1}$ above yields
$1+\Delta_q=(J_p/M)(J_q/M)^{c-1}$. In
$J_i/M=1+ia+i(i-1)d/2+O(R^{3q+3})$, the quadratic moment cancels:
\[
 p+(c-1)q=1,\qquad p(p-1)+(c-1)q(q-1)=0.
\]
Since $d=O(R^{2q+2})$, the product gives
$\Delta_q=a+O(R^{3q+3})$, in particular the required bound throughout
the complex sector.

For real $s$ write $G_s$ for the integral operator $T_s$ with the
factor $W_s$ omitted. Expanding the normalized fixed point once, using
$G_s1=1$ and $\|T_s-G_0\|_{\infty\to\infty}=O(s^{1-\eta})$ above, gives
\[
 b=pG_0 b+(I-\mu)G_s(W_s-1)+O(s^{q+2-\eta}).
\]
Here $\mu G_0b=0$. Both $G_0b$ and the subtracted scalar are radial
and $O(s^{q+1})$. Pairing these two terms with $W_s-1$ therefore contributes $O(s^{2q+4})$,
since $\overline H=0$ gives $\mu|\overline{W_s-1}|=O(s^{q+3})$.
The error in this pairing is $O(s^{2q+4-\eta})$ since $\mu|W_s-1|=O(s^{q+2})$.
Consequently
\[
 \Delta_q=\mu[(W_s-1)G_s(W_s-1)]+o(s^{2q+3})
         =L_0\|b_s\|^2+o(s^{2q+3}).
\]
Indeed $\mu[(W_s-1)G_s(W_s-1)]$ is $L_0$ times the squared dual norm of the functional $\phi\mapsto\int\omega(W_s-1)g^p\phi$ on all tests,
and the component of this functional along $g$ contributes only
$(M-1)^2=O(s^{2q+6})$. All moment bounds follow from the rescaling
above, uniformly for $1\le q\le2$.

The Fourier coefficients of $H$ at $\pm n$ are $2e^{-nx}/n$. Inserting the Dirichlet Green kernel $(e^{-n|x-y|}-e^{-n(x+y)})/(2n)$ and rescaling each mode gives
\[
 \|B_q\|^2=4\zeta(2q+3)\iint (xy)^p[e^{-2\max(x,y)}-e^{-2(x+y)}]\,dx\,dy=N_q.
\]
To evaluate the first integral, split the domain at $y=x$. The second integral factors into two one-dimensional integrals. At $q=1$ define the linear functional $\ell_s(v)=\int\omega W_s v$ on the real energy space. Since $W_s>0$ and $E_s(|v|)=E_s(v)$, taking the supremum over $E_s(v)=1$ gives $1/L_1(s)=\|\ell_s\|_{E_s^*}^2$. The component of $\ell_s$ along $g$ has squared norm $M_s^2/L_0$, and its restriction to $g^\perp$ is $b_s$. Thus $\Delta_1=L_0\|b_s\|^2/M_s^2$. The integrand is positive for $x,y>0$, so $N_q>0$. To justify differentiation, along any $s_i\downarrow0$, $q_i\to Q$ the holomorphic functions $s_i^{-2q_i-3}\Delta_{q_i}(s_i w)$ are bounded on a fixed disk about $w=1$. Their real limit is $L_0(Q)\sigma_Q^{2Q-2}N_Q w^{2Q+3}$. Normal families and the identity theorem give convergence of derivatives on smaller disks.

Finally $\overline{H^2}/2=4\sum_{n\ge1}e^{-2nx}/n^2$ in the boundary coordinate. In $M_\xi$ the linear term has zero mean. Using $\omega g^q=\sigma_q^q t^q+O(t^{q+1})$ and integrating the quadratic term termwise, with $\int_0^\infty x^qe^{-2nx}\,dx=\Gamma(q+1)/(2n)^{q+1}$, gives $M_\xi=1+D_q\xi^{q+3}+O(R^{q+4})$ uniformly in the sector; the cubic remainder is bounded by the same moment estimate. The quotient $L_q=L_0M_\xi^{-c}/(1+\Delta_q)$, Cauchy's formula on a disk of radius proportional to $s$, and the area series give the first expansion and its derivatives.
\end{proof}

\paragraph{How the normalized curvature varies with $q$ at small openings.}\label{par:small-opening-order}
On some interval $0<s<s_0$ common to $1\le q\le2$, the normalized curvature
$F_q''(2\pi s)/F_q(2\pi s)$ strictly decreases with $q$.
To see the sign, set $B=(2/q)A+\log(L_q(s)/L_0(q))$.
Up to a factor independent of $s$, $F_q(2\pi s)=s e^B$, hence
\[
 4\pi^2\frac{F_q''(2\pi s)}{F_q(2\pi s)}
       =2B'/s+B''+(B')^2.
\]
Primes on $B$ denote derivatives in $s$. The holomorphic bound in
Proposition~\ref{S:prop:disk-limit} and Cauchy's formula on a fixed smaller
$q$ disk and an $s$ disk of radius proportional to $s$ bound
$\partial_q\partial_s^j\log(L_q(s)/L_0(q))$ by $O(s^{4-2\eta-j})$
for $j\le2$ (also without $\partial_q$). A finite cover makes the bounds
uniform in $q$. Since $A=4\zeta(3)s^3+O(s^5)$, it follows that
\[
 \partial_q\!\left(\frac{F_q''(2\pi s)}{F_q(2\pi s)}\right)
 =-\frac{24\zeta(3)}{\pi^2q^2}s+O(s^{2-2\eta}).
\]
As $2-2\eta>1$, the leading term fixes the sign throughout that interval.

For example, take a regular $n$-gon $R_n$ and a disk $B$ of the same area. Reflection of the mixed cell of $R_n$ and radial symmetry of the disk minimizer give
\[
 \frac{\lambda_q(R_n)}{\lambda_q(B)}
 = Z(1/n)^{2/q}\left[1-\frac{2D_q}{q}\,n^{-(q+3)}+o(n^{-(q+3)})\right],
\]
where $D_q$ is the coefficient of Proposition~\ref{S:prop:disk-limit}. For $1<q<2$ the first correction to the conformal area factor has order strictly between four and five.

\subsection{Convexity for acute openings}\label{app:real-angular-body}
Set $h=1/4$, $t=s/h$ and $Q=29/20$. We work on
$0<s\le h$, $Q\le q\le2$. The positive area term in the second
variation is explicit. We estimate the change in the minimizer by expanding the source term in a Fourier series in the angular variable.

\subsubsection{Second variation on the cylinder}\label{app:real-variations}
Let $h=1/4$, $s=\alpha/(2\pi)$, $p=q-1$, $r=2-q$, $c=2/q$ and $f(s)=F_q(2\pi s)$. The estimates below hold on $Q\le q\le2$, where $Q=q_\ast-1/4000=29/20$. Use the angular measure $d\theta/(2\pi)$ on the half-cylinder and set
\[
 E_s(v)=\int(|v_y|^2+s^{-2}|v_\theta|^2),\qquad W_s=|1-e^{-y/s+i\theta}|^{-4s},\qquad \omega=e^{-2y}.
\]
Tests vanish at $y=0$. The value $L$ is the minimum energy under $\int\omega W_s |u|^q=1$. The radial minimizer $g>0$ on the disk has $\int\omega g^q=1$ and energy $L_0$. Its equation gives $\sigma^2=4L_0/q$ for $\sigma=g'(0)$ and $\sigma(1-e^{-2y})/2\le g\le\sigma y$. The elementary radial trial $1-e^{-2y}$ yields $\sigma^q\le 2(q+1)(4/q)^{q/2}$. Up to a positive constant, $f=s\exp(cA)L$, where
\[
 A(s)=\sum_{m\ge3,\ m\ {\rm odd}}\frac{(2^{m+1}-4)\zeta(m)}m\,s^m .
\]
Let $D_0=W_s-1$, $D_j=s^j\partial_s^jW_s$ for $j=1,2$, with $y$ fixed. We bound
\[
 b_j=\left|\int\omega D_j g^q\right|,\quad d_j=\sqrt{L_0}\|\omega D_jg^p\|_{E_s^*},\quad
 k_j=L_0\sup_{E_s(v)=1}\int\omega D_j(y,0)g^{-r}v^2 \quad(j=1,2).
\]
The formula for $d_j$ also applies at $j=0$. Write $u=a(g+v)$ with $v\perp_{E_s}g$, $V^2=E_s(v)/L_0$. Section~\ref{app:radial-gap} bounds the quadratic form $B$ defined below, and the terms containing $D_1,D_2$, by forms with radial weights. At $s_0$ let $u$ be the normalized minimizer, $B(\phi)=L\int u^{q-2}\phi^2\,d\mu_{s_0}$,
and $d\mu_s=\omega W_s\,dy\,d\theta/(2\pi)$. Put $K=sZ^c$. Up to the constant factor $2\pi2^{-c}$, which is positive, independent of $s$ and cancels in $f''/f$, the profile is
\[
 f(s)=\inf_{v\ne0} K E_s(v)\Big(\int |v|^q\,d\mu_s\Big)^{-c}.
\]
At $s_0$ the Euler equation gives $E_{s_0}(u,\phi)=L\int u^{q-1}\phi\,d\mu_{s_0}$.
For variations orthogonal to $u$ in energy, the first variation of the mass vanishes. Differentiating the quotient twice in the test direction gives the Hessian $2K\mathcal J_p$, where $\mathcal J_p=E_{s_0}-pB$.

In the following identity $s=s_0$, $P=1+csA'$, $m_j=\int\omega
(\partial_s^j W_s)u^q$, $\beta=s^{-2}\int u_\theta^2/L$ and
$\mathfrak f(\phi)=E_s'(u,\phi)-L\int\omega(\partial_sW_s)u^p\phi$.
Differentiating first with $u$ fixed and then minimizing the quadratic form in the variation of $u$ gives
\[
 \frac{f''}{f}=\frac{K''}{K}
 -c(m_2+2Pm_1/s)+c(c+1)m_1^2
 +(6-4P+4csm_1)\beta/s^2-2\|\mathfrak f\|_{\mathcal J_p^*}^2/L.
\]
We justify differentiating the minimizer with the varying weight. Fix $s_0>0$ in the stated range and let $\tau$ vary in a small neighborhood with upper endpoint $s_+$ such that $4s_++r<2$. Write $a_0(y)=\min(1,y)$. For $j=0,1,2$, the formula for $W_\tau$ gives
\[
 \omega a_0^{-r}|\partial_\tau^j W_\tau|
 \le \Gamma(y):=C e^{-2y}a_0^{-4s_+-r}(1+|\log a_0|^2).
\]
Indeed $|1-e^{-y/\tau+i\theta}|\ge 1-e^{-y/\tau}\ge c_0a_0$, and the first two $\tau$-derivatives of its logarithmic modulus are bounded uniformly in $y,\theta$. The energy space embeds compactly in $L^2(\Gamma\,dy\,d\theta)$. On compact interior strips this is Rellich compactness. At the two ends use $\int|v(y,\theta)|^2\,d\theta/(2\pi)\le y E_{s_0}(v)$ and the integrability of $y\Gamma(y)$. The same trace bound with exponent $q/2$ controls the tails of the $q$-mass, since $q/2-4s_+>-1$. Compactness, uniqueness of the normalized positive minimizer and convergence of the energies now give $u_\tau\to u_{s_0}$ in energy and $L_\tau\to L_{s_0}$.

For $q<2$, set $U_\tau=L_\tau^{-1/r}u_\tau$, so that $E_\tau(U_\tau,\phi)=\int\omega W_\tau U_\tau^p\phi$. There is a uniform lower bound $U_\tau\ge\epsilon g\ge c_1a_0$. Indeed $\epsilon g$ is a subsolution if $\epsilon^r L_0\le2^{-4s_+}$, since $W_\tau\ge2^{-4s_+}$. Testing the difference on $(\epsilon g-U_\tau)_+$, concavity of the power and the ground-state inequality for $U_\tau$ bound the energy of this test by $p$ times itself. Let $S_\tau$ be the divided difference of $x^p$ between $U_\tau$ and $U_{s_0}$, using the derivative at coincident values. Its absolute value is bounded by $C a_0^{-r}$. Along every almost everywhere convergent subsequence the bilinear forms with density $\omega W_\tau S_\tau$ converge in operator norm to the form with density $p\omega W_{s_0}U_{s_0}^{-r}$, by dominated convergence and compactness of the energy unit ball in $L^2(\Gamma\,dy\,d\theta)$. Thus the convergence holds along the full family. The limiting linearized form dominates $(1-p)E_{s_0}$ by the ground-state inequality. Subtracting the two equations for $U$ and dividing by $\tau-s_0$ therefore gives a convergent difference quotient in energy. The right-hand side only differentiates $E_\tau$ and $W_\tau$ at the fixed $U_{s_0}$, and the same majorant justifies this in the dual norm. Normalizing gives the derivative of $u_\tau$. For $q=2$ the weighted mass form is twice differentiable in operator norm by the same compactness argument. Its first eigenvalue is simple, so the normalized eigenfunction is differentiable, by inversion on its energy-orthogonal complement. We may now differentiate the stationary quotient in a normalization with variation orthogonal to $u$. Solving the resulting variational equation produces the last term of the formula.

Set $M_j=b_j+qd_jV+qp k_jV^2/2$ and $S_p=(2V+d_1+pk_1V)^2/J_p$,
where $\mathcal J_p\ge J_p E_s$ on $u^\perp$. Then
$|s^jm_j|\le M_j$, $\beta\le V^2$ and
$s^2\|\mathfrak f\|^2_{\mathcal J_p^*}/L\le S_p$.
Taylor expansion of the mass along the segment from $g$ to $g+v$ gives a second derivative containing
$a^q |D_j|(g+\theta v)^{-r}$. The two endpoint terms are bounded by
$D_j(y,0)g^{-r}$ using $a\le1$, $a^{q+r}=a^2\le L/L_0$ and the
comparison with the weight at $\theta=0$ on p.~\pageref{par:angular-damping}.
Convexity of the inverse power bounds the whole segment. For the projected source, the three contributions to $s\|\mathfrak f\|_{E_s^*}/\sqrt L$ are bounded by $2V$, $d_1$ and $pk_1V$, respectively.
For the last one the corresponding factor is $\sqrt{L/L_0}a^p$,
and $a^{p+r}\sqrt{L/L_0}\le L/L_0$. Here $L\le L_0$ by the radial trial.
When $V<1/4$ the discarded quadratic terms have nonnegative sum.
After multiplication by $s^2$ they equal
\[
 c(c+1)\big(sm_1+2\beta/(c+1)\big)^2
     +\beta\big(6-4P-4c\beta/(c+1)\big).
\]
Indeed $\beta\le V^2$ and $P<10/7$ by the area calculation
in Section~\ref{sec:cylinder-green-sources}. Thus
\begin{equation}\label{real:second}
 \frac{f''}{cs f}\ge\frac{K''}{csK}
       -\frac{M_2+2PM_1+qS_p}{s^3}.
\end{equation}

\subsubsection{A mass bound on the orthogonal complement of the minimizer}\label{app:radial-gap}
Put $r=2-q$, $p=1-r$, $k=20r/11$, $a=1-k/3$,
$N=4(1+3a)=16-4k$ and $b=a(1-a)/(1+3a)$. The radial mass form with density $L_0 e^{-2y}g^{-r}W_s(y,0)$ has second eigenvalue at most $w$ relative to the energy. The same $w$ bounds $B(\phi)/E_s(\phi)$ for $\phi\perp_{E_s}u$. We take
\[
 w=w(r)=\frac{(6+r-4r^2+r^3)(6/5)^k}{16-4k},\qquad
 w\le73/120,\qquad pw\le3/8 .
\]
Indeed, set $x=e^{-2y}$, $t=1-x$. The disk equation gives
$g\ge\sigma t/2$ and the radial trial $t$ gives
$2^r L_0\sigma^{-r}\le 2(3-r)(2-r)^{r/2}\le6+r-4r^2+r^3$.
Since $s\le1/4$, the weight at $\theta=0$ satisfies $W_s(y,0)\le(1-x^2)^{-1}$.
We claim
\[
 t^{1-r}\le(6/5)^k(2-t)\bigl(t+b(1-t)\bigr).
\]
At $r=0$ this is $t\le(2-t)t$. At $r=11/20$ we have $b=2/27$ and
$t^{1-r}\le t^{4/9}$. For $t\le1/8$, weighted AM--GM gives
\[
 \frac{t^{4/9}}{(2+25t)/27}
 \le (12/25)^{4/9}(15/2)^{5/9}<9/4\le(6/5)(2-t).
\]
The middle inequality follows by taking ninth powers, since the ratio
of the ninth powers equals $(128/135)^3$. On $[1/8,1]$ the derivative
of $(2-t)(2+25t)/(27t^{4/9})$ has the sign of $-16+240t-350t^2$.
Its minimum is at an endpoint. Its values exceed $5/6$ at both endpoints:
at $1/8$ use $2^{-4/3}<2/5$. Raise the inequality at $r=0$ to the power $1-k$ and the one at $r=11/20$ to the power $k$. Multiply them and use
$t^{1-k}(t+2(1-t)/27)^k\le t+2k(1-t)/27$ together with
$b\ge2k/27$.

Consequently the radial majorant $L_0e^{-2y}g^{-r}W_s(y,0)$ is at most
$Nw(x/t)(1+bx/t)$. In the angular mean, impose the single constraint
$\int_0^1t^{a-1}\bar v\,dx=0$. Subtract the potential $4a(1-a)x^2/t^2$ and substitute $v=t^a P$.
In the variable $x$ the eigenvalue equation for the mass density $x/t$ becomes the hypergeometric equation
$x(1-x)P''+[1-(1+2a)x]P'+(\lambda/4-a)P=0$. Its Jacobi polynomial solutions, of degrees $j\ge0$, are complete in $L^2((0,1),(1-x)^{2a-1}\,dx)$. They give the radial eigenvalues
$\lambda=4[j(j+2a)+a]$, and the constraint removes $j=0$.
For nonzero angular frequency $m=|n|/s\ge4$, the positive test
$x^{m/2}t^a$ gives the lower bound $4a(1+m)\ge N$ since $a\ge2/3$.
The min--max principle therefore bounds the second eigenvalue of this mass form relative to the energy by $w$
(the density near $y=0$ is $O(y^{-1-r})$ and $r<1$).
The pointwise comparison below gives the same bound for $B$ on $u^\perp$. To check the stated uniform constants, use
$6+r-4r^2+r^3\le6+1/12$, $N\ge12$ and $(6/5)^k\le6/5$.
For $pw$ use $(6/5)^k\le1+4r/11$ and
$(1-r)(6+r)(1+4r/11)\le6-30r/11=(3/8)N$.

\paragraph{Comparison with the weight at $\theta=0$.}\label{par:angular-damping}
Put $t=e^{-y/s}$ and $W^+(y)=W_s(y,0)$. For $r>0$ let $H=(L_0/L)^{1/r}u/g$ and $Kf=(L_0/g)E_s^{-1}(\omega g^p f)$. The disk equation makes $K$ a positive rotation-invariant Markov kernel, and $H=K(W_sH^p)$. We prove
\[
 W_sH^{-r}\le W^+,\qquad |D_j|H^{-r}\le D_j(y,0)\quad(j=1,2),\qquad
 \int\log H(y,\theta)\,\frac{d\theta}{2\pi}\ge0.
\]
Indeed $H\ge c>0$ by restricting the positive Green source to $1\le y\le2$ and bounding it there below by its minimum, which gives $u\ge c'\min(y,1)\ge c''g$. Jensen gives $\log H\ge K\log W_s+pK\log H$. The circular mean of $\log W_s$ vanishes at every height, so averaging and iteration prove the last inequality.

In frequency $n\ge1$ the radial action $K_n$ of $K$ has positive kernel, from the Dirichlet inverse of $-\partial_y^2+(n/s)^2$. One has $K_n t^n\le t^n$ since
\[
 (-\partial_y^2+(n/s)^2)(g t^n)=t^n(-g''+2(n/s)g')\ge t^n(-g'').
\]
Iterating Jensen and using $\log W_s=4s\sum_{n\ge1}t^n\cos(n\theta)/n$ therefore gives
\[
 \log(W_sH^{-r})\le4s\sum_{n\ge1}(t^n-d_n)\cos(n\theta)/n,\quad
 d_n=r\sum_{j\ge1}p^{j-1}K_n^j t^n\in[0,t^n].
\]
To justify the Fourier expansion, truncate the sum and use the radial majorant $-\log(1-t)$. Its mean under $K$ is uniformly bounded: at height $b$ the radial measure is proportional to $\min(b,y)\omega g^pdy$, and the factor $\min(b,y)$ increases in $y$ while $-\log(1-t)$ decreases. The sum displayed is at most $\log W^+(y)$. The derivatives $s\partial_s\log W_s$ and $s^2\partial_s^2\log W_s$ have nonnegative cosine coefficients (the respective multipliers are $1+ny/s$ and $(ny/s)^2$). Their absolute values are bounded by their values at $\theta=0$, and $D_2/W_s=s^2\partial_s^2\log W_s+(s\partial_s\log W_s)^2$. This proves the two derivative bounds as well. At $q=2$ the corresponding potential estimates hold with factor $L/L_0\le1$ by the radial trial and $\overline W_s\ge1$.

\paragraph{Distance from the radial minimizer on the disk.}
Let $V=\|v\|_{E_s}/\sqrt{L_0}$ in $u=A(g+v)$, $v\perp_{E_s}g$,
and put $d=d_0$. The reference mass form $B_0(\phi)=L_0\int g^{-r}\phi^2\,d\mu_s$ restricted to $g^\perp$ has norm at most $w_c=w+d^2/(1-w)$. Indeed $e=g/\sqrt{L_0}$ has unit energy. For $v\perp_{E_s}g$ of unit energy let $M$, $\beta$ and $c$ be the values of $B_0$ at $(e,e)$, $(e,v)$ and $(v,v)$. Then $M=\int g^q\,d\mu_s\ge1$, since the circular mean of $W_s$ is at least one. The disk equation gives $\int\omega g^pv=0$, so $\beta=\sqrt{L_0}\int\omega D_0g^pv$ and $|\beta|\le d$. Since $\log W_s$ has nonnegative cosine coefficients, $W_s\le W_s(y,0)$, so $B_0$ is at most the radial form of Section~\ref{app:radial-gap}. By the min--max principle the smaller eigenvalue of $B_0$ on the span of $e$ and $v$ is at most $w$, and the larger is at least $M>w$. Hence $(M-w)(c-w)\le\beta^2$, which gives $c\le w+d^2/(1-w)$. The mass form at the minimizer has top eigenvector $u/\sqrt L$, normalized in energy. Its norm on $g^\perp$ is at most
$w_T=w+(1-w)V^2/(1+V^2)$ by decomposing each test into its component along this eigenvector and its orthogonal complement. The lower bound for the angular mean of $\log H$ gives $LA^{-r}\le L_0$:
integrate $H$ against $\omega g^q$ and use Jensen (at $r=0$ use the radial
trial). Along the segment the two endpoint mass forms are
$L\omega W_s u^{-r}$ and $LA^{-r}\omega W_s g^{-r}$.
Applying convexity of
the inverse power along $g+\theta v$ to the projected Euler equation gives
\[
 V\le d/[1-p(w_c+w_T)/2].
\]

\subsubsection{The radial contribution and the three source terms}\label{sec:cylinder-green-sources}\paragraph{Estimating the radial integrals $b_1$ and $b_2$.}
Set $\omega=e^{-2y}$, $h=1/4$, $Q=29/20$ and $J(q)=2(q+1)(4/q)^{q/2}$. We estimate $b_j=\int\omega s^j(\partial_s^jW_s)g^q$ using the area expansion $Z=1+\sum_{n\ge1}Z_n$,
where $Z_n=c_n^2s/(n+s)$ and $c_n=(2s)_n/n!$. Write
$d_n=\sum_{i=0}^{n-1}2s/(2s+i)$, $u_n=s\partial_s d_n$,
$l=n/(n+s)$, and put
\[
 P_1(A)=2d_n+A l,\qquad
 P_2(A)=4d_n^2-2d_n+2u_n+(4d_n-2)Al+A(A+1)l^2 .
\]
The $n$-th contribution to $s^j Z^{(j)}$ is $Z_nP_j(1)$.
Because $g(y)\le\sigma y$, its contribution to $b_j$ is at most
\[
 U_{nj}=Z_n M_n P_j(q+1),\qquad
 M_n=\tfrac12 J(q)\Gamma(q+1)[s/(2(n+s))]^q .
\]
For the $n$th term, the density in $y$ before insertion of $g^q$ is proportional to $e^{-2(n+s)y/s}$. Integrating against $y^q$ gives the gamma moments in $U_{nj}$.

The tail needs only the total unweighted area. For $n\ge N$,
\[
 \frac{U_{nj}}{Z_nP_j(1)}
 \le M_N\begin{cases}
 1+q/(2d_N+1),& j=1,\\
 (q+2)(q+3)/6,&j=2 .
 \end{cases}
\]
For the second ratio, the quadratic has coefficients
$4d_n^2-2d_n+2u_n\ge4d_n-2\ge2$ in the unweighted moments
$1,l,l^2$, with the final coefficient two.
With $l\le1$ these put more relative weight on the smaller moment ratios
$1,q+1,(q+1)(q+2)/2$ than three equal weights.

\paragraph{Uniform bounds for $b_2+2P b_1$.}
Put $h=1/4$, $t=s/h$,
$Q=29/20$, $\eta=(q-Q)/(2-Q)$ and $J(q)=2(q+1)(4/q)^{q/2}$.
We claim that
\[
 \frac{b_2+2P b_1}{s^3}<(30.2-14.7\eta)\sqrt t,\qquad
 \frac{b_2+2P b_1}{s^3}<30\sqrt t .
\]
Here $P=1+(2/q)sA'$ and $A=\log Z$. Logarithmic differentiation of the
gamma quotient for $Z$, using the recurrence and duplication formulas and
$\psi_1(1/4)-\psi_1(3/4)=16G$, gives $hA'=1-\log2$ and $h^2A''=2G-1$, where
$G=\sum_{j\ge0}(-1)^j/(2j+1)^2$. In particular $2P<20/7$ throughout
the range by positivity of the power series of $A$. At $h$ that series
gives
\[
 A<1.21/16+2.1/(5\cdot32(1-1/4))<3/32<\log(11/10).
\]
Use $1-\log2<.307$ (three terms in the series for $2\operatorname{arctanh}(1/3)$)
and $G<.92$ (seven alternating terms). Thus $E_1=hZ'<.338$ and
$E_2=h^2 Z''<1.03$ at $h$.

In the notation $U_{nj}=Z_n M_nP_{nj}(q+1)$ of the radial moment
calculation, each term of $U_{nj}/(s^3\sqrt t)$ increases with $s$.
Indeed $c_n/s$, $d_n$, $u_n$ increase, and the remaining logarithmic
slope is at least $q-1/2-(q+3)s/(n+s)\ge3/50$.
At fixed $s$ each term decreases with $q$: the factors
$\Gamma(q+1)$, $(q+1)\Gamma(q+1)$ and
$(q+1)(q+2)\Gamma(q+1)$ have logarithmic derivatives at most
$\log(q+3)$, and $J'/J\le\log2$. The total logarithmic derivative is
therefore at most $\log(5s/(n+s))\le0$.
In addition, the sum of the positive majorants for $b_2+(20/7)b_1$,
divided by $J(q)$, is log-convex at $s=h$: before the gamma integrals
are evaluated each summand is an integral of $y^q$ against a fixed
positive measure. H\"older's inequality applies to their sum.

At $s=h$, $Z_1=1/20$, $Z_2=1/64$ and
\[
 (P_{11}(1),P_{12}(1))=(2.8,4.88),\qquad
 P_{21}(1)>3.55,\quad P_{22}(1)>9.43.
\]
Here $P_{nj}$ denotes the polynomial for harmonic $n$ and derivative
order $j$. For $n=2$ its parameters are
$d_n=4/3$, $u_n=2/9$, $l=8/9$. In the tail $n\ge3$, use the
factor $M_3$ and the two moment ratios
$C_1=1+q/(2d_3+1)$, $C_2=(q+2)(q+3)/6$ from the area comparison.
After subtracting the first two harmonics, the remaining contributions to $hZ^{\prime}$ and $h^2Z^{\prime\prime}$ are bounded respectively by
$R_1=.338-.14-3.55/64$, $R_2=1.03-.244-9.43/64$.

For $q=Q$ we have $M_1<19/80$.
Indeed $J(Q)5^{-Q}<1$ by $2/(5\sqrt Q)<1/3$ and
$3^{29}>(49/10)^{20}$. Log-convexity of $2^{1-q}\Gamma(1+q)$
from $1$ to $3/2$ bounds its value at $Q$ by
$(3\sqrt{2\pi}/8)^{9/10}<(16/17)^{9/10}<19/20$.
Also $M_2/M_1=(5/9)^Q<35/81$ and $M_3/M_1=(5/13)^Q<13/50$:
use $Q-1>4/9$, $5^4 9^5<7^9$, and $Q>10/7$,
$(5/13)^{10}<(13/50)^7$. Hence $M_2<.103$, $M_3<.062$.
For $q=2$ the exact values are $M_n=12/(8n+2)^2$.
Inserting the first two harmonics and $M_3 C_j R_j$ gives the following
upper bounds at $h$:
\begin{align*}
 b_1(Q)&<.047025+.103(4.85)/64+.062(1.40)R_1<.068,\\
 b_2(Q)&<.134663+.103(18.9)/64+.062(2.56)R_2<.268,\\
 b_1(2)&<.0264+.0371(5.34)/64+.0178(1.50)R_1<.0335,\\
 b_2(2)&<.08688+.0371(23.27)/64+.0178(3.34)R_2<.139.
\end{align*}
Dividing by $h^3$ bounds $(b_2+2Pb_1)/h^3$ at $q=Q$ by
$64(.268+(20/7).068)<29.60$ and at $q=2$ by
$64(.139+(20/7).0335)<15.1$. Monotonicity gives the
uniform bound $30\sqrt t$. To obtain a bound that depends on $q$, use
$-(\log J)''=(q+1)^{-2}+(2q)^{-1}<.512$ on $[Q,2]$. Thus $J(q)/(J(Q)^{1-\eta}J(2)^\eta)$ is less than
$\exp(.512(2-Q)^2/8)<1.020$. Weighted AM--GM bounds the result
by $1.020[(1-\eta)29.60+15.1\eta]\sqrt t<(30.2-14.7\eta)\sqrt t$.

\paragraph{Bounds for the three source norms $d_0,d_1,d_2$.}
Write $h=1/4$, $t=s/h$, $Q=29/20$, $p=q-1$ and
$C(q)=\sqrt q(q+1)(4/q)^{q/2}$. The estimates at the two ends $q=Q$ and $q=2$ of the range of exponents are
\[
\begin{array}{c|ccc}
 &d_0/t^{3/2+\delta_0}&d_1/t^{3/2+\delta_1}&d_2/t^{3/2+\delta_2}\\\hline
Q&.080&.215&1.08\\
2&.044&.132&.85
\end{array}
\]
Here $\delta_j=(4q-2j-1)/6$.
The entries bound the dual norm $d_j^*(q,h)$ of the majorant $C(q)\,y^p\omega D_j$ defined below. Although $D_0=W_s-1$ changes sign in $\theta$, the Fourier coefficients of the majorant are nonnegative, and each radial Green kernel $G_n$ is positive, as shown below. At fixed $s=h$, the square of $d_j^*(q,h)/C(q)$ is thus a sum of integrals of $(yy')^{q-1}$ against nonnegative measures independent of $q$. H\"older's inequality gives the log-convexity of $d_j^*(q,h)/C(q)$ used for interpolation. We first prove the entries at $s=h$.

Put $a=2s$, $c_i=(a)_i/i!$, $d_i=\sum_{j<i}a/(a+j)$ and $x=y/s$.
For the pair $i=n+m,m$ with $n,m\ge0$, $k=n+2m>0$, write
$D=d_i+d_m$, $v=D^2-\sum_{j<i}(a/(a+j))^2-\sum_{j<m}(a/(a+j))^2$,
$B=k+a$, $b=B/(B+n)$ and $l=k/B$.
In the variable $x$ the Dirichlet inverse at frequency $n$ has the kernel
\[
 G_n(x,x')=\int_0^{\min(x,x')}e^{-n(x-u)}e^{-n(x'-u)}\,du ,
\]
which is $(e^{-n|x-x'|}-e^{-n(x+x')})/(2n)$ for $n\ge1$ and $\min(x,x')$ for $n=0$. Hence the squared dual norm of a source $f$ at this frequency is $\int_0^\infty|T_nf(u)|^2\,du$, where $T_n f(u)=\int_u^\infty e^{-n(x-u)}f(x)\,dx$.
Since $g\le\sigma y$ and $\sqrt{L_0}\,\sigma^p=\tfrac12\sqrt q\,\sigma^q\le C(q)$, each Fourier coefficient in $\theta$ of the source $\sqrt{L_0}\,\omega D_jg^p$ is at most the corresponding coefficient of the majorant $C(q)\,y^p\omega D_j$. These coefficients are nonnegative, since $D\ge1$ and $v\ge0$ when $k>0$. Each kernel $G_n$ is positive, so $d_j$ is at most the dual norm $d_j^*(q,s)$ of the majorant, and it is this norm that we estimate. In the variable $x$ each Fourier coefficient of the majorant is $x^p$ times a sum of exponentials with polynomial factors. The common
factor is $C(q)s^{q+1/2}$. We compute the remaining norm from the following moments. For an exponential variable $U$ of mean one put
\[
 C_{nm}=\frac{\sqrt{2-\mathbf1_{n=0}}\,c_i c_m}{(B+n)\sqrt{2B}},
\qquad P_j=\sum_{e=0}^j\frac{j!}{e!}(U/2)^e b^{j-e}.
\]
For $p=0$ or $1$ the entries for $d_0,d_1,d_2$ are $C_{nm}B^{-p}$
times the $L^2$ norms of
\[
 P_p,\qquad DP_p+lP_{p+1},\qquad
 vP_p+2(D-1)lP_{p+1}+l^2P_{p+2}.
\]
Indeed direct integration gives $T_n(x^j e^{-Bx})$ with these
polynomials at $U=2Bu$. All moment products use $\mathbb E U^j=j!$.
Add the bounds at each frequency, then take the square root of the sum of their squares over all frequencies.

At $q=Q$ there is a sharper fractional-power bound. Let $H,J$ be any
one of these polynomials at $p=0,1$ respectively. In the last
case $J=vP_1+2(D-1)lP_2+l^2P_3$. The entry at $Q$ is bounded by
\[
 C_{nm}B^{-p}\|H\|\left(\frac{\langle H,J\rangle}{\|H\|^2}\right)^p .
\]
Jensen's inequality for the normalized measure in $T_n f$
bounds $T_n(x^p f)$ by $(T_n f)^{1-p}(T_n(xf))^p$. Integrating the
square and using $2p=9/10\le1$ gives the displayed estimate.
For a single moment $P_j$ denote its bound before $B^{-p}$ by $U_j$.
For $j=0,1$ it increases with $b$; the two means in the bound are
$b+1/2$ and $(2b^3+2b^2+3b/2+3/4)/(b^2+b+1/2)$.
In particular $U_0<1.21$, $U_1<2.40$ at $b\le1$.

\emph{Estimates for $d_0$ and $d_1$.}
Here $a=1/2$ gives $c_i\le(3i+1)^{-1/2}$ by induction and
$d_i\le1+(\log i)/2$ for $i\ge1$ by midpoint integration.
The table bounds the contributions of the indicated pairs $(n,m)$ before multiplication by $C(q)h^{q+1/2}$.
The columns at $Q$ follow from the fractional-power formula, or from
$U_0D+lU_1$ when that expression gives a simpler bound.
\[
\begin{array}{c|cc|cc}
(n,m)&d_0(Q)&d_1(Q)&d_0(2)&d_1(2)\\\hline
(0,1)&.036&.128&.0284&.115\\
(0,2)&.0065&.029&-&-\\
(1,0)&.143&.285&.132&.307\\
(1,1)&.0143&.054&.0089&.039\\
(1,2)&.0042&.019&-&-\\
(2,0)&.036&.091&.025&.073\\
(2,1)&.0064&.031&-&-
\end{array}
\]
For example $B=3/2$, $b=3/5$, $D=1$, $l=2/3$ in the $(1,0)$ row.
For the $d_1$ entry the two polynomials are
$7/5+U/3$ and $27/25+9U/10+U^2/6$.
All fractional comparisons at $Q$ can be raised to the fortieth power.

For $n=0,1,2$ and the remaining $m$, put $y=m+n/2$,
$A_0=1/12$, $A_1=A_2=1/(6\sqrt2)$. We have
$C_{nm}\le A_n(y+1/4)^{-5/2}$ and $D\le2+\log y$.
For the former use $u=n/(2y+1/2)\le1/2$ and
$(1-u)(1+u)^3\ge1$. As $B=2(y+1/4)$, midpoint integration
bounds the two tails at $Q$ by the integrals of
\[
 A_n2^{-p}z^{-5/2-p}(1.21,\ 4.82+1.21\log z)
\]
with lower limits $11/4$, $13/4$ and $11/4$. The bounds are
\[
 (.0054,.030),\qquad(.0055,.031),\qquad(.0075,.042).
\]
At $q=2$ we only need $n=0,1$ in this comparison and start both tails
at $m=2$. Use $C_{nm}/B\le(A_n/2)y^{-7/2}$ and the factors
$1.59,7.25+1.59\log y$, since $\|P_1\|<1.59$, $\|P_2\|<4.07$.
Integrate from $y=3/2,2$ respectively, obtaining bounds
$(.0097,.052),(.0067,.038)$.

At $Q$ and $n\ge3$ the terms $m=0$ have $b\le.54$,
$U_0<1.02$, $U_1<1.54$, and $C_{n0}/B^p\le n^{-2-p}/(2\sqrt3)$.
Use $1.02D+1.54\le2.56 n^{1/5}$ and take the square sum,
integrating from $5/2$, to bound these groups by $.026,.082$.
For $m\ge1$ use
$C_{nm}/B^p\le [6\sqrt m(n+m+1/4)^{2+p}]^{-1}$ and
$1.21D+2.40\le4.84[m(n+m)]^{1/8}$. Midpoint integration first
of the squares in $n$, then in $m$, bounds the $d_1$ group by
\[
 \frac{4.84}{6\sqrt{3.65}}\int_{1/2}^\infty
 z^{-3/8}(z+11/4)^{-1.825}dz<.129 .
\]
Indeed set $v=z/(z+11/4)$ and bound $(1-v)^{1/5}$ by $1-v/5$.
For $d_0$ the bound is $1.21/(6\sqrt{3.9})$ times the integral of
$z^{-1/2}(z+11/4)^{-1.95}$ from $1/2$, below $.033$ by the same
substitution. The contributions of the four groups of frequencies $n=0$, $n=1$, $n=2$ and $n\ge3$ therefore have norms
below $(.049,.168,.051,.060)$ and $(.188,.392,.165,.212)$.

At $q=2$, for $n\ge3$, $m=0$, use
$C_{n0}/B\le n^{-3}/(2\sqrt3)$, $\|P_1\|<1.17$,
$\|P_2\|<2.35$ and $1.17D+2.35\le3.53 n^{1/6}$.
The square sum from $n=3$ together with the $(2,0)$ row gives
bounds $.030,.097$ for the terms with $m=0$, $n\ge2$.
For $n\ge2,m\ge1$ we have $C_{nm}/B\le
[6\sqrt m(n+m)^3]^{-1}$, and $7.25+1.59\log(n+m)\le
7.25(n+m)^{2/9}$. The two group bounds become
\[
 \frac{1.59}{6\sqrt5}\int_{1/2}^\infty
 z^{-1/2}(z+3/2)^{-5/2}dz<.023,\qquad
 \frac{7.25}{6\sqrt{41/9}}\int_{1/2}^\infty
 z^{-1/2}(z+3/2)^{-41/18}dz<.157 .
\]
For the latter set $v=\sqrt{z/(z+3/2)}$ and use $1+v\le2$.
The midpoint bounds above apply because the integrands are positive, decreasing and convex. For $z^{-w}(A+B\log z)$ differentiate twice.
The contributions of the three groups of frequencies $n=0$, $n=1$ and $n\ge2$ have norms below $(.039,.148,.053)$
and $(.168,.386,.254)$. The square sums at $Q,2$, multiplied
by the common factors $C(Q)h^{Q+1/2}<.416$, $C(2)h^{5/2}=3\sqrt2/16$,
give the first two entries in each endpoint row. For the first
factor use $6/5<\sqrt Q<241/200$ and $(5/12)^{9/20}<27/40$.

\emph{Estimate for $d_2$.}
At fixed total degree $k$ the angular polynomial has squared norm
at most $W_k^2=\sum_{i=0}^k c_i^2 c_{k-i}^2$. Partial fractions give
\[
 W_k^2\le\frac2{3k+2}\left(1+\frac{\log(1+3k)}3\right)
 \le (k+1/2)^{-4/5}.
\]
For the second bound use $e^{z/5}\ge 2/3+2(\log3+z)/9$ for $z\ge0$.
Write $H_j=\sum_{i=1}^j1/i$. Discrete concavity places the largest $d_i+d_{k-i}$ at the center:
it is at most $H_k+\log2<\log(k+1/2)+1.3$; for
$k=2j,2j+1$ its excess over $H_k$ is $H_{2j}-H_j<\log2$,
and $H_k-\log(k+1/2)$ decreases.
For the first four total degrees the exact polynomial formula and angular
orthogonality give, before the common factor, the following bounds:
\[
\begin{array}{c|rrrr}
 k&1&2&3&4\\\hline
 p=0&.181&.350&.230&.205\\
 p=1&.357&.370&.150&.110
\end{array}
\]
For example the two pairs at $k=2$ are $(0,1),(2,0)$;
$(v,D,l)=(2,2,4/5),(2/3,4/3,4/5)$. At $Q$ H\"older in the
positive Green integral followed by weighted AM--GM bounds each
entry by $(1-p)$ times the first-row bound plus $p$ times the second-row bound.

For all higher degrees use the trace bound $\|v(x,\cdot)\|_2\le
\sqrt{x}\|v\|_E$. With $u=q+1/2$, integration of the three monomials
in $D_2$ gives, including the common factor, the tail estimate
\[
 C(q)h^u\Gamma(u)\sum B^{-u-2/5}
 [d^2+2u(d-1)+u(u+1)],\qquad d=\log B+1.3,\quad B=k+1/2 .
\]
For $q=2$ the prefactor is $9\sqrt{2\pi}/64<.354$; the polynomial is
bounded by $z^2+7.6z+12$, $z=\log B$.
For $k\ge3$ the midpoint integral from $B=3$ is below $1.81$
before multiplication by the prefactor. Together with degrees one and two this gives
$d_2(2)<.85$. For $q=Q$ the prefactor is below $.416$ by $\Gamma(39/20)\le1$,
and the polynomial is below $z^2+6.5z+8.7$. Sum odd and even
$k\ge5$ separately with step two, using half the midpoint integral from
$9/2$ and $11/2$ respectively. After multiplication by the prefactor, the odd and even sums are below
$.60$ and $.50$, respectively. The contributions for $k=1,2,3,4$ are below $.112,.153,.083,.069$, respectively.
Odd and even frequencies are orthogonal, giving
$d_2(Q)<[(.60+.112+.083)^2+(.50+.153+.069)^2]^{1/2}<1.08$.
For these integrals $\int_L^\infty B^{-1-b}dB=L^{-b}/b$,
and the mean and variance of $\log B$ are $\log L+1/b,1/b^2$.
Convexity of the two integrands on the integration ranges follows
by two differentiations.

\emph{Scaling the sources.}
Let $t=4s$, $\delta_j=(4q-2j-1)/6$.
Recall that $d_j^*(q,s)$ is the dual norm of the majorant $C(q)\,y^p\omega D_j$ at the opening $s$. Then
\[
 d_j(s)\le d_j^*(q,s)\le t^{3/2+\delta_j}d_j^*(q,h),\qquad j=0,1,2 .
\]
To see this put $\lambda=(1+2h)/(1+2s)\le t^{-1/3}$ by weighted AM--GM.
In each Fourier coefficient the factors $c_i c_m$ and the coefficient
polynomials for the derivatives can be expanded in nonnegative powers of $s$.
Indeed each monomial $s^d e^{-ky/s}$, $d\ge1$, acquires the factor
$d+kx$ or $d(d-1)+2(d-1)kx+(kx)^2$, where $x=y/s$.
Because $0<t\le1$, each coefficient at $s$ is at most $t$ times its value at $h$.
The exponential rates satisfy $k+2s\ge(k+2h)/\lambda$ at every degree $k\ge1$.
Rescale the radial variable by $x=\lambda X$. The half-line Dirichlet kernel
satisfies $G_n(\lambda X,\lambda Y)=\lambda G_{\lambda n}(X,Y)
\le\lambda G_n(X,Y)$ by its positive integral representation.
Hence rescaling a source containing $x^{q-1+l}$, $0\le l\le j$, increases
its dual norm by a factor of at most $\lambda^{q+1/2+j}$. The same comparison holds for
the full sum, since the kernels and the majorizing coefficients are positive.
Including the prefactor $t^{q+1/2}$ gives the stated power.
The endpoint bounds above and H\"older interpolation therefore apply to
the same majorant for $0<s\le 1/4$.

\subsubsection{Bounds for \texorpdfstring{$k_1$ and $k_2$}{k1 and k2}}\label{sec:cylinder-potentials}
Write $h=1/4$, $Q=29/20$, $r=2-q$, $x=y/s$ and
\[
 A(x)=-\log(1-e^{-x})+x/(e^x-1),\qquad
 B(x)=\{x/[2\sinh(x/2)]\}^2 .
\]
On the axis $D_1=W_s(4s)A$ and $D_2=W_s[(4sA)^2+4sB]$.
Both $A$ and $B$ decrease ($xA'=-B$). At fixed $y$,
$W_s$ increases with $s$, so each $k_j(s)$ is bounded by
$4s\,k_j(h)$. At $s=h$ let
\[
 B_*(q)=\tfrac12q(q+1)(2\sqrt q)^{-q},\qquad
 U_q(x)=\frac{x^q e^{-qx/4}}{1-e^{-x}},\qquad
 H(v)=\frac{(1+v)^2}{1+4v}.
\]
The radial bound $L_0e^{-2y}g^{-r}\le C_*(q)y^{q-2}e^{-qy}$,
$C_*(q)=q(q+1)(4/q)^{q/2}/2$, and the positive trial
$\phi_\lambda=\sqrt{x/(1+\lambda x)}$, $\lambda\ge0$, give
\[
 k_1(h)\le 4B_*(q)\sup_x H(\lambda x)U_q A,\qquad
 k_2(h)\le 4B_*(q)\sup_x U_q(A^2+B).
\]
Indeed $-\phi_\lambda''/\phi_\lambda=[4x^2H(\lambda x)]^{-1}$
and completing the square with $\phi_\lambda$ gives the one-dimensional
energy bound on each ray. The prefactors are $B_*(2)=3/8$ and $B_*(Q)<1/2$,
since $1/(2\sqrt Q)<5/12$ and $(5/12)^{29}<(9/32)^{20}$.

For use below,
\[
 \begin{gathered}
 A\le L:=21/20-\log x\quad(0<x\le1),\qquad
 A\le1-\log x+x^2/24\quad(0<x\le2),\\
 U_q=x^{q-1}e^{rx/4}\sqrt B .
 \end{gathered}
\]
In fact $A+\log x\to1$ at zero and $x(A+\log x)'=1-B\ge0$.
For $z=x/2<\sqrt3$, the series of $\sinh z/z$ is bounded termwise by
that of $(1-z^2/3)^{-1/2}$, by comparison of successive coefficients.
Thus $1-B\le x^2/12$, proving the second bound and $A(1)\le25/24<21/20$.
We also have $\sqrt B\le(1+x^2/24)^{-1}$.
On $x\le1$, $e^{rx/4}\le1+3x/20$ by convexity and
$e^{11/80}<23/20$.

The product $U_q A$ decreases for $x\ge1$ when $q=Q$ and for $x\ge3/2$ when
$q=2$. Its logarithmic derivative with respect to $\log x$ is
$q(1-x/4)-x/(e^x-1)-B/A$. Using $A\le(1+x)/(e^x-1)$ and
$(1-e^{-x})^{-1}\ge1/x+1/2$ gives
$x/(e^x-1)+B/A\ge1+x/[2(1+x)]$. Also, for $x\ge2$,
$H(x)\le(1+x)/3$ gives
\[
 H(x)U_q A\le\frac{3}{4(1-e^{-2})^2}x^{q+2}e^{-(1+q/4)x}
 <1.02\,x^{q+2}e^{-(1+q/4)x}.
\]
The power times the exponential is at most one since
$(q+2)/(1+q/4)\le8/3<e$.

For $k_1(Q)$ take $\lambda=1$ and write $p=9/20$.
The upper bound $(1+.15x)x^p L$ is below $1.1$ for $x\le1/20$
and below $1.22$ for $x\le1/10$, by monotonicity there. In
$[1/20,1/10]$ one has $H\le.92$. On $[1/10,1/5]$ use
$\sup x^pL=e^{(21/20)p-1}/p<4/3$, $1+.15x\le1.03$ and
$H\le121/140$. On $[1/5,1]$, $H\le.8$ and
$U_Q A\le4/3+.15(21/20)<1.5$, since $x^{p+1}L\le21/20$.
Thus $H U_Q A<1.2$ up to one. On $[1,2]$ use $H\le1$ and
$U_Q(1)A(1)<1.15(24/25)(21/20)<1.2$; the tail bound covers $x\ge2$.

For $k_1(2)$ take $\lambda=1/2$. On $x\le.4$, $U_2 A\le xL<.79$,
and on $[.4,.6]$ one has $xL<1$, $H(x/2)\le.8$.
On $[.6,1.5]$ use $H(x/2)<.77$ and $U_2 A<1.035$.
For the latter bound, $xL<1.02$ on $[.6,.8]$, while
$xL\sqrt B\le1.05/(1+.8^2/24)$ on $[.8,1]$.
For $1\le x\le1.5$, $x\log x\ge x-1+(x-1)^2/3$ by its second
derivative, whence
$xA\le25/24+(x-1)/8-7(x-1)^2/48<1.07$ by the bound on $A$ above.
Multiply by $\sqrt B\le24/25$. On $[1.5,2]$ use
$U_2 A\le1.07/(1+2.25/24)<1$, $H(x/2)\le.8$.
On $[2,4]$ use $H(x/2)\le1$ and $U_2(2)A(2)<.8$:
indeed $A(2)<.46$ by
$e^{-2}<.1354$ and $-\log(1-z)\le z+z^2/[2(1-z)]$.
For $x\ge4$, $H(x/2)\le H(x)$ and the tail bound is at most
$1.02\cdot4^4e^{-6}<.65$. Hence $H(x/2)U_2 A<.8$.

Finally consider $k_2$. For $q=Q$ and $x\le1$, bound
$U_Q(A^2+B)\le(1+.15x)x^p(L^2+1)$. The maximum of $x^p(L^2+1)$
occurs at $L=(1+\sqrt{1-p^2})/p\in(4.2,4.25)$ and is
below $e^{-.45(4.2-1.05)}(4.25^2+1)<14/3$.
Also $x^{p+1}L^2\le4e^{(p+1)1.05-2}/(p+1)^2<6/5$,
so $U_Q(A^2+B)<14/3+.15(11/5)<5$.
On $[1,2]$ use $U_Q(A^2+B)\le2^{.45}e^{.275}(1+1.05^2)<5$.

For $q=2$, $x(L^2+1)$ is increasing on $(0,1]$ and is at most $1+1.05^2$.
On $[1,2]$ use $xA^2\le A(1)^2$ since $B/A\ge x(1+x/2)/(1+x)\ge3/4$,
and $xB\le2/(1+4/24)^2=72/49$. This gives $U_2(A^2+B)<8/3$.
For $x\ge2$,
\[
 \frac{A^2+B}{1-e^{-x}}
 \le\frac{(1+x)^2e^{-2x}+x^2e^{-x}}{(1-e^{-x})^3}
 <2.1 x^2e^{-x},
\]
so $U_q(A^2+B)<2.1$ by the preceding bound on $x^{q+2}e^{-(1+q/4)x}$.
We have proved
\[
 k_1(Q,h)<12/5,\qquad k_1(2,h)<6/5,\qquad
 k_2(Q,h)<10,\qquad k_2(2,h)<4 .
\]
To interpolate in $q$, first bound the radial weight on each ray by
$C_*(q)y^{q-2}e^{-qy}$. After factoring out $C_*(q)$, H\"older's inequality gives log-convexity in $q$ for each test function. Taking
the supremum preserves log-convexity. Since
$-(\log C_*)''=q^{-2}+(q+1)^{-2}+(2q)^{-1}<1$ on $[Q,2]$,
the factor by which $C_*(q)$ can exceed its geometric interpolant is below $\exp((2-Q)^2/8)<1.04$.

\subsubsection{Uniform curvature on the cylinder}\label{sec:cylinder-curvature-bound}
\paragraph{Combining the estimates.}
Put $Q=29/20$, $p=q-1$, $r=2-q$, $\eta=(q-Q)/(2-Q)$,
$t=4s$ and $[a,b]_\eta=a^{1-\eta}b^\eta$.
For $j=0,1,2$ put $\delta_j=(4q-2j-1)/6$. Let $d,e,h_2$ bound $d_j/t^{3/2+\delta_j}$ for these three indices in order. Let $k,l$ bound $k_1,k_2$, respectively.
Rounding the preceding estimates upwards gives the following pairs of bounds at $q=Q$ and $q=2$.
\[
 (.080,.044),\quad(.215,.132),\quad(1.12,.85),\quad
 (2.4,1.2),\quad(10,4).
\]
For a pair $(a,b_2)$, the upper bound at $q$ is
$[a,b_2]_\eta\exp(b\eta(1-\eta))$. Here $b=\tfrac38(2-Q)^2$ for the first three pairs and $b=\tfrac12(2-Q)^2$ for the last two. These bounds follow by H\"older after removing
the respective prefactors $C,C_*$, using
$-(\log C)''\le3/4$ and $-(\log C_*)''\le1$.

For the coercivity bound write
\[
 w=\max\left(3/8,\frac{(6+r-4r^2+r^3)(1+4r/11)}
                         {16-80r/11}\right).
\]
Section~\ref{app:radial-gap} bounds the mass form on the orthogonal complement by $w$ times the energy.
Throughout $0\le r\le11/20$,
\[
 w\le .61,\qquad pw\le3/8-.395r^2 .
\]
Indeed the bound on $w$ follows from $6+r-4r^2+r^3\le 6+1/12$,
$1+4r/11\le6/5$ and $16-80r/11\ge12$. For the second fraction
in the maximum the difference between $3/8$ and $pw$ is
$[r+75r^2-35r^3-9r^4+4r^5]/(176-80r)\ge70r^2/176$,
since $75-35r-9r^2\ge70(1-5r/11)$ on this interval.
The same bound holds with $w=3/8$. Also $p/(1-w)<1.61$, since $pw\le3/8$ and
$w(1-w)\ge\min\{(3/8)(5/8),.61(.39)\}$.
The displayed interpolant for $d_0$ is decreasing in $q$, hence $d_0\le .080$.
To bound the distance $V$ from the radial minimizer on the disk, use
$p w_c\le pw+1.61d_0^2$ and
$p w_T\le\min(1,pw+5V^2/8)$. Successive substitutions in
$V\le d_0/[1-p(w_c+w_T)/2]$ give $V<.261$, then $V<.134$ and hence
\[
 J_p\ge .625+.395r^2,\qquad
 V\le D t^{3/2+\delta_0},\qquad
 D=\frac{d}{.614+.395r^2}.
\]

In the curvature inequality the nonradial terms divided by $s^3$ are
bounded by $Nt^{14/15}$, since $\delta_0+\delta_2=2\delta_1
=(4q-3)/3\ge14/15$ and $\delta_0\ge\delta_1$. Here
\[
 N=64q\left[D(h_2+20e/7)+\frac p2(l+20k/7)D^2+
                 \frac{(2D+e+pkD)^2}{J_p}\right].
\]
We prove $N\le\min(61+16\eta,64)$ by interpolation between the
exponents in the table. For each of $d,e,h_2,k,l$, let $a,b_2$
be the bounds at $Q,2$. At each tabulated $\eta$ we round
$[a,b_2]_\eta\exp(b[\eta(1-\eta)+9/484])$ upwards.
Between adjacent entries we interpolate their logarithms linearly.
Indeed, on an interval $[u,v]$, the difference between
$\eta(1-\eta)$ and its linear interpolant is
$(\eta-u)(v-\eta)\le(v-u)^2/4$. Here $v-u\le3/11$, so the
added $9/484$ ensures an upper bound throughout the interval.

For $q$ and $p$, the ratio of the linear interpolant to the geometric
interpolant with endpoint values $a_0,b_0$ is at most
$\exp(\log^2(b_0/a_0)/8)$. We therefore multiply their geometric
interpolants by $1.002$ and $1.008$, respectively.
For the denominator of $D$ and for $J_p$ we instead need lower bounds.
On each interval $[r_0,r_1]$, the linear interpolant of $r^2$
exceeds $r^2$ by $(r-r_0)(r_1-r)\le .15^2/4$. Since
$.395(.15)^2/4<.0025$, linear interpolation of
$.6115+.395r^2$ and $.6225+.395r^2$ gives the required lower bounds
for $.614+.395r^2$ and $.625+.395r^2$, respectively.
The reciprocal of each positive affine denominator is log-convex.
Thus these substitutions bound $N$ from above by a sum of log-convex
terms. The sum is convex on each interval, so it suffices to compare
its endpoint values with $61+16\eta$ and $64$. The endpoint bounds are as follows; by this convexity they control every exponent between adjacent rows.
\[
\begin{array}{c|rrrrr|rr}
 q&d&e&h_2&k&l&D&N\\\hline
 1.45&.0802&.2155&1.123&2.408&10.030&.110&61\\
 1.55&.0732&.201 &1.087&2.18 &8.69  &.106&63.9\\
 1.70&.063 &.178 &1.019&1.825&6.87  &.0974&63\\
 1.85&.0532&.155 &.941 &1.50 &5.31  &.0858&55\\
 2   &.0441&.1323&.852 &1.204&4.012&.0722&43
\end{array}
\]
For example, at $q=1.55$ the upper bound for $N$ is $64(1.002)(1.55)$ times
\[
 \begin{aligned}
 &.106(1.087+20(.201)/7)
 +\tfrac12(1.008)(.55)(8.69+20(2.18)/7)(.106)^2\\
 &\hspace{12mm}+\frac{(2(.106)+.201+(1.008)(.55)(2.18)(.106))^2}
 {.6225+.395(.45)^2}.
 \end{aligned}
\]
Multiplying by the prefactor gives less than $63.9$. Every entry in the $N$ column is at most both $61+16\eta$ and $64$.

By Section~\ref{sec:cylinder-green-sources} the radial error term is less than
$(30.2-14.7\eta)\sqrt t$.

After division by $s^3$, the radial and nonradial error terms together are at most
$30.2\sqrt t+61t^{14/15}+1/4$. Indeed, the excess over the first two terms is bounded by
$\min(16\eta,3)-14.7\eta$ times $\sqrt t$ whenever this is positive,
and its maximum is $1.3(3/16)<1/4$. By concavity of the two powers,
this is at most $19.42+72.04t$.
Keeping degrees $3,5,7,9$ in the positive area series gives
\[
 \frac{K''}{csK}\ge A_0(t):=57.6+23.2t^2+7.875t^4+2.49t^6+
           \frac{t^3}{64}(14.4+3.75t^2)^2 .
\]
The polynomial $A_0$ is convex. Direct substitution gives
$A_0(13/16)>79.4$ and $70<A_0'(13/16)<72$. Its tangent there
exceeds $20.42+72.04t$ at both zero and one. Thus $f''/(csf)>1$
on $0<s\le 1/4$. Since
$\mathcal C_q=\cos^2(\pi s) f''/(2\pi^2f)$, this gives
$\mathcal C_q>s/(4\pi^2)$ throughout. In particular, this gives $\mathcal C_q>1/(20\pi^2)>1/200$ for $1/5\le s\le1/4$.

\phantomsection\label{par:angular-estimates-end}

\subsection{Positive curvature at acute openings}
\begin{table}[htbp]
\centering\small
\caption{Curvature bounds for acute openings. Here $q_\ast=5801/4000$,
$Q=29/20$ and $s=\alpha/(2\pi)$.}
\label{tab:acute-coverage}
\begin{tabularx}{\linewidth}{@{}lllX@{}}
\toprule
Exponent & Opening $\alpha$ & Curvature bound & Proof\\
\midrule
$1\le q\le q_\ast$ & $(0,\pi/2]$ & $F_q''>0$ & Disk estimate, Proposition~\ref{S:prop:lowq}.\\[2pt]
$Q\le q\le2$ & $(0,\pi/2]$ & $\mathcal C_q>s/(4\pi^2)$ & Cylinder second variation, Section~\ref{sec:cylinder-curvature-bound}.\\
\bottomrule
\end{tabularx}
\end{table}
\begin{proposition}[Positive curvature at acute openings]\label{prop:real-acute-curvature}
For every $1\le q\le2$ one has $F_q''(\alpha)>0$ for $0<\alpha\le\pi/2$.
\end{proposition}
\begin{proof}
Proposition~\ref{S:prop:lowq} gives the first row of Table~\ref{tab:acute-coverage}.
For $Q\le q\le2$, Section~\ref{sec:cylinder-curvature-bound} gives
$f''/(csf)>1$. Since $\mathcal C_q=\cos^2(\pi s) f''/(2\pi^2f)$,
$c=2/q\ge1$ and $\cos^2(\pi s)\ge1/2$ up to $s=1/4$, this gives the second row.
In particular $\mathcal C_q>1/200$ on $1/5\le s\le1/4$.
The two exponent ranges overlap since $Q<q_\ast$.
\end{proof}

\subsection{Proof of Theorem~A}
Differentiability follows from Lemma~\ref{S:lem:exponent-continuity}.
For $q\le q_\ast$, Proposition~\ref{S:prop:lowq} covers every opening.
For $q_\ast\le q\le2$, Propositions~\ref{S:prop:spectral-wide}
and~\ref{S:prop:central-sum} give $\mathcal C_q>0$ on $z\le1$.
Proposition~\ref{prop:real-acute-curvature} covers $z\ge1$. \qed

\paragraph{Uniformity for exponents near $[1,q_\ast]$ and near $2$.}\label{par:exponent-collars}
Put $K(q,\alpha)=\mathcal C_q(\cot^2(\alpha/2))/\sin(\alpha/2)$.
It extends continuously to $[1,2]\times[0,\pi]$, with
\[
 K(q,0)=24(2/q)\zeta(3)/\pi^3,\qquad K(q,\pi)=1 .
\]
Indeed, writing $s=\alpha/(2\pi)$, $f=F_q(2\pi s)/(2\pi c_q)$, one has
$K=\cos^2(\pi s)f''/[2\pi^2\sin(\pi s)f]$.
Proposition~\ref{S:prop:disk-limit} gives $f/s\to1$ and
$f''/s^2\to48(2/q)\zeta(3)$ uniformly. Proposition~\ref{S:prop:uniform-flat}
gives the other endpoint, and Lemma~\ref{S:lem:exponent-continuity} covers the interior.
For $q\le q_\ast$, Proposition~\ref{S:prop:lowq} shows that $K>0$ at every opening. For $q=2$ this follows from
Proposition~\ref{S:prop:wide-linear}, the estimate on p.~\pageref{par:spectral-obtuse}, and
Proposition~\ref{prop:spectral-endpoint-curvature}. Hence the continuous function $K$ has a positive minimum on the compact set
$([1,q_\ast]\cup\{2\})\times[0,\pi]$. By uniform continuity there are
$\eta,\kappa>0$ with $K\ge\kappa$ for every opening whenever
$q\in[1,q_\ast+\eta]\cup[2-\eta,2]$. Joint continuity in the interior and the uniform endpoint limits extend the positive lower bounds for $q\le q_\ast$ and $q=2$ to these neighborhoods. No curvature lower bound on the whole interval $q_\ast<q<2$ is needed for this step. The argument gives no explicit value for $\eta$.

\part{Polygons with a bounded number of convex vertices}\label{part:convex-vertices}
\section{Counting convex vertices instead of sides}\label{C:sec:intro}
A vertex of a simple polygon is \emph{convex} if its interior angle is less than $\pi$ and \emph{reflex} if its interior angle exceeds $\pi$. As before, vertices at which adjacent sides are collinear are removed. The exterior turn at a vertex is $\pi$ minus its interior angle, so it is positive exactly at the convex vertices; we also speak of positive turns. A polygon with at most $n$ sides has at most $n$ convex vertices, but a polygon with few convex vertices may have any number of sides. In this part we bound only the number of convex vertices. Write
\[
 \Phi_q(P)=|P|^{2/q}\lambda_q^{\rm S}(P),
\]
and let $R_k$ be a regular $k$-gon. Thus Theorem~\ref{thm:main} states that $\Phi_q(P)\ge\Phi_q(R_n)$ when $P$ has at most $n$ sides. We obtain the following.

\begin{itemize}
\item \emph{Three and four convex vertices.} If $P$ has at most three convex vertices, then $\Phi_q(P)\ge\Phi_q(R_3)$, and if it has at most four, then $\Phi_q(P)\ge\Phi_q(R_4)$, for every $1\le q\le2$, with equality only for an equilateral triangle and a square, respectively (Theorems~\ref{C:thm:three-convex-turns-sharp} and~\ref{C:thm:four-convex-sharp}).
\item \emph{Comparison with $R_n$ for at most $n-1$ convex vertices.} If $n\ge5$ and $P$ has at most $n-1$ convex vertices, then at equal area $T(P)<T(R_n)$ and $\lambda_1(P)>\lambda_1(R_n)$, with explicit gaps (Theorem~\ref{C:thm:convex-vertices}), and $\Phi_q(P)>\bigl(1+1/(1000qn^4)\bigr)\Phi_q(R_n)$ for every $1\le q\le2$ (Corollary~\ref{C:cor:all-q-gap}).
\item \emph{Symmetric polygons.} If $P$ has at most $m$ convex vertices and is invariant under rotation by $2\pi/m$, then $\Phi_q(P)\ge\Phi_q(R_m)$ for every $1\le q\le2$, with equality only for $R_m$ (Theorem~\ref{C:thm:cyclic-sharp}). The same holds under an equal-area condition on sectors (Corollary~\ref{C:cor:equal-sector-turns}).
\item \emph{Near the regular polygon.} Normalize $P$ and $R_n$ to unit area. If $P$ has at most $n$ convex vertices and is close to $R_n$ in Hausdorff distance, then $\Phi_q(P)\ge\Phi_q(R_n)$, with equality only for $R_n$, for every $n\ge3$ and every $1\le q\le2$ (Theorem~\ref{C:thm:local}). The deficit is at least a constant times the square of the area of the symmetric difference with the nearest copy of $R_n$ (Corollary~\ref{C:cor:local-stability}); unlike for polygons with at most $n$ sides, the Hausdorff distance cannot be used here. Pushing a side of $R_n$ inward along a concave polygonal arc increases $\Phi_q$ at an explicit first-order rate (Proposition~\ref{C:prop:first-order}), and filling in a small corner of a reflex vertex decreases $\Phi_q$ (Lemma~\ref{C:lem:fill-reflex}).
\item \emph{Many convex vertices.} If $m\ge5$ and $P$ has area $\pi$ and at most $m$ convex vertices, then $1-8T(P)/\pi>8\zeta(3)m^{-3}-16\zeta(4)m^{-4}-160m^{-5}$. Hence the least torsion deficit in this class equals that of $R_m$ up to $O(m^{-5})$, uniformly in the number of reflex vertices (Theorem~\ref{C:thm:sharp-asymptotic}).
\end{itemize}
For $n\ge5$ we do not know whether every polygon with at most $n$ convex vertices satisfies $\Phi_q(P)\ge\Phi_q(R_n)$ (Question~\ref{C:que:sharp}).

The proof of the quantitative comparison with $R_n$ revises the argument in Sections~3--4 of the first version \cite{wangv1}. In particular Lemma~\ref{C:lem:shifted-hardy} is \cite[Lemma~3.3]{wangv1}. The proofs in this part use Part~\ref{part:endpoints} only through Theorem~\ref{thm:main}, the Steiner symmetrization argument for triangles and, in Corollary~\ref{C:cor:local-stability}, Theorem~\ref{S:thm:full-stability}. They use Part~\ref{part:intermediate} only through Proposition~\ref{S:prop:exponent-round}.

Sections~\ref{C:sec:hankel}--\ref{C:sec:certification} work at area $\pi$ and use the following quantities.
\begin{center}
\begin{tabular}{@{}p{.19\textwidth}p{.73\textwidth}@{}}
\toprule
Quantity & Meaning at area $\pi$ \\
\midrule
$D(P)$ & The torsion deficit $1-8T(P)/\pi$, which is zero only for a disk. \\
$\rho_n$ & An explicit strict upper bound for $\lambda_1(R_n)/j_{0,1}^2$. \\
$d_n$ & The number $1-\rho_n^{-2}+1/(1000n^4)$; it is not $D(R_n)$. \\
\bottomrule
\end{tabular}
\end{center}

\section{Three convex vertices}\label{C:sec:three}
In this section we prove that the equilateral triangle is optimal among polygons with exactly three convex vertices. We use the triangle case of Theorem~\ref{thm:main} and the Steiner symmetrization argument on p.~\pageref{par:small-polygons}.
\subsection{One side of the hull on the boundary}
\begin{theorem}[One side of the hull on the boundary]\label{C:thm:three-turn-straight-hull-side}
Let $P$ be a simple polygon with exactly three convex vertices, and suppose that it contains an entire side of its convex hull in its boundary. For every $1\le q\le2$,
\[
 |P|^{2/q}\lambda_q^{\rm S}(P)\ge |R_3|^{2/q}\lambda_q^{\rm S}(R_3),
\]
with equality only for an equilateral triangle.
\end{theorem}
\begin{proof}
Delete flat vertices. The three hull vertices are the three convex vertices, since every extreme vertex is convex. Place the straight hull side on $y=0$ and the opposite hull vertex at height $H$. The chain on each other hull side bounds a convex pocket with that side:
the chain vertices are reflex as seen from $P$, and the endpoints lie on
the supporting chord. For every $0<y<H$, each pocket cuts off an interval
attached to the corresponding end of the horizontal hull section.
The rightmost point of the left pocket is a concave function of height
and the leftmost point of the right pocket is convex, by convexity of
the pockets. A straight side uses just the affine endpoint function.
The pockets are disjoint in each open section since the two chains
of a simple polygon cannot meet. Horizontal Steiner symmetrization preserves area and the $L^q$ norm of the rearranged test and decreases its Dirichlet energy, as in the argument on p.~\pageref{par:small-polygons}. The resulting polygon is
\[
 D=\{(x,y):0<y<H,\ |x|< f(y)\},
\]
where $f$ is positive, piecewise affine and convex, $f(0)=b>0$ and $f(H)=0$. Each slope is strictly negative (the last one is negative and convexity orders the slopes). Every genuine kink of either original chain yields a genuine kink of $f$, since both contributions to the jump in $f'$ are nonnegative. It suffices to lower the quotient strictly from each kinked $D$ to some triangle of the same area. If there are no kinks the original polygon is a triangle and Theorem~\ref{thm:main} applies directly.

Let $E=[C,J]$ be a side to be moved, with convex endpoint $C$ and midpoint $M$. Suppose the perpendicular bisector of $E$ cuts off a triangle $T$ containing $C$. Let $S$ be reflection across the bisector. Assume $S$ maps the interior of $T$ into the polygon and sends an open segment of $(\partial T\cap\partial D)\setminus E$ into the interior of $D$. We claim that $u\circ S>u$ inside $T$ for the positive minimizer $u$ of unit $L^q$ mass. To verify this put $p=q-1$ and $v=u\circ S$. The test $w=(u-v)_+$ has zero trace on the whole boundary of the triangle. If $p<1$, the equations, concavity and the positive solution $v$ give
\[
 \int|\nabla w|^2\le p\lambda\int v^{p-1} w^2\le p\int|\nabla w|^2,
\]
where the second bound follows by testing the equation of $v$ with $w^2/(v+\eta)$ and taking the limit (for $p=0$ simply subtract the equations). For $p=1$, the first Dirichlet eigenvalue of the triangle, a proper subdomain, strictly exceeds that of $D$, giving the same conclusion. The nonnegative difference $v-u$ is superharmonic and is positive on an open portion of the boundary, hence positive inside. Hopf's lemma on the relative interior of $CM$ gives
\[
 g(SX)>g(X)\qquad (X\in (C,M)),\qquad g=|\partial_\nu u|^2 .
\]
The reflection preserves the normal direction to $E$. If the normal velocity of $E$ is affine along it, with negative sign at $C$, positive at $J$ and zero mean, it follows by pairing reflected points that $\int_E g V_\nu>0$.

Each deformation below preserves area. Fix a time before two neighbouring vertices collide, and triangulate the polygon at that time. At all sufficiently close times the corresponding triangles stay nondegenerate. The piecewise affine maps between the triangulations and their inverses have uniformly bounded Lipschitz constants, and the velocity field is continuous and piecewise affine. The minimum $\lambda_q^{\rm S}$ is differentiable with respect to time. Indeed, using the minimizer at one time as a test function at a nearby time, in both orders, bounds the difference quotient from above and from below. The minimizers, carried back to the reference polygon, converge strongly in $H^1$ by compactness and uniqueness, so the two bounds have the same limit. Integrating by parts in the resulting formula on the reference polygon gives
\[
 (\lambda_q^{\rm S})'=-\int_{\partial D}|\partial_\nu u|^2 V_\nu.
\]
This is the same formula as in the proof of Proposition~\ref{C:prop:dihedral-sharp}. The integration by parts is carried out after removing a disk of radius $\rho$ about each vertex. At a vertex of angle $\omega<2\pi$ one has $|\nabla u|=O(\rho^{s-1})$ on the circular arc of radius $\rho$ inside the polygon, for any $1/2<s<\min(1,\pi/\omega)$, so the boundary terms on the removed arcs are $O(\rho^{2s-1})$ and vanish as $\rho\to0$. To prove the gradient bound, extend the source $\lambda u^p$ by zero to a containing disk. It lies in $L^t$ for every finite $t$ by planar Sobolev embedding. Its Dirichlet potential on the disk is bounded and dominates $u$ by comparison. Hence $\lambda u^p$ is bounded. In a truncated sector at the vertex, with polar coordinates $(r,\theta)$, compare $u$ with $C r^s\cos(s(\theta-\omega/2))-c r^2$, where $c$ is so large that $-\Delta$ of this function exceeds the right-hand side. Since the cosine has a positive minimum, this gives $u=O(r^s)$. Interior gradient estimates on rescaled annuli, after odd reflection across their straight sides, then give the stated bound for $|\nabla u|$. On a segment of the triangulation inside the polygon, the boundary terms from its two sides cancel.

Finally, the quotient is continuous up to a collision. All the domains in the motions below have the form $\{0<y<H,\ |x|<f(y)\}$ with $f$ positive on $(0,H)$. Their boundaries consist of fixed straight segments and the moving facets, and their areas are equal. The limiting domain at a collision has the same form and is a Lipschitz domain. Every compact subset of its interior is contained in the domains shortly before the collision, so testing with compactly supported functions bounds the upper limit of the quotients by the quotient of the limiting domain. For the reverse inequality extend the minimizers by zero to a common disk. By compactness a subsequence converges weakly in $H^1$ and strongly in $L^q$. The limit vanishes outside the limiting polygon, so it has zero trace on its boundary, and lower semicontinuity of the energy gives the bound.

Write $A=|D|$. First suppose $A\ge\sqrt3\,b^2$. We eliminate a kink by moving the last facet and its mirror, keeping the base fixed. Its endpoints are $(X,a)$ and $C=(0,H)$, with $X=t(H-a)>0$. Extend the preceding facet as the line $d(H_0-y)$; then $d>t$ and $H>H_0$. Keep
\[
  X=d(H_0-a),\qquad C_0=X_{\rm old}(H_{\rm old}-H_0)>0,\qquad
  H=H_0+C_0/X,\qquad t=X/(H-a).
\]
Move $a$ down along the preceding facet, continuing the unchanged lower boundary up to $y=a$ and using $t(H-y)$ above it. The excess of half-area above the continuation $[d(H_0-y)]_+$ on $[a,H]$ is $X(H-H_0)/2=C_0/2$. Thus the total area is constant and $0<t=d/(1+C_0d/X^2)<d$ preserves convexity. The velocity is outward near $(X,a)$, since the last line lies strictly to the right of the old preceding line beyond the new intersection, and inward at the top since $H$ decreases. Its normal component is affine on the facet by affine interpolation of endpoint velocities. Its integral is zero because only this facet and its mirror move normally and the area is fixed.

For this facet use the coordinate $s=y-tx$, which is constant on lines parallel to its perpendicular bisector. The bisector is the line $s=s_m$, where $s_m=a+(1-t^2)(H-a)/2$, and it cuts off the part $s>s_m$ of $D$, which contains the top vertex $C$. The chord from $(b,0)$ to $(0,H)$ lies to the right of the graph, so $A\le bH$. The last supporting line gives $tH\le b$. Therefore $H\ge\sqrt3\,b$ and $t\le1/\sqrt3$, strictly when the graph is kinked. On $0\le y\le a$ the maximal $s$ at each level equals $y+t f(y)$, a convex function. At $0$ it is $tb\le H/3\le s_m$. At $a$ it is $a+t^2(H-a)\le s_m$. The part cut off thus lies in the isosceles triangle with vertices $(\pm X,a),(0,H)$ at the top of $D$. It is a triangle with one vertex at $C$ and two vertices on the bisector. The reflection fixes the two vertices on the bisector and sends $C$ to $(X,a)$. These three points lie in the same isosceles triangle, so the interior of the reflected triangle lies in $D$. The apex angle $2\arctan t$ is strictly smaller than the base angle $\pi/2-\arctan t$, so a part of its left side near the apex reflects strictly inside the isosceles triangle near $(X,a)$. The comparison at reflected points gives $\int_E|\partial_\nu u|^2V_\nu>0$, so the quotient strictly decreases until the next collision. Stop when $a$ reaches the lower endpoint of the preceding facet. If that endpoint is on the base, the new graph has only one facet on the right and the domain is a triangle. Otherwise one kink is removed. Repeat until reaching a triangle.

If $A<\sqrt3\,b^2$, move the first right facet and its mirror. Write its endpoints as $C=(b,0)$ and $J=(X,a)$ on the following line $c-ty$, with old first slope $-d$, $d>t>0$. Keep
\[
 C_0=a_{\rm old}(b_{\rm old}-c)>0,\qquad X=c-ta,\qquad
 b=c+C_0/a,\qquad d=t+C_0/a^2.
\]
Increase $a$ along the following facet. Below $a$ use the first line $b-dy$ and above $a$ retain the old graph. The excess half-area over the extension of the following line on $[0,a]$ is $a(b-c)/2=C_0/2$. Thus area is fixed and the graph stays convex. The base slides tangentially, $b$ decreases, and the normal velocity is inward near $C$ and outward near $J$ (the new first line lies strictly to the right of the following line below their new intersection). It is affine and has zero mean for the same reason as above.

While $A<\sqrt3\,b^2$ we have $d\ge b^2/A>1/\sqrt3$, since $f(y)\ge(b-dy)_+$ throughout the polygon and $H\ge b/d$. Use $s=dx-y$, so $s_m=[d(b+X)-a]/2$. Above $y=a$, the largest $s$ occurs at the right edge and is at most $dX-a<s_m$. On the left side below $a$ both endpoints have $s\le s_m$: the upper one has $s=-dX-a$, and the base endpoint has $s=-db\le s_m$ since $a=(b-X)/d\le d(3b+X)$. Hence the bisector cuts off only a triangle at the corner $(b,0)$. Its two vertices on the bisector are fixed by the reflection, and the corner is sent to $J$. These three points lie in the triangle $0\le y\le b/d$, $|x|\le b-dy$, which is contained in $\overline D$, so the interior of the reflected triangle lies in $D$. The part of its side on the base near $(b,0)$ reflects strictly into $D$ near $J$, since $J$ lies in the relative interior of the right side of that triangle, and the reflection preserves the inward half-plane along that side. The comparison at reflected points again gives strict decrease. Stop either when the next facet disappears or when $A=\sqrt3\,b^2$. In the first event the number of kinks drops (and the domain is a triangle if $X=0$). In the second event apply the preceding top-facet construction, which keeps $b$ fixed. There are finitely many facets, and the motion passes from the case $A<\sqrt3\,b^2$ to the case $A\ge\sqrt3\,b^2$ at most once. This proves the assertion.
\end{proof}

\subsection{All three sides of the hull indented}
Throughout, polygons are open interiors, flat vertices may be suppressed, and an ear at $A$ with successive neighbours $D,F$ is called a weak ear if $\operatorname{conv}\{A,D,F\}\subset\overline P$.

\begin{lemma}[The pockets of the hull]\label{C:lem:three-pockets}
Let $P$ have exactly three convex vertices $A,B,C$. Its hull $H$ is their triangle. Each boundary chain from one hull vertex to the next, together with the intervening hull side, encloses a convex polygonal pocket $K_{AB},K_{BC},K_{CA}$. For a straight chain take the pocket to be the closed hull side. The interiors of the pockets are pairwise disjoint, and removing the pockets from the interior of $H$ gives $P$.
\end{lemma}
\begin{proof}
Every extreme vertex of $H$ is a convex polygon vertex. There are at least three, so there are precisely three. The boundary chain cannot touch the relative interior of any hull side at a genuine vertex, since an interior angle greater than $\pi$ cannot meet a supporting line there. A nonstraight chain and its chord therefore give a Jordan polygon. Its turns at the interior vertices of the chain are all convex as seen from the side outside $P$; at the endpoints the pocket angles are strictly smaller than the angles of the triangular hull. The pocket is convex. The description of the complement now follows from the Jordan curve theorem, or by following cuts along the three hull sides.
\end{proof}

\begin{lemma}[A reflected triangle]\label{C:lem:three-reflected-cap}
Suppose all three hull sides are dented. There is a convex vertex $A$ with incident edges $AD,AF$, both ending at reflex vertices, and a point $E\in(A,F]$ such that $|AE|=|DE|$ and $\operatorname{conv}\{A,D,E\}\subset\overline P$. The open triangle lies in $P$, and an open subinterval of $DE$ is in $P$.
\end{lemma}
\begin{proof}
A simple polygon with more than three vertices has at least two weak ears. For points in general position, triangulate the polygon. Its dual tree has at least two leaves, each defining an ear at a convex vertex. For arbitrary vertices, perturb them into general position, preserving simplicity and the strict signs of the finitely many turns. Pass to a subsequence having the same indices for two ears. The limit ear triangles belong to $\overline P$: if an interior point of a limit triangle were outside $\overline P$, its winding number would be zero before and after small perturbations by uniform convergence of the boundary. Closure gives the result for the whole triangle. If a nondegenerate closed triangle lies in $\overline P$, its interior lies in $P$ because a simple polygon is the interior of its closure.

If a weak ear vertex $A$ has angle $\theta\le\pi/3$, choose the names so that $l=|AD|\le |AF|$. The perpendicular bisector of $AD$ meets $AF$ at $E$ with $|AE|=l/(2\cos\theta)\le |AF|$. Its triangle $ADE$ is contained in the weak ear. The assertion about an open part of $DE$ holds when $E$ is interior to $AF$, since the segment near $E$ enters the interior transversally. If $E=F$, then near $D$ the ray $DF$ is strictly inside the local interior sector of $P$: the angle $ADF$ of the ear is less than $\pi$ and lies on the interior side of $DA$, whereas the interior angle at $D$ exceeds $\pi$. This proves the strictness statement in the borderline case.

It remains to consider two weak ear vertices $B,C$ both having actual angles $\theta_B,\theta_C>\pi/3$. Denote the hull angles by $\alpha_A,\alpha_B,\alpha_C$. Then $\alpha_A<\pi/3$, and $\alpha_B-\theta_B<\pi/3$. Name $B,C$ so that the immediate neighbour $D$ of $A$ on the chain towards $B$ satisfies $l=|AD|\le |AF|$, where $F$ is the other neighbour. Set $R=|AB|$, $\eta=\angle BAD$, $\theta=\angle DAF$. The pocket $K_{AB}$, viewed from $B$, is confined between $BA$ and its first chain edge; the angular width of this cone is at most $\alpha_B-\theta_B$. In particular $\angle ABD<\pi/3$, and $\eta+\theta\le\alpha_A<\pi/3$. In the triangle $ABD$ the angle at $D$ exceeds $\pi/3$, so $l<R$.

Choose $E$ on $AF$ with $|AE|=l/(2\cos\theta)<|AF|$. We claim $\angle ABE\le\pi/3$. Measure angles from $AB$ towards the inside of $H$. We have $|AE|<R$, and $E$ has coordinates $(|AE|\cos(\eta+\theta),|AE|\sin(\eta+\theta))$ relative to $A$, $AB$ horizontal. The claim is equivalent to
\[
 (l/R)\sin(\pi/3+\eta+\theta)\le\sqrt3\cos\theta .
\]
When $\cos\theta\ge1/\sqrt3$ the bound is immediate. Otherwise $\pi/3+\theta>\pi/2$ and $\pi/3+\eta+\theta<2\pi/3$, whence
\[
 \sin(\pi/3+\eta+\theta)\le\sin(\pi/3+\theta)
 = (\sqrt3\cos\theta+\sin\theta)/2\le\sqrt3\cos\theta .
\]
The last inequality uses $\theta<\pi/3$.

At $B$, the opposite pocket $K_{BC}$ is confined between its first chain edge and the ray $BC$. The first of these rays makes an angle at least $\theta_B>\pi/3$ with $BA$. Every point of $ADE$ is in the closed cone of width $\pi/3$ from $BA$ towards $BC$, and $B\notin ADE$. Thus $ADE$ misses the opposite pocket. The pockets adjacent to $A$ lie outside the open sector between $AD$ and $AF$, by the supporting lines of their first facets. The open triangle is therefore in $P$. The relative interior of $DE$ is inside this open sector and inside the hull, and also misses $K_{BC}$, as required.
\end{proof}

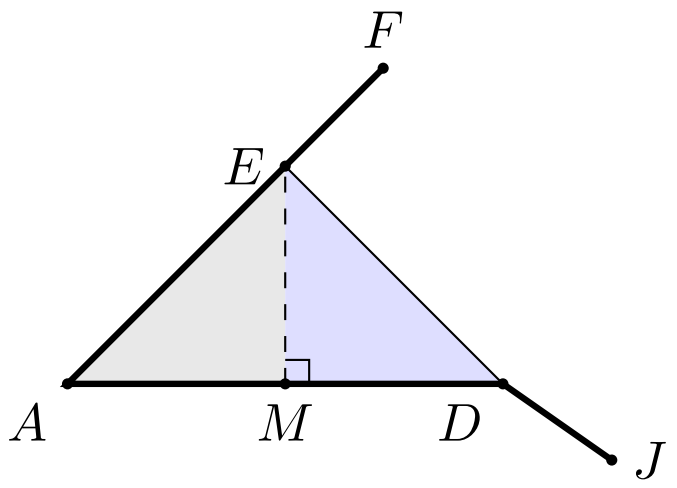
\begin{center}
\begin{tikzpicture}[scale=.72,every node/.style={font=\small}]
\coordinate (A) at (0,0); \coordinate (D) at (4,0);
\coordinate (M) at (2,0); \coordinate (E) at (2,2);
\coordinate (F) at (2.9,2.9); \coordinate (J) at (5,-.7);
\fill[black!9] (A)--(M)--(E)--cycle;
\fill[blue!13] (D)--(M)--(E)--cycle;
\draw[very thick] (F)--(A)--(D)--(J);
\draw[dashed] (M)--(E);
\draw (E)--(D);
\draw (2,.22)--(2.22,.22)--(2.22,0);
\foreach \name/\where in {A/below left,D/below left,M/below,E/left,F/above,J/right}{
  \fill (\name) circle (1.45pt) node[\where] {$\name$};}
\end{tikzpicture}

{\small The sides $AM$ and $AE$ are on the boundary.
Reflection across $ME$ sends $A$ to $D$ and $AE$ to $DE$.}
\end{center}

\begin{lemma}[A descent step]\label{C:lem:good-facet-descent}
When all three hull sides are indented there is a local deformation
that preserves the area, the vertex count and the signs of all turns,
stays in the old closed convex hull, and strictly lowers the quotient.
\end{lemma}
\begin{proof}
Use $A,D,E,F$ from Lemma~\ref{C:lem:three-reflected-cap}, and let $M$
be the midpoint of $AD$. Reflection in $ME$ maps $AME$ to $DME$.
The reflected open boundary segment has an open portion in $P$.
The comparison in the proof of
Theorem~\ref{C:thm:three-turn-straight-hull-side} applies on this
triangle alone. At points that correspond under the reflection, $|\partial_\nu u|^2$ on $DM$ is strictly
larger than on $AM$.

Let $J$ be the other neighbour of $D$. Move $A$ towards $F$ along
$AF$ and $D$ towards $J$ along $DJ$, with all other vertices fixed.
Relative to the outer normal along $AD$ the respective velocities
have strictly negative and strictly positive components: $A$ is
convex and $D$ reflex. The two adjacent sides move only tangentially.
Differentiating the area with respect to either travel distance gives
one half the length of $AD$ times the corresponding unit normal
component. The two derivatives are nonzero with opposite signs.
The implicit function theorem therefore gives an exactly
area-preserving motion with both initial travel speeds positive.
The normal speed on $AD$ is affine, has zero mean, and is positive
at $D$. The shape derivative proved above, paired across $M$, gives
a strictly negative derivative of $\lambda_q^{\rm S}$. Small positive
motions preserve the strict turns and simplicity. Both moving vertices
remain on closed segments in the old hull, establishing the constraint.
\end{proof}

\begin{lemma}[Existence of a minimizer]\label{C:lem:three-pockets-compact}
Fix a positive area $A$, a bounded closed triangle $H_0$ and an integer $N$. The infimum of $\lambda_q^{\rm S}$ among simple polygons with three convex vertices, at most $N$ sides, area $A$ and closure contained in $H_0$, if the class is nonempty, is attained.
\end{lemma}
\begin{proof}
Choose a minimizing sequence $P_j$. All polygons lie in a common ball. Pass to a subsequence with fixed combinatorics, with all vertices convergent and with triangular hulls $H_j=A_jB_jC_j$ convergent to a triangle $H=ABC$. It is nondegenerate since its area is at least $A$. The three convex pockets converge to compact convex sets $K_{AB},K_{BC},K_{CA}$. Each contains its hull side and is a convex polygon or that side. The interiors of any two pockets are disjoint, since every compact subset of the interior of a Hausdorff limit of convex bodies is eventually contained in those bodies. Continuity of polygonal area gives
\[
 |H|-|K_{AB}|-|K_{BC}|-|K_{CA}|=A>0.
\]
In particular the three pockets have no common point: if $X$ belonged to all, the three triangles $ABX,BCX,CAX$ would cover $H$ inside their respective pockets.

For a positive-area pocket let $J_{AB}$ be its boundary arc complementary to the relative interior of the hull side; for a degenerate pocket put $J_{AB}=AB$. These arcs have no bends of the convex sign seen from inside $H$ away from the pockets. For the number of facets, the limiting pocket is the convex hull of the limiting vertices on the corresponding chain, so removing redundant or coincident vertices can only lower the count.

The intersection of two adjacent pockets is either their common hull vertex $A$, or a segment from $A$. Indeed the intersection is convex and has empty interior. The point $A$ is an extreme point of $H$, hence an endpoint of any nontrivial intersection segment. Away from $A$ the segment lies on the boundaries of both pockets: for two nondegenerate pockets this follows from their disjoint interiors; if one is a hull chord it follows from the supporting line of $H$. Such a segment cannot reach the opposite hull endpoint of either pocket, since this would give a triple intersection. If the direction of the segment is along the complementary arc at both pockets, the arcs share exactly a common initial straight segment seen from $A$. If the intersection direction is along the base of a positive-area pocket, its complementary arc has no other point in the relative interior of that base and the complementary arcs meet only at $A$. These alternatives also apply with a degenerate pocket, whose complementary arc is its base. It follows that the two boundary arcs meet in a common initial straight segment (possibly a single vertex).

At each hull vertex cancel this common initial segment from the two adjacent arcs, retaining its further endpoint as their new common end. The three cancellation segments are pairwise disjoint since a triple intersection is impossible. In each arc there remains a connected nonconstant polygonal arc: the cancelled closed portions at the two endpoints are disjoint. Their three remainders meet only at successive endpoints and give a simple Jordan polygon $D$.

Orient $H$ counterclockwise and orient the three arcs accordingly. Away from all limiting boundary segments, the winding number of their concatenation equals the winding number of $H$ minus the sum of those of the three positively oriented pockets, by cancelling their hull-side segments. It equals one on the part of $H$ outside the pockets and zero outside that part, since their interiors are disjoint. Erasing opposite traversals of a segment does not change winding numbers. Consequently $D$ is that part up to a set of area zero, and $|D|=A$. In this orientation every bend away from the three junctions is reflex. A simple polygon has at least three convex vertices, so these junctions are exactly the convex vertices of $D$. Its closure lies in $H_0$. Its side count is at most $N$, since each remaining bend other than a junction comes from an interior vertex of an original chain.

Finally extend the mass-one minimizers by zero to a larger common ball. Since the minimizing sequence has bounded energy, extract a weak $H^1_0$ limit with strong $L^q$ convergence. Its limit vanishes a.e. outside $D$: outside $H$ and in each positive-area pocket's interior every point has a neighbourhood eventually disjoint from $P_j$, and the remaining segments have area zero. The bounded Lipschitz Jordan polygon $D$ has the usual zero-extension characterization for $H^1_0(D)$, so this limit is a unit-mass admissible function there. Lower semicontinuity of the Dirichlet norm proves attainment.
\end{proof}

\begin{theorem}[Three convex vertices]\label{C:thm:three-convex-turns-sharp}
Every simple polygon $P$ with precisely three convex vertices satisfies
\[
 |P|^{2/q}\lambda_q^{\rm S}(P)\ge |R_3|^{2/q}\lambda_q^{\rm S}(R_3)
 \qquad(1\le q\le2),
\]
and equality forces an equilateral triangle.
\end{theorem}
\begin{proof}
Fix the area of $P$, take its closed convex hull as $H_0$, and take its total number of sides as $N$ in Lemma~\ref{C:lem:three-pockets-compact}. If the minimizer had all three hull sides dented, Lemmas~\ref{C:lem:three-reflected-cap} and~\ref{C:lem:good-facet-descent} would give a strictly better member of its class. Thus a minimizer has an entire hull side in its boundary. Theorem~\ref{C:thm:three-turn-straight-hull-side} now gives the desired lower bound. If $P$ attains equality, a strictly improving area-preserving perturbation of it would contradict the now proved bound, so $P$ itself has a straight hull side and the equality assertion in that theorem applies.
\end{proof}

The same comparison holds for curved boundaries.
\begin{corollary}\label{C:cor:three-arcs}
Let $\Omega$ be a bounded Jordan domain whose positively oriented boundary consists of finitely many regular $C^2$ arcs with nonpositive signed curvature. Suppose that the one-sided tangents at a junction are never opposite and that at most three junctions have positive exterior turn. Then $\Phi_q(\Omega)\ge\Phi_q(R_3)$ for every $1\le q\le2$.
\end{corollary}
\begin{proof}
Subdivide each arc finely, keeping all junctions, and replace each short piece by its chord as in Section~\ref{C:sec:arc-proof}. In local coordinates the domain lies above a concave graph, so each chord lies outside $\Omega$. The resulting polygons contain $\Omega$, their areas converge to $|\Omega|$, and their only positive turns are at the original junctions with positive turn. A simple polygon has at least three convex vertices, so these polygons have exactly three. Domain monotonicity and Theorem~\ref{C:thm:three-convex-turns-sharp} give the comparison in the limit.
\end{proof}
This approximation does not decide the equality case for curved boundaries.

\section{Torsional rigidity and the conformal map}\label{C:sec:hankel}
The classical conformal representation of torsional rigidity goes back to
P\'olya and Szeg\H{o} \cite[p.~120]{ps1951}; see also \cite[Theorem~4]{dittmar2011}.
The proposition below expresses it in terms of the singular values of a
Hardy-space Hankel matrix. For the four-index coefficient sum see also
\cite[proof of Proposition~3.5]{blesa2022}. The last display of that proof
includes the factor $\pi/2$ in the convention $\rho=4T$.
A coefficient proof of Saint--Venant's inequality and its equality case
appears in \cite[Theorem~3.6]{blesa2022}. Janson, Upmeier and Wallst\'en give a
boundary integral for the fourth Schatten norm of a Hardy-space Hankel
operator \cite{juw1994}. On Bergman space Olsen and Reguera derive
Saint--Venant's inequality from a sharp operator norm bound
\cite{olsenreguera}, and Fleeman and Lundberg identify the squared
Bergman analytic content with $4T$ in our normalization
\cite[Theorem~1.2]{fleemanlundberg}. We use the second-order minors of
the matrix to estimate the deficit from the exterior turns of a polygon.
The corresponding sum of squares appears in \cite[Lemma~3.1]{wangv1}.

\begin{proposition}[Hankel formula for torsional rigidity]\label{C:prop:torsion-hankel}
Let $\Omega$ be a simply connected planar domain of finite area $A$.
For any conformal map $F(z)=a_0+\sum_{j\ge1}a_jz^j$ from the unit
disk onto $\Omega$, let $\sigma_j$ be the singular values of the
Hankel matrix $C=(a_{i+k+1})_{i,k\ge0}$ on $\ell^2$. Then
\begin{equation}\label{C:eq:torsion-hankel}
 \frac A\pi=\sum_j\sigma_j^2,\qquad
 \frac{8T(\Omega)}\pi=\sum_j\sigma_j^4 .
\end{equation}
In particular, at area $\pi$ the torsion deficit satisfies
\[
 D(\Omega):=1-\frac{8T(\Omega)}{\pi}
 =1-\sum_j\sigma_j^4=2\sum_{i<j}\sigma_i^2\sigma_j^2.
\]
\end{proposition}
\begin{proof}
The area formula gives $A/\pi=\sum_{j\ge1}j|a_j|^2
=\sum_{i,k\ge0}|C_{ik}|^2$, so $C$ is Hilbert--Schmidt.
When $F$ is analytic across the unit circle and maps it to a Jordan curve, expand the conformal representation of torsion
in the normalization $-\Delta u=1$. The resulting coefficient identity is

\[
 \frac{8T(\Omega)}\pi=\sum_{p+q=r+s}
 \min(p,q,r,s)\,a_pa_q\overline{a_ra_s}\qquad(p,q,r,s\ge1).
\]
For the classical conformal representation see \cite[p.~120]{ps1951}.
Lemma~\ref{C:lem:sos} derives the displayed coefficient identity directly;
see also \cite[proof of Proposition~3.5]{blesa2022}. The minimum counts the ways of writing
$(p,q,r,s)=(i+k+1,j+l+1,i+l+1,j+k+1)$ with nonnegative indices. One may choose $i$ from $\max(0,r-q)$ through $\min(p,r)-1$.
Thus the right side is
$\sum_{i,j}|\sum_k C_{ik}\overline{C_{jk}}|^2
=\|CC^*\|_{\rm HS}^2=\sum_j\sigma_j^4$.
The coefficient sums are absolutely convergent whenever the area is finite:
replace $C$ entrywise by $|C|$ and bound the corresponding fourth-order
sum by $\|\,|C|\,\|_{\rm HS}^4$. For a general $\Omega$, apply this identity to $F(tz)$ for $0<t<1$.
As $t\uparrow1$, its Hankel matrix tends to $C$ in Hilbert--Schmidt norm
by the area formula, hence also in the fourth Schatten norm, and the
torsion of $F(t\mathbb D)$ increases to $T(\Omega)$ by exhaustion in
the variational principle: every compactly supported smooth test lies in
one of these increasing domains, and such tests are dense in $H^1_0(\Omega)$.
Squaring the area formula now gives the deficit.
\end{proof}

Consequently $8\pi T(\Omega)\le A^2$, the Saint--Venant inequality
for simply connected domains. Equality holds precisely when $C$ has
rank one. Indeed, $a_1=F'(0)\ne0$, so the first two columns must then
satisfy $a_{k+2}=q a_{k+1}$ for every $k\ge0$, where $q=a_2/a_1$.
Finite area forces $|q|<1$. Thus $F(z)=a_0+a_1z/(1-qz)$ maps the unit
disk onto a disk. Conversely a disk attains equality, as its centered
map is linear.
\section{The gap for fewer convex vertices}\label{C:sec:reproduce}
Here the number of sides may be arbitrarily large. A bound on positive
exterior turns supplies the moment inequality used to compare with the
regular polygon.

Let \(j=j_{0,1}\) denote the first positive zero of \(J_0\).

\begin{theorem}[A bound in terms of convex vertices]\label{C:thm:convex-vertices}
Let \(n\ge5\) and let \(P\) be a simple polygon of area \(A\) with at
most \(n-1\) convex vertices. For the regular \(n\)-gon \(R_n\) of area \(A\),
\[
 \frac{T(R_n)-T(P)}{A^2}>\frac1{8000\pi n^4},\qquad
 A\bigl(\lambda_1(P)-\lambda_1(R_n)\bigr)>\frac{\pi j^2}{2000n^4}.
\]
\end{theorem}

The comparison extends to domains with the same bound on convex
corners whose boundary arcs bend inward (Corollary~\ref{C:cor:concave-arcs}).
Section~\ref{C:sec:signed-proof} gives a sharp inequality for the Fourier coefficients of a signed turning measure (Theorem~\ref{C:thm:signed-moments}). The proof of the quantitative Theorem~\ref{C:thm:convex-vertices} uses the bound in Theorem~\ref{C:thm:turn-gap}. The quantitative gap does not require the sharp signed moment theorem.

Normalize the area to $\pi$, and define $D(P)=1-8T(P)/\pi$. The Kohler--Jobin inequality
\cite[Eq.~(1.1) and Theorem~1.1]{brasco2014} states that
\(\lambda_1(\Omega)T(\Omega)^{1/2}\ge\lambda_1(\mathbb D)T(\mathbb D)^{1/2}
=j^2(\pi/8)^{1/2}\) for every open set \(\Omega\) of finite area. For
a domain of area \(\pi\) we have \(T=(\pi/8)(1-D)\), so the inequality
reads
\begin{equation}\label{C:eq:kj-uniform}
 \frac{\lambda_1(P)}{j^2}\ge (1-D(P))^{-1/2}.
\end{equation}
The torsional rigidity in \cite{brasco2014} is
\(\max_v(\int|v|)^2/\int|\nabla v|^2\). It agrees with ours by
optimizing the scalar multiple of a test function in
\eqref{C:eq:torsion-variation}.

Take the disk's first eigenfunction as a trial function on the
regular polygon by transporting it with the centered disk map.
Section~\ref{C:sec:regular-radial} computes an explicit upper bound
\(\lambda_1(R_n)/j^2<\rho_n\) from this test.
The function \(D\mapsto(1-D)^{-1/2}\) reaches \(\rho_n\) when
\(D=1-\rho_n^{-2}\). Thus a single bound on torsion can give both
comparisons, provided the deficit crosses that threshold. We prove
it crosses by more than \(1/(1000n^4)\). Set
\[
                   d_n=1-\rho_n^{-2}+\frac{1}{1000n^4}.
\]

\begin{proposition}\label{C:prop:nonconvex}
Let \(n\ge5\) and let \(P\) be a simple polygon of area \(\pi\).
If \(P\) has at most \(n-1\) convex vertices, then \(D(P)>d_n\). For \(n=5\) one has the stronger bound \(D(P)\ge1/10\).
In particular,
\[
 D(P)-D(R_n)>\frac{1}{1000n^4},\qquad
 \lambda_1(P)-\lambda_1(R_n)>\frac{j^2}{2000n^4}.
\]
For $n=5$ the same comparisons give the area-invariant ratios
$T(P)/T(R_5)<54/55$ and $\lambda_1(P)/\lambda_1(R_5)>\sqrt{55/54}$.
\end{proposition}

\begin{proof}

Apply Theorem~\ref{C:thm:turn-gap} with $m=n-1$. Its threshold $d_{m+1}$
is defined by the same disk test function as $d_n$ above. This gives $D(P)>d_n$
for every $n\ge5$. The four-vertex bound is Theorem~\ref{C:thm:four-turn-sweep}.

To read off the gaps, write \(\Psi(D)=(1-D)^{-1/2}\) and
\(d_*=1-\rho_n^{-2}\), so \(\Psi(d_*)=\rho_n\). Applying
\eqref{C:eq:kj-uniform} and the trial bound to \(R_n\) gives
\(\Psi(D(R_n))\le\lambda_1(R_n)/j^2<\Psi(d_*)\), hence \(D(R_n)<d_*\).
For any area-\(\pi\) domain with \(D\ge d_n\), the difference of
deficits is therefore greater than \(d_n-d_*=1/(1000n^4)\).
Also \(d_*\ge0\) by Faber--Krahn and \(\Psi'(D)\ge1/2\) on
\([0,1)\). Applying \eqref{C:eq:kj-uniform} to this domain gives
\[
 \frac{\lambda_1}{j^2}\ge\Psi(d_n)
 \ge\rho_n+\tfrac12(d_n-d_*)
 >\frac{\lambda_1(R_n)}{j^2}+\frac1{2000n^4}.
\]
Indeed $d_5<1/12$ by the explicit bound in Proposition~\ref{C:prop:four-axis-turns},
so $1-D(R_5)>\rho_5^{-2}>11/12$, while $1-D(P)\le9/10$.
The torsion ratio follows immediately, and \eqref{C:eq:kj-uniform}
gives the eigenvalue ratio strictly larger than $(10/9)^{1/2}/\rho_5$.
\end{proof}

\begin{proof}[Proof of Theorem~\ref{C:thm:convex-vertices}]
Both gaps are Proposition~\ref{C:prop:nonconvex} rescaled from area
\(\pi\) to area \(A\): torsion scales quadratically in area, while
the eigenvalue scales inversely in area.
\end{proof}

The comparison of exponents in Part~\ref{part:intermediate} extends the two gaps to every exponent.
\begin{corollary}\label{C:cor:all-q-gap}
Let $n\ge5$ and let $P$ be a simple polygon with at most $n-1$ convex vertices. Then $\Phi_q(P)>\bigl(1+1/(1000qn^4)\bigr)\Phi_q(R_n)$ for every $1\le q\le2$.
\end{corollary}
\begin{proof}
Scale to area $\pi$, let $\mathbb D$ be the unit disk and put $r_q(\Omega)=\lambda_q^{\rm S}(\Omega)/\lambda_q^{\rm S}(\mathbb D)$. By Proposition~\ref{S:prop:exponent-round}, $q\mapsto r_q(\Omega)^q$ increases on $[1,2]$ for every polygon $\Omega$. At $q=1$, $r_1(P)=T(\mathbb D)/T(P)=(1-D(P))^{-1}$, and Proposition~\ref{C:prop:nonconvex} gives $D(P)>d_n>1-\rho_n^{-2}$, so $r_1(P)>\rho_n^2$. At $q=2$, $r_2(R_n)=\lambda_1(R_n)/j^2<\rho_n$. Therefore
\[
 r_q(P)^q\ge r_1(P)>\rho_n^2>r_2(R_n)^2\ge r_q(R_n)^q .
\]
For the explicit factor, Proposition~\ref{C:prop:nonconvex} gives $D(P)>d_n=1-\rho_n^{-2}+1/(1000n^4)$, so $r_1(P)>\rho_n^2/(1-y)$ with $y=\rho_n^2/(1000n^4)\ge1/(1000n^4)$, since $\rho_n\ge1$ by the Faber--Krahn inequality. Hence $(\Phi_q(P)/\Phi_q(R_n))^q=r_q(P)^q/r_q(R_n)^q>(1-y)^{-1}$, and convexity of $y\mapsto(1-y)^{-1/q}$ gives $\Phi_q(P)/\Phi_q(R_n)>1+y/q$.\qedhere
\end{proof}

A nonconvex polygon with at most $n$ sides has at most $n-1$ convex vertices. Thus for $n\ge5$ Theorem~\ref{C:thm:convex-vertices} gives a second proof of the torsion and eigenvalue cases of Theorem~\ref{thm:main} for nonconvex polygons, with explicit gaps.

Approximation by polygons gives a version for boundaries with inward-curving arcs. Its proof is in Section~\ref{C:sec:arc-proof}.
\begin{corollary}\label{C:cor:concave-arcs}
Let \(\Omega\) be a bounded Jordan domain of area \(A\), whose positively
oriented boundary consists of finitely many regular \(C^2\) arcs with
nonpositive signed curvature. Assume the one-sided tangents at the
junctions are never opposite and at most \(n-1\) junctions have
positive exterior turn, where \(n\ge5\).
For the regular \(n\)-gon \(R_n\) of area \(A\),
\[
 \frac{T(R_n)-T(\Omega)}{A^2}>\frac1{8000\pi n^4},
 \qquad A\bigl(\lambda_1(\Omega)-\lambda_1(R_n)\bigr)>
 \frac{\pi j^2}{2000n^4}.
\]
Moreover $\Phi_q(\Omega)>\Phi_q(R_n)$ for every $1\le q\le2$.
\end{corollary}

\section{Turning moments and the quantitative gap}\label{C:sec:certification}
This section proves Theorem~\ref{C:thm:turn-gap}, the bound used in Proposition~\ref{C:prop:nonconvex}. Section~\ref{C:sec:arc-proof} extends the gap to boundaries with inward-curving arcs, and Section~\ref{C:sec:signed-proof} proves a sharp form of the moment inequality behind it. Finite decimals used in inequalities denote exact rationals. Write $j=j_{0,1}$.
\subsection{Exterior angles and the torsion deficit}
\label{C:sec:conformal}

Scale \(P\) to the area \(\pi\) of the unit disk \(\mathbb D\).
The disk has torsion function \((1-|z|^2)/4\) and torsional rigidity
\(\pi/8\). Thus \(D(P)=1-8T(P)/\pi\) is the fraction of the disk's
torsion lost on passing to \(P\). We want a lower bound for this
loss in terms of the number of convex vertices. The disk map lets us
express both quantities in the same coordinates. We use the
variational formula
\begin{equation}\label{C:eq:torsion-variation}
 T(P)=\max_{v\in H^1_0(P)}
       \left(2\int_P v-\int_P|\nabla v|^2\right).
\end{equation}
For conformal formulas for torsion, see \cite{dittmar2011,fleeman2019}.
We are free to choose which point of \(P\) corresponds to the center
of the disk. We choose a point of maximum conformal radius so that
the quadratic coefficient of the map vanishes. In detail, take a
Riemann map \(F:\mathbb D\to P\) with this choice of center. By the Koebe quarter theorem the conformal radius tends to zero at the boundary, so this maximum occurs inside \(P\). Translate and rotate to arrange \(F(0)=0\) and \(F'(0)=a>0\). The conformal radius at \(F(z)\) is \((1-|z|^2)|F'(z)|\). Its first derivatives at \(z=0\) vanish, which gives \(F''(0)=0\). Write
\[
 F(z)=\sum_{j\ge1}a_jz^j,\quad a_1=a,\quad a_2=0,\qquad
 b=a^2,\quad B=1-b,\quad
 h=F'/a=1+\sum_{k\ge2}g_kz^k.
\]
The area formula gives \(\sum_{j\ge1}j|a_j|^2=1\). Set
\begin{equation}
 S=\frac{B}{b}=\sum_{k\ge2}\frac{|g_k|^2}{k+1},\qquad
 Q=\sum_{k\ge2}\frac{(k-1)|g_k|^2}{(k+1)^2}.
 \label{C:eq:SQ}
\end{equation}

For the disk, \(F(z)=z\), all \(g_k\) vanish, and \(S=Q=0\).
For a general polygon, \(S\) measures the total size of the higher
coefficients of its normalized derivative. The factor
\((k-1)/(k+1)\) in \(Q\) gives relatively more weight to higher
orders. We first express \(D\) as a sum of nonnegative terms, then
keep the terms needed for a lower bound.

\begin{lemma}[Coefficient form of the torsion deficit]
\label{C:lem:sos}
Let \(A_{p,s}=a_pa_{s-p}\) and \(K_{p,q,s}=\min\{p,q,s-p,s-q\}\). Then
\begin{equation}
 D(P)=\frac12\sum_{s\ge2}\sum_{p,q=1}^{s-1}
       K_{p,q,s}|A_{p,s}-A_{q,s}|^2.
 \label{C:eq:sos}
\end{equation}
Equivalently, for the Hankel matrix $C_{ik}=a_{i+k+1}$, $i,k\ge0$,
\begin{equation}\label{C:eq:hankel-minors}
 D(P)=2\sum_{i<j,\ k<l}|C_{ik}C_{jl}-C_{il}C_{jk}|^2 .
\end{equation}
\end{lemma}

\begin{proof}
For completeness, the coefficient formula in
Proposition~\ref{C:prop:torsion-hankel} follows directly from the disk
torsion problem. First use $F(tz)$ for $t<1$, so the boundary is
analytic, and translate the image to make the constant coefficient
zero. Write the resulting map as $F$ and its coefficients as $a_p$.
The function $4u_P\circ F+|F|^2$ is the harmonic extension of
$|F|^2$ from the circle. Thus the $(p,q)$ term of $4u_P\circ F$
is $a_p\overline{a_q}(r^{|p-q|}-r^{p+q})e^{i(p-q)\theta}$.
Integration against $|F'|^2\,dA$ pairs it with derivative indices
$i=s-p$ and $j=s-q$. The radial factor in $4T(P)/\pi$ is
\[
 2ij\int_0^1 (r^{|p-q|}-r^{p+q})r^{i+j-1}\,dr
 =\min(p,q)-pq/s=:G_{pq}.
\]
As $A_{s-q,s}=A_{q,s}$ and $G_{pq}+G_{p,s-q}=K_{p,q,s}$, this
gives $8T(P)/\pi=\sum_{s,p,q}K_{p,q,s}A_{p,s}\overline{A_{q,s}}$.
Pass to $t\uparrow1$ as in Proposition~\ref{C:prop:torsion-hankel}.

At area $\pi$, Lagrange's identity gives
\[
 D(P)=\|C\|_{\rm HS}^4-\|CC^*\|_{\rm HS}^2
 =\tfrac12\sum_{i,j,k,l\ge0}|C_{ik}C_{jl}-C_{il}C_{jk}|^2.
\]
Ordering the two row and column indices gives \eqref{C:eq:hankel-minors}.
Grouping by $s=i+j+k+l+2$, $p=i+k+1$ and $q=i+l+1$ gives
\eqref{C:eq:sos} by the same multiplicity count as in the proposition.
\end{proof}

\begin{lemma}[Bounds on the conformal coefficients]
\label{C:lem:coercivity}
Put $z_k=a_{k+1}/a$, $Y=\sum_{k\ge2}|z_k|^2=(S-Q)/2$ and
$Q^\sharp=\sum_{k\ge4}(k-3)|z_k|^2$. For the Hankel matrix
$E=C/a-e_0e_0^*$ let $\mathcal D(E)=\|E\|_{\rm HS}^4-\|EE^*\|_{\rm HS}^2$.
Then
\begin{equation}\label{C:eq:homogeneous-hankel}
 \frac{D}{b^2}=2Q-4\operatorname{Re}\sum_{i,j\ge2}z_{i+j}\overline{z_i z_j}
     +\mathcal D(E),\qquad 6Y^2\le\mathcal D(E)\le S^2 .
\end{equation}
With $V=2Y$,
\begin{align}
 D/b^2&\ge 2Q+\tfrac32 V^2-2V\sqrt{Q^\sharp},
 \label{C:eq:pooled-deficit}\\
 D&\ge B^2/4,\qquad D\ge 2bB/3 .
 \label{C:eq:B-deficit}
\end{align}
More generally, whenever $F(z)=a_1z+O(z^{d+1})$ with $d\ge2$ and
$|P|=\pi$, one has
\begin{equation}\label{C:eq:initial-mode-deficit}
 D\ge\frac{2(d-1)}{d+1}b(1-b),\qquad b=|a_1|^2 .
\end{equation}

\end{lemma}

\begin{proof}
\emph{The bound $D\ge B^2/4$.}
Let $M=\max u_P$ and $\mu(t)=|\{u_P>t\}|$. At almost every
$t\in(0,M)$ the flux identity gives
$\int_{u_P=t}|\nabla u_P|=\mu(t)$. The isoperimetric inequality,
Cauchy--Schwarz on the level curve, and coarea therefore give
\[
 4\pi\mu(t)\le {\rm Per}(\{u_P>t\})^2\le \mu(t)(-\mu'(t)).
\]
As $\mu(0)=\pi$, integration yields $\mu(t)\le\pi(1-4t)$, hence
$T(P)\le\pi(M-2M^2)$ and $D\ge(1-4M)^2$.

Recenter a disk map $f$ at a maximum point of $u_P$, translate so
$f(0)=0$, and write $f(z)=\sum c_jz^j$. The harmonic function
$4u_P\circ f+|f|^2$ has boundary values $|f|^2$, since $P$ is a
Jordan polygon and the torsion function vanishes continuously on its
boundary. The mean-value property and Parseval give
\[
 4M=\sum_j|c_j|^2\le\tfrac12\bigl(|c_1|^2+\sum_j j|c_j|^2\bigr)
 \le(1+b)/2 .
\]
The last step uses the maximality of the conformal radius $\sqrt b$
and the area normalization. Thus $1-4M\ge B/2$, proving the bound.
\emph{The identity \eqref{C:eq:homogeneous-hankel}.}
Adding $e_0e_0^*$ to $E$ changes only the ordered minors with rows
$(0,i)$ and columns $(0,j)$, $i,j\ge1$. The old minor is $-z_i z_j$
and its increment is $z_{i+j}$. Its increments have total squared norm
$\sum_{k\ge2}(k-1)|z_k|^2=Q$. Expanding \eqref{C:eq:hankel-minors}
gives the identity. The disjoint minor families of $E$ with smallest row and column
indices $(0,0),(0,1),(1,0)$ each contribute $2Y^2$.
For example, the first family has minors $-z_i z_j$, $i,j\ge1$,
with $z_1=0$; the next has $-z_lz_{j+1}$, $l\ge2$, $j\ge1$,
and the third is its transpose. Thus no minor is counted twice; the upper bound follows from $\|E\|_{\rm HS}^2=S$.
On $i,j\ge2$ there are $k-3$ pairs with $i+j=k$, so Cauchy--Schwarz
bounds the cubic sum by $Y\sqrt{Q^\sharp}$, proving \eqref{C:eq:pooled-deficit}.
\emph{The bound \eqref{C:eq:initial-mode-deficit}.}
For $k\ge d$ the antidiagonal $s=k+2$ in \eqref{C:eq:sos} has two
end entries $E=a_1a_{k+1}$. All other entries vanish for $k<2d$.
For $k\ge2d$ its $2(d-1)$ forced zeros are at positions $2,\ldots,d$
and their reflections. Each of the $k-2d+1$ inner entries $w$
pairs with those zeros with total weight $J=2\sum_{j=2}^d j$ and with
the two ends with weight one each. Thus it contributes at least
$[2J/(2+J)]|E|^2$ by completing the square.
Including the end--zero pairs, the coefficient of $|E|^2$ is at least
\[
 2(k-1)-\frac{4(k-2d+1)_+}{d(d+1)}
 \ge \frac{2(d-1)}{d+1}(k+1).
\]
Summing proves \eqref{C:eq:initial-mode-deficit} and its case $d=2$
gives the second bound in \eqref{C:eq:B-deficit}.
\end{proof}

\subsection{Fourier moments of the turning measure}
\label{C:sec:finite}

Let \(\xi_i\in\partial\mathbb D\) be the boundary point that the
disk map sends to the \(i\)-th polygon vertex. These are the
Schwarz--Christoffel prevertices \cite[Ch.~4]{dt2002}. At a vertex of interior angle
\(\theta_i\) put \(\nu_i=(\pi-\theta_i)/(2\pi)\). Write \(\nu=\sum_i\nu_i\delta_{\xi_i}\)
and \(\widehat\nu(k)=\sum_i\nu_i\overline{\xi_i}^{\,k}\). The turning formula is
\begin{equation}
 \sum_i\nu_i=1,\qquad
 h(z)=\prod_i(1-\overline{\xi_i}z)^{-2\nu_i},
 \qquad
 z(\log h)'=2\sum_{k\ge1}\widehat\nu(k)z^k .
 \label{C:eq:sc}
\end{equation}
In particular, \(\widehat\nu(1)=0\). Positive weights correspond to
convex vertices. For the regular \(n\)-gon we may take the
prevertices to be the \(n\)-th roots of unity. They give the product
\(\prod_j(1-\overline{\xi_j}z)=1-z^n\). The turning formula
\eqref{C:eq:sc} therefore gives the normalized derivative $h_R=F_R'/F_R'(0)$ of the centered disk map $F_R$ as
\begin{equation}\label{C:eq:regular-series}
 h_R(z)=(1-z^n)^{-2/n}=1+\frac2n z^n+O(z^{2n}).
\end{equation}
For general polygons, the link between these turning data and
torsion also involves the normalized derivative \(h\).

\subsection{Counting positive turns on the disk}
\label{C:sec:shifted-hardy}

The generating function \(H\) below has the turning moments as its
Taylor coefficients. Multiplication by \(h\) brings back the
coefficients \(g_k\) that occur in the torsion bounds. For
\(s>0\), set \(r=e^{-s/2}\) and write
\[
 H(z)=\sum_i\frac{\nu_i}{1-\overline{\xi_i}z}-1,\qquad
 h_s(z)=h(rz),\qquad H_s(z)=H(rz).
\]
Use normalized circle measure and put
\(G(s)=\|h_s-1\|_2^2=\sum_{k\ge2}|g_k|^2e^{-ks}\).
Since \(2h_sH_s=zh_s'\) by \eqref{C:eq:sc}, Parseval gives
\begin{equation}\label{C:eq:shifted-parseval}
 \|h_s(H_s-c)\|_2^2
 =c^2+\sum_{k\ge2}(k/2-c)^2 |g_k|^2e^{-ks}
 =\tfrac14 G''(s)+cG'(s)+c^2(1+G(s)).
\end{equation}

\begin{lemma}\label{C:lem:shifted-hardy}
If the turning measure has at most \(m\) positive weights, then for
every \(s>0\) and \(c\ge0\),
\begin{equation}\label{C:eq:shifted-jensen}
 \|h_s(H_s-c)\|_2^2\ge(1+c)^2e^{-ms}.
\end{equation}
In particular, with \(M_m(s)=(e^{ms}-1)^{-1}\),
\begin{equation}\label{C:eq:hardy-differential}
 \tfrac14G''(s)\ge M_m(s)(1-G'(s))-M_m(s)^2G(s) .
\end{equation}
\end{lemma}

\begin{proof}
Omit zero weights and let $N_+$ denote the number of positive atoms.
For real $c\ge0$, cancel the poles in $H-c$ by multiplying by
$\prod_i(1-z/\xi_i)$. The resulting polynomial never vanishes on the unit
circle. At an atom $\xi_i$ its value is
$\nu_i\prod_{j\ne i}(1-\xi_i/\xi_j)\ne0$, and away from
the atoms $\operatorname{Re}H=-1/2$, since each
$1/(1-z/\xi_i)$ has real part $1/2$ there. Its degree remains $N$, since $H(\infty)-c=-(1+c)\ne0$.
The argument principle therefore shows that the number of roots in the disk,
with multiplicity, is independent of $c$ on $[0,\infty)$.

To find this number, put $\varepsilon=(1+c)^{-1}$ and let it tend to
zero. Near $\xi_i$ the root equation is
\[
 (1-z/\xi_i)\left(1-\varepsilon\sum_{j\ne i}
                  \frac{\nu_j}{1-z/\xi_j}\right)=\varepsilon\nu_i .
\]
At $\varepsilon=0$ its root $\xi_i$ is simple, so the implicit function
theorem gives $z_i=\xi_i(1-\varepsilon\nu_i+O(\varepsilon^2))$.
These are all the roots, one per atom, and precisely the positive
atoms give roots inside the disk. Hence the count is $N_+$ at every $c\ge0$.
The numerator over $\prod_i(z-\xi_i)$ has leading coefficient $-(1+c)$
because $H(\infty)=-1$.
For related zero counts for polygonal maps see \cite[Theorem~2]{cp2014}
and \cite[Theorem~4(a)]{cdo2012}.

Now apply the circle mean identity
\(\int_{\mathbb T}\log|r\zeta-z|=\log\max(r,|z|)\) to each root.
The denominator has logarithmic mean zero. This identity remains valid
when a root lies on the circle of radius $r$, since its logarithmic
singularity is integrable. Each of the \(N_+\) interior
roots contributes at least \(\log r\) and the other terms are
nonnegative. Thus
\begin{equation}\label{C:eq:log-mean}
 \int_{\mathbb T}\log|H(r\zeta)-c|\ge\log(1+c)+N_+\log r .
\end{equation}
Since \(h\) is zero-free with \(h(0)=1\), its logarithmic mean is
zero. Exponentiating and using Jensen's inequality gives
\(\|h_s(H_s-c)\|_2^2\ge(1+c)^2r^{2N_+}\ge(1+c)^2e^{-ms}\).
Taking \(c=M_m(s)\), so that \((1+c)^2e^{-ms}=c(1+c)\), gives
\eqref{C:eq:hardy-differential} by \eqref{C:eq:shifted-parseval}.
\end{proof}

\subsection{Combining bounds at different radii}\label{C:sec:split-moment}

\begin{lemma}[Integrated turning estimate]\label{C:lem:integrated}
Suppose the normalized disk map has at most $m$ positive turns, $m\ge4$.
Let $0\le\sigma\le4/3$ and let $w$ be measurable with $0\le w(s)\le s^2e^{-\sigma s}/2$.
Put $p=4/m$, $S_k=|g_k|^2/(k+1)$ and
\[
 P_m(s)=e^{-ms}(1-e^{-ms})^{-1-p},\qquad
 \ell_m=\frac{2}{m(m-2)},\qquad
 \delta_\sigma=m+2-\left(2-\frac{3\sigma}{2}\right)\frac{2-p}{3-p}.
\]
Then
\begin{equation}\label{C:eq:shifted-moment}
 4\int_0^\infty w(s)P_m(s)\,ds
 \le\sum_{k\ge2}\left[k^2(k+1)\int_0^\infty w(s)e^{-ks}\,ds
              +\ell_m\frac{k+1}{k+\delta_\sigma}\right]S_k .
\end{equation}
\end{lemma}
\begin{proof}
Write $M=M_m$, $\phi=(1-e^{-ms})^p$, and $f=1-G'$.
Since $\phi'=4M\phi$, \eqref{C:eq:hardy-differential} gives
$(\phi f)'\le4\phi M^2G$. Integrating to infinity, where $\phi f\to1$,
yields $f(s)\ge[1-4\int_s^\infty\phi M^2G]/\phi(s)$.
Insert this bound in $G''/4\ge Mf-M^2G$, multiply by $4w$ and interchange
the order of integration in the resulting double integral, whose integrand is nonnegative. As $M/\phi=P_m$, the result is
\[
 4\int_0^\infty wP_m\le\int_0^\infty wG''
       +\int_0^\infty G(y)M(y)^2J_w(y)\,dy,\qquad
 J_w(y)=4w(y)+4\phi(y)\int_0^y\frac{w(s)\phi'(s)}{\phi(s)^2}\,ds .
\]

The kernel $J_w$ increases with $w$. Take $w_0(y)=y^2e^{-\sigma y}/2$.
By integration by parts, $J_{w_0}=4\phi\int_0^y w_0'/\phi$. Hence
$K=M^2J_{w_0}$ solves the scalar equation
\[
 K'+[2m+(2m-4)M]K=4M^2w_0'.
\]
Its limit at zero is $\ell_m$, by $\phi(y)\sim(my)^p$ and $w_0'(y)\sim y$.
Put $a=3\sigma/2$, $B=2m-4$ and $\theta=(2m-\delta_\sigma)/B$.
Then $0\le\theta\le1$, $4/\ell_m=mB$ and
$\delta_\sigma-a=m(1/2+\theta)$. With $x=my$, weighted AM--GM gives
\[
 (e^x-1)(1-\theta+\theta e^x)\ge x e^{(1/2+\theta)x},
\]
since $2\sinh(x/2)\ge x$. Multiplication by $B e^{-\delta_\sigma y}M^2$
shows that $\ell_m e^{-\delta_\sigma y}$ is a supersolution: here
$w_0'(y)\le y e^{-ay}$. The integrating factor vanishes at zero like
$y^{2-p}$, so scalar comparison gives $M^2J_w\le K\le\ell_m e^{-\delta_\sigma y}$.
Expanding $G$ and $G''$ into their nonnegative modes now gives
\eqref{C:eq:shifted-moment}. All exchanges can first be made on
$[\varepsilon,R]$ and for finite modal sums, then passed to the limit by
Tonelli's theorem. The coefficient of $|g_k|^2$ in the final upper sum is
$O((k+1)^{-1})$, since $w\le s^2/2$; the area formula therefore
makes that sum finite. The turning integral is finite, since
$wP_m=O(s^{1-p})$ at zero and decays exponentially at infinity.
\end{proof}
The comparison function $\ell_me^{-\delta_\sigma y}$ agrees with $K$ to first order at zero: $K(y)=\ell_m(1-\delta_\sigma y+O(y^2))$.
Indeed the exponential comparison matches through first order after division by $x$,
and $w_0'/(ye^{-ay})=1+O(y^2)$. The residual in the differential equation
is therefore $O(y)$; its integrating factor gives $K-\ell_m e^{-\delta_\sigma y}=O(y^2)$.
The same applies to $w(s)=e^{-s}(1-\cos(\sqrt6\,s))/6$ at $\sigma=1$,
since $w-w_0=O(s^4)$ changes $M^2J_w$ by $O(s^2)$.

The binomial series gives
\[
 P_m(s)=\sum_{j\ge1}c_j e^{-mjs},\qquad
 c_1=1,\qquad c_{j+1}=(1+4/(mj))c_j .
\]
Its terms are nonnegative, which allows termwise integration.

\subsection{Testing the regular polygon with disk minimizers}\label{C:sec:regular-radial}
\label{C:subsec:regular-tail}

\paragraph{The area factor of the regular polygon.}
For the centered regular map, \(h_R\) is given by \eqref{C:eq:regular-series}.
With \(s=1/n\), the map with derivative \(h_R\) and value zero at
the origin sends \(1\) to a vertex at radius \(s\,\mathrm B(s,1-2s)\).
Its image has area \(\pi Z(s)\), where the \(n\) radial triangles give
\[
 Z(s)=\frac{\sin(2\pi s)}{2\pi s}
       \bigl(s\,\mathrm B(s,1-2s)\bigr)^2 .
\]
The beta identity and reflection formula yield, for \(0<s<1/2\),
\begin{equation}\label{C:eq:regular-area}
 \begin{aligned}
 Z(s)&=\frac{\Gamma(1+s)^2\Gamma(1-2s)}
                {\Gamma(1-s)^2\Gamma(1+2s)},\\
 \log Z(s)&=\sum_{\substack{m\ge3\\m\ {\rm odd}}}
       \frac{(2^{m+1}-4)\zeta(m)}{m}s^m .
 \end{aligned}
\end{equation}
The series follows by expanding \(\log\Gamma(1\pm s)\) and
\(\log\Gamma(1\pm2s)\); it converges for \(|s|<1/2\).
Both \(\log Z\) and \(Z-1\) have nonnegative coefficients, beginning
at degree three.

For \(t=1/n\), Euler's product for the gamma function cancels the
linear exponential factors in \eqref{C:eq:regular-area} and gives
\begin{equation}\label{C:eq:regular-product}
 Z(t)=\prod_{j=1}^\infty\frac{(jn-1)^2(jn+2)}{(jn+1)^2(jn-2)}
     =\prod_{j=1}^\infty\left(1+\frac4{(jn+1)^2(jn-2)}\right).
\end{equation}
The second equality follows by subtracting the two cubic polynomials.

\paragraph{The eigenvalue bound.}
Keeping the exact area factor in the disk test function, put
\begin{equation}\label{C:eq:spectral-comparator-uniform}
\begin{gathered}
 L_n=\frac{1+3/(n+4)+3/((n+4)(n+5))}{2(n+1)(n+2)(n+3)},\\
 \rho_n=\frac{Z(t)}{1+23t^2 L_n},\qquad
 d_n=1-\rho_n^{-2}+t^4/1000 .
\end{gathered}
\end{equation}

The disk ground state \(\varphi(r)=J_0(jr)/(j\sqrt\pi J_1(j))\)
has unit Dirichlet energy and squared \(L^2\) norm \(j^{-2}\).
Set $g(x)=\sqrt\pi\,\varphi(\sqrt{x})$ for $0\le x\le1$. The radial
equation $(xg')'=-j^2g/4$ gives
\[
 -g'(x)=\frac{j^2}{4}\int_0^1 g(tx)\,dt .
\]
Starting with $g\ge0$ and differentiating this identity shows
inductively that $(-1)^r g^{(r)}\ge0$. At the boundary $g(1)=0$
and $g'(1)=-1/2$ by the normalization above, so the equation gives
$g''(1)=1/2$. Taylor's formula now yields
$g(1-v)\ge v/2+v^2/4$. Squaring and integrating against
$(1-v)^n\,dv$ gives
\[
 \int_{\mathbb D}|z|^{2n}\varphi^2\,dA\ge L_n
\]
by the beta integrals for \(v^2/4+v^3/4+v^4/16\).

In conformal coordinates the Dirichlet energy is unchanged.
In the mass of \(\varphi\), angular orthogonality and the first
nonconstant term \((2/n)z^n\) in \(h_R\) give the lower bound
\(Z(t)^{-1}[j^{-2}+(2/n)^2L_n]\).
Hence the Rayleigh quotient gives \(\lambda_1(R_n)/j^2<\rho_n\), since
\(j^2>23/4\) and \(L_n>0\).
For completeness, $J_0(x)>0$ on $0\le x\le12/5$. Indeed put $y=x^2/4\le36/25$. The alternating Bessel series after its linear term gives the lower bound $1-y+y^2/4-y^3/36+y^4/576-y^5/14400$. Its derivative is negative there, since $-1+y/2+y^3/144<0$, and its value at $36/25$ is positive by substitution. Thus $4j^2>4(12/5)^2>23$.

\subsection{The moment and deficit comparisons}\label{C:sec:uniform-modes}\label{C:sec:uniform-deficit}
For a simple polygon $P$ of area $\pi$, write $S_k=|g_k|^2/(k+1)$
and $d_k=(k-1)/(k+1)$, so $Q=\sum d_k S_k$ and $V=\sum(1-d_k)S_k$.
In the comparison with $R_n$ put $m=n-1$ and $t=1/n$, and assume
that $P$ has at most $m$ positive turns.
For four positive turns we choose the direction of two Steiner symmetrizations. For five or more positive turns we estimate the centered moments directly.

We first record how the bound from the disk test function depends on the side count. Its formula
\eqref{C:eq:spectral-comparator-uniform} makes sense for real $n\ge5$.
With $t=1/n$, let
\[
 R(t)=2L_n/t^3=(1+12t+38t^2)\Big/\prod_{i=1}^5(1+it).
\]
Then $R$ decreases from $1$. Indeed, the ratios of the numerator and
denominator coefficients, in ascending order, are
$1,12/15,38/85,0,0,0$. In the numerator of $R'$ the terms of
each pair of degrees $i<j$ therefore have nonpositive sum, strictly
negative for $i=0$, $j=1$.

\begin{theorem}[A torsion gap from four convex vertices]\label{C:thm:four-turn-sweep}
Every simple polygon $P$ with at most four convex vertices satisfies
\[
                 1-\frac{8\pi T(P)}{|P|^2}\ge\frac1{10}>\frac1{12}>d_5 .
\]
\end{theorem}
\begin{proof}
For a direction $e=(\cos\theta,\sin\theta)$, write $L(t)$ for the
total length of $P$ on the line $x\cdot e^\perp=t$. If $e$ is parallel
to no side, $L$ is continuous and piecewise affine. It is positive
between the extreme levels, by connectedness, and zero at the ends.
At a vertex whose incoming and outgoing edges have the same sign of
projection on $e^\perp$, a convex turn gives a downward jump of $L'$:
locally the upper graph is the minimum of two affine functions or the
lower graph is their maximum. A reflex turn reverses the sign.
When both neighbours lie on the same side of the vertex level, a convex
vertex adds a triangular slice of length a positive multiple of
$(t-t_0)_+$ or $(t_0-t)_+$, while a reflex vertex subtracts it.
Mark a vertex if its contribution to the jump is downward. Contributions
at a common level add.

We may first perturb the vertices to general position so that no diagonal
is parallel to a side and no three vertices are collinear. To count the
marks, orient the boundary counterclockwise and traverse it with a line
parallel to each side. At a vertex with signed exterior turn
$\gamma\in(-\pi,\pi)\setminus\{0\}$, rotate the unoriented line
clockwise to the next side. For a convex vertex this sweeps an interval
of length $\pi-\gamma$ in which the projections of the incoming and
outgoing edges have equal signs. For a reflex vertex the length is
$-\gamma$ and the signs are opposite. These are exactly the directions
that mark the vertex. The line turns clockwise through
$m\pi-\sum\gamma=(m-2)\pi$ on the full circuit, where $m$ is the number
of convex vertices. Each nonside direction therefore marks exactly
$m-2$ vertices.

A polygon with three convex vertices has just one mark in every such
direction. A polygon with four is a \emph{pseudo-quadrilateral}.
By the bitangent lemma of Rote, Santos and Streinu
\cite[Lemma~2.3 and Section~2.2]{rote-santos-streinu}, it has an interior
diagonal which splits it into two polygons with three convex vertices
each. Slice parallel to this diagonal. The two endpoints contribute
a total increase of $6-4=2$ in the number of convex vertices after the
cut. Each endpoint contributes at least one. A convex endpoint therefore
has its neighbours on opposite sides of the diagonal line. At a reflex
endpoint one of the split angles is still reflex, so both neighbours lie
on the same side of the line. Both endpoints are marked. Since there
are exactly two marks, all downward jumps occur on the diagonal's level.

\emph{Two perpendicular symmetrizations.}
In the direction just obtained, or in any nonside direction when $m=3$,
$L$ is convex on either side of a single level $t_0$. It increases
strictly to $t_0$ and then decreases strictly, since it is positive
inside its support and vanishes at the ends. The first Steiner image
is $|y|<L(t)/2$. At height $0<y<L(t_0)/2$ its left endpoint is the
inverse of an increasing convex function, hence concave in $y$.
Its right endpoint is the inverse of a decreasing convex function,
hence convex in $y$. Their difference $W(y)$ is decreasing, convex
and piecewise affine, and vanishes at $L(t_0)/2$. The second Steiner
image, $|t|<W(|y|)/2$, consequently has just four convex vertices,
one on each half-axis. Both symmetrizations preserve area and increase
torsion in the weak sense by the variational principle and the
Steiner rearrangement inequality \cite{ps1951}.
Proposition~\ref{C:prop:four-axis-turns} gives the asserted gap.

Finally any simple polygon with effective sides can be approximated by
ones in the stated general position, preserving simplicity and the signs
of all turns. A fixed triangulation extends the vertex perturbations
to piecewise affine homeomorphisms whose derivatives tend uniformly
to the identity on its triangles. Areas and torsions thus converge by
change of variables in \eqref{C:eq:torsion-variation}. This proves the
nonstrict bound for every polygon in the statement.
\end{proof}

\begin{proposition}[Four convex vertices on two reflection axes]\label{C:prop:four-axis-turns}
Let $P$ be a simple polygon invariant under reflection in two perpendicular
lines. Suppose its four convex vertices lie one on each of the four rays
of these lines. Then $D(P)>1/10>1/12>d_5$, irrespective of the number of
reflex vertices.
\end{proposition}
\begin{proof}
Normalize $|P|=\pi$ and take the axes as coordinate axes.
Their intersection belongs to $P$, since the half-turn induces a
nonidentity holomorphic involution of the disk and has an interior fixed
point. Normalize the disk map by $F(0)=0$ and $a_1=F'(0)>0$. Uniqueness
gives $F(-z)=-F(z)$ and $F(\bar z)=\overline{F(z)}$. Thus $h=F'/a_1$ has real even coefficients. Since $P$ is Jordan, $F$
extends homeomorphically to the boundary. Injectivity and the two
reflection identities show that the boundary points on the real axis
have preimages $1,-1$, and those on the imaginary axis have preimages
$i,-i$. Oddness puts one point on each half-axis. By hypothesis these
four points are the convex vertices.

Put $X_k=|a_{k+1}/a_1|^2$ and $S=\sum_{k\ge2}(k+1)X_k$,
so $a_1^2=(1+S)^{-1}$ and only even $k$ occur. On the antidiagonal
$s=k+2$ in \eqref{C:eq:sos}, pair its two ends $E=a_1a_{k+1}$ with
the $k/2$ even positions, all zero, to obtain $k|E|^2$.
For each interior odd position $n$ put $l=s-n$, $w=a_na_l$.
Its pairings with the zeros have total weight $(nl-1)/2$:
sum the number of even positions in $[d,s-d]$ for
$1\le d\le\min(n,l)$. Set aside $nl|w|^2/10$; over all odd
$n,l\ge3$ this gives $a_1^4S^2/10$. The remaining weight
$J=2nl/5-1/2$ and the two end pairings contribute at least
$2J|E|^2/(2+J)$, by completing the square in $w$.
Thus, discarding other nonnegative pairs,
\[
 \begin{gathered}
 D(1+S)^2\ge(S^2+2S)/10+\sum_{k\ge2,\ k\ {\rm even}}L_kX_k,\\
 L_k=k-(k+1)/5+(k/2-1)\frac{2j_k}{2+j_k},
 \qquad j_k=6(k-1)/5-1/2 .
 \end{gathered}
\]
Here $nl\ge3(k-1)$ whenever an interior entry is present. All $L_k$
are positive. For $R_k=(k+1)^2/[(k+5)L_k]$ one has $R_4=85/129<2/3$,
$R_8=13365/23569<4/7$. For $k\ge12$, $j_k\ge k$ gives
$L_k\ge(9k-21)/5$, hence $R_k\le4/7$ since
$4(k+5)(9k-21)-35(k+1)^2=(k-12)(k+38)+1$.

Let $c_j=(1/2)_j/j!$, $s_*=Z(1/4)-1=\sum_{j\ge1}c_j^2/(4j+1)$
and $s_5=\sum_{j\ge1}c_j^2/(4j+5)$.
Every $\widehat\nu(4j)$ is real and at least one: the positive atoms
have phase one and the real part of each negative atom's phase is at most
one. Both $h(r)$ and $h(ir)$ are positive for $0\le r<1$ by symmetry,
nonvanishing and $h(0)=1$. Averaging their logarithms in \eqref{C:eq:sc}
and using AM--GM gives $[h(r)+h(ir)]/2\ge(1-r^4)^{-1/2}$.
Integrate at $r^4=XU^4$, with independent $X$ arcsine on $(0,1)$ and
$U$ of density $5u^4$. Since $\mathbb E X^j=c_j$, division by five
after subtracting one gives $s_5\le\sum_{j\ge1}c_jg_{4j}/(4j+5)$.
Absolute integration follows from Cauchy--Schwarz and the area formula.
Another Cauchy--Schwarz and the bounds on $R_{4j}$ give
\[
 \sum_{k} L_kX_k\ge
 \frac{s_5^2}{\sum_{j\ge1}R_{4j}c_j^2/(4j+5)}
 \ge\frac{s_5^2}{4s_5/7+1/378},
\]
where the first mode has mass $c_1^2/9=1/36$.

To bound the shifted moment, put $M(z)=\sum_{j\ge0}c_j^2/(j+z)$.
The recurrence $(j+1)c_{j+1}=(j+1/2)c_j$ telescopes:
\[
 (z+1/2)^2 M(z+1)-z^2M(z)
 =\lim_{N\to\infty}\frac{(N+1)^2c_{N+1}^2}{N+1+z}=1/\pi
\]
by Wallis' formula. Indeed the $j$th summand is the increment of
$j^2c_j^2/(j+z)$. At $z=1/4$ this gives
$s_5=(1+s_*)/9+4/(9\pi)-1/5>1/16$: the first three product factors
in \eqref{C:eq:regular-product} give $s_*>9/100$, and $\pi<22/7$.
The last lower bound for $\sum L_kX_k$ increases in $s_5$ and equals
$189/1856>1/10$ at $1/16$. It follows that $D>1/10$.

It remains to prove $d_5<1/12$. At $n=5$ the first factor of
\eqref{C:eq:regular-product} is $28/27$. For $j\ge2$ its denominator
$(5j+1)^2(5j-2)=125j^3-15j-2$ is at least $120j^3$. Since
$\zeta(3)-1\le1/8+\int_{5/2}^\infty x^{-3}\,dx<5/24$ by the
midpoint bound, $Z(1/5)<(28/27)e^{1/144}<(28/27)(144/143)$.
The mass term in \eqref{C:eq:spectral-comparator-uniform} satisfies
$23L_5/25>1/600$, since $L_5=41/20160$. It follows that
\[
 \rho_5^{-1}>\frac{27}{28}\frac{143}{144}(1+1/600)>
 \frac{27}{28}\frac{161}{162}=\frac{23}{24}.
\]
For the second inequality use $1/144-1/162=1/1296<1/1200$.
Thus $d_5<47/576+1/625000<1/12$.
\end{proof}

\begin{proposition}[The threshold for the comparison with the regular polygon]\label{C:prop:regular-threshold}
For integers $n\ge6$, set $t=1/n$. The bound
\eqref{C:eq:spectral-comparator-uniform} satisfies
\begin{equation}\label{C:eq:continuous-regular}
 d_n<t^3(8\zeta(3)+5t/4)<10t^3.
\end{equation}
\end{proposition}
\begin{proof}
Put $a=8\zeta(3)$. The ratio $R(t)=2L_n/t^3$ decreases to $219/385$ at $1/6$.
Put $x=23t^2 L_n\le23/(2\cdot6^5)$. Since
$2\log(1+x)\ge 2x/(1+x)$, the lower bound on $R$ gives
\[
       2\log(1+x)\ge
       \frac{23(219/385)}{1+23/(2\cdot6^5)}t^5>13t^5 .
\]
Indeed $23(219/385)-13=32/385>13\cdot23/(2\cdot6^5)$.
The quintic coefficient of $2\log Z$ is $24\zeta(5)<25$, since
$\zeta(5)\le1+2^{-5}+\int_{5/2}^\infty x^{-5}\,dx<25/24$
by the midpoint rule. For the tail in \eqref{C:eq:regular-area} use
$1/k\le1/7$ and $\zeta(k)\le\zeta(5)$ for odd $k\ge7$. Dropping
the subtraction in $2^{k+1}-4$ leaves a geometric series with ratio
$4t^2$. Its bound for the tail of $2\log Z$ is
$(512/7)(25/24)t^7/(1-4t^2)<90t^7$. Subtracting the mass correction gives
\[
             2\log\rho_n<H=at^3+12t^5+90t^7 .
\]
The integral bounds
$1+1/8+\int_3^\infty x^{-3}\,dx<\zeta(3)
 \le1+1/8+\int_{5/2}^\infty x^{-3}\,dx$
give $85/9<a<29/3$. Hence $H<1/20$.
Taylor's inequality gives
$1-e^{-H}\le H-(1/2-H/6)H^2\le H-(49/100)H^2$.
The quadratic correction absorbs the entire tail:
$(49/100)a^2>43$ and $90t\le15$. Therefore
\[
 \frac{d_n-at^3}{t^4}<\frac1{1000}+12t-28t^2
 \le\frac1{1000}+\frac{11}{9}<\frac54 .
\]
The parabola is increasing up to $1/6$; the bound $a<29/3$ gives
the second inequality in \eqref{C:eq:continuous-regular}.
\end{proof}

\paragraph{One weight for five or more positive turns.}
The integrated turning estimate gives a lower bound on the quadratic moment $Q$.
We need a positive excess in every mode and a little more above the second
mode, to absorb its interaction with the other modes in the Hankel minors.
For real $m\ge5$ use the following notation.
Set
\[
 v=\frac1m,\quad t=\frac1{m+1},\quad p=4v,\quad
 \ell=\frac2{m(m-2)},\quad A=1+\ell,\quad \delta=\delta_1=m+2-\frac{m-2}{3m-4} .
\]
For the disk map of a simple polygon of area $\pi$, centered at a point of maximal conformal radius as in Section~\ref{C:sec:conformal}, write
\(f'/f'(0)=1+\sum_{k\ge2}g_k z^k\). Define
\[
 X_k=\frac{|g_k|^2}{(k+1)^2},\quad Y=\sum_{k\ge2}X_k,\quad
 Q=\sum_{k\ge2}(k-1)X_k,\quad S=Q+2Y.
\]
Here $Y>0$, since an affine map would have a disk as its image.

To retain a third moment, start with the kernel $s^2e^{-\sigma s}/2$,
whose Laplace transform is $(k+\sigma)^{-3}$. Its mode coefficient
$k^2(k+1)^2/(k+\sigma)^3=k+2-3\sigma+O(k^{-1})$ agrees with the coefficient
$k-1$ of $Q$ through constant order at $\sigma=1$. For this reference
weight $w_0(s)=e^{-s}s^2/2$ the coefficient $k^2/(k+1)$ still exceeds
$k-1$ by $1/(k+1)$. The excess is exactly $1/k^2$ of the coefficient.
A relative correction of this order preserves the first two terms of
the moment expansion. We therefore try a single reciprocal factor
$[1+c/(k+1)^2]^{-1}$, $c\ge0$:
\[
 W_c(k)=\frac{1}{(k+1)((k+1)^2+c)},\qquad
 w_c(s)=w_0(s)\left(\frac{\sin(\sqrt c\,s/2)}{\sqrt c\,s/2}\right)^2 .
\]
Here $W_c$ is the Laplace transform of $w_c$, with the quotient
interpreted by continuity. In particular $0\le w_c\le w_0$.

The turning integral will retain the third and fourth moments, so we
require $k^3W_c(k)\ge1-3/k$. With $x=1/k$, the numerator of the
difference over the positive denominator $(1+x)((1+x)^2+cx^2)$ is
\[
 x^2\big[(6-c)+(8+2c)x+3(1+c)x^2\big].
\]
Thus this bound holds at every frequency exactly for $0\le c\le6$.
In the other direction the second-mode excess before the turning
correction, $1-36W_c(2)=(c-3)/(c+9)$, increases strictly with $c$.
Within this family $c=6$ retains the largest such excess compatible
with the moment bound. Write $W=W_6$ and $w=w_6$.
Let \(c_j=\prod_{i=1}^{j-1}(1+p/i)\) and \(I=4\sum_{j\ge1}c_jW(mj)\).

\stepcounter{theorem}
\begin{turnlemma}[The second mode and the remaining modes]\label{C:lem:turn-modal}
For real $m\ge5$, let
\[
 R_k=A(k-1)-k^2(k+1)^2W(k)-\ell(k+1)^2/(k+\delta).
\]
Then $R_2\ge At$ and $R_k\ge2At$ for $k\ge3$.
\end{turnlemma}
\begin{proof}
Put $s_k=(5k-7)/[(k+1)^2+6]$ and $\tau=m+1$. Since $\delta>m+3/2$,
$R_2>1/5-\ell/17\ge1/5-2/255>17/90\ge At$.
For $k\ge3$, $(k+\tau)s_k\ge(k+6)s_k>16/5$: clearing the last
denominator leaves $8+(k-3)(9k+110)>0$. Decrease $\delta$ to $\tau$ and
divide $(k+1)^2$ by $k+\tau$. This gives
\[
 (k+\tau)(R_k-2At)\ge (k+\tau)s_k-\ell m^2+(k+\tau)(2/m-2At).
\]
The last coefficient $2/m-2At$ is nonnegative. Thus it suffices that
$\ell m^2-\tau(2/m-2At)=2+2(m+4)/[m(m-2)]\le 16/5$, which follows
at once from $m\ge5$.
\end{proof}

\begin{turnlemma}[Size of the turning integral]\label{C:lem:turn-size}
For real $m\ge5$,
\begin{equation}\label{C:eq:all-cosine-intercept}
             I/A>4t^3\bigl(\zeta(3)+t\bigr).
\end{equation}
\end{turnlemma}
\begin{proof}
Set $a=\zeta(4)-1$, $b=\zeta(3)-1$, $v=1/m$, $p=4v$ and
$K(x)=[(1+x)(1+2x+7x^2)]^{-1}$. Then $J=m^3I/4=\sum c_j B_j$,
where $B_j=j^{-3}K(v/j)$. For the constant part of $c_j$ keep $j=1$
and use $K(x)\ge1-3x$ on the tail. For $j\ge2$ the increment $c_j-1$
is nonnegative and $B_j\ge1/[j^2(j+p)]$: the difference of the two
denominators is $v[(j-10v)(j+v)+3v^2]\ge0$. Thus
\[
 J\ge K(v)+b-3av+U_p,\qquad U_p=\sum_{j\ge1}\frac{c_j-1}{j^2(j+p)}.
\]
The generating series $\sum c_jx^j/j=[(1-x)^{-p}-1]/p$ evaluates
this correction. Write $M(a)=\sum c_j/[j(j+a)]$ and let $M_0$ be the
same sum with $c_j=1$. Integrate the generating series against $x^{a-1}$:
$pM(a)=\mathrm B(a,1-p)-1/a$. At $a=0$ its limit is
$-\gamma-\psi(1-p)$, while $pM_0(p)=\psi(1+p)+\gamma$.
The gamma reflection and recurrence formulas therefore give
\[
 pU_p=M(0)+M_0(p)-M(p)-\zeta(2)
 =\frac2{p^2}-\frac\pi p\cot(\pi p/2)-\zeta(2).
\]
Logarithmically differentiate Euler's product for the sine. The resulting
positive partial fractions give
\[
 U_p=v\sum_{k\ge1}\frac1{k^2(k^2-4v^2)}
 \ge \frac{v}{1-4v^2}+av.
\]
Consequently
\begin{equation}\label{C:eq:J-lower}
 J\ge K(v)+\zeta(3)-1-2(\zeta(4)-1)v+\frac{v}{1-4v^2}.
\end{equation}
For the comparison at a given $m\ge5$ put $H=(1-2v)(1+v)^3$, $E=1-2v+2v^2$.
We need $HJ>E(1+b+t)$. In their difference the coefficient $H-E$ of $b$
is positive, since $(H-E)/v=3-5v-5v^2-2v^3>0$. Use $b>1/5$, $a<1/12$:
sum the first six terms of $\zeta(3)$ and integrate its tail from 7
for the first, and sum the first three terms of $\zeta(4)$ and use the
midpoint tail from $7/2$ for the second. After multiplying
$[HJ-E(1+b+t)]/v$ by $30(1+v)(1+2v)(1+2v+7v^2)$ we obtain the lower bound
\[
 13-90v+5v^2+878v^3+87v^4-1400v^5+475v^6+432v^7+140v^8.
\]
As $v\le1/5$, this exceeds $13-90v+750v^3=1+(1-5v)^2(12+30v)>0$.

\end{proof}

\begin{lemma}[An ordered Hankel sum]\label{C:lem:ordered-hankel}
For $X_k\ge0$, $Q=\sum_{k\ge2}(k-1)X_k<\infty$ and $Y=\sum_{k\ge2}X_k$,
put $Z=2Y-X_2$ and $H=\sum_{i,j\ge2}\sqrt{X_iX_jX_{i+j}}$. For all $a>0$,
\[
 4H\le a Z+(2Y/a)(Q-Y-X_3).
\]
Minimizing in $a$ gives $H^2\le (Q-Y-X_3)Y(Y-X_2/2)$.
\end{lemma}
\begin{proof}
Suppose $Y>0$. Write $h_i=1+\mathbf1_{\{i=2\}}$. Apply Young's inequality
termwise with summands $2aX_iX_j/(Yh_i)$ and $2Y h_iX_{i+j}/a$.
The first sum is $aZ$. In the second, each $k\ge4$ occurs $k-2$ times,
counting the pair with $i=2$ twice.
\end{proof}

The comparison has to hold for every $m\ge5$, not only in the limit of large $m$: it must place every competitor above the explicit threshold $d_{m+1}$, including at the smallest allowed $m$. Set $t=1/(m+1)$. Proposition~\ref{C:prop:regular-threshold} gives $d_{m+1}<8\zeta(3)t^3+\frac54t^4$. Lemmas~\ref{C:lem:turn-modal} and~\ref{C:lem:turn-size} give the lower bound $Q>e+tZ$, with $e$ and $Z$ as in the next proof. That proof uses this bound to control the cubic interaction of the Hankel coefficients and shows $D(P)>8\zeta(3)t^3+\frac32t^4$. The difference of the two bounds is the gap $t^4/4$.

\begin{theorem}[Deficit from the number of convex vertices]\label{C:thm:turn-gap}
Let $m\ge4$ be an integer. Every simple polygon of area $\pi$ with at most
$m$ convex vertices satisfies $D(P)>d_{m+1}$.
For $m\ge5$, putting $t=1/(m+1)$, one has the stronger estimate
\[
 D(P)>8\zeta(3)t^3+\tfrac32t^4>d_{m+1}+\tfrac14t^4 .
\]
\end{theorem}
\begin{proof}
For $m=4$, Theorem~\ref{C:thm:four-turn-sweep} gives $D(P)\ge1/10>d_5$.
Let $m\ge5$, $t=1/(m+1)$, $e=4t^3(\zeta(3)+t)$ and
$d=t^3(8\zeta(3)+3t/2)$. With $Z=2Y-X_2$, the integrated estimate
at $\sigma=1$ gives $AQ\ge I+\sum R_kX_k$. Lemmas~\ref{C:lem:turn-modal}
and~\ref{C:lem:turn-size} give the lower bound $Q>e+tZ$.
It suffices to take $D\le10t^3<1/20$, as $8\zeta(3)<29/3$.
By \eqref{C:eq:B-deficit}, $B^2\le4D$ first gives $B<1/2$, then
$B(1-B)\le3D/2<3/40$ gives $B<1/8$. Thus $S\le3D/[2(1-B)^2]<2D$.
Set $\gamma=D(2+S)<21t^3$ and $c=1-\gamma-2Y/t>0$,
using $Y\le S/3$ and $21t^3+(40/3)t^2<1$.

In \eqref{C:eq:homogeneous-hankel}, $D/b^2=D+\gamma S$.
Lemma~\ref{C:lem:ordered-hankel} with $a=t$ gives
\[
 D+\gamma S\ge 2Q+6Y^2-tZ-(2Y/t)(Q-Y).
\]
Inserting $S=Q+2Y$ and collecting the nonnegative terms yields
\[
 \begin{split}
 D\ge {}&(2-\gamma)e+[t-2e/t-\gamma(2+t)]Y+(4+2/t)Y^2\\
         &+(1+c)(Q-e-tZ)+ct(Z-Y).
 \end{split}
\]
Now $e<(11/2)t^3$, $\gamma<21t^3$, and $1/t\ge6$. With $w=Y/t^2$,
the first line alone gives
\[
 (D-2e)/t^4>-13/4-13w+16w^2\ge-13/4-(13/8)^2>-6,
\]
since $21(11/2)t^2<13/4$ and $11+21t(2+t)-1/t<13$.
As $2e-d=13t^4/2$, we obtain $D>d$.

Proposition~\ref{C:prop:regular-threshold} gives
$d_{m+1}<8\zeta(3)t^3+5t^4/4$, proving the last inequality.
\end{proof}

\subsection{Boundaries with inward-curving arcs}\label{C:sec:arc-proof}
\begin{proof}[Proof of Corollary~\ref{C:cor:concave-arcs}]
Orient the boundary positively. Around every point, including a junction,
the two one-sided tangents have positive projection on a common direction.
In coordinates with this direction as horizontal the boundary is a
graph, with $\Omega$ locally above it. Each smooth piece is concave,
since its signed curvature is nonpositive. Short chords therefore lie
below the graph. Take successively finer inscribed polygonal curves,
including every junction. For fine mesh they are simple: nearby pieces
remain graphs in a common chart, and separated parameter intervals
remain disjoint by uniform convergence and compactness. Replacing a
graph piece by its chord enlarges the domain by the strip between the
two graphs. Thus the polygonal domains $P_j$ contain $\Omega$, and
$|P_j|\to|\Omega|$, since these strips lie in a shrinking neighbourhood
of the boundary.

Chord slopes on each smooth graph are nonincreasing. At a junction
the incoming chord direction is at least the incoming tangent direction
and the outgoing one is at most its tangent direction, in the same local
angle coordinate. A positive turn can therefore occur only at a
junction with positive turn. Proposition~\ref{C:prop:nonconvex} gives
$1-8\pi T(P_j)/|P_j|^2>d_n$. By torsion monotonicity and convergence of
the areas, $1-8\pi T(\Omega)/|\Omega|^2\ge d_n$.
The strict gaps for torsion and the first eigenvalue now follow as in Proposition~\ref{C:prop:nonconvex}.
For the last assertion, scale to area $\pi$ and use $r_q$ from the proof of Corollary~\ref{C:cor:all-q-gap}. That proof gives $r_q(P_j)^q\ge(1-D(P_j))^{-1}>(1-d_n)^{-1}$ for the polygons $P_j$ scaled to area $\pi$. Since $\Omega\subset P_j$ and $|P_j|\to|\Omega|$, domain monotonicity gives $r_q(\Omega)^q\ge(1-d_n)^{-1}$. Finally $(1-d_n)^{-1}>\rho_n^2>r_q(R_n)^q$, as $d_n>1-\rho_n^{-2}$.
\end{proof}

\subsection{Sharp moments of a signed turning measure}\label{C:sec:signed-proof}
The zero count of Lemma~\ref{C:lem:shifted-hardy} gives a sharp moment bound for signed measures on the circle.

\begin{theorem}[Signed turning moments]\label{C:thm:signed-moments}
Let \(\nu=\sum_{i=1}^N\nu_i\delta_{\xi_i}\) have distinct nodes
\(\xi_i\in\partial\mathbb D\), nonzero real weights and total mass
one. Write \(N_+=\#\{i:\nu_i>0\}\) and
\(\widehat\nu(k)=\sum_i\nu_i\overline{\xi_i}^{\,k}\).
Then, for \(0<r<1\) and \(\sigma>1\),
\[
 \sum_{k\ge1}|\widehat\nu(k)|^2r^{2k}\ge\frac{r^{2N_+}}{1-r^{2N_+}},
 \qquad
 \sum_{k\ge1}\frac{|\widehat\nu(k)|^2}{k^\sigma}
 \ge\frac{\zeta(\sigma)}{N_+^\sigma}.
\]
Equality in either inequality holds exactly when \(\nu\) is the uniform
measure on a rotated set of \(N_+\)-th roots of unity.
\end{theorem}

For positive weights both inequalities in Theorem~\ref{C:thm:signed-moments}, including
equality, are due to Cheng, Gui, Hu and Li
\cite[Theorem~4.3(i)]{cghl2026}. The signed version counts only the
positive weights. We give the zero-count proof and its equality analysis below.

\begin{proof}[Proof of Theorem~\ref{C:thm:signed-moments}]
For the measure \(\nu\) in the theorem, set
\[
 H(z)=\sum_{i=1}^N\frac{\nu_i}{1-\overline\xi_i z}-1.
\]
The zero count in Lemma~\ref{C:lem:shifted-hardy} uses only
the total mass and the number of positive weights. It gives \eqref{C:eq:log-mean} for every \(c\ge0\). The mean of
\(|H(r\zeta)-c|^2\) is at least the exponential of the mean of
\(2\log|H(r\zeta)-c|\), and
\(H(r\zeta)-c=-c+\sum_{k\ge1}\widehat\nu(k)r^k\zeta^k\). Parseval's
identity therefore yields
\[
 \sum_{k\ge1}|\widehat\nu(k)|^2r^{2k}+c^2
 \ge (1+c)^2r^{2N_+}.
\]
Write \(R=r^{2N_+}\). The lower bound for the moment sum is the
quadratic \(R(1+c)^2-c^2\), whose maximum over \(c\ge0\) is
\(R/(1-R)\) at \(c=R/(1-R)\). Equality requires
\(|H(r\zeta)-c|=r^{N_+}/(1-r^{2N_+})\) for all \(\zeta\in\mathbb T\).
Squaring and using the value of $c$ gives $|H|^2=R|1+H|^2$ on this circle. In particular $1+H$ never vanishes there. Thus \(B=H/(1+H)\) satisfies \(|B(r\zeta)|=r^{N_+}\).
On the unit circle, \(|B(\zeta)|=1\). Reflection across the two circles
gives
\[
 B(z)\overline{B(1/\overline z)}=1,\qquad
 B(z)\overline{B(r^2/\overline z)}=r^{2N_+},
\]
as rational identities. Thus \(B(r^2z)=r^{2N_+}B(z)\); its power series
at zero gives \(B(z)=\eta z^{N_+}\), \(|\eta|=1\).
The poles of \(H=B/(1-B)\) are the solutions of \(\eta z^{N_+}=1\),
and their weights are \(1/N_+\). Conversely, for this measure
\(|\widehat\nu(k)|\) equals one when \(N_+\) divides \(k\) and zero
otherwise, so it attains equality. For the second inequality, put
\(r=e^{-s/2}\) and integrate against \(s^{\sigma-1}/\Gamma(\sigma)\).
Equality after integration forces equality at every radius by
continuity.
\end{proof}

\begin{corollary}\label{C:cor:log-derivative}
For a conformal map \(f:\mathbb D\to P\), take the branch of
\(\log(f'/f'(0))\) that vanishes at zero. If \(P\) has \(N_+\) convex
vertices, then
\[
 \int_{\mathbb D}\left|\log\frac{f'(z)}{f'(0)}\right|^2\,dA(z)
 \ge 4\pi\sum_{\ell\ge1}\frac{1}{(N_+\ell)^2(N_+\ell+1)}.
\]
Equality holds exactly when \(P\) is regular and \(f(0)\) is its center.
\end{corollary}

\begin{proof}[Proof of Corollary~\ref{C:cor:log-derivative}]
Schwarz--Christoffel gives the Taylor coefficients
\(2\widehat\nu(k)/k\) for every choice of \(f(0)\). Disk orthogonality
gives the squared norm as \(4\pi\sum_{k\ge1}|\widehat\nu(k)|^2/[k^2(k+1)]\).
Integrate the radius inequality at \(r=e^{-s/2}\) against
\((s-1+e^{-s})\,ds\), whose Laplace transform at \(k\) is \(1/[k^2(k+1)]\).
The density is positive for \(s>0\), so equality forces the uniform
turning measure from the radius inequality. Its Schwarz--Christoffel
derivative is proportional to \((1-\eta z^{N_+})^{-2/N_+}\), \(|\eta|=1\),
which gives a regular polygon centered at \(f(0)\). Conversely a centered
map onto a regular polygon has this turning measure.
\end{proof}

For the centered regular pentagon we may take
\(\xi_j=e^{2\pi i j/5}\) and \(\nu_j=1/5\). Summing a geometric
series gives \(\widehat\nu(k)=0\) unless \(5\mid k\), in which case
it is one. The squared moments thus sum to
\(r^{10}/(1-r^{10})\) with weights \(r^{2k}\). For a pentagon
with a reflex vertex, Theorem~\ref{C:thm:signed-moments} gives at least
\(r^8/(1-r^8)\), a strictly larger value for \(0<r<1\).

\section{Symmetric polygons}\label{C:sec:sharp-sectors}
We first show that the normal derivative of the minimizer is monotone along the last facet in a sector between two mirror rays, and then move that facet at fixed area. Throughout these sector arguments write $\lambda_q^{\rm S}$ for the
Dirichlet quotient $\inf\int|\nabla u|^2/\|u\|_q^2$, $1\le q\le2$.
At $q=1$ it equals $1/T$, and at $q=2$ it equals $\lambda_1$.
Circular rearrangement will then remove the reflection hypothesis.
\begin{proposition}[The dihedral case]\label{C:prop:dihedral-sharp}
Let $m\ge3$ be an integer. Suppose that a simple polygon $P$ is invariant under a
rotation of order $m$ and reflection in a line through the rotation center, and has at
most $m$ convex vertices. At the area of $P$, one has
$\lambda_q^{\rm S}(P)\ge\lambda_q^{\rm S}(R_m)$ for $1\le q\le2$,
with equality only for the regular $m$-gon.
\end{proposition}

\begin{proof}
Write $O$ for the rotation center and $\gamma=\pi/m$. The point $O$
lies in $P$. Indeed, if it did not, take a path in $P$ from a point
to its rotation by $2\pi/m$ and concatenate its $m$ rotated copies.
Its winding number about $O$ is $1$ modulo $m$, contrary to simple
connectivity. The convex vertices are therefore one orbit of size $m$.
Reflection permutes this regular orbit, so each convex vertex lies on a
mirror ray. A mirror line meets the Jordan boundary in just two points:
its restriction to the boundary is an orientation-reversing circle
involution and hence has two fixed points. As $O\in P$, there is one
boundary point on each mirror ray.

Consider the sector from a ray without a convex vertex to the
neighboring ray through a convex vertex, with polar angles $0$ and
$\gamma$. Write $x,y$ for coordinates along and perpendicular to
the first ray. The boundary arc goes from the first ray to the second
inside this sector. Its interior turns are all nonpositive. If $2h$ is
the positive turn at the tip, the negative turning in one sector,
including half the turn at its initial point, has magnitude
$h-\gamma$: all $2m$ sectors make the same contribution to the total
turn. At the initial mirror, the outgoing tangent has angle $\pi/2$
minus half the magnitude of the turn there, by reflection. Its angle
then decreases to $\pi/2+\gamma-h>\gamma$. Thus, with $d=\cot\gamma$,
the sector of $P$ is the part of the wedge $W=\{0<y<x/d\}$ satisfying
\[
 x<f(y),\qquad f(y)=\max_{1\le j\le k}(p_jy+c_j),\qquad
 0\le p_1<\cdots<p_k<d .
\]
Here $x=f(y)$ is the chain of boundary facets in increasing $y$.
The slopes increase because the tangent angles decrease; we extend
each affine piece to its whole line. All $c_j$ are positive, since
$f(y)>dy$ before the tip and $p_j<d$. Beyond the tip $f(y)<dy$.

A positive mass-one minimizer is unique. Indeed the root mean square of
two minimizers has no greater energy, by Cauchy--Schwarz, and has
strictly greater $q$-mass for $q<2$ unless the two agree, by strict
concavity of $t^{q/2}$. For $q=2$ use simplicity of the first eigenvalue.
It follows by symmetry and reflection that a minimizer on the polygon
restricts, up to scaling, to the mixed minimizer on its sector.
We claim its inward normal derivative on the last facet is strictly
decreasing towards the tip. Use orthonormal coordinates along
\[
 e_1=(1,-p_k)/\sqrt{1+p_k^2},\qquad
 e_2=(p_k,1)/\sqrt{1+p_k^2}.
\]
These are the last facet's outward normal and its tangent toward the
tip. Both wedge rays have nonnegative coordinates in this basis, as do
all the facet normals $(1,-p_j)$. The part of the wedge outside the cell, which will be penalized, is
\[
 W\setminus\{x<f(y)\}
   =W\cap\bigcap_j\{x-p_j y\ge c_j\}.
\]
Choose a smooth nonnegative increasing cutoff $\eta$ with bounded
derivative, zero on $(-\infty,0]$ and positive on $(0,\infty)$.
Then $V=\prod_j\eta(x-p_jy-c_j)$ vanishes on the cell, is positive
on the interior of its complement in the wedge, and has nonnegative
derivatives in both new coordinates.

Intersect the wedge with a disk $B_R(O)$ containing the closed cell
in its interior relative to the wedge. On the resulting sector impose
zero trace on the outer arc and natural conditions on the two rays.
Let $U_M>0$ minimize $E_M(u)=\int(|\nabla u|^2+MVu^2)$ at
$\int u^q=1$, and write $\lambda_M=E_M(U_M)$, $p=q-1$.
Then $(-\Delta+MV)U_M=\lambda_M U_M^p$. For any admissible $w$,
testing first with $w^2/(U_M+\epsilon)$ for bounded $w$, then by truncation and $\epsilon\downarrow0$, gives
\[
 E_M(w)\ge\lambda_M\int U_M^{p-1}w^2 .
\]
Let $v=-\nabla U_M$ in the new frame, and take negative parts componentwise.
On a ray, $v$ and $v^-$ are tangent to the ray since its outgoing
direction has nonnegative components; $\partial_\nu v$ is normal
by tangential differentiation of the Neumann condition. On the
outer arc, $v$ is a positive multiple of the radial vector by Hopf's
lemma, so $v^-=0$. Testing the differentiated equation by $-v^-$ yields
\[
 \sum_i E_M(v_i^-)\le p\lambda_M\int U_M^{p-1}|v^-|^2 .
\]
Indeed its source $MU_M\nabla V$ is componentwise nonnegative.
For each fixed $M$, even reflection tiles the disk. The reflected
potential is Lipschitz; elliptic regularity up to its smooth circle
gives $U_M\in C^{2,\delta}$ for some $\delta>0$.
At $0<p<1$ the differentiated power is locally integrable up to the
circle, since $U_M$ vanishes there to first order. The tests above
follow by approximation, with the derivative of $U_M^0$
taken as zero. For $p<1$ the two energy inequalities force $v^-=0$.
For $p=1$ each negative part would have to be a multiple of the
positive first eigenfunction; each vanishes on an open set next to
the outer arc where both components of the radial vector are positive.
Thus $v^-=0$ here as well.

A cell test extended by zero bounds $\lambda_M$ independently of $M$.
Compact embedding gives a mass-one limit as $M\to\infty$ and
$\int VU_M^2\to0$. The limit, extended by zero off the cell,
has zero trace on the Dirichlet chain. Lower semicontinuity and the
same tests identify it as the cell minimizer $U$. The inequalities
$-\partial_{e_i}U\ge0$ pass to the limit in distributions.
In particular $w=-\partial_{e_2}U\ge0$ is superharmonic in the cell
since $-\Delta w=p\lambda U^{p-1}w$. It cannot vanish identically:
the Neumann condition on the upper ray would force a function of
the first coordinate alone to be constant on an interval, contrary to
$-\Delta U=\lambda U^p>0$. Hence $w>0$. On the open last facet
$w=0$ and Hopf's lemma gives
$\partial_{e_2}(\partial_{\rm in}U)=-\partial_{\rm in}w<0$.

We now vary only that facet, preserving area. If $k\ge2$, its
predecessor is $x=py+c$. The predecessor's line meets $x=dy$ at height
$L_0=c/(d-p)$. The last facet crosses the predecessor at height
$Y<L_0$ and meets the ray at height $L>L_0$. Extending the predecessor
to the ray would remove a triangle of area
\[
             \tfrac12(d-p)(L_0-Y)(L-L_0).
\]
Keep this area fixed and decrease $Y$ until it reaches the preceding
breakpoint (height zero if $k=2$). Then $L$ decreases, and the new
terminal slope
\[
       s=\frac{d(L-L_0)+p(L_0-Y)}{L-Y}
\]
remains strictly between $p$ and $d$. At the inner endpoint the outward
normal velocity has the sign of $(p-s)\dot Y>0$, and at the tip the
sign of $(d-s)\dot L<0$. On the segment the normal velocity $V$ is
affine and has integral zero by area preservation. Since $g^2=(\partial_{\rm in}U)^2$ strictly
decreases along it, $\int g^2 V>0$. The other facets move only
tangentially. For the normalized minimizer on the full polygon the shape derivative
of the quotient is $-\int_{\partial P}(\partial_\nu U)^2V_n$.
It is strictly negative. Indeed, on each open interval of motion a
triangulation from $O$ transports the domain by piecewise affine maps.
Testing in both directions and using uniqueness and strong energy
convergence of the transported minimizers gives the derivative of the
pulled-back quotient at its minimizer. Integrating the two bulk terms
by parts using $-\Delta U=\lambda_q^{\rm S} U^{q-1}$ gives the displayed
boundary derivative. The minimizer $U$ is bounded: extend the positive right-hand side by zero to a
containing disk. It belongs to $L^2$ by the planar Sobolev embedding,
and by comparison $U$ is at most the bounded solution of the Dirichlet problem in the disk with this right-hand side. There is no
corner term: at a Dirichlet vertex of angle $\omega<2\pi$ choose
$1/2<\sigma<\tau<\min(1,\pi/\omega)$. In its sector
$r^\sigma\cos(\tau(\theta-\omega/2))$ bounds $U$ up to a
constant by comparison, since its negative Laplacian is positive of
order $r^{\sigma-2}$. The gradient estimate on rescaled annuli, using
odd reflection across a straight side, is $O(r^{\sigma-1})$.
The integration boundary on a circle of radius $r$ therefore contributes
$O(r^{2\sigma-1})$. At the end of the motion domains converge between
concentric dilates with ratios tending to one; indeed their radial
functions in the fixed sector are
$\max_j c_j/(\cos\theta-p_j\sin\theta)$ and converge uniformly.
Inclusion and scaling give continuity of the quotient.

The collision deletes the penultimate facet. Repeating leaves one
facet between points at distances $a,b$ on the initial and tip rays.
Its nonnegative slope is equivalent to $a\le b\cos\gamma$. Keeping
$ab$ fixed preserves area. Increase $a$ and decrease $b$ until
$a=b\cos\gamma$. The outward normal velocity is again positive at the
initial point and negative at the tip, so the monotonicity of the normal derivative along the facet again
gives strict improvement whenever $a<b\cos\gamma$. The final facet
is perpendicular to the first ray; the reflected polygon is $R_m$.
\end{proof}

\begin{theorem}[The cyclic case]\label{C:thm:cyclic-sharp}
Let $m\ge3$ be an integer and let $P$ be a simple polygon invariant under rotation by
$2\pi/m$ with at most $m$ convex vertices. The regular $m$-gon of the
same area satisfies $\lambda_q^{\rm S}(P)\ge\lambda_q^{\rm S}(R_m)$ for
every $1\le q\le2$, with equality only for the regular polygon.
\end{theorem}
\begin{proof}
As in Proposition~\ref{C:prop:dihedral-sharp}, the rotation center $O$
lies in $P$ and the convex vertices form one orbit.
They are the vertices of the convex hull: each extreme vertex is convex,
and its rotations give the entire orbit. Put $\gamma=\pi/m$ and let
$A,B$ be two consecutive hull vertices at arguments $0,2\gamma$ on
the circle of radius $b$ centered at $O$. The boundary visits the hull
vertices in cyclic order. Indeed two disjoint boundary arcs cannot join alternating points of the hull.
If their intervening chain is straight, $P$ is
regular. Otherwise the chain and $AB$ bound a convex polygon $K$:
its interior is on the exterior side of the chain, so every
intermediate reflex vertex of $P$ is convex in $K$. At the endpoints
convexity follows from the supporting line $AB$. The pockets obtained
by rotations have disjoint interiors, by simplicity of the boundary.

\emph{The radius along one chain.}
For every interior vertex $X$ of the chain,
\begin{equation}\label{C:eq:cyclic-pocket}
                  \max(A\cdot X,B\cdot X)>|X|^2 .
\end{equation}
To see this, write $X=r(\cos\theta,\sin\theta)$ with origin $O$.
The triangle $ABX$ lies in $K$. Replacing each pocket by the
corresponding smaller triangle leaves $A$ an extreme point of the
same regular convex hull, hence a convex vertex of the resulting
alternating polygon. For a point in $\operatorname{int}K$ the
replacement is a simple polygon, since the new open segments lie in
the pocket interiors. Passing to $X$ by a limit thus gives
$(A-U^{-1}X)\mathbin{\times}(X-A)\ge0$, where $U$ is rotation by
$2\gamma$ and $\times$ denotes the oriented plane determinant. Expanding gives
\[
  0\le br[\sin\theta-\sin(\theta-2\gamma)]-r^2\sin(2\gamma)
    =2r\sin\gamma\,[b\cos(\theta-\gamma)-r\cos\gamma].
\]
Thus $r\le b\cos(\theta-\gamma)/\cos\gamma$.
Take the lift with $|\theta-\gamma|<\pi/2$. If $\theta<0$, the
last ratio is $\cos\theta+\sin\theta\tan\gamma<\cos\theta$.
If $\theta>2\gamma$, reflection in the bisector gives the same
conclusion with $B$. For $0\le\theta\le2\gamma$ the support of the
regular hull in direction $X$ is attained at $A$ or $B$, and $X$
lies strictly inside that hull. This proves \eqref{C:eq:cyclic-pocket}.

On a straight segment the squared radius is strictly convex, so an interior local maximum can occur only at a vertex.
There the two adjacent edge vectors pointing away from $X$ would
both have nonpositive dot product with $X$. By convexity these
vectors span a cone containing $K-X$, contradicting \eqref{C:eq:cyclic-pocket}. Thus the
radius decreases strictly from $A$ to a unique minimum $b_0>0$
and increases strictly to $B$; it cannot be constant on a line segment.
On either branch oriented from the minimum to a tip, the radial
derivative stays positive on each open segment and on both sides of
every intermediate breakpoint. Indeed the vector from a breakpoint toward the
preceding point has strictly negative dot product with $X$ (a tangent
line to the circle immediately goes to larger radius); the vector
toward the next point must have positive dot product by \eqref{C:eq:cyclic-pocket} and the
same cone argument.

\emph{Circular symmetrization.}
Use $t=\log r$ and angular period $2\gamma$. For
$\log b_0<t<\log b$ the section of $P$ on this circle is one proper
arc in the period: it has exactly two boundary crossings. Lift its
endpoints to $l(t)<h(t)<l(t)+2\gamma$. Its width tends to $2\gamma$
at $b_0$ and to zero at $b$. At radii below $b_0$ the circle is full.
Center every such arc at argument zero, giving $P^*$. Its left
boundary has argument $\Theta=(l-h)/2$. The two functions
$\beta_1=l$, $\beta_2=-h$ parametrize the original left branch and the
reflected right branch, both oriented with the domain on their left.
For $z(t)=e^t(\cos\beta,\sin\beta)$,
\[
 \det(z',z'')=e^{2t}\bigl(\beta''+\beta'+(\beta')^3\bigr).
\]
Here primes denote derivatives in $t$. The determinant vanishes on each straight piece. At a reflex vertex the tangent
angle $\beta+\arctan\beta'$ decreases, so $\beta'$ jumps downward.
The preceding radial argument ensures that both one-sided derivatives
are finite at intermediate breakpoints. With $\psi=\Theta'$ and
$\delta=(\beta'_1-\beta'_2)/2$, averaging the two identities gives
\[
             \psi'+\psi(1+\psi^2+3\delta^2)=0
\]
between breakpoints and downward jumps of $\psi$ at the breakpoints. On each piece the sign of $\psi$ is preserved by the linear equation
$\psi'=-(1+\psi^2+3\delta^2)\psi$. Since jumps are downward, a nonpositive value
cannot be followed by a positive one. But
$\int_t^{\log b}\psi=(h(t)-l(t))/2>0$, hence $\psi>0$ on the
pieces. Consequently $\psi'+\psi+\psi^3=-3\psi\delta^2\le0$. Thus the centered arc turns
only clockwise. At its initial mirror ray the tangent makes angle at
most $\pi/2$ with that ray, since the radius increases. On the terminal pieces the two angular derivatives remain bounded up to
the tip because each nontrivial edge leaves the outer circle
transversely. The coefficient in the linear equation is bounded there, so its positive solution has a
strictly positive terminal value of $\psi$. The tangent there lies strictly
ahead of the tip ray. Consequently the centered boundary is a convex graph $x=f(y)$
in the coordinates of Proposition~\ref{C:prop:dihedral-sharp}.
The chordal interpolants of this convex graph lie above it, have
nondecreasing slopes in $[0,\cot\gamma)$ and converge uniformly.
Reflecting them gives outer approximating dihedral polygons with at most $m$ convex vertices, with areas tending
to $|P^*|=|P|$. Inclusion, scaling and that proposition give
$\lambda_q^{\rm S}(P^*)\ge\lambda_q^{\rm S}(R_m)$.

A nonnegative test function $u$ on $P$ that is invariant under the rotation, extended by zero, may
be rearranged decreasingly about zero on the quotient circle of length
$2\gamma$, with the result repeated in each of the $m$ periods.
This preserves its $q$-mass. Angular Dirichlet
energy decreases by one-dimensional symmetric rearrangement on the
circle. The radial energy decreases by $L^2$ contraction of
rearrangement applied to difference quotients. More precisely, if $v$
is the rearrangement of $u$, then $\|v_t(t,\cdot)\|_2\le
\|u_t(t,\cdot)\|_2$ for almost every $t$, including radii where the
section is a full circle. In $(t,\theta)$
coordinates the energy is $\int(u_t^2+u_\theta^2)\,dt\,d\theta$,
and planar measure is $e^{2t}\,dt\,d\theta$. The rearranged test
is supported on $\overline{P^*}$ and has zero trace: $P^*$ is
Lipschitz, since the sector graph has bounded slope and meets the tip
ray transversely. At the initial ray the graph $x=f(y)$ and its
reflection form the single Lipschitz graph $x=f(|y|)$. Take $u$ to be the positive minimizer on $P$ with $\int_Pu^q=1$. Uniqueness makes it invariant under the rotation, so its rearrangement is an admissible test and $\lambda_q^{\rm S}(P^*)\le\lambda_q^{\rm S}(P)$.

It remains to check strictness, since the outer approximation by
itself gives only non-strict bounds. Equality forces equality in both energy comparisons for the positive
minimizer on $P$. For almost every $t$ with a proper arc, its endpoints cross
straight polygon edges transversely. Hopf's lemma and boundary
$C^1$ regularity give nonzero one-sided angular derivatives there.
The positive slice is analytic inside because $u^{q-1}$ is analytic
where $u>0$. Its critical points are isolated away from its endpoints,
hence finite. For
a regular positive level with superlevel length $L(s)$, coarea gives
$-L'(s)=\sum |u_\theta|^{-1}$ over its boundary points, while the
angular energy density is $\sum |u_\theta|$. For a centered interval
the energy density is $4/(-L')$. Equality in Cauchy--Schwarz forces
two boundary points with equal derivative magnitudes at almost every
regular level. The midpoint of this interval consequently has
derivative zero in $s$. The intervals are nested and their endpoints
are continuous across the finitely many critical values, since a
nonconstant analytic slice has no level interval. Therefore the slice
is a translate of its centered rearrangement by $\phi(t)=(h+l)/2$.

Write $v(t,\theta)$ for the centered rearrangement. Wherever $\phi$ is
smooth, $u(t,\theta)=v(t,\theta-\phi(t))$ gives
\[
 \int |u_t|^2\,d\theta=\int |v_t|^2\,d\theta
                  +|\phi'|^2\int |v_\theta|^2\,d\theta ,
\]
as $v_t$ is even and $v_\theta$ odd. This Sobolev chain rule applies on
compact subintervals between breakpoints. The last angular integral is
positive on every proper nonconstant slice. Equality of the radial
energies forces $\phi'=0$ there, and continuity across the breakpoints
gives a common midpoint for all proper arcs. Thus $P$ itself has a
reflection axis, and Proposition~\ref{C:prop:dihedral-sharp} gives the
claimed rigidity.
\end{proof}

\subsection{Equal-area sectors without a symmetry assumption}
\begin{corollary}\label{C:cor:equal-sector-turns}
Let $m\ge3$ and $\gamma=\pi/m$. Suppose the rays
$\theta=j\gamma$ about an interior point divide a simple polygon $P$
into $2m$ cells of equal area. In each cell take the even-indexed ray
as the positive $x$-axis and point the $y$-axis into the cell.
Suppose its outer boundary is a continuous convex polygonal graph
$x=f(y)$ joining the two rays, with all slopes in $[0,\cot\gamma)$.
Then $\lambda_q^{\rm S}(P)\ge\lambda_q^{\rm S}(R_m)$ at fixed area
for all $1\le q\le2$, with equality only for $P\cong R_m$.
\end{corollary}
\begin{proof}
Reflect a single cell successively across its radial sides. The graph
condition gives one boundary point on every intermediate ray:
$f(y)-y\cot\theta$ decreases strictly for $0<\theta\le\gamma$.
The copies form a simple $D_m$-invariant polygon. There is one convex
tip on each odd ray, since the terminal slope is less than
$\cot\gamma$. The nonnegative initial slope gives a flat or reflex
corner on each even ray, and convexity of $f$ makes its other turns
nonpositive. This reflected polygon has area $|P|$. For $N=2m$, uniqueness and reflection give
$\lambda_q^{\rm S}(P_{\rm refl})=N^{1-2/q}\lambda_{q,{\rm DN}}(C)$
for this cell $C$. By Proposition~\ref{C:prop:dihedral-sharp}, its mixed
constant is therefore at least $N^{2/q-1}\lambda_q^{\rm S}(R_m)$.
On the actual polygon, for each cell let $m_i>0$ be its share of the
$q$-mass of the positive unit minimizer. Summing its energy restrictions
and using $\sum_i m_i^{2/q}\ge N^{1-2/q}$ proves the bound.
At equality every reflected
polygon is regular. Each cell is thus a right triangle with its outer
side perpendicular to its even ray. Equal areas give the same altitude
in all cells; their union is $R_m$.
\end{proof}

The cells need not be congruent. For example, start with the same
convex graph in each cell, with three distinct slopes strictly between
$0$ and $\cot\gamma$. Fix its endpoints and the heights of its two
interior breakpoints. Small changes to their $x$-coordinates preserve
strict ordering of the slopes, and the cell area is affine in those
coordinates. Its level line therefore gives distinct graphs with the
same endpoints and area, which can be chosen independently in the
cells.

\section{Polygons near the regular polygon}\label{C:sec:near-regular}
This section shows that $R_n$ is a strict local minimizer among polygons with at most $n$ convex vertices (Theorem~\ref{C:thm:local}). We first show that filling in a small corner at a reflex vertex lowers $\Phi_q$.

The gain at a reflex corner comes from the first Dirichlet harmonic
of a sector. For related self-similar corner perturbations of
Dirichlet problems see Dauge, Tordeux and Vial \cite{dauge-corner}.

\begin{lemma}[Filling a reflex corner]\label{C:lem:fill-reflex}
Let $P$ be a simple polygon with a reflex corner of angle
$\omega>\pi$. Cut the tip off the exterior notch by joining the points
at distance $\varepsilon$ from the corner on its two adjacent sides,
and add the resulting triangle to $P$. There are constants $c,\varepsilon_0>0$ depending only on $P$ and the chosen corner
such that $P_\varepsilon$ has the same number of convex vertices and
\[
 |P|^{2/q}\lambda_q^{\rm S}(P)
 -|P_\varepsilon|^{2/q}\lambda_q^{\rm S}(P_\varepsilon)
 \ge c\varepsilon^{2\pi/\omega}
 \qquad (0<\varepsilon<\varepsilon_0,\ 1\le q\le2).
\]
Here $\lambda_q^{\rm S}$ is the infimum of $\int|\nabla v|^2/\|v\|_q^2$
over nonzero Dirichlet tests. In particular filling increases normalized torsion
and decreases the area-normalized first eigenvalue.
\end{lemma}
\begin{proof}
The exterior notch has angle $2\pi-\omega<\pi$. The added triangle
replaces one reflex corner by two reflex corners and has area
$O(\varepsilon^2)$. Write $\sigma=\pi/\omega<1$. For
$1\le q\le2$, let $u>0$ be an $L^q$-normalized minimizer for
$\Lambda=\lambda_q^{\rm S}(P)$, extended by zero to $P_\varepsilon$.
It solves $-\Delta u=\Lambda u^{q-1}$ in $P$.

In polar coordinates at the vertex, some fixed sector
$0<r<R$, $0<\theta<\omega$ agrees with $P$ inside $B_R$.
The comparison constant on its circular arc can be chosen independently of $q$.
Indeed, a fixed smooth trial bounds $\Lambda$ uniformly above. Every sequence
of normalized minimizers with $q_j\to q_*$ thus has a subsequence converging weakly in
$H^1_0(P)$ and strongly in each finite $L^p(P)$. Uniform integrability gives
$\int u_*^{q_*}=1$. Testing against a fixed $q_*$-minimizer and using weak lower
semicontinuity identifies $u_*$ as a nonnegative minimizer. On a neighborhood
of the arc the convergence is also $C^1$: the sources $\Lambda u^{q-1}$ are
bounded in every finite $L^p$ by the planar Sobolev embedding, and interior
$W^{2,p}$ estimates apply up to the endpoints after odd reflection across the
straight sides. Positivity in the interior and Hopf's lemma at the endpoints
give a positive minimum for $u_*(R,\theta)/(R^\sigma\sin(\sigma\theta))$,
using the limiting ratios at $0,\omega$. The $C^1$ convergence and the
common zero traces give the same lower bound up to a factor two along
the subsequence. Compactness therefore gives a single $a>0$ for all $q$.
Comparison with the harmonic function $a r^\sigma\sin(\sigma\theta)$
therefore gives $-\partial_\nu u\ge a\sigma r^{\sigma-1}$ on each
side for $0<r<R$, where $\nu$ is the exterior normal.

Center a ball at distance $\varepsilon/3$ along one of the old sides.
By dilation of the filled notch, a radius $c_\omega\varepsilon>0$
can be chosen so the ball lies in $P_\varepsilon$ and meets only that
old side, away from the vertex. Take a nonnegative smooth bump
$\eta_\varepsilon$ supported in the ball, equal to one on its
concentric half-ball, by scaling a fixed bump. Its Dirichlet energy
$E$ is independent of $\varepsilon$. Integration by parts against
this bump gives
\[
 \int_P\nabla u\cdot\nabla\eta_\varepsilon
 -\Lambda\int_Pu^{q-1}\eta_\varepsilon
 =\int_{\partial P}\eta_\varepsilon\partial_\nu u
 \le -c\varepsilon^\sigma .
\]
For $h>0$, H\"older's inequality on $P$ gives
$\|u+h\eta_\varepsilon\|_{L^q(P_\varepsilon)}^2
\ge 1+2h\int_Pu^{q-1}\eta_\varepsilon$.
Thus the energy of this trial function minus $\Lambda$ times
its squared $L^q$ norm is at most $-2ch\varepsilon^\sigma+h^2E$.
Take $h=c\varepsilon^\sigma/E$. The norm is at most $1+h\|\eta_\varepsilon\|_q\le1+Ch$
uniformly, so $\lambda_q^{\rm S}(P_\varepsilon)\le\Lambda-c'\varepsilon^{2\sigma}$
with $c'>0$ independent of $q$. The change in $|P|^{2/q}$ is $O(\varepsilon^2)$
uniformly on $[1,2]$. Since $2\sigma<2$, reducing $\varepsilon_0$ gives the claim.
\end{proof}

\subsection{Inward perturbations of the regular polygon}\label{C:sec:first-order}
The next two statements concern polygons close to $R_n$. The first is the monotonicity of the normal derivative along a side, established inside the proof of Proposition~\ref{C:prop:dihedral-sharp}.

\begin{lemma}[The normal derivative along a side of $R_n$]\label{C:lem:regular-pressure}
Let $n\ge3$, $1\le q\le2$, and let $u$ be a positive minimizer for $\lambda_q^{\rm S}(R_n)$. On each side of $R_n$, $|\partial_\nu u|$ is symmetric about the midpoint and strictly increases from each endpoint to the midpoint.
\end{lemma}
\begin{proof}
Normalize $u$ by $\int_{R_n}u^q=1$. The positive minimizer with this normalization is unique, as shown in the proof of Proposition~\ref{C:prop:dihedral-sharp}, so it has the symmetries of $R_n$. In particular it is symmetric across the perpendicular bisector of each side. Apply that proof with $m=n$ to the sector between the bisector of a side and the ray through an adjacent vertex. Its boundary chain is a single facet perpendicular to the first ray, so $k=1$ and $p_1=0<\cot(\pi/n)$. The proof shows that the inward normal derivative strictly decreases along this facet towards the vertex.
\end{proof}

\begin{proposition}[Inward perturbations of $R_n$]\label{C:prop:first-order}
Let $n\ge3$ and $1\le q\le2$. Parametrize a side $E$ of $R_n$ by arclength $s\in[0,L]$, and let $y$ be the distance from $E$ towards the interior. Let $g\colon[0,L]\to[0,\infty)$ be concave and piecewise affine, with $g(0)=g(L)=0$ and $g\not\equiv0$. For small $\varepsilon>0$ put
\[
 P_\varepsilon=R_n\setminus\{(s,y):0\le s\le L,\ 0\le y\le \varepsilon g(s)\}.
\]
Then $P_\varepsilon$ is a simple polygon with exactly $n$ convex vertices, and
\[
 \lim_{\varepsilon\downarrow0}\frac{\log\Phi_q(P_\varepsilon)-\log\Phi_q(R_n)}{\varepsilon}
 =\sum_{i}\Pi(s_i)\bigl(g'(s_i^+)-g'(s_i^-)\bigr)>0 .
\]
Here $s_i$ runs over the interior breakpoints of $g$, and $\Pi<0$ on $(0,L)$ is defined in \eqref{C:eq:pressure-primitive} below. The same holds if several sides are pushed in at once; the first-order changes add.
\end{proposition}
\begin{proof}
For small $\varepsilon$ the vertices of $R_n$ remain convex, since the slopes of $g$ are finite. At a breakpoint of $g$ the boundary turns away from the interior, because $g$ is concave, so these vertices are reflex. Hence $P_\varepsilon$ has exactly $n$ convex vertices.

Let $u$ be the positive minimizer on $R_n$ with $\int u^q=1$, and write $\Lambda=\lambda_q^{\rm S}(R_n)$ and $p=|\partial_\nu u|^2/\Lambda$ on $E$. The vertices of $P_\varepsilon$ move affinely in $\varepsilon$. As in the proof of Proposition~\ref{C:prop:dihedral-sharp}, a fixed triangulation transports the domains by piecewise affine maps, the quotient is differentiable, and there is no contribution from the corners. The outward normal velocity on $E$ is $-g$, so
\[
 \frac{d}{d\varepsilon}\Big|_{0^+}\log\lambda_q^{\rm S}(P_\varepsilon)=\int_0^Lp\,g\,ds,\qquad
 \frac{d}{d\varepsilon}\Big|_{0^+}\log|P_\varepsilon|=-\frac{1}{|R_n|}\int_0^Lg\,ds .
\]
The same formula for the dilation $x\cdot\nu$, which equals the inradius $r$ on every side, gives $\int_{\partial R_n}p\,ds\cdot r=4/q$ because $\Phi_q$ does not change under dilation. Since $r|\partial R_n|=2|R_n|$ and all sides are alike, the mean of $p$ over $E$ is $\bar p=2/(q|R_n|)$. Hence the derivative of $\log\Phi_q(P_\varepsilon)$ is $\int_0^L(p-\bar p)g\,ds$.

Put
\begin{equation}\label{C:eq:pressure-primitive}
 \Pi(s)=\int_0^s(s-t)\bigl(p(t)-\bar p\bigr)\,dt .
\end{equation}
Then $\Pi''=p-\bar p$ and $\Pi(0)=\Pi'(0)=0$. By Lemma~\ref{C:lem:regular-pressure}, $p$ is symmetric about $L/2$, so $\int_0^Lt(p-\bar p)\,dt=\tfrac L2\int_0^L(p-\bar p)\,dt=0$; thus $\Pi(L)=\Pi'(L)=0$. Integrating by parts on each interval where $g$ is affine, and using $g(0)=g(L)=0$ and $\Pi(0)=\Pi(L)=0$,
\[
 \int_0^L(p-\bar p)g\,ds=-\int_0^L\Pi'g'\,ds=\sum_i\Pi(s_i)\bigl(g'(s_i^+)-g'(s_i^-)\bigr).
\]
Each jump of $g'$ is nonpositive, and at least one is negative because $g$ vanishes at both ends without vanishing identically.

It remains to show $\Pi<0$ on $(0,L)$. By Lemma~\ref{C:lem:regular-pressure}, $p-\bar p$ strictly increases on $[0,L/2]$ and is symmetric about $L/2$. Since it has mean zero, there is $s_1\in(0,L/2)$ with $p-\bar p<0$ on $[0,s_1)\cup(L-s_1,L]$ and $p-\bar p>0$ on $(s_1,L-s_1)$. Hence $\Pi'(s)=\int_0^s(p-\bar p)$ decreases on $[0,s_1]$ and increases on $[s_1,L-s_1]$. By symmetry $\Pi'(L-s)=-\Pi'(s)$, so $\Pi'(L/2)=0$. Therefore $\Pi'<0$ on $(0,L/2)$ and $\Pi'>0$ on $(L/2,L)$. As $\Pi(0)=\Pi(L)=0$, we get $\Pi<0$ on $(0,L)$.
\end{proof}

The first-order gain is not uniform. Since $p$ vanishes at the vertices, $\Pi(s)=-\bar ps^2/2+o(s^2)$ as $s\downarrow0$. For a triangular pocket of height $H$ with apex at $s=a$, the derivative is $-H\Pi(a)(1/a+1/(L-a))$, which is asymptotic to $H\bar p a/2$ as $a\downarrow0$. Thus the positive first variation has no lower bound uniform in the apex position. Theorem~\ref{C:thm:local} below proves local minimality by a different comparison.

\subsection{Local minimality of the regular polygon}\label{C:sec:local}
Proposition~\ref{C:prop:first-order} compares a perturbation directly with $R_n$. To treat all nearby polygons at once, we compare $P$ instead with a convex $n$-gon inside $P$, and then use Theorem~\ref{thm:main}. Among convex polygons with at most $n$ sides, local optimality of $R_n$ was proved earlier for torsion when $n\ge5$ \cite{bbf2026}, by validated computation for torsion when $5\le n\le25$ \cite{bogosel2026}, and for the first eigenvalue when $n=5,6$ \cite{bb2026}; Theorem~\ref{thm:main} now gives the global statement in that class. The class here is larger, since the number of reflex vertices is unrestricted.

\begin{theorem}[Local minimality of $R_n$]\label{C:thm:local}
Let $n\ge3$ and $1\le q\le2$. There are $\varepsilon,c>0$, depending only on $n$ and $q$, with the following property. Let $P$ be a simple polygon with at most $n$ convex vertices, and let $R_n$ be the regular $n$-gon of the same area $A=|P|$. If $d_H(\overline P,g\overline{R_n})<\varepsilon\sqrt A$ for some rigid motion $g$, then $P$ contains a convex $n$-gon $K$ with
\[
 \Phi_q(P)\ge\Phi_q(K)+c\,\frac{|P\setminus K|}{A} .
\]
Consequently $\Phi_q(P)\ge\Phi_q(R_n)$, with equality only if $P$ is congruent to $R_n$. The number of reflex vertices of $P$ is unrestricted.
\end{theorem}

The second statement follows from the first. If $|P\setminus K|=0$, then $P=K$, because $P$ is open and $K$ is convex; in that case Theorem~\ref{thm:main} applies to $P$. Otherwise $\Phi_q(P)>\Phi_q(K)\ge\Phi_q(R_n)$ by Theorem~\ref{thm:main}.

The proof uses a dual form of the variational principle. For a convex polygon $K$, the nearest-point map onto $\overline K$ divides the exterior of $K$ into a strip over each open side and a sector at each vertex.

\begin{lemma}[Extending the flux outside a convex polygon]\label{C:lem:flux}
Let $K\subset P$ be bounded polygons with $K$ convex, and let $f\ge0$ be bounded on $P$, extended by zero outside $P$. Let $w_K\in H^1_0(K)$ and $w_P\in H^1_0(P)$ solve $-\Delta w=f$, and put $\sigma=|\partial_\nu w_K|$ on the sides of $K$. Define $J=\nabla w_K$ in $K$,
\[
 J(y+r\nu)=-\Bigl(\sigma(y)+\int_0^rf(y+\rho\nu)\,d\rho\Bigr)\nu
\]
in the strip over a side with outer normal $\nu$ ($y$ on the side, $r>0$), and
\[
 J(x+r\omega)=-\frac{\omega}{r}\int_0^r\rho f(x+\rho\omega)\,d\rho
\]
in the sector at a vertex $x$ ($|\omega|=1$, $r>0$). Then
\[
 \int_Pfw_P\le\int_K|\nabla w_K|^2+\int_{P\setminus K}|J|^2 .
\]
\end{lemma}
\begin{proof}
In a strip $\operatorname{div}J=-f$ by differentiation in $r$, and in a sector by $\operatorname{div}(\phi\,\omega)=r^{-1}\partial_r(r\phi)$. Since $w_K\ge0$, the outer normal derivative of $w_K$ on a side is $-\sigma$, which is the normal component of $J$ on the outside. On a ray separating a strip from a sector, $J$ is parallel to the ray on both sides, so its normal component vanishes there. Hence $\operatorname{div}J=-f$ in the sense of distributions in the plane. The gradient of $w_K$ is bounded because the angles of $K$ are less than $\pi$, so $J\in L^2(P)$. Testing with $w_P\in H^1_0(P)$ gives
\[
 \int_P|\nabla w_P|^2=\int_Pfw_P=\int_PJ\cdot\nabla w_P\le\|J\|_{L^2(P)}\|\nabla w_P\|_{L^2(P)},
\]
and the lemma follows.
\end{proof}

\begin{proof}[Proof of Theorem~\ref{C:thm:local}]
Dilate by $A^{-1/2}$ and apply a rigid motion so that $R=R_n$ is a fixed regular polygon of area one. This dilation divides Hausdorff distances by $\sqrt A$ and areas by $A$, and leaves $\Phi_q$ unchanged. It therefore suffices to prove the estimate at unit area.
Let $l$ be half the side length of $R$, and let $u$ be the positive minimizer for $\Lambda=\lambda_q^{\rm S}(R)$ with $\int_Ru^q=1$. Thus $-\Delta u=\Lambda u^{q-1}$, read as $-\Delta u=\Lambda$ when $q=1$. Put $w_R=u/\Lambda$ and let $p_R=|\partial_\nu w_R|^2$ on the sides of $R$. We use $s$ for arclength measured from the midpoint of a side.

\emph{(a) The function $p_R$.} By Lemma~\ref{C:lem:regular-pressure}, $p_R$ is symmetric about $s=0$ and strictly decreasing in $|s|$; it tends to zero at the vertices by the corner bound in (d) below. The dilation identity in the proof of Proposition~\ref{C:prop:first-order} gives $\int_{\partial R}|\partial_\nu u|^2=4\Lambda/(qr)$, where $r$ is the inradius. As $r|\partial R|=2|R|$ and all half-sides are alike, the mean of $p_R$ over $0\le s\le l$ is
\[
 M=\frac{2}{q\Lambda|R|} .
\]
Since $s$ increases and $p_R$ strictly decreases on $[0,l]$, Chebyshev's inequality gives $\frac{2}{l^2}\int_0^lsp_R\,ds<M$. Moreover, if $0<l_0\le l$ and $g\ge0$ is convex on $[0,l_0]$ with $g(0)=0$, then $g(s)/s$ is nondecreasing, and Chebyshev's inequality for the measure $s\,ds$ gives
\begin{equation}\label{C:eq:local-chebyshev}
 \int_0^{l_0}p_Rg\,ds\le\Bigl(\frac{2}{l_0^2}\int_0^{l_0}sp_R\,ds\Bigr)\int_0^{l_0}g\,ds .
\end{equation}
Fix $\gamma>0$ and $l_0<l$ such that $\frac{2}{l_0^2}\int_0^{l_0}sp_R\,ds<M-4\gamma$ and $p_R(s)<M-4\gamma$ for $l_0\le|s|\le l$. They depend only on $n$ and $q$.

It suffices to show the following. If $P_j$ are simple polygons with $|P_j|=|R|$ and at most $n$ convex vertices, and $\overline{P_j}\to\overline R$ in Hausdorff distance, then for all large $j$ there are convex $n$-gons $K_j\subset P_j$ with $\Phi_q(P_j)-\Phi_q(K_j)\ge\gamma\Lambda^2|R|^{2/q}|P_j\setminus K_j|$. Indeed, if the theorem failed with $c=\gamma\Lambda^2|R|^{2/q}$, there would be a sequence of polygons with $\varepsilon\to0$ for which no such $K$ exists.

\emph{(b) An inscribed convex $n$-gon.} Let $H_j$ be the convex hull of $\overline{P_j}$. An extreme point of $H_j$ is a vertex of $P_j$ with angle less than $\pi$, that is, a convex vertex. Since $H_j\to\overline R$ and each vertex of $R$ is an exposed point of $\overline R$, for large $j$ there is an extreme point of $H_j$ near each vertex of $R$. As $P_j$ has at most $n$ convex vertices, $H_j$ is a convex $n$-gon whose vertices are exactly the convex vertices of $P_j$. The boundary of $P_j$ passes through the vertices of $H_j$ in their cyclic order, since two boundary arcs joining alternate hull vertices would cross. So between consecutive hull vertices $V,V'$ every vertex of $\partial P_j$ is reflex. The boundary arc from $V$ to $V'$ and the segment $VV'$ therefore bound a convex polygon outside $P_j$, which we call a pocket. Its angles at the intermediate vertices are $2\pi$ minus reflex angles, and its angles at $V,V'$ are at most $\pi$ because it lies on one side of the line $VV'$. If the arc is the segment itself, the pocket is that segment. Thus $P_j$ is the interior of $H_j$ with the $n$ pockets removed.

The pockets become thin. A pocket of depth $h$ contains a triangle with base $VV'$ and height $h$, hence a disk whose radius depends only on $h$ and on bounds for the side lengths of $H_j$. This disk lies outside $P_j$ and within distance $o(1)$ of $\overline R$, which is impossible for a fixed $h>0$ and large $j$, since every point of $\overline R$ is within $o(1)$ of $\overline{P_j}$.

Take coordinates in which the side $VV'$ is $\{0\le x\le L,\ y=0\}$ and $y$ increases into $H_j$. The upper boundary of the pocket is the graph of a concave function $h\ge0$ on the projection of the pocket. Let $y=\ell(x)=h(L/2)+\sigma(x-L/2)$ be a supporting line of this graph at $x=L/2$. Then the pocket lies in $\{y\le\ell(x)\}$, $\ell(0)\ge0$, $\ell(L)\ge0$, and $|\sigma|\le2\max h/L\to0$. Let $K_j$ be the intersection of the $n$ open half-planes $\{y>\ell(x)\}$, one for each side of $H_j$. The bounding lines converge to the side lines of $R$, so for large $j$ this intersection is an open convex $n$-gon. It is disjoint from every closed pocket, although its boundary can meet a pocket at a point of contact.

To prove $K_j\subset P_j$ it remains to check that $\overline{K_j}\subset H_j$, since $K_j$ is open and avoids the pockets. It suffices to check the vertices of $\overline{K_j}$. Each such vertex is close to a vertex $V$ of $H_j$, so it satisfies the inequalities of the sides of $H_j$ not containing $V$. Let $a,b$ be its signed distances to the two lines of $H_j$ through $V$, positive towards $H_j$, and let $\theta<\pi$ be the angle of $H_j$ at $V$. Its coordinates along the two sides from $V$ are $x_1=(b+a\cos\theta)/\sin\theta$ and $x_2=(a+b\cos\theta)/\sin\theta$, both $o(1)$. The vertex lies on both new lines, so $a=\delta_1+\sigma_1x_1$ and $b=\delta_2+\sigma_2x_2$, where $\delta_i\ge0$ is the value of the new line at $V$, $\sigma_i\to0$ is its slope, and $\delta_i+\sigma_iL_i\ge0$. Suppose $a<0$. If $\sigma_1\le0$, then $a\ge\delta_1$ for $x_1\le0$ and $a\ge\delta_1+\sigma_1L_1$ for $0<x_1<L_1$, both nonnegative; so $\sigma_1>0$. Now $a(1-\sigma_1\cot\theta)=\delta_1+\sigma_1b/\sin\theta$ with $1-\sigma_1\cot\theta>0$, so $b<0$ and $|a|(1-\sigma_1\cot\theta)\le\sigma_1|b|/\sin\theta$. The same argument applied to $b$ gives $|b|(1-\sigma_2\cot\theta)\le\sigma_2|a|/\sin\theta$. Multiplying, $(1-\sigma_1\cot\theta)(1-\sigma_2\cot\theta)\le\sigma_1\sigma_2/\sin^2\theta$, which is false for large $j$. Hence $a\ge0$, and likewise $b\ge0$.

Put $a_j=|P_j\setminus K_j|$ and $\eta_j=\sup_{x\in P_j}\operatorname{dist}(x,K_j)$. Both tend to zero. On each side of $K_j$ measure arclength $s$ from the contact point $(L/2,h(L/2))$ of the supporting line; it tends to the midpoint of the corresponding side of $R$. For large $j$ and $|s|\le l_0$, the points of $P_j\setminus\overline{K_j}$ on the outer normal through $s$ form the segment $0<r<g(s)$, where $g\ge0$ is convex and $g(0)=0$. Indeed $g$ is the gap between the line and the concave graph $h$, measured along the normal of the line; the other pockets do not reach this region because their depths tend to zero.

\emph{(c) The minimizers $u_j$ and the functions $w_j$.} Let $u_j$ be the positive minimizer for $\Lambda_j=\lambda_q^{\rm S}(P_j)$ with $\int u_j^q=1$, and let $f_j=u_j^{q-1}$ on $P_j$ ($f_j=1$ on $P_j$ when $q=1$), extended by zero. Let $w_j\in H^1_0(K_j)$ solve $-\Delta w_j=f_j$ in $K_j$, and let $p_j=|\partial_\nu w_j|^2$ on the sides of $K_j$. Since $K_j\subset P_j\subset\operatorname{int}H_j$, domain monotonicity gives $\lambda_q^{\rm S}(\operatorname{int}H_j)\le\Lambda_j\le\lambda_q^{\rm S}(K_j)$. The convex polygons $K_j$ and $H_j$ converge to $R$, so the piecewise affine maps used in the proof of Proposition~\ref{C:prop:first-order} show that both bounds tend to $\Lambda$; also $|K_j|\to|R|$.

Fix a disk $B$ containing all the polygons. Let $z_j$ solve $-\Delta z_j=\Lambda_jf_j$ in $B$ with zero boundary values. Since $\|f_j\|_{L^{q/(q-1)}}=1$ and $q/(q-1)\ge2$ (and $f_j\le1$ when $q=1$), the Green potential of $B$ gives $\sup_j\|z_j\|_\infty<\infty$, and comparison on $P_j$ gives $0\le u_j\le z_j$. The functions $u_j$ are bounded in $H^1_0(B)$. Every limit point $\tilde u$ in $L^q$ lies in $H^1_0(R)$, since it vanishes outside $\overline R$ and $R$ is a Lipschitz domain, and it satisfies $\int\tilde u^q=1$ and $\int|\nabla\tilde u|^2\le\Lambda$. The positive minimizer on $R$ is unique (see the proof of Proposition~\ref{C:prop:dihedral-sharp}), so $u_j\to u$. For $q>1$ the uniform bound gives $f_j\to f_R$ in every $L^\rho(B)$, $1\le\rho<\infty$, where $f_R=u^{q-1}$ on $R$ and is zero outside. For $q=1$ the sources are the indicator functions of $P_j$ and $R$. They have the same convergence because $|P_j|=|R|$ and every compact subset of $R$ eventually lies in $K_j\subset P_j$. Finally $w_j\to w_R$ in $H^1(B)$: the $w_j$ are bounded in $H^1_0(B)$, every compact subset of $R$ lies in $K_j$ for large $j$, so every weak limit solves the equation of $w_R$, and $\int|\nabla w_j|^2=\int f_jw_j\to\int f_Rw_R=\int|\nabla w_R|^2$.

\emph{(d) Convergence of $p_j$ to $p_R$.} Match each side of $K_j$ with the corresponding side of $R$ so that the contact point corresponds to the midpoint, and extend $p_R$ by zero beyond the endpoints of the side. We claim that $p_j(s)-p_R(s)\to0$ uniformly. On a closed subsegment of an open side, move the side of $K_j$ onto the side of $R$ by a rigid motion close to the identity. The boundary $W^{2,\rho}$ estimate on half-disks, with $\rho>2$, then gives $C^1$ convergence up to the subsegment. Near a vertex, let $\theta_*<\pi$ bound all angles of the polygons $K_j$ and $R$, and fix $1<\alpha<\beta<\min\{2,\pi/\theta_*\}$. In the sector of angle $\theta$ at a vertex, with polar coordinates $(\rho,\phi)$, the function $\rho^\alpha\cos(\beta(\phi-\theta/2))$ is positive and satisfies $-\Delta(\rho^\alpha\cos(\beta(\phi-\theta/2)))=(\beta^2-\alpha^2)\rho^{\alpha-2}\cos(\beta(\phi-\theta/2))$. Since $w_j$ and $f_j$ are uniformly bounded, a large multiple of it dominates $w_j$ in a fixed small sector at the vertex, so $w_j\le C\rho^\alpha$ there. Rescaled gradient estimates on annuli give $|\nabla w_j|\le C\rho^{\alpha-1}$ uniformly in $j$. The same bound holds for $w_R$, so $p_R\to0$ at the vertices, and the claim follows.

\emph{(e) The cost outside $K_j$.} Let $J_j$ be the field of Lemma~\ref{C:lem:flux} for $K_j\subset P_j$ and $f_j$, and let $M_0$ bound $f_j$ and $p_j^{1/2}$. Every point of $P_j\setminus K_j$ is within $\eta_j$ of $K_j$. In a strip, therefore, $|J_j|^2\le(p_j^{1/2}+M_0\eta_j)^2\le p_R(s)+o(1)$, uniformly; in a sector, $|J_j|\le M_0\eta_j/2$. In the strips with $|s|\le l_0$, the set $P_j\setminus\overline{K_j}$ is $\{0<r<g(s)\}$. The polygonal boundary $\partial K_j$ has area zero, so its inclusion in the integration region does not affect the integral. Thus \eqref{C:eq:local-chebyshev} on $[0,l_0]$ and $[-l_0,0]$ bounds its contribution by $(M-4\gamma+o(1))\int g\,ds$, that is, by $(M-4\gamma+o(1))$ times its area. On the rest of the strips and on the sectors, $|J_j|^2\le M-4\gamma+o(1)$ pointwise. Hence, for large $j$,
\[
 \int_{P_j\setminus K_j}|J_j|^2\le(M-3\gamma)\,a_j .
\]
No bound on the number of reflex vertices has been used.

\emph{(f) Conclusion.} The function $u_j/\Lambda_j$ solves $-\Delta w=f_j$ in $P_j$, and $\int_{P_j}f_ju_j/\Lambda_j=1/\Lambda_j$. On $K_j$, Hölder's inequality and the definition of $\lambda_q^{\rm S}(K_j)$ give
\[
 \int_{K_j}|\nabla w_j|^2=\int_{K_j}f_jw_j\le\|f_j\|_{L^{q'}(K_j)}\|w_j\|_{L^q(K_j)}
 \le\frac{\|\nabla w_j\|_{L^2}}{\lambda_q^{\rm S}(K_j)^{1/2}},
\]
since $\|f_j\|_{L^{q'}(K_j)}^{q'}=\int_{K_j}u_j^q\le1$ (and $f_j\le1$ when $q=1$). Thus $\int_{K_j}|\nabla w_j|^2\le1/\lambda_q^{\rm S}(K_j)$. Write $\lambda_K=\lambda_q^{\rm S}(K_j)$. By Lemma~\ref{C:lem:flux}, (e), and $2/(q\lambda_K|K_j|)\to M$,
\[
 \frac1{\Lambda_j}\le\frac1{\lambda_K}+(M-3\gamma)a_j
 \le\frac1{\lambda_K}\Bigl(1+\frac{2a_j}{q|K_j|}\Bigr)-2\gamma a_j
 \le\frac1{\lambda_K}\Bigl(\frac{|P_j|}{|K_j|}\Bigr)^{2/q}-2\gamma a_j
\]
for large $j$, by Bernoulli's inequality and $2/q\ge1$. Multiplying by $\Lambda_j\lambda_K|K_j|^{2/q}$ gives
\[
 \Phi_q(P_j)-\Phi_q(K_j)\ge2\gamma\,\Lambda_j\lambda_K|K_j|^{2/q}a_j\ge\gamma\Lambda^2|R|^{2/q}a_j
\]
for large $j$. This proves the theorem.
\end{proof}

In particular, after normalization to unit area, within a sufficiently small Hausdorff neighbourhood of $R_n$ every polygon with at most $n$ convex vertices has exactly $n$ of them, and Question~\ref{C:que:sharp} below has a positive answer there. The approximation in the proof of Corollary~\ref{C:cor:three-arcs} gives $\Phi_q(\Omega)\ge\Phi_q(R_n)$ also for the curved domains $\Omega$ described there with at most $n$ junctions of positive turn, if $\Omega$ is close enough to $R_n$, though without the equality case.

The theorem also gives quantitative stability. For a domain $P$ let $\alpha(P)=\min_g|P\,\triangle\,gR|/|R|$, where $R$ is the regular $n$-gon of area $|P|$ and $g$ runs over plane isometries.

\begin{corollary}[Local stability]\label{C:cor:local-stability}
Let $n\ge3$ and $1\le q\le2$. There is $c'>0$, depending only on $n$ and $q$, such that every polygon $P$ as in Theorem~\ref{C:thm:local} satisfies
\[
 \Phi_q(P)-\Phi_q(R_n)\ge c'\alpha(P)^2 .
\]
The exponent two is optimal.
\end{corollary}
\begin{proof}
Both the deficit and $\alpha(P)$ are invariant under dilation, so normalize $|P|=|R_n|=1$. Let $K$ be the convex $n$-gon of Theorem~\ref{C:thm:local}, put $a=|P\setminus K|$, and let $R'$ be the regular $n$-gon with $|R'|=|K|$. Theorem~\ref{S:thm:full-stability} and the scale invariance of $\Phi_q$ give $\Phi_q(K)-\Phi_q(R_n)\ge c_n\min\{1,d^2/|K|\}$, where $d=d(K,R')$. Choose an isometry $g$ with $d_H(\overline K,g\overline{R'})=d$, and let $R''\supset gR'$ be the concentric regular $n$-gon of area $|R_n|$. Two convex sets in a fixed disk at Hausdorff distance $d$ differ by an area at most a constant times $d$. Hence
\[
 |R_n|\,\alpha(P)\le|P\,\triangle\,R''|\le|P\setminus K|+|K\,\triangle\,gR'|+|R''\setminus gR'|\le2a+Cd .
\]
After shrinking $\varepsilon$, both $a$ and $d$ are small, and Theorem~\ref{C:thm:local} gives
\[
 \Phi_q(P)-\Phi_q(R_n)\ge ca+\frac{c_nd^2}{|K|}\ge c''(2a+Cd)^2\ge c''|R_n|^2\alpha(P)^2 .
\]
For optimality, center $R_n$ at the origin, put $A_t=\operatorname{diag}(e^t,e^{-t})$ and $P_t=A_tR_n$, a convex $n$-gon of area $|R_n|$. Let $u$ be a minimizer for $\lambda_q^{\rm S}(R_n)$ with $\int_{R_n}|u|^q=1$, and let $M=\int_{R_n}\nabla u\otimes\nabla u$. The function $u\circ A_t^{-1}$ on $P_t$ has the same $L^q$ mass and the energy $e^{-2t}M_{11}+e^{2t}M_{22}$. By Theorem~\ref{thm:main} this energy is at least $\lambda_q^{\rm S}(R_n)=M_{11}+M_{22}$ for every $t$, with equality at $t=0$, so its derivative vanishes there. Hence $M_{11}=M_{22}$ and
\[
 \Phi_q(R_n)\le\Phi_q(P_t)\le\cosh(2t)\,\Phi_q(R_n)=\Phi_q(R_n)+O(t^2).
\]
On the other hand $\alpha(P_t)\ge c|t|$ for small $|t|$. Indeed, for a set $E$ of finite area let $J(E)$ be the difference of the two diagonal entries of $\int_E(x-\bar x_E)\otimes(x-\bar x_E)\,dx$, where $\bar x_E$ is the centroid of $E$. For $R_n$ this matrix commutes with the rotation by $2\pi/n$, $n\ge3$, so it is $jI$ for some $j>0$. Hence $J(gR_n)=0$ for every isometry $g$, while $J(P_t)=2j\sinh(2t)$. For small $|t|$ the minimum defining $\alpha(P_t)$ is attained at some $gR_n$ that meets $P_t$, and both sets then lie in a fixed disk. For two sets of equal area in a fixed disk, $J$ differs by at most a constant times the area of their symmetric difference. Thus $|R_n|\,\alpha(P_t)\ge c'|\sinh(2t)|$, and the exponent two cannot be lowered.
\end{proof}

Here $\alpha(P)$ cannot be replaced by the Hausdorff distance. Start with a regular $n$-gon $R^0$ of circumradius $\rho$.
Near each vertex choose two points, one on each adjacent side, symmetric about the angle bisector and $\delta^2$ apart. Replace the boundary between them by two segments meeting at a tip a distance $\delta$ radially beyond the old vertex. The resulting polygon $P$ contains $R^0$ and has $n$ convex vertices at the tips and $2n$ reflex vertices. At fixed $n$ the added area is $O(\delta^3)$. Let $R$ be the regular $n$-gon of area $|P|$. Its circumradius is $\rho\sqrt{|P|/|R^0|}=\rho+O(\delta^3)$, while the smallest disk enclosing $P$ has radius $\rho+\delta$, as its boundary contains the $n$ equally spaced tips. Since this radius is invariant under rigid motions and is $1$-Lipschitz with respect to Hausdorff distance,
\[
 d(P,R)\ge (\rho+\delta)-\rho\sqrt{|P|/|R^0|}=\delta-O(\delta^3).
\]
On the other hand, domain monotonicity and scale invariance give $\Phi_q(P)\le(|P|/|R^0|)^{2/q}\Phi_q(R^0)=(1+O(\delta^3))\Phi_q(R)$. Among polygons with at most $n$ sides this cannot happen, by Theorem~\ref{S:thm:full-stability}.

The local theorem also covers polygons that are nearly convex in area, whatever their shape.

\begin{corollary}[Nearly convex polygons]\label{C:cor:small-hull-loss}
Let $n\ge3$ and $1\le q\le2$. There is $\eta>0$, depending only on $n$ and $q$, such that every simple polygon $P$ with at most $n$ convex vertices and $|\operatorname{conv}P|<(1+\eta)|P|$ satisfies $\Phi_q(P)\ge\Phi_q(R_n)$, with equality only for $R_n$.
\end{corollary}
\begin{proof}
Otherwise there are polygons $P_j$ of area one, not congruent to $R_n$, with $\Phi_q(P_j)\le\Phi_q(R_n)$ and $|H_j|\to1$, where $H_j=\operatorname{conv}P_j$. If $w_j$ is the width of $H_j$ perpendicular to a diameter of length $d_j$, then $d_jw_j\le2|H_j|$. Since $P_j$ lies in a strip of width $w_j$, H\"older's inequality and the one-dimensional Poincar\'e inequality give $\lambda_q^{\rm S}(P_j)\ge\lambda_1(P_j)\ge\pi^2/w_j^2$. Hence $w_j$ is bounded below and $d_j$ above. The extreme points of $H_j$ are convex vertices of $P_j$, so after translation a subsequence of the $H_j$ converges to a convex polygon $H$ of area one with at most $n$ sides. By domain monotonicity $\Phi_q(H_j)\le|H_j|^{2/q}\lambda_q^{\rm S}(P_j)\le|H_j|^{2/q}\Phi_q(R_n)$, so $\Phi_q(H)\le\Phi_q(R_n)$, and Theorem~\ref{thm:main} shows that $H$ is congruent to $R_n$.

The polygons $P_j$ also converge to $H$ in Hausdorff distance. Otherwise there are $\delta>0$ and points $x_j\in H_j$ with $\operatorname{dist}(x_j,P_j)\ge\delta$ along a subsequence. The disk $B(x_j,\delta)$ then misses $P_j$, and its intersection with the convex set $H_j$, whose inradius is bounded below, has area at least a positive constant times $\delta^2$. This contradicts $|H_j\setminus P_j|\to0$. Theorem~\ref{C:thm:local} now gives a contradiction for large $j$.
\end{proof}

\section{Four convex vertices}\label{C:sec:four}
\begin{theorem}[Four convex vertices]\label{C:thm:four-convex-sharp}
Let $P$ be a simple polygon with at most four convex vertices. Then $\Phi_q(P)\ge\Phi_q(R_4)$ for every $1\le q\le2$, with equality only for squares. The number of reflex vertices is unrestricted.
\end{theorem}

Two Steiner symmetrizations first make $P$ symmetric in two perpendicular axes. In one quadrant we then straighten the boundary by an area-preserving motion that lowers $\lambda_q^{\rm S}$. As in Section~\ref{C:sec:three}, the sign of the change comes from a reflection across the perpendicular bisector of a side; here the axes of symmetry supply the Neumann condition that makes the reflection work.

\subsection{Reduction to polygons with two axes of symmetry}
Let $\mathcal Q$ be the class of simple polygons $Q$ that are symmetric in both coordinate axes and whose boundary in the closed first quadrant is the graph of a convex, strictly decreasing, piecewise affine function joining a point $(a,0)$ to a point $(0,b)$. Such a polygon has exactly four convex vertices, one on each half-axis; its other vertices are reflex.

The first part of the proof of Theorem~\ref{C:thm:four-turn-sweep} shows that, after a small perturbation into general position, two Steiner symmetrizations in perpendicular directions take $P$ to a member of $\mathcal Q$. Steiner symmetrization preserves the area and the $L^q$ norm of a nonnegative function and does not increase its Dirichlet energy, so it does not increase $\Phi_q$. The perturbation changes $\Phi_q$ arbitrarily little, by the piecewise affine transports used at the end of that proof. Hence it suffices to show $\Phi_q(Q)\ge\Phi_q(R_4)$ for every $Q\in\mathcal Q$. The equality case is treated at the end.

\subsection{A reflection on the first side}
Fix $Q\in\mathcal Q$. Interchanging the axes if necessary, we assume $a\ge b$. Let $u$ be the positive minimizer for $\Lambda=\lambda_q^{\rm S}(Q)$ with $\int_Qu^q=1$, so $-\Delta u=\Lambda u^{q-1}$ (with right side $\Lambda$ when $q=1$). It is unique (see the proof of Proposition~\ref{C:prop:dihedral-sharp}), hence even in both coordinates, and it satisfies Neumann conditions on the axes inside $Q$.

We first show $u_x<0$ and $u_y<0$ in the open first quadrant. Let $k>0$. Since every horizontal section of $Q$ is an interval centred on the $y$-axis, the reflection of $Q\cap\{x>k\}$ in the line $x=k$ lies in $Q$. Let $v$ be the reflected function on $Q\cap\{x>k\}$ and put $w=(u-v)_+$ there, extended by zero; then $w\in H^1_0(Q)$. With $p=q-1$, the inequality $u^p-v^p\le u^{p-1}(u-v)$ for $u>v\ge0$ gives
\[
 \int|\nabla w|^2=\Lambda\int(u^p-v^p)w\le\Lambda\int u^{q-2}w^2 .
\]
Testing the equation of $u$ with $w^2/(u+\delta)$ and letting $\delta\downarrow0$ gives the reverse inequality, with difference $\int|\nabla w-(w/u)\nabla u|^2$. Hence $w/u$ is constant, and since $w$ vanishes on an open set, $w=0$. Thus $u_x\le0$. The function $-u_x\ge0$ satisfies $-\Delta(-u_x)=\Lambda(q-1)u^{q-2}(-u_x)\ge0$ and is not identically zero, so it is positive. The same argument applies to $y$.

Use the coordinates $s=a-x$ and $y$. The boundary of $Q$ in the quadrant is the graph $y=h(s)$, $0\le s\le a$, of a convex increasing function with $h(0)=0$ and $h(a)=b$, and $Q$ lies below it. In these coordinates $u_s>0$ and $u_y<0$. Let $A=(0,0)$ be the tip and $J=(d,td)$ the other end of the first side. If the graph has at least two sides, then $0<t<b/a\le1$ by convexity.

\begin{lemma}\label{C:lem:four-reflection}
Suppose the boundary of $Q$ in the quadrant has at least two sides, and let $S$ be the reflection in the perpendicular bisector of $AJ$. Then $|\partial_\nu u(SX)|>|\partial_\nu u(X)|$ for every $X$ on the open segment from $A$ to the midpoint $M$ of $AJ$.
\end{lemma}
\begin{proof}
The bisector is $s+ty=c$ with $c=d(1+t^2)/2$, and it meets the axis at $E=(c,0)$, where $0<c<d$. Let $C$ be the open triangle $AME$. The reflection $S$ fixes $M$ and $E$ and maps $A$ to $J$, so $S(C)$ is the triangle $JME$, which lies in $\{0<y<ts,\ 0<s<d\}\subset Q$. The open side $AE$ lies on the axis, and its image, the open segment $JE$, lies in $Q$. Put $\theta=2\arctan t\in(0,\pi/2)$. For $X$ on $AE$ and the outer normal $\nu=(0,-1)$ of $C$ there,
\[
 \partial_\nu(u\circ S)(X)=\sin\theta\,u_s(SX)-\cos\theta\,u_y(SX)>0,
 \qquad \partial_\nu u(X)=0 .
\]
Let $v=u\circ S$ on $C$. Then $-\Delta v=\Lambda v^{q-1}$, $v=u=0$ on $AM$, and $v=u$ on $ME$. Put $w=(u-v)_+$ in $C$ and extend it by zero to the part $Q^+$ of $Q$ in the open quadrant; $w$ vanishes on $\partial Q$. Green's formula on $C$ gives
\[
 \int|\nabla w|^2=\Lambda\int(u^p-v^p)w-\int_{AE}w\,\partial_\nu v\,ds\le\Lambda\int u^{q-2}w^2 .
\]
On $Q^+$, $u$ solves its equation with Neumann conditions on the axes, so the argument above again gives the reverse inequality and $w=0$. Thus $v\ge u$ in $C$, and $v\not\equiv u$ because their normal derivatives differ on $AE$. Since $-\Delta(v-u)=\Lambda(v^p-u^p)\ge0$, the strong maximum principle gives $v>u$ in $C$. Both functions vanish on $AM$, so Hopf's lemma gives $|\partial_\nu v|>|\partial_\nu u|$ on the open segment $AM$. As $S$ maps the line $AJ$ to itself, $|\partial_\nu v(X)|=|\partial_\nu u(SX)|$.
\end{proof}

\subsection{Straightening the boundary}
Suppose the boundary in the quadrant has at least two sides, and let $JK$ be the second side, of slope $t_1>t$. In the coordinates $(s,y)$, move $A$ along the axis with velocity $(1,0)$ and $J$ along $JK$ with velocity $t(1,t_1)/(t_1-t)$, keep the other vertices of the quadrant fixed, and move the other quadrants symmetrically. The normal velocity $V$ on $AJ$ is affine, equal to $-t/\sqrt{1+t^2}$ at $A$ and $t/\sqrt{1+t^2}$ at $J$, and it vanishes on all other sides. Its integral is zero, and adjusting the speed of $J$ by the implicit function theorem gives a motion that preserves the area exactly. For small times the polygon stays in $\mathcal Q$ with the same number of sides, and it stays inside any square that contains it, since $A$ moves inwards and $J$ moves towards $K$.

Along this motion
\[
 \frac{d}{d\varepsilon}\lambda_q^{\rm S}=-\int_{\partial Q}|\partial_\nu u|^2V\,ds .
\]
This follows as in the proof of Proposition~\ref{C:prop:first-order}: piecewise affine transports carry the minimizers, and the corners contribute nothing because $|\nabla u|=O(r^{\gamma-1})$ with $\gamma>1/2$ near a corner of angle less than $2\pi$. Since $V(SX)=-V(X)$ and $V<0$ on $AM$, Lemma~\ref{C:lem:four-reflection} gives
\[
 \int_{AJ}|\partial_\nu u|^2V\,ds=\int_{AM}V(X)\bigl(|\partial_\nu u(X)|^2-|\partial_\nu u(SX)|^2\bigr)\,ds>0 .
\]
Hence $\lambda_q^{\rm S}$ strictly decreases at fixed area.

\begin{proof}[Proof of Theorem~\ref{C:thm:four-convex-sharp}]
Let $Q\in\mathcal Q$ have $K$ sides in the quadrant, and fix a closed square centred at the origin that contains $Q$. Consider the members of $\mathcal Q$ with the area of $Q$, at most $K$ sides in the quadrant, contained in this square, with the axes allowed to be interchanged. In a minimizing sequence for $\Phi_q$, pass to a subsequence along which the number of sides is constant and the vertices converge. Remove sides of zero length, and remove limit sides lying on the axes; they can occur only at the two ends of the chain, since the slopes stay ordered. Merging collinear sides gives a member of the class with the same area. The minimizers, extended by zero, converge weakly in $H^1$ and strongly in $L^q$ along a further subsequence; the limit vanishes outside the limit polygon and lies in its Dirichlet space. Lower semicontinuity shows that the minimum is attained.

By the motion above, a minimizer has a single side in each quadrant, so it is a rhombus. Theorem~\ref{thm:main} with $n=4$ gives $\Phi_q\ge\Phi_q(R_4)$ for rhombi, and hence for every member of $\mathcal Q$. By the reduction, $\Phi_q(P)\ge\Phi_q(R_4)$ for every simple polygon with at most four convex vertices.

If a nonconvex $P$ attained equality, Lemma~\ref{C:lem:fill-reflex} would give a polygon with the same number of convex vertices and smaller $\Phi_q$, contradicting the inequality just proved. A convex $P$ has at most four sides, and the equality case of Theorem~\ref{thm:main} shows that it is a square.
\end{proof}

At $q=1$ the theorem sharpens the torsion bound of Theorem~\ref{C:thm:four-turn-sweep}, whose geometric part it uses. The approximation in the proof of Corollary~\ref{C:cor:three-arcs} gives $\Phi_q(\Omega)\ge\Phi_q(R_4)$ for the curved domains described there with at most four junctions of positive turn, though without the equality case.

The reflection argument of Section~\ref{C:sec:three} does not extend directly to four convex vertices. In the square $[0,8]^2$ push the four sides in to the points $(4,1),(7,4),(4,7),(1,4)$, keeping the corners. At $(0,0)$ the incident edge vectors $D=(4,1)$ and $F=(1,4)$ satisfy $|D|^2=|F|^2=17$ and $2D\cdot F=16$, so the perpendicular bisector of either incident edge misses the other edge. Even the whole component of the domain cut off by the bisector $4x+y=17/2$ near the origin fails the reflection test of Section~\ref{C:sec:three}: points of the domain near the path $(0,0),(1,4),(0,8)$ lie in this component, and the reflection of $(0,8)$ is $(4/17,137/17)$, outside the square.

\section{The general case}\label{C:sec:general}
\begin{question}\label{C:que:sharp}
Let $n\ge5$. Does every simple polygon $P$ with at most $n$ convex vertices satisfy $\Phi_q(P)\ge\Phi_q(R_n)$ for $1\le q\le2$, with equality only for $R_n$?
\end{question}
The number of reflex vertices is unrestricted. For $n=3$ and $n=4$ the inequality holds (Theorems~\ref{C:thm:three-convex-turns-sharp} and~\ref{C:thm:four-convex-sharp}). For every $n$ it is yes under rotational symmetry (Theorem~\ref{C:thm:cyclic-sharp}) or the equal-area sector condition (Corollary~\ref{C:cor:equal-sector-turns}). Without symmetry, Theorem~\ref{C:thm:convex-vertices} and Corollary~\ref{C:cor:all-q-gap} give the comparison with $R_{n+1}$ instead of $R_n$, and the answer is yes for all polygons close to $R_n$ (Theorem~\ref{C:thm:local}). Hence a polygon violating the inequality must have exactly $n$ convex vertices, by Corollary~\ref{C:cor:all-q-gap}, must stay a fixed Hausdorff distance away from $R_n$ after normalization to unit area and optimal rigid alignment, and must have a convex hull larger than itself by a fixed fraction of its area (Corollary~\ref{C:cor:small-hull-loss}). For many convex vertices, Theorem~\ref{C:thm:sharp-asymptotic} below shows that the least torsion deficit agrees with that of $R_n$ up to $O(n^{-5})$.

\subsection{Many convex vertices}\label{C:sec:many}
The weight of Section~\ref{C:sec:uniform-modes} also gives the two leading terms of the least torsion deficit, with an explicit remainder.

\begin{theorem}[Two-term torsion deficit]\label{C:thm:sharp-asymptotic}
Let $m\ge5$ and let $P$ be a simple polygon of area $\pi$ with at most $m$ convex vertices. Then
\[
 D(P)>\frac{8\zeta(3)}{m^3}-\frac{16\zeta(4)}{m^4}-\frac{160}{m^5}.
\]
Consequently, if $d_m^*$ denotes the infimum of $D(P)$ over this class, then as $m\to\infty$
\begin{gather*}
 d_m^*=\frac{8\zeta(3)}{m^3}-\frac{16\zeta(4)}{m^4}+O(m^{-5}),\\
 0\le D(R_m)-d_m^*=O(m^{-5}),\qquad
 d_m^*-D(R_{m+1})=\frac{24\zeta(3)}{m^4}+O(m^{-5}).
\end{gather*}
\end{theorem}

The lower bound is uniform in the number of reflex vertices, and no minimizing polygon is needed. The last relation shows that the gaps of order $n^{-4}$ in Theorem~\ref{C:thm:convex-vertices} have the correct order for torsion.

\begin{proof}
Put $v=1/m$ and $t=1/(m+1)$, and use the notation of Section~\ref{C:sec:uniform-modes}, with $Z=2Y-X_2$, $b=1/(1+S)$ and $B=1-b$.

\emph{A lower bound for $Q$.} As in the proof of Theorem~\ref{C:thm:turn-gap}, Lemmas~\ref{C:lem:integrated} and~\ref{C:lem:turn-modal} give $Q\ge I/A+tZ$. By \eqref{C:eq:J-lower}, $K(v)\ge1-3v$ and $v/(1-4v^2)\ge v$, we have $J\ge\zeta(3)-2\zeta(4)v$. Write $\kappa=1-2v+2v^2$, so that $1/A=1-2v^2/\kappa$ and $\kappa\ge17/25$ for $v\le1/5$. Since $I=4v^3J$ and $\zeta(3)<121/100$,
\[
 \frac IA\ge4v^3\bigl(\zeta(3)-2\zeta(4)v\bigr)\Bigl(1-\frac{2v^2}{\kappa}\Bigr)
 >4\zeta(3)v^3-8\zeta(4)v^4-\frac{242}{17}v^5 .
\]

\emph{The cubic term.} Put $R=Q-I/A\ge tZ$, and let $\mathcal C=\sum_{i,j\ge2}z_{i+j}\overline{z_iz_j}$ be the cubic sum in \eqref{C:eq:homogeneous-hankel}. Then $|\mathcal C|\le H$ in the notation of Lemma~\ref{C:lem:ordered-hankel}, and that lemma gives $H^2\le(Q-Y)\,YZ/2\le Q^2Z/8$. Hence $|\mathcal C|\le Q\sqrt{R/(8t)}$ and
\[
 2R-4\operatorname{Re}\mathcal C\ge 2R-Q\sqrt{2R/t}
 =2\Bigl(\sqrt R-\frac{Q}{2\sqrt{2t}}\Bigr)^2-\frac{Q^2}{4t}\ge-\frac{Q^2}{4t}.
\]
Since $\mathcal D(E)\ge0$ and $D/b^2=D(1+S)^2$, \eqref{C:eq:homogeneous-hankel} gives
\begin{equation}\label{C:eq:asym-hankel}
 D(1+S)^2=2Q-4\operatorname{Re}\mathcal C+\mathcal D(E)\ge\frac{2I}{A}-\frac{Q^2}{4t}.
\end{equation}

\emph{The remainders.} We may assume $D\le10v^3\le2/25$, since $8\zeta(3)<10$. By \eqref{C:eq:B-deficit}, $B\le2\sqrt D<6/7$ and $B(1-B)\le3D/2\le3/25<6/49$; hence $B<1/7$ and $b>6/7$. Then $S=B/b\le3D/(2b^2)<21D/10\le21v^3$. By \eqref{C:eq:pooled-deficit} with $V=2Y$ and $Q^\sharp\le Q$,
\[
 D(1+S)^2\ge2Q+6Y^2-4Y\sqrt Q=2Q+6\Bigl(Y-\frac{\sqrt Q}{3}\Bigr)^2-\frac{2Q}{3}\ge\frac{4Q}{3},
\]
so $Q\le\frac34(1+S)^2D<\frac{21}{20}D\le\frac{21}{2}v^3$. Therefore
\[
 \frac{Q^2}{4t}\le\frac{441}{16}(1+v)v^5\le\frac{1323}{40}v^5,\qquad
 D(2S+S^2)\le420v^6+4410v^9\le\frac{11382}{125}v^5 .
\]
The last inequality is an equality exactly at $v=1/5$, since $420v+4410v^4$ is strictly increasing on $(0,1/5]$.
Combining these remainder bounds with \eqref{C:eq:asym-hankel} and the strict lower bound for $I/A$ gives
\[
 D>8\zeta(3)v^3-16\zeta(4)v^4-Cv^5,\qquad
 C=\frac{484}{17}+\frac{1323}{40}+\frac{11382}{125}=\frac{2594227}{17000}<153 .
\]

\emph{The regular polygon.} For the regular map,
$S_R=Z(v)-1=4\zeta(3)v^3+O(v^5)$ by \eqref{C:eq:regular-area}, and
\[
 Y_R=4v^4\sum_{j\ge1}\frac{\prod_{l<j}(1+2v/l)^2}{j^4(1+v/j)^2}
       =4\zeta(4)v^4+O(v^5).
\]
Here the products start at $l=1$, and the numerator differs from
one by at most $4vH_{j-1}e^{4vH_{j-1}}$, where $H_{j-1}$ is the harmonic number.
For $v\le1/8$ this is bounded by a constant times
$v(1+\log j)j^{1/2}$, so the summable majorant
$(1+\log j)j^{-7/2}$ gives the uniform $O(v^5)$ remainder.
Thus $Q_R=S_R-2Y_R=4\zeta(3)v^3-8\zeta(4)v^4+O(v^5)$. The regular coefficients in \eqref{C:eq:homogeneous-hankel} are nonnegative, so the cubic term has the favorable sign, and $\mathcal D(E)\le S^2$ gives $D(R_m)\le2Q_R+S_R^2=8\zeta(3)v^3-16\zeta(4)v^4+O(v^5)$.

Together with the lower bound, $d_m^*\le D(R_m)$ gives the first two relations. Applying them with $m+1$ in place of $m$ gives $D(R_{m+1})=8\zeta(3)v^3-(24\zeta(3)+16\zeta(4))v^4+O(v^5)$, and the third relation follows.
\end{proof}

\subsection{Attainment}

Fix the area and an upper bound $m$ on the number of convex vertices. By Lemma~\ref{C:lem:fill-reflex}, a polygon attaining the least value of $\Phi_q$ in this class, for some $1\le q\le2$, has no reflex vertex. It is then convex, so it has at most $m$ sides, and Theorem~\ref{thm:main} identifies it as $R_m$. Thus each case of Question~\ref{C:que:sharp} is equivalent to the attainment of the infimum by a polygon. Filling a corner increases the number of reflex vertices, and a minimizing sequence may have unbounded side count, with a limit whose boundary is not polygonal. The lemma therefore does not prove attainment.

\appendix
\renewcommand{\thesection}{\Alph{section}}
\renewcommand{\thesubsection}{\thesection.\arabic{subsection}}
\renewcommand{\theHsection}{appendix.\arabic{section}}
\renewcommand{\theHsubsection}{\theHsection.\arabic{subsection}}

\section{Numerical estimates used in the proofs}\label{app:estimate-register}
This appendix collects the numerical inequalities used in the proofs, so that they can be checked in one place. Each row lists an estimate and, where applicable, its margin over the value required in the proof. All decimals denote exact rational numbers. When the proof clears a denominator, the stated difference is for the resulting numerator.
Table~\ref{tab:acute-coverage} gives the ranges of the estimates at acute openings.

\begingroup
\small
\setlength{\tabcolsep}{4pt}
\setlength{\LTpre}{3pt}
\setlength{\LTpost}{3pt}
\begin{longtable}{@{}>{\raggedright\arraybackslash}p{.205\linewidth}>{\raggedright\arraybackslash}p{.392\linewidth}>{\raggedright\arraybackslash}p{.346\linewidth}@{}}
\toprule
Proof location & Inequality and use & Margin or scale\\
\midrule
\endfirsthead
\toprule
Proof location & Inequality and use & Margin or scale\\
\midrule
\endhead
\bottomrule
\endfoot

\multicolumn{3}{@{}l}{\textbf{The cases $q=1$ and $q=2$; supporting lines}}\\*
\nameref{par:independent-real-projection} (Lemma~\ref{lem:acute-eliminate})
& At $h=1/6$ the complement has coercivity at least
$5/29-1/24-1/72=61/522$.
The two frequency ranges give forcing bounds
$(.0087,.0284,.0911)t^{3/2}$ for the three derivatives.
& The difference between the coercivity bound and $.1168$ is $19/326250$.
The margins in the three forcing bounds, in units of $t^{3/2}$, are at least
$(.0003,.0006,.0039)$.\\
\addlinespace[3pt]
Proposition~\ref{prop:acute-curvature}, $s\le1/6$
& The projected derivatives satisfy $P<.0174$, $Q<.142$, so
$(Q+13P/6)/h^3<38.82<39$ while $48\zeta(3)>57.6$.
& Since $j^2\ge\pi^2/2$ and $1+\gamma<1.003$, the lower bound
for $F_2''(2\pi s)/s^2$ exceeds $7$. The proof requires $4$.\\
\addlinespace[3pt]
Proposition~\ref{prop:endpoint-A5}, $1\le z\le3$
& After subtracting $(1-h^2/4)L$ from the polynomial $P$, the remainder is at least
$h^2[27(y-13/36)^2+19/240]$.
& The lower bound $15/56$ for $F_2''$ exceeds $1/4$ by $1/56$.\\
\addlinespace[3pt]
Proposition~\ref{prop:square-tangent} and \eqref{eq:spectral-quadratic-gap}
 & In (T90), $DL\ge[P^2/(4(P-1))](1+w^3)$. For obtuse openings the tangent at $\pi/2$ reduces the quadratic comparison to its two endpoint values.
 & $DL-1\ge(P-2)^2/[4(P-1)]$. The angular coefficient is $1/1369$ at $n=3,4$, $3/4$ at $n=5$ and $1$ at $n\ge6$, in units of $(\alpha-a)^2/n^2$. For $n\ge5$ stability uses $1/2$.\\
\addlinespace[3pt]
\multicolumn{3}{@{}l}{\textbf{Comparison with the disk; obtuse openings}}\\*
Lemma~\ref{S:lem:shared-disk-barrier};
Proposition~\ref{S:prop:lowq}
& For $p=q-1\le1801/4000$, the disk test gives $X\ge(1801/4000)W$.
The bound \eqref{S:eq:twoform} and integration of \eqref{S:eq:curv} give $1/\beta\ge3z+29/11$.
& The power supersolution gives $k>(99/200)(1-p/5)$
for $0<p\le1/2$.
The comparison covers $p\le p_*$ since $p_*<495/1099$.
This gives $\mathcal C_q\ge[10(1+z)^2]^{-1}$.\\
\addlinespace[3pt]
Section~\ref{S:sec:compact-curvature};
Proposition~\ref{S:prop:wide-linear}
& For $t=z^{1/3}$, $P=\pi^2/4$, the first two unweighted eigenvalues satisfy
$\lambda/P\le1+61t/50$ and $\lambda_2/P\ge(1+t)^3$ ($z\le1/3$).
& The gap is at least $89Pt/50$.
The other term $5P$ in the minimum of Lemma~\ref{S:lem:horizontal-projection} exceeds $P(1+t)^3$ by more than $87P/1000$, since $t<7/10$.\\
\addlinespace[3pt]
Lemma~\ref{S:lem:horizontal-sturm};
Lemma~\ref{lem:obtuse-projection}
& $X\ge(P/\Lambda_2)W$. The Green operator on $L^2((0,1),dx/x)$
has $\operatorname{tr}(L^3)=1/3$.
& In (W) the constant coefficient is at most
$\log(P/\sqrt[3]{3})-4/27<2/5$.
At $z\le1/3$, $q\le27/16$, the direct estimate already gives
$\mathcal C_q\ge109/3840$.\\
\addlinespace[3pt]
Proposition~\ref{S:prop:spectral-wide},
first test in \eqref{S:eq:widechecks}
& For $z\le2/9$, $R_\delta/J_\delta<4t/3<(3-61/50)\log(1+t)$.
For $z\ge2/9$, $R_\delta/J_\delta<1$.
& On the latter interval, $v/u\ge640/231>11/4>e$.\\
\addlinespace[3pt]
Proposition~\ref{S:prop:spectral-wide},
second test in \eqref{S:eq:widechecks}
& The lower bound for the second eigenvalue in the longitudinal variable and the spectral gap at $z=1/3$ give $v\ge(1+21t/20)^3$ for all $z\le1/3$.
AM--GM gives $R_\delta/\tau<8/15$.
& Linear interpolation gives $(7/5)\log(1+21t/20)-\log(1+61t/50)>(3/14)t$, which exceeds the required $(16/75)t$ by $t/1050$. Together the tests give $\mathcal C_q>1/1000$.\\
\addlinespace[3pt]
Lemmas~\ref{S:lem:eqspec} and~\ref{S:lem:square-slope};
\eqref{S:eq:local-c1}
& The bound $\nu_2/\nu_1\ge4$ at $z=1/3$ and the eigenvalues $2P,10P$ at $z=1$ give a lower bound for $\lambda+\lambda_2$ by concavity. Subtracting the trial upper bound for $\lambda$ bounds $\lambda_2$ from below.
Integration from $z=1$ starts at $\beta_q(1)<2/15$.
& The difference $2/15-\beta_q(1)$ is at least
$\pi^2/16-37/60>1/10000$ (use $\pi>3.1415$).
The resulting bound on $\lambda_2$ is exact at $z=1$.\\
\addlinespace[3pt]
Proposition~\ref{S:prop:central-sum},
\nameref{par:middle-uniform}
& On $[1/3,1]$, $B/v\le13/16$ and $BH/S\le16/9$.
& The ratio $(B-Sr)e^{rH}/v$ is at most $(117/256)e^{7/9}<1$.
\\
\addlinespace[3pt]
Proposition~\ref{S:prop:central-sum},
\nameref{par:middle-backward}
& Substitution of $b=14/(105/4+74z)$ gives the quadratic
$(9168b^2-3901b+416)/1480$ as a lower bound for the curvature.
& Its minimum $37751/54274560$ exceeds the required $1/1500$ by more than $1/40000$.\\
\addlinespace[3pt]
Proposition~\ref{S:prop:wide-linear} and \nameref{par:spectral-obtuse}
& At $q=2$ the estimates for $z\le1/3$ and for $1/3\le z\le1$ give
$\mathcal C_2>1/16$ and $\mathcal C_2\ge1/37^2$, respectively.
& For $z\le1/3$ the two terms subtracted from one in the proof are at most $189/256$ and $39/200$.
These are curvature bounds for $q=2$ at obtuse openings.\\

\multicolumn{3}{@{}l}{\textbf{Curvature for $29/20\le q\le2$ and $0<\alpha\le\pi/2$}}\\*
\nameref{app:radial-gap}
 & The radial eigenvalues and the estimates for the nonzero angular modes give $pw\le3/8$. In Section~\ref{sec:cylinder-green-sources}, each source norm has a bound containing
 $(4s)^{3/2+(4q-2j-1)/6}$ for $j=0,1,2$.
 & On the energy-orthogonal complement of the minimizer, $E_s-pB\ge5E_s/8$.\\
\addlinespace[3pt]
Section~\ref{sec:cylinder-curvature-bound}
 & The tangent to the polynomial $A_0$ at $t=13/16$ exceeds the upper bound for the error terms $19.42+72.04t$ by more than one on $0\le t\le1$.
 & $f''/(csf)>1$, giving $\mathcal C_q>s/(4\pi^2)$ and hence $\mathcal C_q>1/200$ for $2\pi/5\le\alpha\le\pi/2$.\\

\multicolumn{3}{@{}l}{\textbf{Stability}}\\*
Sections~\ref{sec:small-deficit-convexity}--\ref{sec:stability-cell-geometry}, area stability
 & With $e=\sqrt{\Delta/L_0}$, the angular errors use
 $1/(\sqrt{18}\,\pi^3)<1/100$ and $\sqrt{256/(72\pi^7)}<1/25$.
 The geometric bound is $d(P)\le16ne$.
 & The margins in these two bounds, in units of $e$, exceed $.002$ and $.005$, respectively.
 $L_0/(16n)^2=(18432\pi^6n^5)^{-1}$.\\
\addlinespace[3pt]
Lemma~\ref{lem:torsion-curvature} and \eqref{eq:stability-angle-torsion}
 & On $|\alpha-a|\le a/8$, the flux estimate gives $-g''/a^3\ge 2/375$ where the proof needs $1/256$. Outside this neighborhood convexity preserves the endpoint gap.
 & The difference is $2/375-1/256$ in units of $a^3$. Integration gives $K_T=\pi/(256n^3)$.\\
\addlinespace[3pt]
Three and four sides, p.~\pageref{par:small-n-area}
 & At $x_0=1/12{,}800{,}000$ the parameter $e$ in \eqref{eq:one-cap-data} is $1/64$, and for $x\le x_0$ the distance is at most $x/x_0$.
 & $x_0^2/(72\cdot73)>10^{-18}$ for the torsion deficit.\\
\addlinespace[3pt]
Section~\ref{sec:perimeter-geometry}, fixed perimeter
 & The sharp quadratic coefficient of the tangent remainder for $\tan(\theta/2)$ is $c_a=(a\sec^2(a/2)/2-\tan(a/2))/a^2$.
 & The coefficient of $L\sqrt\delta$ is less than $1.870$ for $n=3$, $1.785$ for $n=4$, and $1.779$ for $n\ge5$; the theorem uses $2$.\\
\end{longtable}
\endgroup

\section*{AI use declaration}

The author used AI tools to suggest improvements and to help polish the exposition. The author is responsible for all content.

\end{document}